\documentclass[reqno,english]{amsart}

\usepackage{amssymb,latexsym,amsmath,amsthm}
\usepackage{graphicx,color, amsfonts}

\usepackage[colorlinks=true]{hyperref}
\hypersetup{urlcolor=blue, linkcolor=blue, citecolor=red, anchorcolor=blue}

\usepackage{upgreek}
\usepackage{framed}
\usepackage{comment}
\usepackage{import}

\newtheorem{theorem}{Theorem}[section]
\newtheorem{proposition}[theorem]{Proposition}
\newtheorem{lemma}[theorem]{Lemma}
\newtheorem{corollary}[theorem]{Corollary}
\newtheorem{remark}[theorem]{Remark}

\theoremstyle{definition}

\numberwithin{equation}{section}

\newcommand{\R}{\mathbb R}
\newcommand{\N}{\mathbb{N}}
\newcommand{\C}{\mathbb{C}}
\newcommand{\Z}{\mathbb{Z}}
\newcommand{\abs}[2][]{#1\lvert #2 #1\rvert}

\let\Re\relax
\let\Im\relax

\DeclareMathOperator{\Re}{Re}
\DeclareMathOperator{\Im}{Im}

\long\def\hide#1{}

\newcommand{\norm}[1]{\left \| #1 \right \|}
\newcommand{\D}[1]{\mathcal{D}_{#1R}}

\newcommand{\LtU}[1]{\tilde{L}_U\left({#1}\right)}
\newcommand{\piU}[1]{\pi_{U^\perp}\big({#1}\big)}

\newcommand{\normZ}{\|Z_0\|_{0}}
\renewcommand{\Re}[1]{\text{Re}\left({#1}\right)}
\renewcommand{\Im}[1]{\text{Im}\left({#1}\right)}

\def\smath#1{\text{\scalebox{.8}{$#1$}}}
\def\sfrac#1#2{\smath{\frac{#1}{#2}}}

\includecomment{commentA}

\begin{document}

\title[Stability HMF]{Stability of least-energy solutions of the harmonic map flow from the plane into the 2-sphere} 

\begin{abstract}

We prove the asymptotic stability of the least-energy harmonic map $W:\mathbb R^2\to S^2$ for the harmonic map flow under small, sufficiently decaying, non-symmetric perturbations. Solutions exist globally and are asymptotically described by a modulation of $W$ plus a remainder that converges to zero as $t\to\infty$. The modulation parameters associated with  scaling, rotation, and translation remain close to their initial values and converge to finite limits. The proof relies on an inner-outer gluing decomposition and the analysis of a non-local equation governing the modulation parameters.

\end{abstract}

\author[J.~D\'avila]{Juan D\'avila}
\address{\noindent J.~D\'avila: Department of Mathematical Sciences University of Bath, Bath BA2 7AY, United Kingdom}
\email{jddb22@bath.ac.uk}

\author[E.~Van Steirteghem]{Eliott Van Steirteghem}
\address{\noindent E.~Van Steirteghem: Department of Mathematical Sciences University of Bath, Bath BA2 7AY, United Kingdom}
\email{eliottvans@hotmail.com}

\maketitle

\section{Introduction}

We consider the harmonic map flow from $\mathbb{R}^2$ into the sphere
$S^2\subset\mathbb{R}^3$:
\begin{equation}
\label{eq:main-HMF}
\begin{aligned}
\partial_t u &= \Delta u + |\nabla u|^2 u
&& \text{in } \mathbb{R}^2\times(0,T),\\
u(\cdot,0) &= u_0
&& \text{in } \mathbb{R}^2,
\end{aligned}
\end{equation}
where $u_0:\mathbb{R}^2\to S^2$ is a sufficiently regular initial condition.
Here $T\in(0,+\infty]$ denotes the maximal time of existence of a smooth
solution. The equation is the $L^2$-gradient flow of the Dirichlet energy
\[
E(u)=\frac12\int_{\mathbb{R}^2}|\nabla u|^2\,dx
\]
for maps $u:\R^2\to S^2$.

The harmonic map flow was first introduced by Eells and Sampson \cite{Eells} as a
method for constructing harmonic maps in a prescribed homotopy class.
When the target manifold has non-positive sectional curvature, they
established global existence and convergence of the flow to a harmonic map.
For positively curved targets, however, finite-time singularity formation
may occur.
The first example of a finite-time blow-up solution in dimension two was
constructed by Chang, Ding and Ye \cite{Chang} for $1$-corotational maps
from a disk into $S^2$. In dimension two, a fundamental result of Struwe
\cite{Struwe} shows that, for finite-energy initial data, solutions exist
globally in a weak sense and are smooth away from finitely many singular
points. Moreover, singularities occur through the bubbling of finite-energy
harmonic maps. Finite-energy harmonic maps from $\mathbb{R}^2$ into $S^2$
are classified by their degree, and their energy is quantised according to
\[
E(B)=4\pi |k|,
\]
where $k\in\mathbb{Z}$ is the degree of $B$. The simplest non-trivial example
is the degree-one harmonic map
\[
W(x)=\frac{1}{1+|x|^2}
\begin{pmatrix}
2x\\
|x|^2-1
\end{pmatrix},
\]
which has energy $4\pi$ and is obtained from the inverse stereographic
projection.
There are several further works on the analysis of blow-up scenarios
\cite{Ding-Tian,lin-wang-1998,Qing,Qing-Tian,topping_anals,wang-1996,jendrej-lawrie-schlag}.

Two central questions in the study of harmonic map flow concern the dynamics
near singularities and the dynamics near harmonic maps. In the first
direction, for $1$-corotational maps, formal asymptotics for finite-time
blow-up solutions were derived in \cite{vdBBHK}, and the corresponding
blow-up rate was later rigorously established in
\cite{Raphael} under the same symmetry assumption. Finite-time blow-up solutions without imposing corotational symmetry were later constructed in \cite{HMF_DdPW}; near each blow-up point, the leading-order profile is given by $W$. Infinite-time blow-up solutions with a profile of degree 2 were constructed in \cite{ddwz-deg2}.

A second direction concerns the dynamics near harmonic maps. A fundamental question is whether a harmonic map is dynamically stable under the flow: given initial data close to a harmonic map, does the corresponding solution remain close to the family of harmonic maps for all time and eventually converge to another harmonic map of the same energy? 
The stability problem for the harmonic map heat flow from $S^2$ to $S^2$ is well understood. In this setting, the Łojasiewicz-Simon inequality introduced in \cite{Simon,Lojasiewicz1965}, together with the energy dissipation along the flow, allows the control of the evolution of solutions near a harmonic map, see also \cite{topping_almost_harm}. In particular, if the initial datum is sufficiently close to a harmonic map $W$ in a suitably strong topology, then the corresponding solution exists globally and converges to a harmonic map as $t \to \infty$. The situation for maps defined on $\R^2$ is more delicate as the Łojasiewicz-Simon inequality fails, and the stability problem requires a different approach.

Some aspects of the stability problem on $\R^2$ are nevertheless understood under $k$-corotational symmetry, with the behaviour depending on the degree $k$. For $|k|\geq 3$, Guan, Gustafson and Tsai \cite{GuanGusTsai}, and Gustafson, Nakanishi and Tsai \cite{GusTsaiNak}, proved that no blow-up occurs even under the weaker assumption $E(u_0)\leq 4\pi |k| + \delta$.  For $|k|=2$, more complicated phenomena arise: there are initial conditions with energy close to the threshold energy $8\pi$ such that the solution is global but oscillates between different harmonic maps, or may converge pointwise but not uniformly to a constant map (this is infinite-time blow-up or concentration), see  \cite{GusTsaiNak}.
For $|k|=1$, Guan, Gustafson and Tsai \cite{GuanGusTsai} constructed solutions with energy arbitrarily close to $4\pi$ that nevertheless develop finite-time singularity. Therefore, closeness of the energy to its minimal value alone does not imply the stability of degree-one harmonic maps.

Connected to this dynamical question is the problem of \emph{rigidity of minimisers}, that is, the extent to which the proximity of the Dirichlet energy to these discrete values $4\pi |k|$, $k \in \mathbb{Z}$, implies proximity to the set of harmonic maps. In the case of maps of degree $k = \pm 1$, Bernand-Mantel, Muratov and Simon \cite{skyrmions} proved that for each $u \in W^{1, 2}(\R^2, S^2)$ of degree $\pm 1$, there exists a harmonic map $v:\R^2 \to S^2$ of degree $\pm 1$ such that 
$$\int_{\R^2} |\nabla (u - v)|^2 \le C \big(E(u) - 4\pi\big).$$
Their proof was later simplified by Hirsch and Zemas \cite{HZ} and Topping \cite{Topping_rigidity} using tools from geometric analysis. Deng, Sun and Wei \cite{DSW_rigidity} established a similar estimate for $|k|\ge 2$ while assuming an a priori closeness to a harmonic map of degree $k$. Without this assumption, the situation is more complex, and this estimate fails in general. An optimal quantitative estimate for $|k| \ge 2$ was established by Rupflin \cite{Rupflin}.

In this paper, we establish the asymptotic stability of the degree-one harmonic map under general non-symmetric
perturbations.
We consider initial data of the form
\[
u_0=W+Z_0,
\]
where $Z_0$ is sufficiently small in a $C^1$ weighted norm
and decays at infinity.
To make this precise, we introduce the following norm
\begin{align}
\label{eq:normZ0}
\|Z_0\|_0= \sup_{x\in\R^2}
t_0^{d}\Bigl(1 + \frac{\abs{x}}{\sqrt{t_0}}\Bigr)^{2d}
\left[|Z_0(x)| + \sqrt{t_0}|\nabla Z_0(x)|
\right]
\end{align}
depending on the parameters $d>0$ and $t_0>0$.
Our main result
shows that the corresponding solution exists globally and converges, in a
suitable sense, towards a modulated degree-one harmonic map.

\begin{theorem}
\label{thm:main}
Let $d>\frac12$. Then, for all sufficiently large
$t_0$, there exists $\delta>0$ such that if
\begin{align}
\label{eq:ini-Z0}
\|Z_0\|_0<\delta
\qquad\text{and}\qquad
|W(x)+Z_0(x)|=1
\quad\text{for all }x\in\mathbb{R}^2,
\end{align}
the solution $u$ of \eqref{eq:main-HMF} with initial condition
$u_0=W+Z_0$ exists globally in time. Moreover, there exist $C,\sigma>0$
and modulation parameters $(\lambda(t),\omega(t),\xi(t))$ such that
\begin{align}
\nonumber
\Big\|u(t)-Q_{\omega(t)}
W\Big(\frac{\cdot-\xi(t)}{\lambda(t)}\Big)\Big\|_{L^\infty}
+
\Big\|\nabla\Big(
u(t)-Q_{\omega(t)}
W\Big(\frac{\cdot-\xi(t)}{\lambda(t)}\Big)
\Big)\Big\|_{L^2}
\leq
C\frac{\|Z_0\|_0}{(t_0+t)^\sigma}.
\end{align}
Here, $Q_\omega$ denotes the rotation around the $z$-axis by an angle
$\omega$. The parameters satisfy
\[
\lambda(0)=1,\qquad
\omega(0)=0,\qquad
\xi(0)=0,
\]
and
\[
\lambda(t)\to\lambda_\infty,\quad
\omega(t)\to\omega_\infty,\quad
\xi(t)\to\xi_\infty,
\]
as $t\to\infty$, for some $\lambda_\infty$ close to 1, $\omega_\infty$ close 
to 0, and $\xi_\infty$ close to 0.
\end{theorem}

The proof of Theorem~\ref{thm:main} is based on a decomposition
of the solution into inner and outer components
developed in \cite{HMF_DdPW} for the construction of
solutions blowing up in finite time.
In contrast, here we use a similar decomposition
of the solution into a modulated harmonic map, an inner correction near the
harmonic map core, and an outer component 
and find estimates for these components.
The modulation parameters correspond to the symmetries of the problem:
scaling, rotation and translation. The main difficulty is caused by the
scaling and rotational modes, whose associated errors decay slowly.
To address this difficulty, we introduce a non-local correction term $\Phi^0$, inspired by previous gluing constructions for critical flows \cite{HMF_DdPW,delpino-musso-wei-3d-crit,Keller_segel_infinite_time}.

The orthogonality conditions associated with the kernel of the linearized
operator lead to an integro-differential equation for the complex modulation
parameter
\[
p(t)=\lambda(t)e^{i\omega(t)}.
\]
The analysis of this non-local equation is one of the main technical
ingredients of the paper. Rather than solving the equation exactly, we
decompose it into a principal part and a remainder with improved decay.
This decomposition allows us to recover the required orthogonality
conditions through an additional inner correction and ultimately to close
the argument.

The paper is organised as follows. In Section~\ref{sec:Decomposition of the
solution}, we derive the inner--outer gluing system and introduce the
non-local correction. Section~\ref{Sect:orth_mode01} contains the modulation
equations obtained from the orthogonality conditions. The linear theory for
the inner and outer problems is presented in
Sections~\ref{sec:Estimates_inner} and~\ref{subSect:est_HE}, respectively.
The non-local equation governing the parameter $p$ is analyzed in
Section~\ref{Sect:eq_for_p}. Finally, these estimates are combined through
a fixed-point argument to prove Theorem~\ref{thm:main}.
Several definitions and computations are taken from \cite{HMF_DdPW} and reproduced here without proofs.

\section{Decomposition of the solution and the inner-outer gluing system}
\label{sec:Decomposition of the solution}
\subsection{Decomposition of the solution}

Consider the harmonic map heat flow from $\R^2$ into $S^2$ given by 
\begin{equation}
\label{HMF}
\left\{
\begin{aligned}
\partial_t u &= \Delta u + \abs{\nabla u}^2 u & & \text{in } \R^2 \times (t_0,\infty),\\
u(\cdot, t_0) &= u_0,& &\text{in } \R^2,
\end{aligned}
\right.
\end{equation}
starting at $t_0 > 0$ sufficiently large. This formulation is equivalent to \eqref{eq:main-HMF} and leads to simpler expressions for the norms used throughout the paper.

From now on we write 
\begin{equation}
\label{eq:W}
W(x) = \frac{1}{1 + \abs{x}^2} \begin{pmatrix}
2x\\
\abs{x}^2 - 1
\end{pmatrix}, \quad x \in \R^2,
\end{equation}
instead of $W_1$.
We consider $u_0=W+Z_0$ with $\|Z_0\|_0\ll1$, where
$\|\ \|_0$ is defined in \eqref{eq:normZ0} for some $\frac12<d<1$.
We also assume that $|u_0(x)|=1$ for all $x\in\R^2$.
We decompose the solution $u$ into an explicit part plus remainders.
The explicit part is a modulation of $W$ defined by  
\begin{equation}
\nonumber
U(x, t) = Q_{\omega(t)} W\left(y\right), \quad y = \frac{x - \xi(t)}{\lambda(t)},
\end{equation}
for the parameters $(\lambda(t), \omega(t), \xi(t)) \in \R  \times \R \times \R^2$ and $Q_\omega$ is the rotation around the $z$-axis given by the matrix
\[
Q_\omega = \begin{pmatrix}
\cos \omega & -\sin \omega & 0\\
\sin \omega & \cos \omega & 0\\
0 & 0 & 1
\end{pmatrix}.
\]
Moreover, we impose
\begin{align}
\label{eq:ini-param}
\lambda(t_0) = 1, \quad \omega(t_0) = 0 \quad \text{and} \quad \xi(t_0) = 0
\end{align}
so that $U(x, t_0) = W(x)$.

Define the orthogonal projection onto the plane orthogonal to $U$. 
\begin{equation}
\label{eq:def-pi}
\pi_{U^\perp}(v) = v - (v \cdot U) U, \quad v\in\R^3,
\end{equation}
and
\begin{equation}
\label{eq:def-a}
a(v) = \sqrt{1 - \abs{v}^2} - 1, \quad v\in\R^3 \text{ with } |v|<1.
\end{equation}
so that if $|v|<1$, then
\[
|U+\pi_{U^\perp}(v)+ a( \pi_{U^\perp}(v) ) U|=1.
\]
We note that if $|v|<1$ and  $|U+v|=1$, then
\begin{align}
\nonumber
v=\pi_{U^\perp}(v)+ a( \pi_{U^\perp}(v) ) U.
\end{align}

The aim of this section is to rewrite \eqref{HMF} as an equation for the perturbation $\varphi$. 
Let us define 
\begin{align}
S(u) = -u_t + \Delta u + |\nabla u|^2 u. \label{eq_S}
\end{align}
Then $u$ is a solution of \eqref{HMF} if and only if 
$$S(u) = 0.$$
We write $u$ as
\begin{align}
\label{eq:decomp-u}
u =  U + \pi_{U^\perp}(\varphi) + a\left( \pi_{U^\perp}(\varphi) \right) U.
\end{align}
where $\varphi(x,t)\in\R^3$ and satisfies $|\varphi|<1$ (not necessarily satisfying $|U+\varphi|=1$).
Then we derive an equation for $\varphi$, which is equivalent to $S(u)=c(x,t)U$ (so up to terms collinear with $U$). 
The next lemma ensures that $c(x,t)=0$.
\begin{lemma}
\label{lem:colinear_U}
Suppose $u$ is a solution of 
\begin{align}
\label{colinear_U}
S(u) = c(x,t) U
\end{align}
with $|u|=1$ and $u \cdot U\not=0$. Then $c(x,t)=0$.
\end{lemma}
\begin{proof} 
Indeed, taking the scalar product of \eqref{colinear_U} with $u$ yields
\begin{align*}
c(U\cdot u) &= S(u)\cdot u= \left(-u_t + \Delta u + |\nabla u|^2 u\right) \cdot u
\\
&= -\frac{1}{2} \partial_t \abs{u}^2 + \frac{1}{2} \Delta \abs{u}^2 
= 0.
\end{align*}
This implies that $c(x, t) = 0$.
\end{proof}
Since $Z_0$ satisfies \eqref{eq:ini-Z0} and $\varphi(t_0)=Z_0$, we get that 
\[
u(t_0)=W+Z_0.
\]
Substituting $u = U + \tilde\varphi$ in \eqref{eq_S} gives us 
\begin{align*}
S(U + \tilde\varphi) 
&= -U_t - \tilde\varphi_t + L_U(\tilde\varphi) + \abs{\nabla \tilde\varphi}^2 U + 2(\nabla U : \nabla \tilde\varphi)\tilde\varphi + \abs{\nabla \tilde\varphi }^2 \tilde\varphi
\end{align*}
where $L_U$ is the linearized operator around $U$ given by 
\begin{equation}
\label{linearised_op_U}
L_U(\tilde\varphi) = \Delta \tilde\varphi  + 2(\nabla U : \nabla \tilde\varphi)U + |\nabla U|^2\tilde\varphi. 
\end{equation}
Here $A: B$ is the Frobenius inner product
\[
A:B = \sum_{i,j} A_{ij} B_{ij}.
\]
Substituting $\tilde\varphi = \pi_{U^\perp}(\varphi) + a U$, with $a = a\left(\pi_{U^\perp}(\varphi) \right)$ as defined in \eqref{eq:def-a}, leads to 
\begin{align*}
&
S(U +  \pi_{U^\perp}(\varphi)  + a U) 
\\
&= -U_t - \left( \pi_{U^\perp}(\varphi)  + a U\right)_t + L_U(\pi_{U^\perp}(\varphi) + a U)
\\
&\quad + |\nabla \left(\pi_{U^\perp}(\varphi) + a U\right)|^2 U + 2(\nabla U \cdot \nabla \left(\pi_{U^\perp}(\varphi) + a U\right))\left( \pi_{U^\perp}(\varphi) + a U\right) 
\\
&\quad+ \abs{\nabla \left( \pi_{U^\perp}(\varphi) + a U\right) }^2 \left( \pi_{U^\perp}(\varphi) + a U\right).
\end{align*}
Expanding term by term, one gets
\begin{align*}
\tilde\varphi_t 
&= \left( \pi_{U^\perp}(\varphi)\right)_t + a_t U + aU_t,
\\
L_U(\tilde\varphi) &= L_U(\pi_{U^\perp}(\varphi)) + (\Delta a) U + 2 \nabla a \cdot \nabla U + 2(\nabla U \cdot \nabla (a U)).
\\
2(\nabla U \cdot \nabla \tilde\varphi)\tilde\varphi
&= 2(\nabla U \cdot \nabla \left(\pi_{U^\perp}(\varphi) + a U\right))\pi_{U^\perp}(\varphi) + 2(\nabla U \cdot \nabla \left(\pi_{U^\perp}(\varphi) + a U\right)) a U
\\
\abs{\nabla \tilde\varphi }^2 \tilde\varphi &= \abs{\nabla \left(\pi_{U^\perp}(\varphi) + a U\right) }^2 \pi_{U^\perp}(\varphi) + \abs{\nabla \left(\pi_{U^\perp}(\varphi) + a U\right) }^2aU.
\end{align*}
Gathering all the lower order terms and the terms collinear to $U$, the equation \eqref{HMF} becomes 
\begin{align}
\nonumber
0=S(U +  \pi_{U^\perp}(\varphi) + a U) 
& = -U_t - (\pi_{U^\perp}(\varphi))_t + L_U(\pi_{U^\perp}(\varphi))
\\
\label{S(U)}
&\qquad + N_U(\pi_{U^\perp}(\varphi)) + c(\pi_{U^\perp}(\varphi))U 
\end{align}
where
\begin{align}
\nonumber
N_U(\pi_{U^\perp}(\varphi)) 
&= \Big\{2\nabla(U + \pi_{U^\perp}(\varphi) ) \cdot \nabla (aU) + 2\nabla U \cdot \nabla \pi_{U^\perp}(\varphi) + 
|\nabla \pi_{U^\perp}(\varphi)|^2 
\\
\nonumber
&\qquad + \abs{\nabla (aU)}^2\Big\} \pi_{U^\perp}(\varphi) - a U_t + 2\nabla a \cdot \nabla U
\end{align}

and $c$ is a scalar function defined as follows 
\begin{align*}
c(\pi_{U^\perp}(\varphi)) 
& =  \Delta a - a_t - 2\nabla U \cdot \nabla \pi_{U^\perp}(\varphi)
\\
&\qquad + (1 + a)\Big\{ |\nabla(U + \pi_{U^\perp}(\varphi) + a U)|^2 
- |\nabla U|^2\Big\} .
\end{align*}
The perturbation $\varphi$ will later be decomposed into the inner and
outer components mentioned in the introduction, and we will derive the
equations they satisfy. Before doing so, we recall
some basic properties of the
linearized operator $L_U$.

\subsection{Linearized operator around \texorpdfstring{$U$}{U}}
The linearized operator of the harmonic map equation around the modulated harmonic map $U$ is the elliptic operator given by \eqref{linearised_op_U}.
The least energy harmonic map $W$, defined in \eqref{eq:W} and its modulation $U$ can be written in polar coordinates
\begin{align}
\nonumber
W(y) = \begin{pmatrix}
e^{i\theta} \sin w(\rho)\\
\cos w(\rho)
\end{pmatrix} \quad \text{and} \quad U(x, t) = \begin{pmatrix}
e^{i(\theta + \omega)} \sin w(\rho)\\
\cos w(\rho)
\end{pmatrix}. 
\end{align}
where
$$w(\rho) = \pi - 2\arctan \rho, \quad y = \rho e^{i\theta}.$$
The function $w(\rho)$ satisfies
$$w''(\rho) + \frac{1}{\rho} w'(\rho) - \frac{\sin 2w}{2\rho^2} = 0$$
and
\begin{align}
\label{w(rho)}
w_\rho(\rho) = - \frac{2}{1 + \rho^2}, \quad 
\sin w(\rho) = \frac{2\rho}{1 + \rho^2}, \quad \text{and} \quad 
\cos w(\rho) = \frac{\rho^2 - 1}{\rho^2+1}. 
\end{align}
Let 
\begin{align}
\label{def_Ei}
E_1(y) = \begin{pmatrix}
e^{i\theta} \cos w(\rho)\\
-\sin w(\rho)
\end{pmatrix} \quad \text{and} \quad E_2(y) = \begin{pmatrix}
ie^{i\theta} \\
0
\end{pmatrix}
\end{align}
which form an orthonormal basis of the tangent plane to the sphere $S^2$ at $W(y)$.
Differentiating the map $U$ with respect to its parameters $(\lambda, \omega, \xi)$ yields linearly independent elements of the kernel of $L_U$:
\begin{align}
\label{partial_lambda_U}
\partial_\lambda U (x) &= -\frac{1}{\lambda}Q_\omega Z_{01}(y), 
& Z_{01}(y) &= \rho w_\rho E_1(y), 
\\
\label{partial_lambda_U_2}
\partial_\omega U (x) &= -Q_\omega Z_{02}(y), & Z_{02}(y) &= \rho w_\rho E_2(y), 
\\
\label{partial_xi1_U}
\partial_{\xi_1} U (x) &= -\frac{1}{\lambda}Q_\omega Z_{11}(y), 
& Z_{11}(y) &= w_\rho \big( \cos \theta E_1(y) + \sin \theta E_2(y)) , 
\\
\label{partial_xi2_U}
\partial_{\xi_2} U (x) &= \frac{1}{\lambda}Q_\omega Z_{12}(y), 
& Z_{12}(y) &= w_\rho \big( \sin \theta E_1(y) - \cos \theta E_2(y)) .
\end{align}
The functions $Z_{k j}$, $k = 0, 1$ and $ j= 1, 2$, are in the kernel of the operator $L_W$, that is 
\begin{align}
\label{def:LW}
L_W(Z_{k j}) = 0, \quad L_W(\phi) = \Delta \phi + \abs{\nabla W}^2 \phi + 2(\nabla W : \nabla \phi) W.
\end{align}
Moreover, the operator $L_W$ has two additional bounded elements in its kernel given by 
\begin{align}
\label{Z_{-11}}
Z_{-11}(y) &= \rho^2 w_\rho (\cos \theta E_1 - \sin \theta E_2), 
\\
\label{Z_{-12}}
Z_{-12}(y) &= \rho^2 w_\rho (\sin \theta E_1 + \cos \theta E_2). 
\end{align}
As proved in \cite{CLW}, harmonic maps from $\R^2$ to $S^2$ are
non-degenerate.
That is, for the 1-corotational harmonic map $W$, any bounded element in 
the kernel of the linearized operator $L_W$ is linear combination of 
these six elements $Z_{k j}$, $k = -1, 0, 1$,  $j = 1, 2$. 
For a function $h(y, t) \in \R^3$ defined in $\R^2 \times (t_0, +\infty)$ 
such that $h \cdot W=0$, we say that $h$ is in mode $k$, with $k \in \Z$, 
if $h$ can be written as 
\begin{align*}
h(y, t) = \Re{\tilde{h}(\abs{y}, t)e^{ik\theta}}E_1 
+ \Im{\tilde{h}(\abs{y}, t)e^{ik\theta}}E_2 
\end{align*}
for $\tilde{h}$ a complex valued function. With this terminology, it 
follows that the functions $Z_{-1j}, Z_{0j}$ and $ Z_{1j}$, $j = 1, 2$, 
are in modes $-1, 0$ and $1$ respectively. This notion is motivated by 
the approach used to analyze the inner equation and will be discussed further in Section~\ref{sec:Estimates_inner}.

As the inner component of the perturbation $\varphi$ will be expressed in the variable $y$, it is worth noticing that
\begin{align}\nonumber
L_U(\varphi) = \frac{1}{\lambda^2} Q_\omega L_W(\phi), \quad \text{for } \varphi(x) = Q_\omega\phi(y), \quad y = \frac{x - \xi}{\lambda},
\end{align}
where $L_W$ is defined in \eqref{def:LW}.

\subsection{Linearized operator on functions orthogonal to \texorpdfstring{$U$}{U}}

For a function $\Phi(x) \in \R^3$, the orthogonal projector on the plane orthogonal to $U$ is given by $\pi_{U^\perp}(\Phi)$ as defined in \eqref{eq:def-pi}. The linearized operator $L_U$ applied on the projection of $\Phi$ takes the following form
\begin{align}
\label{formula-L-Pi-U-perp}
L_U(\pi_{U^\perp}(\Phi)) = \pi_{U^\perp}(\Delta \Phi) + \tilde L_U(\Phi)
\end{align}
where
\begin{align}
\label{def_tildeL}
\tilde L_U(\Phi) = \abs{\nabla U}^2 \pi_{U^\perp}(\Phi) - [\nabla(U \cdot \Phi) \cdot \nabla] U, 
\end{align}
with
\begin{align*}
[\nabla(U \cdot \Phi) \cdot \nabla] U = 
\sum_{i=1}^2  \partial_{x_i} (U \cdot \Phi) \partial_{x_i}U.
\end{align*}
In polar coordinates $x = \xi + r e^{i\theta}$, $\Phi(x) = \Phi(r, \theta)$, the operator $\tilde{L}_U$ can be rewritten as 
\begin{align}
\label{form_Ltilde_polar}
\tilde L_U(\Phi) = -\frac{2}{\lambda} w_\rho \Big((\partial_r\Phi \cdot U) Q_\omega E_1 - \frac{1}{r}(\partial_\theta \Phi \cdot U) Q_\omega E_2\Big), \quad \rho = \abs{y} = \frac{r}{\lambda}.
\end{align}
A corollary of formula \eqref{form_Ltilde_polar} when $\Phi$ has the particular form
$
\Phi(x) = \begin{pmatrix}
\phi(r) e^{i\theta}\\
0
\end{pmatrix}, 
$
where $\phi(r) \in \C$, is 
\begin{align}
\label{eq:Ltilde2}
\tilde L_U(\Phi) = \frac{2}{\lambda} \rho w_\rho^2 \Big(\Re{e^{-i\omega} \partial_r\phi(r)} Q_\omega E_1 + \frac{1}{r}\Im{e^{-i\omega} \phi(r)}Q_\omega E_2\Big).
\end{align}

\subsection{Inner and outer solution}
We decompose the perturbation $\varphi$ in 
\eqref{eq:decomp-u} into inner and outer parts as follows
\begin{align*}
\varphi(x, t) = \varphi^{\text{in}}(y, t) + \varphi^{\text{out}}(x, t),
\quad y = \frac{x - \xi(t)}{\lambda(t)}
\end{align*}
where
\begin{align}
\label{phi_in}
\varphi^{\text{in}}(y, t) = \eta(y, t) Q_{\omega} \phi(y, t).
\end{align}
The function $\eta$ is defined by
\begin{align}
\label{def_eta}
\eta(y, t) = \eta_0\Bigl(\frac{\abs{y}}{R(t)}\Bigr), \quad y = \frac{x - \xi(t)}{\lambda(t)}, 
\end{align}
where $\eta_0$ is a smooth cut-off satisfying
\begin{align*}
\eta_0(s) = 1 \text{ for } s \leq 1,\quad
\eta_0(s) = 0 \text{ for } s \geq 2,
\end{align*}
The function $\phi$ is referred to as the inner solution and is expressed in terms of the variable $y$. We assume that $\phi(y,t)$ is pointwise orthogonal to $W(y)$, that is $\phi(\cdot, t) \cdot W = 0$ for all $(y,t)$, and therefore
\[
\pi_{U^\perp}(\varphi^{\text{in}}) = \varphi^{\text{in}}.
\]
The radius $R(t)$ has the form
\begin{align}
\nonumber
R(t)=t^\beta
\end{align}
for some $\beta \in(0,\frac12)$ to be specified later on.
Therefore, in the variable $x$, the term $\eta Q_\omega \phi$ is supported in a ball of radius $2\lambda R$ around $x = \xi$. This region will be referred to as the inner region.
The outer piece $\varphi^{\text{out}}$, expressed in the variable $x$ and supported on $\R^2 \times (t_0, \infty)$, is further decomposed into three components
\begin{align}
\label{phi_out}
\varphi^{\text{out}}(x, t) = \psi(x, t) + Z(x, t) + \Phi^0(x, t). 
\end{align}
The function $\psi$ is called the outer solution and accounts for the pieces of the system supported outside of the inner region.
The function $Z: \R^2 \times (t_0, +\infty) \to \R^3$ solves the homogeneous heat equation with initial datum $Z_0$ such that
\begin{align}
\label{eq:heat-eq-Z}
\left\{
\begin{aligned}
\partial_t Z &= \Delta Z && \text{in } \R^2 \times (t_0, +\infty),
\\
Z(\cdot, t_0) &= Z_0 && \text{in } \R^2. 
\end{aligned}
\right.
\end{align}
We have the following lemma (proof is in Section~\ref{App:est_Z0}).
\begin{lemma}
\label{Lem:est_Z}
Let $0<d<1$.
If $\|Z_0\|_{0}<+\infty$ the solution of \eqref{eq:heat-eq-Z} given by Duhamel's formula satisfies
\begin{align}
\nonumber
\abs{Z(x, t)} + \sqrt{t}\abs{\nabla Z(x, t)} 
\le C\frac{\|Z_0\|_{0}}{t^d\big(1 + \frac{\abs{x}}{\sqrt{t}}\big)^{2d}}
\end{align}
in $\R^2 \times (t_0, \infty)$, for some constant $C = C(d) > 0$.
\end{lemma}
The final term in the outer piece \eqref{phi_out} is the correction $\Phi^0 = \Phi^0[\lambda, \omega, \xi]$, introduced to offset the slow decay coming from $U_t$, and will be discussed in the next section.
Substituting the $\varphi = \varphi^{\text{in}} + \varphi^{\text{out}}$ in \eqref{S(U)} gives us 
\begin{align}
\nonumber
S(U +  \pi_{U^\perp}(\varphi) + a U) 
&= -U_t - \partial_t\varphi^{\text{in}} 
- \partial_t\pi_{U^\perp}(\varphi^{\text{out}})  
+ L_U(\varphi^{\text{in}})
\\
\nonumber
&\qquad  
+ \tilde L_U(\varphi^{\text{out}}) 
+ \pi_{U^\perp}(\Delta_x(\varphi^{\text{out}}))
\\
\nonumber
&\qquad +  N_U\left(\varphi^{\text{in}}
+\pi_{U^\perp}(\varphi^{\text{out}})
\right) + cU 
\end{align}
where we used the formula \eqref{formula-L-Pi-U-perp} for
$L_U(\pi_{U^\perp}(\varphi^{\text{out}}))$. Now substituting $\varphi^{\text{in}}$ and $\varphi^{\text{out}}$ as defined in \eqref{phi_in} and \eqref{phi_out} respectively in each piece of this equation leads to 
\begin{align*}
Q_{-\omega}\partial_t\varphi^{\text{in}} &= \eta_t \phi + \eta \big(\phi_t -  \lambda^{-1}\dot \lambda y \cdot \nabla_y \phi - \lambda^{-1} \dot \xi \cdot \nabla_y \phi + \dot \omega Q_{-\omega} \partial_\omega Q_{\omega} \phi\big)
\\
Q_{-\omega}L_U(\varphi^{\text{in}}) &= \lambda^{-2} \big(\eta L_W(\phi) + (\Delta_y \eta)\phi + 2 \nabla_y \eta \cdot \nabla_y \phi + 2\big(((\nabla_y \eta)\phi \cdot \nabla_y W) W\big)\big)
\\
\pi_{U^\perp}(\Delta(\varphi^{\text{out}})) &= \Delta \psi - \big( \Delta \psi  \cdot U \big)U  + \pi_{U^\perp}(\Delta \Phi^0 + \Delta Z) \\
    \partial_t\piU{\varphi^{\text{out}}} &= \partial_t \psi + \piU{\partial_t\Phi^0 + \partial_t Z} -  \big(\big(\psi + \Phi^0 + Z \big)\cdot U \big)U_t + c(x,t)U
\end{align*}
where $c(x,t)U$ collects terms collinear to $U$,
which can be discarded (see the discussion after Lemma~\ref{lem:colinear_U}).
Assuming $u\cdot U\not=0$,
a function $u$ of the form \eqref{eq:decomp-u}, is a solution to \eqref{HMF} if and only if
\begin{align}
\nonumber
0 &= \lambda^{-2}\eta Q_\omega\Big(-\lambda^2 \phi_t + L_W(\phi) + \lambda^2Q_{-\omega}\LtU{\psi + Z} \Big)
\\
\nonumber
&\qquad +\LtU{\Phi^0} - \piU{\Phi^0_t - \Delta \Phi^0 + U_t} 
\\
\nonumber
&\qquad - \psi_t + \Delta \psi + (1 - \eta) \LtU{\psi + Z} 
\\
\nonumber
& \qquad + \eta Q_\omega \big(\lambda^{-1}\dot \lambda y \cdot \nabla_y \phi + \lambda^{-1} \dot \xi \cdot \nabla_y \phi - \dot \omega J\phi\big) 
\\
\nonumber
&\qquad + Q_\omega \big((\Delta \eta)\phi + 2\lambda^{-1}\nabla \eta \cdot \nabla_y \phi - \eta_t \phi\big)
\\
\label{Full_before_splitting}
&\qquad + N_U\left(\eta \phi +  \piU{\psi +\Phi^0 + Z}\right) + ((\psi +\Phi^0 + Z) \cdot U)U_t  .
\end{align}

\subsection{The non-local correction \texorpdfstring{$\Phi^0$}{Phi0}}\label{sec:Phi0}
In this section we define precisely the operator $\Phi^0 = \Phi^0[\lambda, \omega, \xi]$ introduced in \eqref{phi_out}. 
This term is added to \eqref{phi_out} to compensate the slow space decay coming from the error of approximation $S(U)=-U_t$. 
We have
\begin{align*}
S(U)=-U_t = -\dot \lambda \partial_\lambda U - \dot \omega \partial_\omega U - \dot \xi \cdot \partial_\xi U = \mathcal{E}_0 + \mathcal{E}_1 
\end{align*}
where, using \eqref{partial_lambda_U}--\eqref{partial_xi2_U}, 
\begin{align}
\nonumber
\mathcal{E}_0(x, t) &= \rho w_\rho(\rho) Q_\omega \Big(\frac{\dot \lambda}{\lambda}  E_1(y)  + \dot \omega E_2(y) \Big)
\\
\nonumber
\mathcal{E}_1 (x, t) &= \frac{\dot \xi_1}{\lambda} w_\rho(\rho) Q_\omega\big(\cos \theta~ E_1(y)  + \sin \theta~ E_2(y) \big)
\\
\label{eq:E1}
&\quad -\frac{\dot \xi_2}{\lambda} w_\rho(\rho) Q_\omega\big(\sin \theta~ E_1(y)  - \cos \theta~ E_2(y) \big).
\end{align}
Here $y = \frac{x - \xi}{\lambda} = \rho e^{i\theta}$ and $w_\rho$ is defined in \eqref{w(rho)}. We observe that $\mathcal{E}_0$ decays slower than $\mathcal{E}_1$. More precisely, for large values of $r = \abs{x - \xi}$, we find
\begin{align*}
\mathcal{E}_0(x, t) =  -\frac{2r}{\lambda^2 + r^2} \Big(\dot \lambda Q_\omega E_1(y)  + \lambda \dot \omega Q_\omega E_2(y) \Big) 
\approx \tilde{\mathcal{E}}_0(x, t) 
\end{align*}
where
\begin{align*}
\tilde{\mathcal{E}}_0(x, t) 
=-\frac{2r}{\lambda^2 + r^2}
\begin{pmatrix}
(\dot \lambda + i \lambda \dot \omega) e^{i(\theta + \omega)}\\ 0
\end{pmatrix}.
\end{align*}
We define the parameter $p$ by
\begin{align}\nonumber
p(t) = \lambda(t) e^{i\omega(t)},
\end{align}
so that 
\begin{align}\nonumber
\tilde{\mathcal{E}}_0(x, t) 
= -\frac{2r}{\lambda^2 + r^2} \begin{pmatrix}
\dot p(t) e^{i\theta}\\ 0
\end{pmatrix}.
\end{align}
We define $\Phi^0(x, t)$
\begin{align}
\label{eq:def-Phi0}
\Phi^0(x, t) = \begin{pmatrix}
\tilde{\varphi}^0(x, t)\\ 
0
\end{pmatrix}
\end{align}
where
\begin{align}
\nonumber
\tilde{\varphi}^0(x, t) &= re^{i\theta} \varphi^0(z(r,t), t), \quad z(r,t) = \sqrt{r^2 + \lambda^2}, \quad r = |x - \xi|,
\end{align}
and
\begin{align}
\label{eq:varphi0}
\varphi^0(z, t) = -\int_{t_0}^t \dot p(s) k(z, t - s) ds, \quad k(z, t) = 2\frac{1 - e^{-\frac{z^2}{4t}}}{z^2}.
\end{align}
Hence for $\Phi^0$ \eqref{eq:def-Phi0} we get the expression
\begin{align}
\nonumber
\Phi^0[p,\xi](x, t) = -\int_{t_0}^t \dot p(s) k(z(r,t), t - s) ds \begin{pmatrix}
r e^{i\theta}\\ 0
\end{pmatrix}, 
\quad x=\xi + re^{i\theta}.
\end{align}

Note that $\varphi^0$ satisfies
\begin{align*}
\varphi^0_t = \varphi^0_{zz} + \frac{3}{z}\varphi^0_z - \frac{2\dot p}{z^2},
\quad \varphi^0(z,t_0)=0.
\end{align*}
and so 
\begin{align}\nonumber
\tilde{\varphi}^0_t \approx \Delta \tilde{\varphi}^0 - \frac{2re^{i\theta} }{\lambda^2 + r^2}\dot p .
\end{align}
Therefore  $\Phi^0[p, \xi]$ 
\[
\partial_t \Phi^0 \approx \Delta_x \Phi^0  
+\tilde{\mathcal E}_0,
\]
with remainders that we compute next.
We have
\begin{align}
\nonumber
- \Phi^0_t + \Delta \Phi^0 + \tilde{\mathcal{E}}_0 = \tilde{\mathcal{R}}_0 +\tilde{\mathcal{R}}_1, \quad 
\tilde{\mathcal{R}}_0 = \begin{pmatrix}
{\mathcal{R}}_0\\ 0
\end{pmatrix}, \quad \tilde{\mathcal{R}}_1 = \begin{pmatrix}
{\mathcal{R}}_1\\ 0
\end{pmatrix},
\end{align}
where 
\begin{align}
\nonumber
\mathcal{R}_0
&= r e^{i \theta}\frac{\lambda^2}{z^4} \int_{t_0}^{t} \dot p(s) 
(z^2k_{zz}-zk_z)(z, t - s) \,d s,
\\
\nonumber
\mathcal{R}_1
&=\lambda \dot \lambda   \frac{r}{z^2} e^{i\theta}  \int_{t_0}^t \dot p(s) z k_z(z, t - s)\,d s
+\dot \xi\int_{t_0}^{t} \dot p(s) k(z, t - s)ds 
\\
\label{def_R1}
&\quad + \frac{r}{z^2} e^{i\theta}  \Re{r e^{-i \theta }\dot \xi} \int_{t_0}^t \dot p(s) z k_z(z, t - s)\,ds .
\end{align}
In this splitting $\mathcal{R}_0$ contains the main terms and $\mathcal{R}_1$ is smaller. 
In the previous formulas $\xi(t)$ is interpreted as a 
complex number.
Let us now rewrite the term in the second line of \eqref{Full_before_splitting}:
$\tilde L_U(\Phi^0) + \pi_{U^\perp}(-\partial_t\Phi^0 + \Delta \Phi^0 - U_t)$. 

\begin{lemma}
We have
\begin{align}
\label{eq:error-Phi0}
\tilde L_U(\Phi^0) - \pi_{U^\perp}(\partial_t\Phi^0 - \Delta \Phi^0 + U_t)
&= \mathcal{K}_0[p, \xi] 
+ \mathcal E_1[p,\xi]
+ \piU{\tilde{\mathcal{R}}_1}.
\end{align}
where
\begin{align}
\label{eq:def-K0}
\mathcal{K}_0[p, \xi] = \mathcal{K}_{01}[p, \xi] + \mathcal{K}_{02}[p, \xi]
\end{align}
\begin{align}
\nonumber
\mathcal{K}_{01}[p, \xi] 
=-\frac{2}{\lambda} \rho w_\rho^2 
\int_{t_0}^t 
&\left\{
\Re{e^{-i\omega(t)}\dot p(s)} Q_\omega E_1 
+ \Im{e^{-i\omega(t)}\dot p(s)} Q_\omega E_2\right\}
\\
\label{eq:K01}
&\cdot
k(z, t - s) \,d s .
\end{align}
\begin{align}
\label{eq:K02}
& \mathcal{K}_{02}[p, \xi] 
\\
\nonumber
& =\frac{1}{\lambda} \rho w_\rho^2\bigg\{ -\dot{\lambda} - 2\int_{t_0}^t \Re{\dot p(s) e^{-i \omega(t)}} r k_z(z, t - s) z_r  ds \bigg\}Q_\omega E_1 
\\
\nonumber
& - \frac{1}{4\lambda} \rho w_\rho^2 \cos w \bigg\{\int_{t_0}^t \Re{\dot p(s) e^{-i \omega(t)}}\bigl(zk_z(z, t - s) - z^2 k_{zz}(z, t - s)\bigr)ds\biggr\}Q_\omega E_1 \\
\nonumber
&- \frac{1}{4\lambda} \rho w_\rho^2 \bigg\{\int_{t_0}^t \Im{\dot p(s) e^{-i \omega(t)}}\bigl(zk_z(z, t - s) - z^2 k_{zz}(z, t - s)\bigr)ds\bigg\}Q_\omega E_2 .
\end{align}
\end{lemma}

\begin{proof}
We have
\begin{align}
\nonumber
&
\tilde L_U(\Phi^0) + \pi_{U^\perp}(-\partial_t\Phi^0 + \Delta \Phi^0 - U_t)
\\
\nonumber
&= \tilde L_U(\Phi^0)  
+ \pi_{U^\perp}(\tilde{\mathcal{R}}_0) 
+ \pi_{U^\perp}(\tilde{\mathcal{R}}_1) 
- \pi_{U^\perp}(\tilde{\mathcal{E}}_0) 
+ \mathcal{E}_0 
+ \mathcal{E}_1 .
\end{align}
But using \eqref{eq:Ltilde2} we get
\begin{align*}
\tilde L_U(\Phi^0)
&=
\mathcal{K}_{01}[p, \xi](x,t)
-\frac{2}{\lambda} \rho w_\rho^2 \frac{r^2}{z}
\int_{t_0}^t \Re{e^{-i\omega} \dot p(s)} k_z(z,t-s)\,d s
\, Q_\omega E_1 
\end{align*}
with $\mathcal{K}_{01}[p, \xi](x,t)$ defined in \eqref{eq:K01}.
We also compute
\begin{align*}
\pi_{U^\perp}(\tilde{\mathcal E}_0)
&=-\frac{2r}{\lambda^2+r^2}
\left(
\dot\lambda \cos(w) Q_\omega E_1
+ ( \lambda \dot \omega) Q_\omega E_2
\right)
\end{align*}
and therefore
\begin{align*}
\pi_{U^\perp}(\tilde{\mathcal E}_0)
-\mathcal{E}_0
=\frac{\dot\lambda}{\lambda}\frac{4\rho}{(\rho^2+1)^2}Q_\omega E_1
= \frac{\dot\lambda}{\lambda}\rho w_\rho^2 Q_\omega E_1.
\end{align*}
Similar computations show that 
\begin{align*}
\pi_{U^\perp}(\tilde{\mathcal{R}}_0) 
&= - \frac1{4\lambda} \rho w_\rho^2 \cos(w(\rho)) 
\int_{t_0}^{t} \Re{\dot p(s)e^{-i\omega(t)}} (zk_z - z^2k_{zz})(z, t - s) \,ds \, Q_\omega E_1
\\
& \quad - \frac1{4\lambda} \rho w_\rho^2  \int_{t_0}^{t} 
\Im{\dot p(s)e^{-i\omega(t)}} (zk_z - z^2k_{zz})(z, t - s) \,ds 
\, Q_\omega E_2.
\end{align*}
Then \eqref{eq:error-Phi0} follows.
\end{proof}

Note that $\mathcal E_1$ \eqref{eq:E1} can also be written as 
\begin{align}
\label{eq:K1}
\mathcal E_1[p,\xi] = \frac{1}{\lambda} w_\rho 
\left[
\Re{(\dot \xi_1 - i \dot \xi_2)e^{i\theta}} Q_\omega E_1 + \Im{(\dot \xi_1 - i \dot \xi_2)e^{i\theta}} Q_\omega E_2\right]. 
\end{align}

\subsection{Inner-outer gluing system}
We now return to \eqref{Full_before_splitting}, incorporating the considerations regarding the non-local correction $\Phi^0$ from the previous section, and derive the inner-outer gluing system that will be analyzed subsequently.
Replacing the decomposition \eqref{eq:decomp-u} with 
\[
\varphi(x, t) = \eta(x, t) Q_{\omega} \phi(y, t) + \psi(x, t) + Z(x, t) + \Phi^0[p, \xi](x, t), \quad y = \frac{x - \xi}{\lambda},
\]
into \eqref{HMF} yields a system of equations for  $\phi$, $\psi$, $p$, $\omega$.
For any $R(t)>0$ we use the notation
\begin{align}
\label{eq:def-DR}
\mathcal D_{R} &= \big\{(y, t) \in \R^2 \times (t_0, +\infty)~\big|~ y \in B_{R(t)}(0)\big\}.
\end{align}
The inner equation is 
\begin{align}
\label{Inner_eq}
\left\{
\begin{aligned}
\lambda^2\partial_t \phi &= L_W(\phi) + h[p, \xi, \psi, Z] && \text{in } \mathcal D_{2R},\\
\phi(y) \cdot W(y) &= 0 && \text{in } \mathcal D_{2R},\\
\phi( \cdot, t_0) &= 0 && \text{in } B_{2R(t_0)}(0),
\end{aligned}
\right.
\end{align}
where
\begin{align}
\nonumber
h[p, \xi, \psi, Z] &= \lambda^2Q_{-\omega}\left(\tilde{L}_{U}(\psi + Z)\right)\chi_{\D{2}} 
+\lambda^2Q_{-\omega}\mathcal{K}_0[p, \xi] 
\\
\label{def_rhs_inner}
&\quad + \lambda^2Q_{-\omega}\mathcal E_1[p, \xi]\chi_{\D{2}} .
\end{align}
The outer equation is then defined by 
\begin{align}
\label{Outer_eq}
\left\{
\begin{aligned}
\partial_t\psi &= \Delta \psi + g[\psi, \phi, p, \xi, Z] && \text{in } \R^2 \times (t_0, \infty),\\
\psi( \cdot, t_0) &= 0 && \text{in } \mathbb R^2,
\end{aligned}
\right.
\end{align}
where
\begin{align}
\nonumber
g[\psi, \phi, p, \xi, Z] =& ~(1 - \eta) \LtU{\psi + Z} + \left(\{\psi + Z\} \cdot U\right)U_t 
\\
\nonumber
&+ Q_\omega \big((\Delta \eta)\phi + 2\lambda^{-1}\nabla \eta \cdot \nabla_y \phi - \eta_t \phi\big)
\\
\nonumber
&+ \eta Q_\omega \big(\lambda^{-1}\dot \lambda y \cdot \nabla_y \phi + \lambda^{-1} \dot \xi \cdot \nabla_y \phi - \dot \omega J\phi\big) 
\\
\nonumber
&+ (1 - \eta) \left(\mathcal{K}_0[p,\xi] + \mathcal E_1[p,\xi]\right) + \piU{\tilde{\mathcal R}_1} + \left(\Phi^0 \cdot U\right)U_t
\\
\label{def_g}
&+ N_U\left(\eta  \phi + \pi_{U^\perp}(\psi +\Phi^0 + Z)\right),
\end{align}
and, for $\zeta = \piU \varphi$ and $a = a(\zeta)$ as defined in \eqref{eq:def-a}, $N_U$ is given by
\begin{align*}
N_U(\zeta) &= \left(2\nabla(aU) \cdot \nabla (U + \zeta) + 2 \nabla U \cdot \nabla \zeta + \abs{\nabla \zeta}^2 + \abs{\nabla(aU)}^2 \right) \zeta - aU_t
\\
&\quad + 2\nabla a \cdot \nabla U.
\end{align*}
As discussed earlier, the function $Z: \R^2 \times (t_0, +\infty) \to \R^3$  solves the homogeneous heat equation with initial data $Z_0$ \eqref{eq:heat-eq-Z}.

The initial conditions on the parameters $(\lambda, \omega, \xi)$ \eqref{eq:ini-param} and the inner and outer solutions $(\psi, \phi)$ yield
\begin{align*}
u(t_0) = W + \pi_{W^\perp} \big(Z_0\big) + a\left(\pi_{W^\perp} \big(Z_0\big)\right) W
= W(x) + Z_0(x) .
\end{align*}

\section{Orthogonality conditions} 
\label{Sect:orth_mode01}
Solutions of the inner equation \eqref{Inner_eq} exhibit a better decay 
when the right-hand side $h$ satisfies suitable orthogonality conditions
with respect to the bounded elements of the kernel of $L_W$, denoted by 
$Z_{k j}$, $k = -1, 0, 1$, $j = 1, 2$, and defined in 
\eqref{partial_lambda_U}--\eqref{partial_xi2_U} and 
\eqref{Z_{-11}}--\eqref{Z_{-12}}. In the present section, we concentrate on the 
orthogonality conditions in modes $0$ and $1$.

\subsection{Orthogonality in mode 0}
The parameter $p(t) = \lambda(t)e^{i\omega(t)}$ is used so that the
orthogonality condition in mode $0$ holds, namely with respect to the 
elements $Z_{0j}$, $j = 1, 2$, given in 
\eqref{partial_lambda_U}--\eqref{partial_lambda_U_2}. 
As will become apparent, the equation arising 
from mode $0$ is more delicate to analyze than its counterpart in 
mode $1$, primarily because of the presence of the non-local operator $
\Phi^0[p, \xi]$ introduced in Section \ref{sec:Phi0}. The resulting 
equation for $p$ takes the form of an integro-differential equation. To 
obtain a solution in the appropriate space, we do not solve this equation 
exactly. Instead, we introduce a remainder term, denoted by 
$\mathcal{R}_\alpha$ with $\alpha \in (\frac{1}{2}, 1)$, which will 
exhibit faster time decay. We now describe how to reach this equation and 
informally discuss our strategy to solve it. The statements regarding the 
properties of the solution $p$ are given in Section \ref{Sect:eq_for_p} 
and the proofs are given in Section \ref{sec:eq-p}.

For the inner equation \eqref{Inner_eq},
with the right-hand side $h[p, \xi, \psi, Z]$ 
given by \eqref{def_rhs_inner},
the orthogonality conditions in mode 0 are:
\begin{align}
\label{orth_mode0}
\int_{\R^2} h[p, \xi, \psi, Z](y, t) \cdot Z_{0j}(y)dy = 0, \quad j = 1, 2,  \quad \text{for all } t \ge t_0.
\end{align}
We aim to solve this equation for a parameter $p \in C^1([t_0, +\infty), \C)$ such that
\begin{align}
\label{assumption_p}
p(t_0) = 1, \quad |p(t)-1|<\tfrac12
\quad \text{and}\quad \lim_{t \to \infty} \dot p(t) = 0.
\end{align}
One of the terms arising in \eqref{orth_mode0}
is the operator
\begin{align}
\label{def_B0j}
\mathcal{B}_{0j}[p,\xi](t) = \frac{\lambda}{2\pi} \int_{\R^2} Q_{-\omega}\big(\mathcal{K}_0[p, \xi](y, t) + \mathcal E_1[p, \xi](y, t)\chi_{\D{2}} \big) \cdot Z_{0j}(y) dy, 
\end{align}
$j=1,2$, 
where $\mathcal{K}_0$ and $\mathcal E_1$ are defined in 
\eqref{eq:K01}--\eqref{eq:K1}, and are the remainder of the 
equation solved by $\Phi^0$ discussed in 
Section~\ref{sec:Phi0}. 
These integrals can be rewritten as 
\begin{align*}
\mathcal{B}_{01}[p](t) &= 2 \int_{t_0}^t \Re{\dot p(s) e^{-i\omega(t)}} \Gamma_1\left(\frac{\lambda(t)^2}{t - s}\right)\frac{ds}{t - s} 
+ 2\dot \lambda (t),\\
\mathcal{B}_{02}[p](t) &= 2 \int_{t_0}^t \Im{\dot p(s) e^{-i\omega(t)}} \Gamma_2\left(\frac{\lambda(t)^2}{t - s}\right)\frac{ds}{t - s},
\end{align*}
where the kernels $\Gamma_j$, $j = 1, 2$, are given by 
\begin{align*}
\Gamma_1(\tau) &= -\int_{0}^\infty \rho^3 w_\rho^3 \left[K(\zeta)
+ 4 \frac{\rho^2}{1 + \rho^2} \zeta K_\zeta(\zeta)  
- \cos w(\rho) \zeta^2 K_{\zeta\zeta}(\zeta)\right]_{\zeta = \tau( 1 + \rho^2)}d\rho
\\
\Gamma_2(\tau) &= -\int_{0}^\infty \rho^3 w_\rho^3 \left[K(\zeta) - \zeta^2 K_{\zeta\zeta}(\zeta)\right]_{\zeta = \tau( 1 + \rho^2)}d\rho
\\
K(\zeta) &= 2\frac{1 - e^{-\frac{\zeta}{4}}}{\zeta}.    
\end{align*}
Moreover, these kernels satisfy the expansions
\begin{align}
\label{est_Gamma_tau<1}
|\Gamma_j(\tau) - 1| &\le C\tau \left(1 + \abs{\log \tau}\right)  & \text{for } \tau < 1, 
\\
\nonumber
|\Gamma_j(\tau)| &\le \frac{C}{\tau}  & \text{for } \tau > 1, 
\end{align}
for $j = 1, 2$. 
Note that $\mathcal B_{0j}$ \eqref{def_B0j} does not depend on $\xi$.
Let us set
\begin{align}
\label{def:B0}
\mathcal{B}_{0}[p](t) = \frac{e^{i\omega(t)}}{2} \big(\mathcal{B}_{01}[p](t) + i \mathcal{B}_{02}[p](t)\big). 
\end{align}
Similarly, we let
\begin{align}
\label{def_a0}
a_0[p, \xi, \psi, Z](t) = \frac{e^{i\omega(t)}}{2} 
\big(a_{01}[p, \xi, \psi, Z](t) + i a_{02}[p, \xi, \psi, Z](t)\big) 
\end{align}
where
\begin{align}
\label{eq:def-a0j}
a_{0j}[p, \xi, \psi, Z](t) = -\frac{\lambda}{2\pi} \int_{B_{2R(t)}(0)} Q_{-\omega}\tilde{L}_{U}(\psi + Z) \cdot Z_{0j}(y) \, dy, \quad j = 1, 2. 
\end{align}
Consequently, the orthogonality condition in mode 0
 \eqref{orth_mode0} reduces to the equation 
\begin{align}
\label{eq:B01}
\mathcal{B}_0[p](t) = a_{0}[p, \xi, \psi, Z](t) \quad \text{for all } t \ge t_0,
\end{align}
in the unknown $p(t) = \lambda(t)e^{i\omega(t)}$.
Using formula \eqref{form_Ltilde_polar} for $\tilde{L}_U$,
one can approximate
\begin{align*}
a_{0}[p, \xi, \psi, Z](t) \approx \int_{B_{2R(t)}} \rho w_\rho^2 \big\{\nabla_{x} \psi(\xi(t) + \lambda(t) y, t) + \nabla_{x} Z(\xi(t) + \lambda(t) y, t)\big\}dy,
\end{align*}
where $w_\rho$ is given in \eqref{w(rho)}.
Therefore, the right-hand side $a_0(t)$ exhibits the same time decay as $\nabla \psi(\cdot, t) + \nabla Z(\cdot, t)$. 
We will set up things so that the main time decay is the one coming from $\nabla \psi(\cdot, t)$, and satisfies
\begin{align*}
\abs{\nabla \psi(x, t)} &\le \frac{C}{R(t)t^{\nu_1}} \log t, 
\end{align*}
for all $t \ge t_0$ where
\begin{align*}
0<\nu_1 < 1 , \quad \nu_1 + \beta > 1, \quad R(t)= t^\beta, ~\beta \in (\tfrac{1}{4}, \tfrac{1}{2}),
\end{align*}
and therefore it is expected that
\begin{align}
\label{eq:expected-a0}
|a_{0}[p, \xi, \psi, Z](t) | &\le \frac{C}{R(t)t^{\nu_1}} \log t.
\end{align}

\subsection{Decomposition of the nonlocal operator}
\label{sec:decomp-nonlocal}
Motivated by \eqref{eq:B01}, we consider
\begin{align}
\label{non_local_generic_a0}
\mathcal{B}_0[p](t) = a_{0}(t) \quad \text{for } t \ge t_0 
\end{align}
for any $a_0$ satisfying \eqref{eq:expected-a0},
where $\mathcal B_0$ is defined in \eqref{def:B0}. 
This equation is delicate and appears difficult to solve
exactly.
To overcome this difficulty, we modify the equation by introducing a 
remainder $\mathcal{R}_\alpha$, which exhibits faster time decay.

Let us fix $\frac12<\alpha<1$.
Write
\begin{align*}
\frac{e^{i\omega(t)}}{2} \mathcal{B}_{01}[p](t)
&= e^{i\omega(t)}
\int_{t_0}^t \Re{\dot p(s) e^{-i\omega(t)}} \Gamma_1\left(\frac{|p(t)|^2}{t - s}\right)
\frac{ds}{t - s} 
+ e^{i\omega(t)}\tfrac d{dt}|p(t)|
\\
&= e^{i\omega(t)}
\int_{t-t^\alpha}^{t} \Re{\dot p(s) e^{-i\omega(t)}} \Gamma_1\left(\frac{|p(t)|^2}{t - s}\right)
\frac{ds}{t - s} 
\\
&\quad
+
e^{i\omega(t)}
\int_{t_0}^{t-t^\alpha} \Re{\dot p(s) e^{-i\omega(t)}} \Gamma_1\left(\frac{|p(t)|^2}{t - s}\right)
\frac{ds}{t - s} 
+ \frac{p(t)}{|p(t)|} \tfrac d{dt}|p(t)|
\\
&= 
e^{i\omega(t)}\Re{\dot p(t) e^{-i\omega(t)}}
\int_{t-t^\alpha}^{t}  \Gamma_1\left(\frac{|p(t)|^2}{t - s}\right)
\frac{ds}{t - s} 
\\
&\quad 
-e^{i\omega(t)}
\int_{t-t^\alpha}^{t} 
\frac{
\Re{(\dot p(t)-\dot p(s)) e^{-i\omega(t)}}}{t-s} \Gamma_1\left(\frac{|p(t)|^2}{t - s}\right)
\,d s
\\
&\quad
+
e^{i\omega(t)}
\int_{t_0}^{t-t^\alpha} \Re{\dot p(s) e^{-i\omega(t)}} \Gamma_1\left(\frac{|p(t)|^2}{t - s}\right)
\frac{ds}{t - s} 
+ \frac{p(t)}{|p(t)|} \tfrac d{dt}|p(t)|.
\end{align*}
This decomposition makes sense if 
$t\geq 2t_0$ and $t_0$ is sufficiently large so that 
$t-t^\alpha\geq t_0$.
We do the same with $\mathcal{B}_{02}[p](t)$ and obtain
\begin{align}
\nonumber
\mathcal{B}_0[p](t)
&=
e^{i\omega(t)}\Re{\dot p(t) e^{-i\omega(t)}}
\int_{t-t^\alpha}^{t}  \Gamma_1\left(\frac{|p(t)|^2}{t - s}\right)
\frac{ds}{t - s} 
\\
\nonumber
&\quad
+ ie^{i\omega(t)}\Im{\dot p(t) e^{-i\omega(t)}}
\int_{t-t^\alpha}^{t}  \Gamma_2\left(\frac{|p(t)|^2}{t - s}\right)
\frac{ds}{t - s} 
\\
\nonumber
&\quad
+
e^{i\omega(t)}
\int_{t_0}^{t-t^\alpha} \Re{\dot p(s) e^{-i\omega(t)}} \Gamma_1\left(\frac{|p(t)|^2}{t - s}\right)
\frac{ds}{t - s} 
+ \tfrac{p(t)}{|p(t)|}  \tfrac d{dt}|p(t)|
\\
\nonumber
&\quad
+ie^{i\omega(t)}
\int_{t_0}^{t-t^\alpha} \Im{\dot p(s) e^{-i\omega(t)}} \Gamma_2\left(\frac{|p(t)|^2}{t - s}\right)
\frac{ds}{t - s}
\\
\label{eq:B0-1}
&\quad
-\mathcal R_{\alpha,0}[p](t)
\end{align}
where
\begin{align}
\nonumber
\mathcal R_{\alpha,0}[p](t)
&= e^{i\omega(t)}
\int_{t-t^\alpha}^{t} 
\frac{
\Re{(\dot p(t)-\dot p(s)) e^{-i\omega(t)}}}{t-s} \Gamma_1\left(\frac{|p(t)|^2}{t - s}\right)
\,d s
\\
\label{def:Ralpha0}
&\quad
+ie^{i\omega(t)}
\int_{t-t^\alpha}^{t} 
\frac{
\Im{(\dot p(t)-\dot p(s)) e^{-i\omega(t)}}}{t-s} \Gamma_2\left(\frac{|p(t)|^2}{t - s}\right)
\,d s.
\end{align}
For the first two terms of \eqref{eq:B0-1}, we compute
\begin{align*}
\int_{t-t^\alpha}^{t}\frac{1}{t-s}
\Gamma_j\Bigl(\frac{|p(t)|^2}{t-s}\Bigr)\,d s
&=
\int_0^{t^\alpha}\frac{1}{u}
\Gamma_j\Bigl(\frac{|p(t)|^2}{u}\Bigr)\,d u
\\
&=
\int_0^{t^\alpha/|p(t)|^2} \frac{1}{v}
\Gamma_j\Bigl(\frac{1}{v}\Bigr)\,d v
\\
&=\int_{|p(t)|^2/t^\alpha}^\infty \frac{\Gamma_j(w)}{w}\,dw
\end{align*}
For $0<u<1$ write 
\begin{align*}
\int_u^\infty \frac{\Gamma_j(w)}{w}\,d w
&=- \log(u) +g_j(u)+C_j
\end{align*}
where
\begin{align*}
g_j(u)&= -\int_0^u \frac 1w(\Gamma_j(w)-1)\,d w
\\
C_j&=\int_0^1 \frac 1w(\Gamma_j(w)-1)\,d w
+\int_1^\infty\frac{\Gamma_j(w)}{w}\,d w.
\end{align*}
Thanks to \eqref{est_Gamma_tau<1} we have
\begin{align}
\label{eq:bound-gj}
|g_j(u)|\leq C u |\log(u)|, \quad 0<u<\tfrac12.
\end{align}
Therefore
\begin{align*}
& e^{i\omega(t)}\Re{\dot p(t) e^{-i\omega(t)}}
\int_{t-t^\alpha}^{t}  \Gamma_1\left(\frac{|p(t)|^2}{t - s}\right)
\frac{ds}{t - s} 
\\
&\quad
+ ie^{i\omega(t)}\Im{\dot p(t) e^{-i\omega(t)}}
\int_{t-t^\alpha}^{t}  \Gamma_2\left(\frac{|p(t)|^2}{t - s}\right)
\frac{ds}{t - s} 
\\
&=
\dot p(t) \log\Bigl(\frac{t^\alpha}{|p(t)|^2}\Bigr)
-\mathcal R_{\alpha,1}[p](t)
+H_1(p(t))\dot p(t)
\end{align*}
where
\begin{align}
\nonumber
\mathcal R_{\alpha,1}[p](t)
&=-e^{i\omega(t)}\Re{\dot p(t) e^{-i\omega(t)}} g_1\Bigl(\frac{|p(t)|^2}{t^\alpha}\Bigr)
\\
\label{def:Ralpha1}
&\quad 
-i e^{i\omega(t)}\Im{\dot p(t) e^{-i\omega(t)}} g_2\Bigl(\frac{|p(t)|^2}{t^\alpha}\Bigr)
\\
\label{def:H1}
H_1(p)
&=\bar C I + \Delta C 
\left[
\begin{matrix}
\cos(2\omega)&\sin(2\omega)
\\
\sin(2\omega)&-\cos(2\omega)
\end{matrix}
\right], 
\end{align}
where $\omega$ is the argument of $p\in \C$. We work under the assumption \eqref{assumption_p} and we take the argument in $\omega\in(-\frac\pi2,\frac\pi2)$.
In \eqref{def:H1}
\[
\bar C = \frac{C_1+C_2}{2},\quad \Delta C = \frac{C_1-C_2}{2}
\]
which are real constants.
Here we interpret $p(t)\in\R^2$ and $H_1(p(t))$ is a real $2\times 2$ matrix.
The expression $H_1(p(t))\dot p(t)$ can also be written as
\[
H_1(p(t))\dot p(t)=\bar C \dot p(t) + \Delta C e^{2i\omega(t)}\overline{\dot p}(t).
\]
Next we rewrite
\begin{align*}
&e^{i\omega(t)}
\int_{t_0}^{t-t^\alpha} \Re{\dot p(s) e^{-i\omega(t)}} \Gamma_1\left(\frac{|p(t)|^2}{t - s}\right)
\frac{ds}{t - s} 
\\
&\quad
+
ie^{i\omega(t)}
\int_{t_0}^{t-t^\alpha} \Im{\dot p(s) e^{-i\omega(t)}} \Gamma_2\left(\frac{|p(t)|^2}{t - s}\right)
\frac{ds}{t - s}
\\
\quad 
&=\int_{t_0}^{t-t^\alpha} \frac{\dot p(s)}{t-s}\,d s
-\mathcal R_{\alpha,2}[p](t)
\end{align*}
where
\begin{align}
\nonumber
\mathcal R_{\alpha,2}[p](t)
&=-e^{i\omega(t)}
\int_{t_0}^{t-t^\alpha} \Re{\dot p(s) e^{-i\omega(t)}} \left(\Gamma_1\left(\frac{|p(t)|^2}{t - s}\right)-1\right)
\frac{ds}{t - s} 
\\
\label{def:Ralpha2}
&\quad
-
ie^{i\omega(t)}
\int_{t_0}^{t-t^\alpha} \Im{\dot p(s) e^{-i\omega(t)}} \left(\Gamma_2\left(\frac{|p(t)|^2}{t - s}\right)-1\right)
\frac{ds}{t - s}.
\end{align}
Therefore
\begin{align}
\label{eq:decompB0}
\mathcal{B}_0[p](t)
&=\dot p(t) \log\Bigl(\frac{t^\alpha}{|p(t)|^2}\Bigr)
+H_1(p(t))\dot p(t)
+ \tfrac{p(t)}{|p(t)|} \tfrac d{dt}|p(t)|
+\int_{t_0}^{t-t^\alpha} \frac{\dot p(s)}{t-s}\,d s
-\mathcal R_\alpha[p](t)
\end{align}
where
\begin{align}
\label{def:Ralpha}
\mathcal R_\alpha=\mathcal R_{\alpha,0}+\mathcal R_{\alpha,1}+\mathcal R_{\alpha,2}.
\end{align}
We recall that the previous computation is valid for $t\geq 2t_0$ and $t_0$ large so that $t-t^\alpha\geq t_0$.
To make sense of this decomposition for all $t\geq t_0$, we introduce 
a cutoff $\eta_{2t_0}(t)$ such that 
\begin{align}
\label{def:eta2t0}
\eta_{2t_0}(t) = 1 \text{ for } t \geq 3t_0,\quad
\eta_{2t_0}(t) = 0 \text{ for } t \leq 2t_0.
\end{align}
So, instead of \eqref{non_local_generic_a0}, we consider the equation
\begin{align}
\label{prob12}
\left\{
\begin{aligned}
&\dot p(t)\log(t^\alpha)
+H(p(t))\dot p(t)
+\eta_{2t_0}(t)\int_{t_0}^{t-t^\alpha}\frac{\dot p(s)}{t-s}
\,d s = a_0(t)
\quad t>t_0
\\
&
p(t_0)=1.
\end{aligned}
\right.
\end{align}
where $H(p)$ is a $2\times 2$ matrix, smooth for $|p-1|<\frac12$,
such that
\[
H(p)v=H_1(p)v-2\log(|p|)v+\frac{p}{|p|^2}\Re{\overline p \, v},
\quad |p-1|<\tfrac12
\] 
with $H_1$ defined in \eqref{def:H1}.
We let $h$ be defined on $\C=\R^2$ 
such that $h$ is a smooth, bounded, and Lipschitz continuous function such that $h(p)=H(p)$ for $|p-1|<\frac12$.
Equation \eqref{prob12} (with the modified nonlinearity $h$) can be solved for any 
continuous function $a_0$ on $[t_0,\infty)$ and the solution is unique and local, in 
the sense that if $a_0$, $a_1$ are continuous functions on $[t_0,\infty)$ and 
$a_0=a_1$ on $[t_0,T]$ then the solutions coincide on $[t_0,T]$. 
A proof of this is in Section~\ref{sec:nonlocal3}.
The solution corresponding to $a_0=0$ is $p\equiv 1$.

\subsection{Orthogonality in mode 1}
Similarly to the parameter $p$ discussed in Section~\ref{sec:decomp-nonlocal}, 
the parameter $\xi$ is determined by imposing the orthogonality of the right-hand 
side of the inner equation \eqref{Inner_eq} with the bounded elements of the kernel 
of $L_W$ in mode $1$, namely $Z_{1j}$, $j = 1, 2$, defined in 
\eqref{partial_xi1_U}--\eqref{partial_xi2_U}. In contrast with the mode $0$ case, the 
equation for $\xi$ arising from the orthogonality condition in mode $1$ can be solved 
exactly. More precisely, there exist $\xi \in C^1([t_0, +\infty), \C)$ solving 
\begin{align}
\label{orth_mode1}
\int_{\R^2} h[p, \xi, \psi, Z](y, t) \cdot Z_{1j}(y)dy = 0, \quad j = 1, 2, \quad \text{for all } t \ge t_0, 
\end{align}
satisfying 
\begin{align}
\nonumber
\xi(t_0) = 0, \quad |\xi(t)| \leq 1 , \quad \lim_{t \to \infty} \dot \xi(t) = 0 .
\end{align}
Just as in \eqref{def_B0j}, we define 
\begin{align*}
\mathcal{B}_{1j}[p,\xi](t) 
= \frac{\lambda}{2\pi} \int_{\R^2} 
Q_{-\omega}
\big(\mathcal{K}_0[p, \xi](y, t) + \mathcal E_1[p, \xi](y, t)\chi_{\D{2}} \big) \cdot Z_{1j}(y) dy, \quad j = 1, 2,
\end{align*}
and 
\begin{align*}
\mathcal{B}_{1} = \mathcal{B}_{11} + i \mathcal{B}_{12}.
\end{align*}
Given the form of $\mathcal{E}_1$ in \eqref{eq:K1} and $Z_{1j}$, $j = 1, 2$, 
and using that 
$\int_0^\infty \rho w_\rho^2d\rho = 2$,
we find that
\begin{align}
\nonumber
\mathcal{B}_{1}[\xi](t) = 2\kappa(t)\big(\dot \xi_1(t) + i \dot \xi_2(t)\big) .
\end{align}
where
\[
\kappa(t)=1-\frac{1}{1+4R(t)^2}.
\]

Additionally, we define $a_1$ as
\begin{align}
\nonumber
a_{1j}[p, \xi, \psi, Z](t) = -\frac{\lambda}{2\pi} \int_{B_{2R(t)}(0)} Q_{-\omega}\tilde{L}_{U}(\psi + Z) \cdot Z_{1j}(y) dy, \quad j = 1, 2,
\end{align}
and
\begin{align}
\label{eq:def-a1}
a_1[p, \xi, \psi, Z](t) = e^{i\omega(t)} \big(a_{11}[p, \xi, \psi, Z](t) + i a_{12}[p, \xi, \psi, Z](t)\big). 
\end{align}
Consequently, imposing the orthogonality in mode 1 as in \eqref{orth_mode1} reduces to
\begin{align}
\label{eq:xi}
\dot \xi(t) = \frac1{2\kappa(t)} a_{1}[p, \xi, \psi, Z](t) \quad \text{for all } t \ge t_0, 
\end{align}
where $\xi(t) = \xi_1(t) + i \xi_2(t)$. In contrast with the equation for $p$, equation \eqref{eq:xi} can be solved exactly. The function $a_1[p, \xi, \psi, Z]$ can be approximated by
\begin{align*}
a_{1}[p, \xi, \psi, Z](t) 
\approx 
\int_{B_{2R(t)}} w_\rho^2 \left[\nabla_{x} \psi(\xi(t) + \lambda(t) y, t) + \nabla_{x} Z(\xi(t) + \lambda(t) y, t)\right]\,d y
\end{align*}
and we shall check that 
\begin{align}
\nonumber
\abs{a_1(t)} &\le \frac{C}{t^{\upnu}} \log t .
\end{align}
for some $\upnu>1$.
We then solve \eqref{eq:xi} with $\xi(t_0)=0$.

\section{A further decomposition and the full system}
\subsection{The full system}
\label{sec:decompinner}
In order to estimate the inner component $\phi$ satisfying  \eqref{Inner_eq},
we split it as $\phi = \phi_1 + \phi_2 + \phi_3 $ and give a suitable estimate of each $\phi_i$, $i = 1, 2, 3$.
The first component $\phi_1$ solves an equation similar to \eqref{Inner_eq}, but with a right-hand-side that has been projected. More precisely,
we define the projector $c_{k j}[h]$ by
\begin{align}
\label{def:ckj}
c_{k j}[h](t) = 
\frac{1}{\norm{w_\rho Z_{k j}}_{L^2(\R^2)}^2} \int_{\R^2} h (y, t)\cdot Z_{k j}(y)\, dy
\end{align}
for $k=-1,0,1$, $j=1,2$, where $Z_{k j}$ are the six bounded 
elements in the kernel of $L_W$ defined in 
\eqref{partial_lambda_U}--\eqref{partial_xi2_U} and 
\eqref{Z_{-11}}--\eqref{Z_{-12}}. 
In the definition \eqref{def:ckj}, $h$ is assumed to be defined in $\R^2$. 
The function $h[p, \xi, \psi, Z]$ defined in \eqref{def_rhs_inner} is indeed defined in 
$\R^2$.
We consider then the equation 
\begin{align}	
\label{eq:phi1}
\left\{
\begin{aligned}
\lambda^2 \partial_t \phi_1 &= L_W(\phi_1) + h[p, \xi, \psi, Z] 
- \sum_{\substack{k = 1, 0,  -1\\ j= 1, 2}} 
c_{k j} (h[p, \xi, \psi, Z])w_\rho^2Z_{k j}& \text{in } \D{2}
\\
\phi_1 \cdot W &= 0 \qquad \text{in } \D{2}
\\
\phi_1(\cdot, t_0) &= 0 \qquad \text{in } B_{2R(t_0)}
\end{aligned}
\right.
\end{align}
so that its RHS is orthogonal to all elements in the kernel.
We put the extra terms of \eqref{eq:phi1} in the RHS of $\phi_2$ and $\phi_3$.
To achieve this we introduce the equation
\begin{align}
\label{eq:phi2}
\left\{
\begin{aligned}
\lambda^2 \partial_t \phi_2 &= L_W(\phi_2) 
+(1-\eta_{3t_0}(t)) \sum_{j= 1, 2} c_{0j} (h[p, \xi, \psi, Z])w_\rho^2Z_{0j}
\\
&\qquad -\eta_{3t_0}(t)\sum_{ j= 1, 2} \widetilde{R}_{\alpha j}[p]w_\rho^2Z_{0 j}
&&\text{in } \D{2}\\
\phi_2 \cdot W &= 0 && \text{in } \D{2}
\\
\phi_2(\cdot, t_0) &= 0 && \text{in } B_{2R(t_0)}
\end{aligned}
\right.
\end{align}
where
\hypertarget{tildeRalpha}{}
\begin{align}
\label{eq:def-tildeRalphaj}
\widetilde{R}_{\alpha 1} 
= -\frac{4\pi \lambda }{\norm{w_\rho Z_{0 1}}_{L^2(\R^2)}^2} 
\text{Re}\big(e^{-i\omega}\mathcal R_\alpha\big), \quad 
\widetilde{R}_{\alpha 2} 
= - \frac{4\pi\lambda }{\norm{w_\rho Z_{0 2}}_{L^2(\R^2)}^2}
\text{Im}\big(e^{-i\omega}\mathcal R_\alpha\big)
\end{align}
with $\mathcal R_\alpha$ defined in \eqref{def:Ralpha}, and where $\eta_{3t_0}$ is a smooth cutoff function such that
\begin{align*}
\eta_{3t_0}(t) = 1 \text{ for } t \geq 4t_0,\quad
\eta_{3t_0}(t) = 0 \text{ for } t \leq 3t_0,
\end{align*}
Finally we consider 
\begin{align}
\label{eq:phi3}
\left\{
\begin{aligned}
\lambda^2 \partial_t \phi_3 &= L_W(\phi_3) 
+ \sum_{ j= 1, 2} c_{-1 j} (h[p, \xi, \psi, Z])w_\rho^2Z_{-1 j}
&& \text{in } \D{2}
\\
\phi_3 \cdot W &= 0 && \text{in } \D{2}
\\
\phi_3(\cdot, t_0) &= 0 && \text{in } B_{2R(t_0)}
\end{aligned}
\right.
\end{align}
We keep the outer equation \eqref{Outer_eq} as it is, that is,
 \begin{align}
\label{Outer_eq2}
\left\{
\begin{aligned}
\partial_t\psi &= \Delta \psi + g[\psi, \phi, p, \xi, Z] && \text{in } \R^2 \times (t_0, \infty),\\
\psi( \cdot, t_0) &= 0 && \text{in } \mathbb R^2,
\end{aligned}
\right.
\end{align}
with $\phi=\phi_1+\phi_2+\phi_3$.
The idea is that adding \eqref{eq:phi1}--\eqref{eq:phi3} gives that
$\phi=\phi_1+\phi_2+\phi_3$ satisfies
\begin{align*}
\lambda^2 \partial_t \phi &= L_W(\phi) + h[p, \xi, \psi, Z] 
- \eta_{3t_0}  \sum_{j=1,2} \left( c_{0 j} (h[p, \xi, \psi, Z])
-\widetilde R_{\alpha j}[p]
\right)
w_\rho^2Z_{0 j}
\\
&\quad 
- \sum_{j=1,2}  c_{1 j} (h[p, \xi, \psi, Z])
w_\rho^2Z_{1 j}.
\end{align*}
We impose that $p$ satisfies
\begin{align}
\label{eq:p-main}
\left\{
\begin{aligned}
&\dot p(t)\log(t^\alpha)
+H(p(t))\dot p(t)
+\eta_{2t_0}(t)\int_{t_0}^{t-t^\alpha}\frac{\dot p(s)}{t-s}
\,d s = a_0(t)
\quad t>t_0
\\
&
p(t_0)=1.
\end{aligned}
\right.
\end{align}
 with $a_0=a_0[p, \xi, \psi, Z]$ defined in 
\eqref{def_a0}.
This implies 
\begin{align}
\label{eq:realB0}
\mathcal B_0[p](t)=a_0(t)-\mathcal R_\alpha[p](t),\quad t\geq 3 t_0
\end{align}
which in turn gives 
\begin{align}
\label{eq:coj-and_tildeRalphaj}
c_{0 j} (h[p, \xi, \psi, Z])(t)-\widetilde R_{\alpha j}[p](t) 
=0\quad \text{for } t\geq 3t_0.      
\end{align}
Summing the equations for $\phi_1,\phi_2,\phi_3$, the remaining mode-zero term in the inner equation is 
\[ 
\eta_{3t_0}(t) \sum_{j=1}^2 \bigl(c_{0j}(h)(t)-\widetilde R_{\alpha j}(t)\bigr)Z_{0j}. 
\] 
This term vanishes identically by \eqref{eq:coj-and_tildeRalphaj} and the choice of $\eta_{3t_0}$.

We also ask that $\xi$ satisfies \eqref{eq:xi}, that is,
\begin{align}
\label{eq:xi-main}
\dot \xi(t) = \frac1{2\kappa(t)} a_{1}[p, \xi, \psi, Z](t) \quad \text{for all } t \ge t_0, \quad \xi(t_0)=0,
\end{align}
so that $ c_{1 j} (h[p, \xi, \psi, Z])(t)=0$ for $t\geq t_0$.
In this way, $\phi$ satisfies the original inner equation \eqref{Inner_eq}.

In summary, we will analyze and get estimates for the solutions $\phi_1$, $\phi_2$, 
$\phi_3$, $\psi$, $p$ and $\xi$ of the coupled system 
\eqref{eq:phi1}--\eqref{eq:xi-main}.

\medskip

The main difficulty in the proof arises from the orthogonality 
conditions in mode $0$. Formally, these conditions require solving 
\[
c_{0j}\big(h[p,\xi,\psi,Z]\big)(t)=0, \qquad j=1,2, 
\]
which is equivalent to the equation 
\[
\mathcal{B}_0[p](t) = a_{0}[p, \xi, \psi, Z](t) \quad \text{for all } t \ge t_0,
\]
(or \eqref{eq:B01}).
As explained in Section~\ref{sec:decomp-nonlocal}, the operator $\mathcal B_0$ does not seem suitable for a direct resolution. 
To overcome this difficulty, we replace the  previous equation by the modified problem \eqref{eq:p-main}. Proposition~\ref{prop:p1} below yields a solution $p$ of \eqref{eq:p-main} satisfying appropriate estimates. 
Equation \eqref{eq:p-main} implies
\eqref{eq:realB0}, where the remainder $\mathcal R_\alpha[p]$ has better time decay. 
The role of the component $\phi_2$ is precisely to 
absorb this defect. More precisely, the equation for $\phi_2$ 
contains the forcing term 
\[
-\eta_{3t_0}(t) \sum_{j=1,2} 
\widetilde R_{\alpha j}[p]\, w_\rho^2 Z_{0j}, 
\]
which corresponds to the remainder appearing in the previous 
identity. As a consequence, for $t\geq 4t_0$, the contribution 
of the mode $0$ is completely cancelled and the orthogonality 
condition is effectively restored. A technical point is that the 
equation \eqref{eq:realB0} is only valid when $t\geq 3t_0$.
For this reason the correction through $\mathcal R_\alpha$ 
is only activated after a transition region. The cutoff $\eta_{3t_0}$ 
is introduced to separate these two regimes. For $t\leq 3t_0$, 
we do not use the decomposition of $\mathcal B_0$. Instead, 
the coefficients $c_{0j}(h[p,\xi,\psi,Z])$ are inserted directly into 
the equation for $\phi_2$ through the term 
\[
(1-\eta_{3t_0}(t)) \sum_{j=1,2} 
c_{0j}\big(h[p,\xi,\psi,Z]\big) w_\rho^2 Z_{0j}. 
\]
Since this contribution is supported in a bounded time interval, 
it can be treated by direct estimates and does not influence 
the asymptotic behaviour of the solution. 
Propositions~\ref{Prop:phi1}--\ref{prop:est-Ralpha} below provide the estimates for each component in order to close the estimates for the system.

\subsection{Estimates for the inner solution} 
\label{sec:Estimates_inner}
Below we state several results concerning operators that associate to a given right-hand side $h$ a function $\phi$ satisfying the linear equation
\begin{align}
\label{eq:phi-basic}
\left\{
\begin{aligned}
\lambda^2 \partial_t \phi
&= L_W(\phi) + h && \text{in } \D{2}
\\
\phi \cdot W &= 0 && \text{in } \D{2}
\\
\phi(\cdot, t_0) &= 0  && \text{in } B_{2R(t_0)} .
\end{aligned}
\right.
\end{align}
The proofs can be found in \cite[Section 7]{HMF_DdPW}.
In some cases, we consider this equation without imposing boundary conditions, and it may happen that the solution's boundary values are not zero.

For the right-hand sides $h$ of \eqref{eq:phi-basic} we use the norm
\[
\|h\|_{**,a, \nu} = \sup_{\R^2 \times (t_0, +\infty)} \,t^{\nu}(1 + |y|)^{a} |h(y, t)| ,
\]
for some parameters $\nu>0$ and $a\in(2,3)$.
The norm for $\phi_1$ is defined as 
\begin{align}
\nonumber
\|\phi_1\|_{\text{in}1} 
= \sup_{\D{2}}\,
t^{\nu_1} \bigg(
& \frac{R(t)^{\frac{5 - a_1}{2}}}{(1 + \abs{y})}\min \bigg\{1, \frac{R(t)^{\frac{5 - a_1}{2}}}{\abs{y}^2}\bigg\}\bigg)^{-1} \Big( \abs{\phi_1(y, t)}
\\
& 
\label{def_norm_phi1}
+ (1 + \abs{y})\abs{\nabla_y \phi_1(y, t)} 
\Big).
\end{align}
The following proposition gives a solution with estimate of 
\eqref{eq:phi-basic} when the right-hand side is orthogonal to 
the 6 bounded elements of the kernel of $L_W$ \eqref{partial_lambda_U_2}--\eqref{partial_xi2_U}, \eqref{Z_{-11}}, \eqref{Z_{-12}}.
\begin{proposition}
\label{Prop:phi1}
Let $a_1 \in (2, 3)$, $\nu_1 > 0$.
There exists a linear operator  $\mathcal{T}_{1}$
such that if $\norm{h_1}_{**,a_1, \nu_1} <+\infty$ and
$c_{k j}[h_1]=0$ for $k=-1,0,1$, $j=1,2$, then
$\phi_1 = \mathcal{T}_{1}[h_1]$, satisfies \eqref{eq:phi-basic} with $h=h_1$ and
\begin{align*}
\|\phi_1\|_{\text{\rm in}1}  \le C \norm{h_1}_{**,a_1, \nu_1}.
\end{align*}
\end{proposition}
In this proposition $h$ is assumed to be defined in $\R^2$. In particular
the projectors $c_{kj}[h]$ \eqref{def:ckj} are defined as integrals on
$\R^2$.

A function $h$, defined in $\R^2\times (t_0, +\infty)$, is said to be in mode $k$, with $k \in \Z$, if $h$ can be written as 
\begin{align*}
h(y, t) = \Re{\tilde{h}(\abs{y}, t)e^{ik\theta}}E_1 + \Im{\tilde{h}(\abs{y}, t)e^{ik\theta}}E_2 
\end{align*}
for a complex valued function $\tilde{h}$ , 
where the vectors $E_1$ and $ E_2$ 
are defined in \eqref{def_Ei}. Note that the index $k = -1, 0, 1$ relates to the mode of $Z_{k j}$,  $j = 1, 2$. Thus the right-hand side of \eqref{eq:phi2} is in mode 0 and the right-hand side of \eqref{eq:phi3} is in mode~$-1$.
We define the norm for $\phi_2$ as 
\begin{align}
\label{def_norm_phi2}
\norm{\phi_2}_{\text{\rm in}2} 
= \sup_{\D{2}} \, \frac{t^{\nu_2}(1 + \abs{y})}{R(t)^2}
\left( \abs{\phi_2(y, t)}  + (1 + \abs{y})\abs{\nabla_y \phi_2(y, t)}
\right) .
\end{align}
The next result gives an estimate when the right-hand side is in mode 0, without orthogonality conditions.
\begin{proposition}

\label{Prop:phi2}
Let $a_2 \in (1, 3)$, $\nu_2 > 0$.
There exists a linear operator  $\mathcal{T}_{2}$
such that if $h_2$ is in mode 0 and $\norm{h_2}_{**,a_2, \nu_2} <+\infty$, then
$\phi_2 = \mathcal{T}_{2}[h_2]$ is also in mode 0, satisfies \eqref{eq:phi-basic} with right-hand side $h_2$, and 
\begin{align*}
\norm{\phi_2}_{\text{\rm in}2} \le C \norm{h_2}_{**,a_2, \nu_2}.
\end{align*}
\end{proposition}
Finally, to estimate $\phi_3$ the following norm is useful.
\begin{align}
\label{def_norm_phi3}
\norm{\phi_3}_{\text{\rm in}3} 
= \sup_{\D{2}} \, t^{\nu_3}(\log R(t))^{-1}
\left( \abs{\phi_3(y, t)}  + (1 + \abs{y})\abs{\nabla_y \phi_3(y, t)}
\right).
\end{align}
\begin{proposition}
\label{Prop:phi3}
Let $a_3 \in (2, 3)$, $\nu_3 > 0$.
There exists a linear operator  $\mathcal{T}_{3}$
such that if $h_3$ is in mode $-1$ and $\norm{h_3}_{**,a_3, \nu_3} <+\infty$, 
then $\phi_3 = \mathcal{T}_{3}[h_3]$ is also in mode $-1$, satisfies \eqref{eq:phi-basic} with right-hand side $h_3$, and 
\begin{align*}
\norm{\phi_3}_{\text{\rm in}3} \le C \norm{h_3}_{**,a_3, \nu_3}.
\end{align*}
\end{proposition}

\subsection{Estimate for the outer solution}
\label{subSect:est_HE}
We state the main estimate for the solution of \eqref{Outer_eq2}. 
For this purpose we  consider the heat equation with a generic right-hand side 
$g: \R^2 \times [t_0, +\infty) \to \R^3$,
\begin{align}
\label{eq:linear-outer}
\left\{
\begin{aligned}
\partial_t \psi &= \Delta \psi + g && \text{in }\R^2 \times (t_0, +\infty),
\\
\psi(\cdot, t_0) &= 0 && \text{in } \R^2.
\end{aligned}
\right.
\end{align}
The right-hand side $g[\psi, \phi, p, \xi, Z]$ of  \eqref{Outer_eq}  (defined in \eqref{def_g})  
consists of terms exhibiting different behaviors. To account for them, we define suitable 
weights as follows
\begin{align}
\label{sup_g1}
\varsigma_1(x, t) &=\frac{\chi_{B_{4R(t)}(0)}}{R(t)^2t^{d_0}},
\\
\label{sup_g2}
\varsigma_2(x, t)&= \frac{\chi_{B_{4R(t)}(0)}}{R(t)^{1 + 2\gamma}t^{d_0}\big(1 + \abs{x}\big)^{a_2}},
\\
\label{sup_g4}
\varsigma_3(x, t) &= \frac{(1 - \chi_{B_{R(t)}(0)})}{t^{d_0+\theta}(1 + \abs{x})^2},
\end{align}
for all $(x, t) \in \R^2 \times (t_0, +\infty)$, where $\chi_A$ is
the characteristic function of the set $A$. The constants appearing in these bounds are chosen so that
\begin{align}
\label{eq:param-linear-outer}
\left\{
\begin{aligned}
& 0<d_0 < 1 ,\quad  d_0 + \beta > 1, \quad \gamma \in (0, \tfrac{1}{2}), \quad a_2>2,
\\
&\theta>0,\quad d_0+\theta<1
\end{aligned}
\right.
\end{align}
where we recall that $R(t) = t^\beta$, $\beta\in(0,\frac12)$.
Further restrictions on these constants are given in Section \ref{Sect:Choice_of_cst01}, but in the present section we work under the more general assumptions \eqref{eq:param-linear-outer}.
Using the weights \eqref{sup_g1}--\eqref{sup_g4}
we define the norm for $g$ as  
\begin{align}
\label{def_norm_g}
\norm{g}_{*} = \sup_{\R^2 \times (t_0, +\infty)} \abs{g(x, t)}\bigg( \sum_{i = 1}^3
\varsigma_i(x, t)\bigg)^{-1}
\end{align}
and we define the norm for the solution $\psi$ as
\begin{align}
\nonumber
\|\psi\|_{\text{\rm out}} = \sup_{\R^2 \times (t_0, +\infty)} &\abs{\psi(x, t)} t^{d_0}\Big(1 + \sfrac{\abs{x}}{\sqrt{t}}\Big)^2(\log t)^{-1} 
\\
\nonumber 
&+ \sup_{\R^2 \times (t_0, +\infty)} \abs{\nabla_x \psi(x, t)} t^{d_0}R(t)\Big(1 + \sfrac{\abs{x}}{\sqrt{t}}\Big)^2(\log t)^{-1} 
\\
\label{norm_psi}
&+  \sup \frac{\abs{\nabla_x \psi(x,t) - \nabla_x \psi (x', t')}}{\big(\abs{x - x'}^2 + \abs{t-t'}\big)^{\gamma}} t^{d_0} R(t)^{1 + 2\gamma} (\log t)^{-1}
\end{align}
where the last supremum is taken on
\begin{align*}
(x, t), (x', t') \in \D{4}, \quad
t' \in \big(\max\{t_0, t - (4R(t))^2\}, t\big).
\end{align*}

\begin{proposition}
\label{prop:Estimate_HE_gen}
Assume that $d_0$, $\gamma$, $a_2$ satisfy \eqref{eq:param-linear-outer}.
Let $g\in C^0(\R^2 \times (t_0, +\infty), \R^3)$ be such that $\norm{g}_* < + \infty $. Then 
the solution $\psi $ to \eqref{eq:linear-outer} given by Duhamel's formula satisfies
\begin{align}
\nonumber
\norm{\psi}_{\text{\rm out}} \le C \norm{g}_*.
\end{align}
\end{proposition}
The proof is given in Section~\ref{Sect:est_HE}.

\subsection{Estimates for \texorpdfstring{$p = \lambda e^{i\omega}$}{p}}\label{Sect:eq_for_p}
Following up from Sections \ref{sec:decomp-nonlocal} and \ref{sec:decompinner}, the equation for $p$ is \eqref{eq:p-main}  with $a_0[p, \xi, \psi, Z](t)$ as in \eqref{def_a0}, but here we treat $a_0$ as a given function. 
Given the expected estimate \eqref{eq:expected-a0} for $a_0$,
we define
\begin{align}
\label{def:norm-gamma-m}
\|f\|_{\gamma,m} &= \sup_{t>t_0} \frac{t^\gamma}{\log^{m}(t)}|f(t)|
\end{align}
where $\gamma,m\in \R$.
In particular if $a_0$ satisfies \eqref{eq:expected-a0} then
$\|a_0\|_{\nu,1}<\infty$ for some $\nu>1$. 
The first result for \eqref{prob12} is the following.
\begin{proposition}
\label{prop:p1}
Assume $\alpha\in(\frac58,1)$ and $\nu>1$.
Then for $t_0$ large enough,
if $a_0\in C([t_0,\infty))$ satisfies $\|a_0\|_{\nu,1}\leq t_0^{\nu-1}\log^{-4}(t_0)$,
then the solution $p$ of \eqref{prob12} satisfies
\begin{align}
\label{eq:estp1}
|p(t)-1|
&\leq C \|a_0\|_{\nu,1} \frac{1}{t_0^{\nu-1}} \frac{\log^2(t_0)}{\log(t)}
\\
\label{eq:estdotp1}
|\dot p(t)|
&\leq C \|a_0\|_{\nu,1} \frac{1}{t_0^{\nu-1}} \frac{\log^2(t_0)}{t\log^2(t)}.
\end{align}
\end{proposition}
The proof of this proposition is given in Section~\ref{sec:eq-p}.
The hypothesis $\|a_0\|_{\nu,1}\leq t_0^{\nu-1}\log^{-4}(t_0)$ is so that
the condition $t_0$ large is uniform with respect to $a_0$.

A corollary of Proposition~\ref{prop:p1} is that if 
$\|a_0\|_{\nu,1}\leq t_0^{\nu-1}\log^{-4}(t_0)$, and $t_0$ is large enough,
then the solution $p$ of \eqref{prob12} satisfies 
\begin{align}
\nonumber
\mathcal B_0[p](t)=a_0(t)-\mathcal R_\alpha[p](t),\quad t\geq 3 t_0.
\end{align}
This equation contains the remainder term $\mathcal R_\alpha$ defined in \eqref{def:Ralpha}, which appears in the right-hand side of \eqref{eq:phi2}.
To bound it we need H\"older estimates for $p$.
For this we need to assume some H\"older control of $a_0$. 
Given the formula \eqref{def_a0} for $a_0[p, \xi, \psi, Z](t)$, and that $\psi$, $Z$ solve 
heat-type equations, a natural assumption is
\[
\|a_0\|_{\nu,1} + [a_0]_{\nu+2\beta\mu,1,\mu}<\infty
\]
where $\|a_0\|_{\nu,1} $ is defined in \eqref{def:norm-gamma-m} and
\begin{align}
\nonumber
[a_0]_{\gamma,1,\mu} &= \sup_{t>t_0}\,\sup_{\max(t_0,\frac{3t}{4})<s<t} t^\gamma (\log t)^{-1}  \frac{|a_0(t) - a_0(s)|}{(t - s)^\mu},
\end{align}
for a constant $0<\mu<1$ and $\gamma>0$.

\begin{proposition}
\label{prop:p2}
Assume $\alpha\in(\frac58,1)$,  $\nu>1$, $\gamma>0$, $0<\mu<1$.
Then for $t_0$ large enough,
if $a_0\in C([t_0,\infty))$ satisfies $\|a_0\|_{\nu,1}\leq t_0^{\nu-1}\log^{-4}(t_0)$ and
$ [a_0]_{\gamma,1,\mu}<\infty$, 
then the solution $p$ of \eqref{prob12} satisfies
\begin{align}
\label{eq:holder-p}
|\dot p(t)-\dot p(s)|
\leq C \|a_0\|_{\nu,1} \frac{\log^2(t_0)}{t_0^{\nu-1}}\frac{t-s}{t^2 \log^2(t)}
+ C [a_0]_{\gamma,1,\mu} \frac{(t-s)^\mu}{t^\gamma}
\end{align}
for $t >t_0$ and $\max(t_0,\frac{3t}{4})<s<t$.
\end{proposition}
The proof is given in Section~\ref{sec:holder-p}.
The importance of a H\"older estimate for $p$ is that it gives time decay for the remainder $\mathcal R_\alpha$.

\begin{proposition}
\label{prop:est-Ralpha}
Suppose $p\in C^1([t_0,\infty))$ satisfies  \eqref{eq:estp1}, \eqref{eq:estdotp1},
and \eqref{eq:holder-p}.
Then
\begin{align*}
|\mathcal R_\alpha[p](t)|
\leq 
C \|a_0\|_{\nu,1} \frac{\log^2(t_0)}{t_0^{\nu-1}}\frac{1}{t^{2-\alpha} \log^2(t)}
+ C[a_0]_{\gamma,1,\mu}\frac{1}{t^{\gamma-\mu\alpha}}.
\end{align*}
\end{proposition}

The proof is postponed to Section~\ref{sec:some-estimates}.

\subsection{Summary of parameters}
\label{sect:summary-param}
In the construction, we have
introduced several parameters and some restrictions among them.

\begin{itemize}
\item 
$d>\frac12$ is the parameter in the initial condition \eqref{eq:normZ0}.

\item
We have $\beta\in (0,\frac12)$ such that $R(t)=t^\beta$.
$R(t)$ is used to define the domain for the inner equations $\D{2}$ \eqref{eq:def-DR}.

\item
$\alpha\in(\frac12,1)$ the parameter introduced in the
decomposition of the operator $B_0$ 
associated to the orthogonality in mode 0
\eqref{eq:decompB0},
and appears in the definition of the remainder operator
$\mathcal R_\alpha$ \eqref{def:Ralpha}.

\item
The norm for $\phi_1$ is $\|\ \|_{\text{in}1}$, defined in \eqref{def_norm_phi1} with $\nu_1>0$ and $a_1\in(2,3)$ to be determined.

\item
The norm for $\phi_2$ is $\|\ \|_{\text{in}2}$, defined in \eqref{def_norm_phi2}, with
$\nu_2 = \nu_1 + \beta$.

\item
The norm for $\phi_3$ is $\|\ \|_{\text{in}3}$, defined in \eqref{def_norm_phi3}, with
$\nu_3 = 1$.

\item 
The norm of the outer solution is defined 
in \eqref{norm_psi}. We will use it with $d_0=\nu_1$, and the parameter $\beta\in(0,\frac12)$.
The norm for the right-hand side is defined in \eqref{def_norm_g}, involving the parameters $d_0=\nu_1$, $a_2>2$, $\gamma>0$ small, $\theta>0$ small.
\end{itemize}

\section{Proof of Theorem~\ref{thm:main}}

\subsection{A first look at the estimates}
\label{Sect:outline_proof}
To motivate the general argument, we first discuss  a simplified version of the system 
\eqref{eq:phi1}--\eqref{eq:xi-main}. 
Specifically, we 
consider
\begin{align}
\label{eq_psi_simp} 
\partial_t \psi &= \Delta \psi + \Delta \eta (\phi_1 + \phi_2)
\\
\nonumber
\lambda^2 \partial_t \phi_1 &= L_W(\phi_1) +  \lambda^2\big(\tilde L_U(\psi + Z) + \mathcal K_0[p] \big) 
\\
\label{eq_phi1_simp}
& \quad - \sum_{j=1,2}
c_{0j}\left[   \lambda^2\big(\tilde L_U(\psi + Z) + \mathcal K_0[p] \big)  \right]
w_\rho^2Z_{0 j}
\\
\nonumber
\lambda^2 \partial_t \phi_2 &= L_W(\phi_2) 
+(1-\eta_{3t_0}(t)) \sum_{j= 1, 2} 
c_{0j}\left[   \lambda^2\big(\tilde L_U(\psi + Z) + \mathcal K_0[p] \big)  \right]
w_\rho^2Z_{0j}
\\
\label{eq_phi2_simp}
&\quad -\eta_{3t_0}(t)\sum_{ j= 1, 2}
\widetilde{R}_{\alpha j} [p]w_\rho^2Z_{0 j}
\\
\nonumber
\dot p(t)\log(t^\alpha)
&+\eta_{2t_0}(t)\int_{t_0}^{t-t^\alpha}\frac{\dot p(s)}{t-s}\,ds
= \int_{B_R} \LtU{\psi + Z} \cdot Z_{01}\, dy
\\
\label{eq_c0j_simp}
&\qquad\qquad\qquad\qquad\qquad\qquad
+i\int_{B_R} \LtU{\psi + Z} \cdot Z_{02}\, dy.
\end{align}
The first equation \eqref{eq_psi_simp} is set in $\R^2$ and represents the outer 
equation  \eqref{Outer_eq2}, where we kept $\Delta \eta(\phi_1 + \phi_2)$ as the 
coupling term. 
Equations \eqref{eq_phi1_simp}--\eqref{eq_phi2_simp} are set in $\D{2}$ and 
represent simplified versions of the 
inner equations \eqref{eq:phi1}, \eqref{eq:phi2}. The inner solutions $\phi_1$ and 
$\phi_2$ are assumed to be in mode $0$, with $\phi_1$ satisfying the corresponding 
orthogonality condition, while $\phi_2$ absorbs remainder terms. Finally, equation 
\eqref{eq_c0j_simp} is a simplified version of the integro-differential equation 
\eqref{eq:p-main}  derived from the orthogonality condition in mode 0. 
The operator $\widetilde{R}_{\alpha j}$ appearing in
\eqref{eq_phi2_simp} is the one defined in \eqref{eq:def-tildeRalphaj}.
As in the original problem, we prescribe zero initial data for the solution 
$(\psi, \phi_1, \phi_2)$.  In this simplified model, we have omitted the parameter $\xi$, the 
third inner equation \eqref{eq:phi3}, and several terms in the
right-hand side of  the outer equation \eqref{Outer_eq2}.
Finally $Z$ solves the heat equation \eqref{eq:heat-eq-Z} with initial condition $Z_0$, and we assumed $\|Z_0\|_0\ll1$, where $\|\ \|_0$ is defined in \eqref{eq:normZ0}.

We recall that the outer solution $\psi$ and $Z$ are 
expressed in the variable $x \in \R^2$ and the inner solutions $\phi_1, \phi_2$ are 
expressed in the variable $y = \frac{x}{\lambda}$. 
As a consequence of the estimate for $p$ from Proposition~\ref{prop:p1}, it is reasonable to assume that the parameter $\lambda$ satisfies 
\[
\frac{1}{2} \le \lambda(t) \le \frac{3}{2},
\]
making these variables comparable. 

We will explain how the estimates for the system
\eqref{eq_psi_simp}--\eqref{eq_c0j_simp} close.
The approach for the full
system \eqref{eq:phi1}--\eqref{eq:xi-main}
is the same and explained in Section~\ref{sec:full-sys}.

The main external source is the function $Z$. 
By Lemma~\ref{Lem:est_Z}, in the self-similar region $|x| \le \sqrt{t}$, the function $Z$ decays as
\begin{align*}
|Z(x, t)| \lesssim \frac{1}{t^d} \|Z_0\|_0,
\quad |\nabla Z(x, t)| \lesssim \frac{1}{t^{d + \frac{1}{2}}} \|Z_0\|_0.
\end{align*}
Therefore, using the expression for $\tilde{L}_U$ in \eqref{def_tildeL}, we have 
\begin{align}
\label{eq:decay-tLUZ-before}
|\tilde{L}_U(Z)(y, t)|
&\le \frac{1}{1 + \rho^2} \frac{1}{t^{d
+ \frac{1}{2}}} \|Z_0\|_{0}, \quad \rho = |y|,
\\
\label{main_decay_rhs_phi1}
&\le \frac{1}{(1 + \rho)^{a_1}} \frac{1}{t^{\nu_1}} \|Z_0\|_{0},
\quad \rho  \le 2R(t) ,
\end{align}
where $R(t) = t^\beta$,
$\nu_1$ satisfies
\begin{align}
\label{restriction_nu_a1b}
0<\nu_1 <  d+\tfrac12-\beta(a_1-2), 
\end{align}
and $a_1 \in (2, 3)$.

This term serves as the driving decay for the right-hand side of
\eqref{eq_phi1_simp}, which in turn determines the decay for $\phi_1, \psi$
and $\phi_2$. The second inequality in \eqref{main_decay_rhs_phi1} ensures
that the right-hand side exhibits a space decay of the form
$\rho^{-{a_1}}$, $a_1 \in (2, 3)$, required by Proposition \ref{Prop:phi1}.
We design a scheme such that all terms on the right-hand side of
\eqref{eq:phi1} are bounded by
\begin{align}
\label{eq:model-decay-rhs-eqphi1}
\frac{C}{t^{\nu_1}(1 + \rho)^{a_1}}.
\end{align}
Under the decay \eqref{eq:model-decay-rhs-eqphi1} for the right-hand side of \eqref{eq_phi1_simp}, subject to the orthogonality condition in mode 0, Proposition \ref{Prop:phi1} implies 
a bound for $\|\phi_1 \|_{\text{in}1}$
defined in \eqref{def_norm_phi1} with $\nu_1$ and $a_1\in (2,3)$.
Later we shall see that it is natural to consider $\|\ \|_{\text{in}2}$ defined in \eqref{def_norm_phi2} with $\nu_2=\nu_1+\beta$, and $a_2\in (2,3)$.

\medskip
\noindent
\textbf{Claim.}
We claim that, provided $t_0$ is sufficiently large and $\|Z_0\|_0$ 
is sufficiently small, there exists a constant $\delta\in(0,1)$ such that
\[
\|\phi_1\|_{\text{in}1} \leq1, \quad \|\phi_2\|_{\text{in}2} \leq1,
\]
then
\begin{align}
\label{eq:contraction-phi12}
\|\phi_1\|_{\text{in}1} \leq1-\delta, \quad \|\phi_2\|_{\text{in}2} \leq1-\delta.
\end{align}

Such an apriori estimate implies that the solution is global in time and that 
$\|\phi_1\|_{\text{in}1} \leq1-\delta$, $\|\phi_2\|_{\text{in}2} \leq1-\delta$
by a standard argument. Indeed, the claim also holds with the norms defined as a supremum up to a time $T>t_0$, which we shall call $\|\phi_1\|_{\text{in}1,T}$. If there is a first time $T>t_0$ such that  
$\|\phi_1\|_{\text{in}1,T}=1$ and 
$\|\phi_2\|_{\text{in}2,T}\leq1$ and then the claim implies that $\|\phi_1\|_{\text{in}1,T} \leq1-\delta$, $\|\phi_2\|_{\text{in}2,T} \leq1-\delta$, which is a contradiction.

We start by estimating $\psi$, and for this we estimate the right-hand side of 
\eqref{eq_psi_simp}. The first term is 
\begin{align*}
|(\Delta \eta)\phi_1(y, t)| 
&\leq C \frac{\chi_{\rho \sim R(t)}}{R(t)^2} \frac{1}{R(t)^{a_1 - 2}t^{\nu_1}} \|\phi\|_{\text{in}1}
\leq C \frac{\chi_{\rho \sim R(t)}}{R(t)^2 t^{\nu_1}} \|\phi\|_{\text{in}1}.
\end{align*}
with the cut-off function $\eta$ in \eqref{def_eta}, and
where $\chi_{\rho \sim R(t)}$ denotes the indicator function of the region $\rho \sim R(t)$. 
We also have 
\begin{align*}
|(\Delta \eta)\phi_2(y, t)| 
&\leq C \frac{\chi_{\rho \sim R(t)}}{R(t)^2} R(t) \frac1{t^{\nu_1+\beta}} \|\phi\|_{\text{in}2}
\leq C \frac{\chi_{\rho \sim R(t)}}{R(t)^2 t^{\nu_1}} \|\phi\|_{\text{in}2}.
\end{align*}
Then by Proposition~\ref{prop:Estimate_HE_gen} the outer solution $\psi$ satisfies the following estimates in the self-similar region $|x| \le \sqrt{t}$ (see Section \ref{subSect:est_HE}),
\begin{align}
\nonumber
|\psi(x, t)| &\lesssim \frac{\log t}{t^{\nu_1}},
\\
\label{model_decay_psi}
|\nabla \psi(x, t)| &\lesssim \frac{\log t}{R(t) t^{\nu_1}}, 
\\
\label{eq:holder-psi}
|\nabla \psi(x, t) - \nabla \psi(x, s)| 
&\lesssim \frac{\log t}{R(t) t^{\nu_1}} \frac{|t - s|^\mu}{R^{2\mu}},
\end{align}
where $\mu\in (0, \frac{1}{2})$ and $s$ is comparable to $t$. This implies 
\begin{align}
\nonumber
|\tilde{L}_U(\psi)(y, t)| &\le \frac{1}{1 + \rho^2} \frac{\log t}{R(t)t^{\nu_1}} 
\\
\nonumber
&\le o_{t_0}(1)  \frac{1}{(1 + \rho)^{a_1}} \frac{1}{t^{\nu_1}} , 
\end{align}
for $\rho = |y| \le 2R(t)$.
Then by Proposition \ref{Prop:phi1} 
\begin{align}
\label{eq:phi1o1}
\|\phi_1\|_{\text{in}1}\leq o_{t_0}(1).
\end{align}

The right-hand side of equation \eqref{eq_phi2_simp} 
can be decomposed as 
\[
h_2=h_{2,a}+h_{2,b}
\]
where
\begin{align}
\label{eq:def-h2a}
h_{2,a}
&=
-\eta_{3t_0}(t)\sum_{ j= 1, 2} \widetilde{R}_{\alpha j}[p]w_\rho^2Z_{0 j}
\\
\label{eq:def-h2b}
h_{2,b}
&=
(1-\eta_{3t_0}(t)) \sum_{j= 1, 2} 
c_{0j}\left[   \lambda^2\big(\tilde L_U(\psi + Z) + \mathcal K_0[p] \big)  \right]
w_\rho^2Z_{0j} .
\end{align}
We first estimate the remainder $\widetilde{R}_{\alpha j}$
\eqref{eq:def-tildeRalphaj}.
Let $a_0$ denote the right-hand side of \eqref{eq_c0j_simp},
\begin{align}
\nonumber
a_0(t) = \int_{B_{2R}} \LtU{\psi + Z} \cdot Z_{0j} \,dy.
\end{align}
Note that from \eqref{eq:decay-tLUZ-before}
\begin{align*}
\left|
\int_{B_{2R}} \tilde L_U(Z) \cdot Z_{0j} \,dy
\right|
&\le C \frac{1}{t^{d + \frac{1}{2}}} \|Z_0\|_0.
\end{align*}
Since $\psi$ satisfies \eqref{model_decay_psi} and \eqref{eq:holder-psi}
we get 
\[
|a_0(t)| \leq C\frac{\log t}{ R(t) t^{\nu_1}}(1+\|Z_0\|_{0}), \quad
|a_0(t) - a_0(s)| \leq C \frac{\log t}{R(t) t^{\nu_1}} \frac{|t - s|^\mu}{R(t)^{2\mu}}(1+\|Z_0\|_{0}).
\]
By Propositions \ref{prop:p1}, the solution $p$ of \eqref{eq_c0j_simp} satisfies
\[
|p(t)-1|
\leq C \frac{1}{t_0^{\nu-1}} \frac{\log^2(t_0)}{\log(t)}(1+\|Z_0\|_{0})
,\quad
|\dot p(t)|
\leq C \frac{1}{t_0^{\nu-1}} \frac{\log^2(t_0)}{t\log^2(t)}(1+\|Z_0\|_{0}),
\]
where $\nu=\nu_1+\beta$. Thus, we require
\begin{align}
\label{eq:d1_plus_beta}
\nu_1+\beta>1.
\end{align}
Using Proposition \ref{prop:p2} we then find that
\begin{align}
\nonumber
|\dot p(t)-\dot p(s)|
\leq \left[ C \frac{\log^2(t_0)}{t_0^{\nu-1}}\frac{t-s}{t^2 \log^2(t)}
+ C \frac{(t-s)^\mu}{t^\gamma}\right](1+\|Z_0\|_{0})
\end{align}
where $\gamma=\nu_1+\beta+2\beta\mu$.
By Proposition~\ref{prop:est-Ralpha},
the remainder $\mathcal R_{\alpha j}$ is then bounded by
\begin{align*}
|\mathcal R_\alpha[p](t)|
& \leq 
\left[ C  \frac{\log^2(t_0)}{t_0^{\nu-1}}\frac{1}{t^{2-\alpha} \log^2(t)}
+ C \frac{1}{t^{\gamma-\mu\alpha}}\right](1+\|Z_0\|_{0})
\\
&=
\left[
C  \frac{\log^2(t_0)}{t_0^{\nu-1}}\frac{1}{t^{2-\alpha} \log^2(t)}
+ C \frac{1}{t^{\nu_1+\beta+\mu(2\beta-\alpha)}}\right]
(1+\|Z_0\|_{0}).
\end{align*}
By requiring
\begin{align}
\label{eq:param1}
2-\alpha \geq \nu_1+\beta, \quad 
\alpha < 2\beta,
\end{align}
we ensure that
\begin{align}
\label{eq:est-hatRalpha}
|\widetilde R_{\alpha j} [p] | \leq o_{t_0}(1) \frac{1}{R(t)t^{\nu_1}}(1+\|Z_0\|_{0}).
\end{align}
Then \eqref{eq:est-hatRalpha} implies that 
\begin{align*}
|h_{2,a}|\leq o_{t_0}(1) \frac{1}{R(t)t^{\nu_1}} \frac{1}{(1 + \rho)^{4}} (1+\|Z_0\|_{0}).
\end{align*}
If we ignore $h_{2,b}$, Proposition \ref{Prop:phi2} would imply that
\begin{align*}
|\phi_2(y, t)| \leq o_{t_0}(1) \frac{1}{R(t) t^{\nu_1}} \frac{R(t)^2}{(1 + \rho)}
(1+\|Z_0\|_{0}),
\end{align*}
that is,
\begin{align}
\label{eq:phi2o1}
\|\phi_2\|_{\text{in}2}\leq o_{t_0}(1).
\end{align}
Therefore, if we ignore the term $h_{2,b}$,
from \eqref{eq:phi1o1} and \eqref{eq:phi2o1} we achieve the objective
\eqref{eq:contraction-phi12}.

\subsection{Proof of Theorem~\ref{thm:main}}
\label{sec:full-sys}

We handle here the full system 
\eqref{eq:phi1}--\eqref{eq:xi-main}
in a similar way as the simplified system in Section~\ref{Sect:outline_proof}.

To simplify the presentation we first carry out the proof under the additional assumption $d>\frac34$. The modifications needed to treat the full range $d>\frac12$ are described in Remark~\ref{rem:1/2} and Remark~\ref{rem:adapt-prop:p1v2}.

We consider $\|\ \|_{\text{in}1}$ defined in \eqref{def_norm_phi1} with $\nu_1$ and $a_1\in (2,3)$ and  $\|\ \|_{\text{in}2}$ defined in \eqref{def_norm_phi2} with $\nu_2=\nu_1+\beta$, and $a_2\in (2,3)$.
We will use also the following norms for $\dot p$ and $\dot \xi$:
\[
\|\dot p\|_{\dot p}
=\sup_{t>t_0} t \log^2(t) |\dot p(t)|
\]
\[
\|\dot \xi\|_{\dot\xi}
=\sup_{t>t_0} t^{\nu_1+\beta} |\dot \xi(t)|.
\]
For the proof we need to consider all norms so far restricted to $[t_0,T]$ (where $T>t_0$) as follows:
\begin{align}
\nonumber
\|\phi_1\|_{\text{in}1,T} 
= \sup_{\mathcal{D}_{2R,T}}
t^{\nu_1} \bigg(
& \frac{R(t)^{\frac{5 - a_1}{2}}}{(1 + \abs{y})}\min \bigg\{1, \frac{R(t)^{\frac{5 - a_1}{2}}}{\abs{y}^2}\bigg\}\bigg)^{-1} \Big( \abs{\phi_1(y, t)}
\\
& 
\nonumber
+ (1 + \abs{y})\abs{\nabla_y \phi_1(y, t)} 
\Big) 
\end{align}
where 
\begin{align}
\nonumber
\mathcal D_{R,T} &= \big\{(y, t) \in \R^2 \times (t_0, T)~\big|~ y \in B_{R(t)}(0)\big\}.
\end{align}
We use similar modifications for the other norms:
$\|\ \|_{\text{in}2,T}$, $\|\ \|_{\text{in}3,T}$, and so on.

The system \eqref{eq:phi1}--\eqref{eq:xi-main} consists of non-linear (non-local) evolution equations. 
It can be solved for a short time after $t_0$, where to construct the solution $\phi_1$ to \eqref{eq:phi1} we use the operator given in Proposition~\ref{Prop:phi1}, and similarly for $\phi_2$, $\phi_3$.
The solution to the system \eqref{eq:phi1}--\eqref{eq:xi-main} may blow up in finite time.
But if a solution exists on $[t_0,T]$ and $\phi_1,\phi_2,\phi_3,\psi$ together with their spatial derivatives are bounded  on $[t_0,T]$, then it can be continued to $[t_0,T+h]$ for some small $h>0$. The proof of Theorem~\ref{thm:main} consists in showing that if $\|Z_0\|_0$ is small, then the solution to the system \eqref{eq:phi1}--\eqref{eq:xi-main} does not blow up, and we have some estimates for the components.
The proof of Theorem~\ref{thm:main} is based on the following two propositions.

Recall that we have split the right-hand side of equation \eqref{eq_phi2_simp} $h_2=h_{2,a}+h_{2,b}$,
defined in \eqref{eq:def-h2a} and \eqref{eq:def-h2b} respectively.

\begin{proposition}
\label{prop:main1}
Once $t_0$ is fixed (as large as needed), there is $\delta=\delta(t_0)>0$ and $C=C(t_0)$ such that if $\|Z_0\|_0<\delta$ then
\begin{align}
\label{eq:claim1}
\left\{
\begin{aligned}
&\|\psi\|_{L^\infty(\R^2\times(t_0,4t_0))}
+
\|\nabla_x \psi\|_{L^\infty(\R^2\times(t_0,4t_0))}
\leq C\|Z_0\|_0
&
\\
& 
t_0^{-\nu_1}\|\phi_1\|_{\mathrm{in}1,4t_0} \leq C\|Z_0\|_0,\quad 
t_0^{-\nu_2}\|\phi_2\|_{\mathrm{in}2,4t_0} \leq C\|Z_0\|_0,
\\
&
t_0^{-\nu_3}\|\phi_3\|_{\mathrm{in}3,4t_0} \leq C\|Z_0\|_0
\\
&
t_0^{-1}\log^{-2}(t_0)\|\dot p\|_{\dot p,4t_0} \leq C\|Z_0\|_0,\quad
t_0^{-\nu_1-\beta}\|\dot \xi\|_{\dot\xi,4t_0}\leq C\|Z_0\|_0
\\
&|h_{2,b}(y,t)|\leq C 
\frac{1}{  (1 + \rho)^{2+\sigma}} 
\|Z_0\|_{0} , \quad t_0<t<4t_0
\end{aligned}
\right.
\end{align}
for any $0<\sigma<1$.
\end{proposition}

For the proof we will need the following estimates.
We recall the function $\Phi^0[p, \xi]$ defined in 
\eqref{eq:def-Phi0}
\begin{align*}
\Phi^0[p, \xi](x, t)
= \begin{pmatrix}
\varphi^0(z,t) re^{i\theta}\\ 0
\end{pmatrix},
\quad  \varphi^0(z, t) &= - \int_{t_0}^{t} \dot p(s) k(z, t - s)ds
\end{align*}
where $z(r,t)=\sqrt{r^2+\lambda(t)^2}$.
\begin{lemma}
\label{lem:phi0}
Assume that $p$ satisfies
\[
|\dot p(t)|\leq \frac{1}{t\log^2(t)},\quad t>t_0.
\]
Let $\varphi^0$ be defined in \eqref{eq:varphi0}.
Then
\begin{align}
\nonumber
|\varphi^0(z, t)|
\leq
\frac{C}{\log(t_0)}\frac1{t(1+\frac{z}{\sqrt t})^2}
\|\dot p\|_{\dot p}
\end{align}
for $z>\frac12$, $t>t_0$.
\end{lemma}
The proof is in Section~\ref{sec:some-estimates}.
The function $\varphi^0$ has a $\log$ singularity when $z\to 0$
but we evaluate $\varphi^0(z,t)$ only for $z\geq\frac12$.
We also need an estimate for $\mathcal{K}_0$ defined by \eqref{eq:def-K0},
proved in Section~\ref{sec:some-estimates}.
\begin{lemma}[Estimate of $\mathcal K_0$]
\label{lem:estK0}
Assume that $p$ satisfies
\[
|\dot p(t)|\leq \frac{1}{t\log^2(t)},\quad t>t_0.
\]
Then
\begin{align*}
|\mathcal{K}_0[p, \xi](y,t)|
&\le C \frac{1}{(1 + |y|)^3}
\frac{1}{\log(t_0)}\frac1{t(1+\frac{|y|}{\sqrt t})^2}\|\dot p\|_{\dot p}.
\end{align*}
\end{lemma}

We also need the following estimate for the heat equation,
whose proof is direct using Duhamel's formula.
\begin{lemma}
\label{lem:est-grad-heateq}
Let $g\in L^\infty(\R^2\times (t_1,t_2))$ and let $v$ be the solution to
\begin{align}
\nonumber
\left\{
\begin{aligned}
\partial_t v&=\Delta v+g &&\text{in }\R^2\times (t_1,t_2)
\\
v(t_1) & = v_1&& \text{in }\R^2
\end{aligned}
\right.
\end{align}
given by Duhamel's formula, where $v_1 \in C^1(\R^2)\cap W^{1,\infty}(\R^2)$.
Then for any $t_1<t<t_2$
\begin{align*}
\|v(t)\|_{L^\infty(\R^2)}
& \leq \|v_1\|_{L^\infty(\R^2)}
+ C (t-t_1)  \sup_{s\in [t_1,t]} \|g(s)\|_{L^\infty(\R^2)}
\\
\|\nabla v(t)\|_{L^\infty(\R^2)}
&\leq \|\nabla v_1\|_{L^\infty(\R^2)}
+ C \sqrt{t-t_1} \sup_{s\in [t_1,t]} \|g(s)\|_{L^\infty(\R^2)}.
\end{align*}
\end{lemma}

\begin{proof}[Proof of Proposition~\ref{prop:main1}]
Let $t_0\leq t_1\leq 4t_0$.
We work with functions $\psi$ which are bounded, continuous, and whose gradient $\nabla_x \psi$ is also bounded and continuous.
We use the notation
\begin{align*}
\|\psi(t)\|_{C^1_x}
& =\|\psi(\cdot,t)\|_{L^\infty(\R^2)}
+
\|\nabla_x\psi(\cdot,t)\|_{L^\infty(\R^2)}
\\
\|\psi\|_{C^1_x,(t_0,t_1)}
&=\|\psi\|_{L^\infty(\R^2\times(t_0,t_1))}
+
\|\nabla \psi\|_{L^\infty(\R^2\times(t_0,t_1))}.
\end{align*}
Recall that $\psi$ satisfies the outer equation \eqref{Outer_eq2}
with right-hand side $g$ \eqref{def_g}.

Let $t_0\leq t\leq 4t_0$.
Assume that
\begin{align}
\nonumber
&\|\psi\|_{C^1_x,(t_0,t)}\leq 1,
\\
\nonumber
& 
t_0^{-\nu_1}\|\phi_1\|_{\mathrm{in}1,t} \leq 1,\quad
t_0^{-\nu_2}\|\phi_2\|_{\mathrm{in}2,t} \leq 1,\quad
t_0^{-\nu_3}\|\phi_3\|_{\mathrm{in}3,t} \leq 1,
\\
\label{eq:assumption-gA}
&
t_0^{-1}\log^{-2}(t_0)\|\dot p\|_{\dot p,t} \leq 1,\quad
t_0^{-\nu_1-\beta}\|\dot \xi\|_{\dot\xi,t}\leq 1.
\end{align}
The factors in front of the norms achieve,
for example in the case  of
$t_0^{-\nu_1}\|\phi_1\|_{\mathrm{in}1,T}$, that this norm
is equivalent to 
\begin{align}
\nonumber
\|\phi_1\|_{\text{in}1,T} 
= \sup_{\mathcal{D}_{2R,T}}
\bigg(
& \frac{R(t)^{\frac{5 - a_1}{2}}}{(1 + \abs{y})}\min \bigg\{1, \frac{R(t)^{\frac{5 - a_1}{2}}}{\abs{y}^2}\bigg\}\bigg)^{-1} \Big( \abs{\phi_1(y, t)}
\\
& 
+ (1 + \abs{y})\abs{\nabla_y \phi_1(y, t)} 
\nonumber
\\
\nonumber
&   + (1 + \abs{y})^{1 + 2\gamma_1} [\nabla_y \phi_1]_{\gamma_1; \mathcal P_\ell(y, t) \cap \D{2}} \Big) 
\end{align}
for $T\leq 4t_0$.

\medskip\noindent
\textbf{Step 1.}
We claim that if \eqref{eq:assumption-gA} holds, then
\begin{align}
\label{eq:claim1-about-g}
\|g[\psi, \phi, p, \xi, Z]\|_{L^\infty(\R^2\times(t_0,t))}
&\leq C_{t_0}
(\|\psi\|_{C^1_x,(t_0,t)} +\|Z_0\|_0).
\end{align}

From the definition of $g[\psi, \phi, p, \xi, Z]$ \eqref{def_g},
\begin{align}
\nonumber
\|g[\psi, \phi, p, \xi, Z]\|_{L^\infty(\R^2\times(t_0,t))}
&\leq C_{t_0}
(\|\psi\|_{C^1_x,(t_0,t)}
+t_0^{-\nu_1}\|\phi_1\|_{\mathrm{in}1,t}
+t_0^{-\nu_2}\|\phi_2\|_{\mathrm{in}2,t}
\\
\nonumber
& \quad 
+t_0^{-\nu_3}\|\phi_3\|_{\mathrm{in}3,t}
+t_0^{-1} \log^{-2}(t_0) \|\dot p\|_{\dot p,t}
\\
\label{eq:est-gA}
&\quad
+ t_0^{-\nu_1-\beta}\|\dot \xi\|_{\dot \xi,t}).
\end{align}

We use Proposition~\ref{prop:nonlocal0}, which gives a rough estimate on the solution $p$ of \eqref{eq:p-main}.
This estimate
depends on a parameter $m$, which we call here $m_1$, such that
$m_1>-\frac1\alpha$, $m_1>\frac{2-3\alpha}{\alpha(2\alpha-1)}$.
For $t_0\leq t\leq 4t_0$ and $\rho = |y| \le 2R(t)$,
we have
\begin{align}
\nonumber
|\tilde{L}_U(\psi)(y, t)|
&\le \frac{C}{1 + \rho^2} \sup_{s\in(t_0,t)} \|\nabla \psi(s)\|_{L^\infty}
\\
\nonumber
&\le \frac{t_0}{\log^{m_1}(t_0)} \frac{\log^{m_1}(t)}{ t (1 + \rho)^2}  
\|\psi\|_{C^1_x,(t_0,t)}
\end{align}
and therefore
\[
\left|
\int_{B_{2R}} \tilde L_U(\psi) \cdot Z_{0j} \,dy
\right|
\leq C \frac{t_0}{\log^{m_1}(t_0)} \frac{\log^{m_1}(t)}{ t }
\|\psi\|_{C^1_x,(t_0,t)}.
\]
We also have from \eqref{eq:decay-tLUZ-before}
\begin{align*}
\left|
\int_{B_{2R}} \tilde L_U(Z) \cdot Z_{0j} \,dy
\right|
&\le C \frac{1}{t^{d + \frac{1}{2}}} \|Z_0\|_{0},
\end{align*}
so
\begin{align}
\label{eq:est-tildeL-psi-Z-upto4t0}
\left|
\int_{B_{2R}} \tilde L_U(\psi+Z) \cdot Z_{0j} \,dy
\right|
\leq C \frac{t_0}{\log^{m_1}(t_0)} \frac{\log^{m_1}(t)}{ t }
(\|\psi\|_{C^1_x,(t_0,t)}+\|Z_0\|_{0} ).
\end{align}
This gives for $a_0[p,\xi,\psi,Z]$ \eqref{def_a0} the estimate
\[
|a_0(t)|\leq C \frac{t_0}{\log^{m_1}(t_0)} \frac{\log^{m_1}(t)}{ t }
(\|\psi\|_{C^1_x,(t_0,t)}+\|Z_0\|_{0} ).
\]
Then Proposition~\ref{prop:nonlocal0} implies that
\begin{align*}
|p(t)-1|
&\leq C \frac{t_0}{\log^{m_1}(t_0)} \frac{\log^{{m_1}+1}(t_0)}{ \log(t) }
(\|\psi\|_{C^1_x,(t_0,t)}+\|Z_0\|_{0} ).
\\
&\leq C \frac{t_0 \log(t_0) }{ \log(t) }
(\|\psi\|_{C^1_x,(t_0,t)}+\|Z_0\|_{0} ).
\end{align*}
\begin{align}
\nonumber
|\dot p|
&\leq  \frac{t_0}{\log^{m_1}(t_0)}
\frac{\log^{{m_1}+1}(t_0)}{t\log^2(t)}
(\|\psi\|_{C^1_x,(t_0,t)}+\|Z_0\|_{0} )
\\
\label{eq:est-dotp-bad}
&\leq
\frac{C}{\log(t_0)}
(\|\psi\|_{C^1_x,(t_0,t)}+\|Z_0\|_{0} ).
\end{align}
Note that in particular \eqref{eq:est-dotp-bad} implies
\[
t_0^{-1}\log^{-2}(t_0)
\|\dot p\|_{\dot p}
\leq 
\frac{C}{\log^3(t_0)}
(\|\psi\|_{C^1_x,(t_0,t)}+\|Z_0\|_{0} ).
\]

From a computation using \eqref{eq:est-dotp-bad}, 
similar to the proof of Proposition~\ref{prop:est-Ralpha} in Section~\ref{sec:some-estimates},
we get
\begin{align*}
|\widetilde R_{\alpha j}[p](t)|
\leq  C (\|\psi\|_{C^1_x,(t_0,t)}+\|Z_0\|_{0} )
\end{align*}
and hence  the term $h_{2,a}$ \eqref{eq:def-h2a} 
in the right-hand side of \eqref{eq:phi2} has the estimate
\begin{align*}
|h_{2,a}(t)|
= \left|\eta_{3t_0}(t)\sum_{ j= 1, 2} \widetilde{R}_{\alpha j}[p]w_\rho^2Z_{0 j}\right|
&\leq C \frac{1}{1+\rho^4} 
(\|\psi\|_{C^1_x,(t_0,t)}+\|Z_0\|_{0} ).
\end{align*}
We also have
\begin{align}
\nonumber
|\tilde{L}_U(\psi)(y, t)| 
&\le \frac{C}{1 + \rho^2} 
(\|\psi\|_{C^1_x,(t_0,t)}+\|Z_0\|_{0} )
\\
\nonumber
&\le \frac{R^\sigma}{(1 + \rho)^{2+\sigma}} 
(\|\psi\|_{C^1_x,(t_0,t)}+\|Z_0\|_{0} )
\\
\label{eq:est-tildeL-psi2}
&\le C  R(t_0)^\sigma \frac{1}{  (1 + \rho)^{2+\sigma}} 
(\|\psi\|_{C^1_x,(t_0,t)}+\|Z_0\|_{0} ).
\end{align}
Therefore the term $h_{2,b}$ \eqref{eq:def-h2b} has the estimate 
\begin{align*}
|h_{2,b}|\leq C  R(t_0)^\sigma 
\frac{1}{  (1 + \rho)^{2+\sigma}} 
(\|\psi\|_{C^1_x,(t_0,t)}+\|Z_0\|_{0} ).
\end{align*}
Applying Proposition~\ref{Prop:phi2} to \eqref{eq:phi2}, we find that, 
\begin{align*}
|\phi_2(y,t)| 
&\leq C\frac{R^2}{1+\rho}  
R(t_0)^\sigma  (\|\psi\|_{C^1_x,(t_0,t)}+\|Z_0\|_{0} ).
\end{align*}
which gives
\begin{align}
\label{eq:est-phi2-bad}
t_0^{-\nu_2}\|\phi_2\|_{\mathrm{in}2,t}
\leq C R(t_0)^\sigma  (\|\psi\|_{C^1_x,(t_0,t)}+\|Z_0\|_{0} ).
\end{align}
We also get from \eqref{eq:est-tildeL-psi2} and the estimate
in Proposition~\ref{Prop:phi1}
\begin{align}
\nonumber
|\phi_1(y, t)| 
\leq C \frac{R(t)^{\frac{5 - a_1}{2}}}{(1 + \rho)} 
\min \left(1, R(t)^{\frac{5 - a_1}{2}} \rho^{-2}\right) 
R(t_0)^\sigma 
(\|\psi\|_{C^1_x,(t_0,t)}+\|Z_0\|_{0} )
\end{align}
so that 
\begin{align}
\label{eq:est-phi1-bad}
t_0^{-\nu_1}\|\phi_1\|_{\mathrm{in}1,t}
\leq C R(t_0)^\sigma  (\|\psi\|_{C^1_x,(t_0,t)}+\|Z_0\|_{0} ).
\end{align}
In a similar way we obtain
\begin{align}
\label{eq:est-phi3-bad}
t_0^{-\nu_3}\|\phi_3\|_{\mathrm{in}3,t}
\leq C  (\|\psi\|_{C^1_x,(t_0,t)}+\|Z_0\|_{0} ).
\end{align}
Similarly to \eqref{eq:est-tildeL-psi-Z-upto4t0} we have
\begin{align}
\nonumber
\left|
\int_{B_{2R}} \tilde L_U(\psi+Z) \cdot Z_{1j} \,dy
\right|
\leq C  (\|\psi\|_{C^1_x,(t_0,t)}+\|Z_0\|_{0} )
\end{align}
so that from the ODE \eqref{eq:xi-main}
and the definition of $a_1$ \eqref{eq:def-a1} we get
\begin{align}
\label{eq:est-dotxi-bad}
t_0^{-\nu_1-\beta}\|\dot \xi\|_{\dot\xi,4t_0}\leq
C  (\|\psi\|_{C^1_x,(t_0,t)}+\|Z_0\|_{0} ).
\end{align}

Inserting the estimates \eqref{eq:est-dotp-bad}, \eqref{eq:est-phi2-bad},
\eqref{eq:est-phi1-bad}, \eqref{eq:est-phi3-bad}, and
\eqref{eq:est-dotxi-bad} into \eqref{eq:est-gA}
we obtain \eqref{eq:claim1-about-g}.

\medskip\noindent
\textbf{Step 2.}
Let $0<h<1$ to be determined.
Suppose \eqref{eq:assumption-gA} holds for a certain $t=t_1+h$.
By Lemma~\ref{lem:est-grad-heateq} 
\[
\|\psi\|_{C^1_x,(t_1,t_1+h)}
\leq \|\psi(t_1)\|_{C^1_x}
+ C \sqrt{h}
\|g[\psi, \phi, p, \xi, Z]\|_{L^\infty(\R^2\times(t_0,t))}
\]
and combining with \eqref{eq:claim1-about-g} we get
\begin{align*}
\|\psi\|_{C^1_x,(t_1,t_1+h)}
&\leq \|\psi(t_1)\|_{C^1_x}
+ C_{t_0} \sqrt{h} ( \|\psi\|_{C^1_x,(t_0,t_1+h)}
+ \|Z_0\|_0)
\\
&\leq \|\psi(t_1)\|_{C^1_x}
+ C_{t_0} \sqrt{h} \|\psi\|_{C^1_x,(t_0,t_1)}
+ C_{t_0} \sqrt{h} \|\psi\|_{C^1_x,(t_1,t_1+h)}
\\
& \quad 
+C_{t_0} \sqrt{h}  \|Z_0\|_0 .
\end{align*}
By taking $h$ small (depending on $t_0$),
\[
\|\psi\|_{C^1_x,(t_1,t_1+h)}
\leq 
C\|\psi(t_1)\|_{C^1_x}
+ C_{t_0} \sqrt{h} 
( \|\psi\|_{C^1_x,(t_0,t_1)}+\|Z_0\|_{0} )
\]
and therefore
\begin{align}
\label{eq:iter-psi}
\|\psi\|_{C^1_x,(t_0,t_1+h)}
\leq  C_1 (\|\psi\|_{C^1_x,(t_0,t_1)}+\|Z_0\|_{0} ),
\end{align}
for some constant $C_1>0$.
Let $t_k = t_0+ k h$, $k\in\N$, and 
\[
a_k = \|\psi\|_{C^1_x,(t_0,t_k)}
\]
Iterating \eqref{eq:iter-psi} gives 
\[ 
a_k\leq C_1^k\|Z_0\|_0 
\] 
as long as the assumptions \eqref{eq:assumption-gA} remain valid.
Let $k_*\in\N$ be such that $t_{k_*}>4t_0$.
By asking that
\[
C_1^{k_*} \|Z_0\|_{0} <1
\]
we deduce that $\psi$ does not blow up before $4 t_0$ and 
\begin{align}
\label{eq:est-psi-up-to-4t0}
\|\psi\|_{C^1_x,(t_0,4t_0)}
\leq C_1^{k_*}\|Z_0\|_{0} .
\end{align}
The estimates in the statement \eqref{eq:claim1} follow from 
\eqref{eq:est-psi-up-to-4t0} and \eqref{eq:est-dotp-bad}--\eqref{eq:est-dotxi-bad}.
\end{proof}

Assuming that $t_0$ is large and that $\|\dot p\|_{\dot p}\leq 1$, $\|\dot\xi\|_{\dot\xi}\leq 1$, we have
\begin{equation}
\label{assumption_param_outer}
|p(t)-1|\leq \tfrac12,\quad
|\xi(t)|\leq 1.
\end{equation}

To carry out the proof of this result, we need to derive an estimate for $\psi$ solving the outer equations \eqref{Outer_eq2}.

\begin{proposition}
Let $g$ be given by \eqref{def_g}, which is the right-hand side of the outer equation \eqref{Outer_eq2}.
Assume 
\begin{align}
\label{eq:claim1b}
\left\{
\begin{aligned}
&\|\phi_1\|_{\mathrm{in}1,t_1} \leq1,\quad 
\|\phi_2\|_{\mathrm{in}2,t_1} \leq1,\quad
\|\phi_3\|_{\mathrm{in}3,t_1} \leq1,
\\
&
\|\dot p\|_{\dot p,t_1}\leq 1,\quad 
\|\dot\xi\|_{\dot\xi,t_1}\leq 1.
\end{aligned}
\right.
\end{align}
Then
\begin{align}
\nonumber
\norm{g[\psi, \phi, p, \xi, Z]}_{*,t_1}
&\leq
\frac C{t_0^\sigma}
\left(
\|\psi\|_{\mathrm{out},t_1}
+\norm{Z_0}_{0}
+\|\psi\|_{\mathrm{out},t_1}^2
\right)
\\
\nonumber
&\quad 
+\frac C{t_0^\sigma}
\left(
\|\phi_1\|_{\mathrm{in}1,t_1}
+\|\phi_3\|_{\mathrm{in}3,t_1}
\right)
+C\|\phi_2\|_{\mathrm{in}2,t_1}    
\\
\label{eq:est-RHS-outer}
&\quad+\frac{C}{t_0^\sigma}(
\|\dot p\|_{\dot p,t_1}
+ \|\dot \xi\|_{\dot \xi,t_1})
\end{align}
for some $\sigma>0$.
\end{proposition}
\begin{proof}
We decompose $g=g_1+g_2+g_3+g_4$ where
\begin{align}
g_1[\phi, p, \xi] &= Q_\omega \big((\Delta_x \eta)\phi + 2\lambda^{-1}\nabla_x \eta \cdot \nabla_y \phi - \eta_t \phi\big)\nonumber
\\
\label{def_g1}
&\quad + \eta Q_\omega \big(\lambda^{-1}\dot \lambda y \cdot \nabla_y \phi + \lambda^{-1} \dot \xi \cdot \nabla_y \phi - \dot \omega J\phi\big) 
\\
\label{def_g2}
g_2[\psi, p, \xi, Z] &=  (1 - \eta) \LtU{\psi + Z} + (1 - \eta) \left(\mathcal{K}_0[p, \xi] + \mathcal E_1[p, \xi]\right)
\\
\label{def_g3}
g_3[\psi, p, \xi, Z] 
&= \left(\big(\psi + \Phi^0[p, \xi] + Z\big)  \cdot U\right)U_t + \piU{\tilde{\mathcal R}_1} 
\\
\label{def_g4}
g_4[\psi, \phi, p, \xi, Z] 
&= N_U\left(\eta  \phi + \piU{\psi +\Phi^0[p, \xi] + Z}\right) .
\end{align}

\medskip\noindent
\textbf{Estimate of $g_1$ \eqref{def_g1}.} 
We claim that
\begin{align}
\label{eq:est-g1}
\|g_1[\phi, p, \xi]\|_{*,t_1} \le 
\frac C{t_0^\sigma}
(\|\phi_1\|_{\text{in}1,t_1}
+\|\phi_3\|_{\text{in}3,t_1})
+C\|\phi_2\|_{\text{in}2,t_1}
\end{align}
for some $\sigma>0$.
In the estimates below we always consider $t<t_1$.
Notice that, thanks to the cut-off function $\eta$, defined in \eqref{def_eta}, and \eqref{assumption_param_outer}, each of these pieces is supported in the ball
\begin{align*}
\text{supp } \eta &\subset B_{2\lambda R(t)}(\xi(t)) \subset B_{4R(t)}(0) ,
\end{align*}
for $t_0$ sufficiently large. 
Let us consider the term $Q_\omega (\Delta_x \eta)\phi = (\Delta_x \eta)Q_\omega \big(\phi_1 + \phi_2 + \phi_3\big)$. The remaining terms involving derivatives of the cut-off function $\eta$, namely 
$2\lambda^{-1}Q_\omega\big(\nabla_x \eta \cdot \nabla_y \phi\big)$ and  $\eta_t Q_\omega\phi$,
can be treated in a similar manner and reduced to this case. Notice that, since 
$$\eta(x, t) = \eta_0\left(\frac{\abs{x - \xi(t)}}{\lambda(t)R(t)}\right)$$
we obtain the estimate $|\Delta_x\eta(x, t)| \leq \frac{C}{(\lambda R)^2}$
and  supported in the region $R(t) \le \abs{y}\le 2R(t)$, for $y = \frac{x - \xi}{\lambda}$. Therefore, $\abs{y}$ behaves like $R$ in this region. We compute for $t<t_1$
\begin{align*}
|Q_\omega (\Delta_x \eta)\phi_1| &\le  |\Delta_x \eta|
\frac{R(t)^{\frac{5 - a_1}{2}}}{t^{\nu_1}( 1 + |y|)}\min\bigg\{1, \frac{R(t)^{\frac{5 - a_1}{2}}}{|y|^2}\bigg\}\|\phi_1\|_{\text{in}1,t_1}
\\
&\le C\frac{\chi_{B_{4R}}}{ R(t)^2t^{\nu_1}}
\frac{1}{R(t_0)^{a_1 - 2}}\|\phi_1\|_{\text{in}1,t_1}
\\
&\le \frac C{t_0^\sigma}\varsigma_1(x, t)\|\phi_1\|_{\text{in}1,t_1}
\end{align*}
\begin{align}  
\nonumber
|Q_\omega (\Delta_x \eta)\phi_2| &\le  |\Delta_x \eta| \frac{R(t)^2}{t^{\nu_2}(1 + \abs{y})} 
\|\phi_2\|_{\text{in}2}
\\
\nonumber
&\le C\frac{\chi_{B_{4R}}}{ R(t)^2t^d} 
\|\phi_2\|_{\text{in}2,t_1}
\\
\nonumber
&\le C\varsigma_1(x, t)
\|\phi_2\|_{\text{in}2,t_1}
\end{align}
\begin{align*}
|Q_\omega (\Delta_x \eta)\phi_3| &\le  |\Delta_x \eta| \frac{\log R}{t^{\nu_3}} \|\phi_3\|_{\text{in}3,t_1}
\\
&\le C\frac{\chi_{B_{4R}}}{ R(t)^2t^{\nu_1}}  \frac{1}{t_0^{1 - \nu_1}}\|\phi_3\|_{\text{in}3,t_1}
\\
&\le \frac C{t_0^\sigma}\varsigma_1(x, t)\|\phi_3\|_{\text{in}3,t_1}.
\end{align*}
These estimates imply that  
\begin{align}
\label{est_Delta_eta_phi}
\norm{Q_\omega (\Delta_x \eta)\phi}_{*,t_1} \le \frac C{t_0^\sigma}
(\|\phi_1\|_{\text{in}1,t_1}
+\|\phi_3\|_{\text{in}3,t_1})
+ C\|\phi_2\|_{\text{in}2,t_1}.
\end{align}
Let us estimate the term $\eta Q_\omega \big(\lambda^{-1}\dot \lambda y \cdot \nabla_y \phi\big)$.
For this we remark that we 
will arrange the parameters so that 
\eqref{nl1} and \eqref{nl2} hold.
Together with $\nu_1<1$ we get
\begin{align}
\label{eq:nl-param2}
\frac{\beta(9-a_1)}2<1,\quad 3\beta<1.
\end{align}
Using \eqref{eq:claim1b}, the term $\eta Q_\omega \big(\lambda^{-1}\dot \lambda y \cdot \nabla_y \phi\big)$ has the estimate
\begin{align*}
|\eta Q_\omega \big(\lambda^{-1}\dot \lambda y \cdot \nabla_y \phi_1\big)| 
&\le  C\chi_{B_{4R}}|\dot \lambda(t)| \frac{R(t)^{\frac{5 - a_1}{2}}}{t^{\nu_1}( 1 + |y|)}\min\bigg\{1, \frac{R(t)^{\frac{5 - a_1}{2}}}{|y|^2}\bigg\} 
\|\phi_1\|_{\text{in}1,t_1}
\\
&\le C\frac{\chi_{B_{4R}}}{t\log^2(t)} 
\frac{R(t)^{\frac{5 - a_1}{2}}}{t^{\nu_1}}
\|\phi_1\|_{\text{in}1,t_1}
\\
&\le C\frac{\chi_{B_{4R}}}{ R(t)^2t^{\nu_1}} \frac{R(t_0)^{\frac{9 - a_1}{2}}}{t_0\log^2(t_0)}
\|\phi_1\|_{\text{in}1,t_1}
\\
&\le \frac C{t_0^\sigma}\varsigma_1(x, t)\|\phi_1\|_{\text{in}1,t_1}
\end{align*}
where we have used \eqref{eq:nl-param2}.
We also have
\begin{align*}    
|\eta Q_\omega \big(\lambda^{-1}\dot \lambda y \cdot \nabla_y \phi_2\big)| &\le  C\chi_{B_{4R}}|\dot \lambda(t)|  \frac{R(t)^2}{t^{\nu_2}(1 + \abs{y})} \|\phi_2\|_{\text{in}2,t_1}
\\
&\le  \frac{\chi_{B_{4R}}}{t\log^2(t)}  \frac{R(t)^2}{t^{\nu_2}} 
\|\phi_2\|_{\text{in}2,t_1}
\\
&\le C\frac{\chi_{B_{4R}}}{ R(t)^2t^{\nu_1}}  \frac{R(t_0)^3}{t_0\log^2(t_0)}
\|\phi_2\|_{\text{in}2,t_1}
\\
&\le \frac C{t_0^\sigma}\varsigma_1(x, t)\|\phi_2\|_{\text{in}2,t_1}
\end{align*}
by \eqref{eq:nl-param2}.
Similarly,
\begin{align*}
|\eta Q_\omega \big(\lambda^{-1}\dot \lambda y \cdot \nabla_y \phi_3\big)| 
&\le 
C\chi_{B_{4R}}|\dot \lambda(t)| \frac{\log R}{t^{\nu_3}} 
\|\phi_3\|_{\text{in}3,t_1}
\\
&\le   C\frac{\chi_{B_{4R}}}{t\log^2(t)} \frac{\log R}{t^{\nu_3}}
\|\phi_3\|_{\text{in}3,t_1}
\\
&\le C\frac{\chi_{B_{4R}}}{ R(t)^2t^{\nu_1}} \frac{R(t_0)^2 }{t_0^2 \log(t_0)}
\|\phi_3\|_{\text{in}3,t_1}
\\
&\le \frac C{t_0^\sigma}\varsigma_1(x, t)\|\phi_3\|_{\text{in}3,t_1}.
\end{align*}
These bounds imply
\begin{align}
\label{est_y_gradphi}
\big\|\eta Q_\omega \big(\lambda^{-1}\dot \lambda y \cdot \nabla_y \phi\big)\big\|_{*,t_1} \le 
\frac C{t_0^\sigma}
(\|\phi_1\|_{\text{in}1,t_1}
+\|\phi_2\|_{\text{in}2,t_1}
+\|\phi_3\|_{\text{in}3,t_1}).
\end{align}
The last two terms in $g_1$, $ \eta Q_\omega\big(\lambda^{-1} \dot \xi \cdot \nabla_y \phi\big)$ and $ \eta Q_\omega\big(\dot \omega J\phi\big)$, are estimated in the same way, it suffices to notice that 
\begin{align*}
|\dot \xi(t)| \le \frac{1}{R(t)t^{d - \theta}} \quad \text{and} \quad \abs{\dot \omega(t)} \le \abs{\dot p(t)} \le  
\frac{1}{t\log^2(t)}.
\end{align*}
This observation, \eqref{est_Delta_eta_phi}, and \eqref{est_y_gradphi} prove \eqref{eq:est-g1}.

\medskip\noindent
\textbf{Estimate of $g_2$ \eqref{def_g2}.}
We claim that
\begin{align}
\label{eq:est-g2}
\|g_2[\psi, p, \xi, Z]\|_{*,t_1} \le 
\frac{C}{t_0^\sigma}\left(
\|\psi\|_{\text{out},t_1} 
+ \|\dot p\|_{\dot p,t_1}
+ \|\dot \xi\|_{\dot \xi,t_1}
+ \norm{Z_0}_{0} 
\right),
\end{align}
where $\sigma>0$ is a small constant.

Let us first estimate $(1 - \eta) \tilde L_U(\psi + Z)$, starting with the term $\tilde L_U(\psi)$. Notice that the support of $(1 - \eta)$ is included in
\begin{align*}
\mathrm{supp}(1 - \eta) &\subset \big(B_{\lambda R(t)}(\xi(t))\big)^c \subset \big(B_{\frac{1}{2}R(t)}(0)\big)^c,
\end{align*}
provided that $t_0$ is sufficiently large. 
This implies, using the definition of the operator $\tilde{L}_U$ in \eqref{def_tildeL}, that 
\begin{align*}
(1 - \eta) |\tilde L_U(\psi)(x, t)| \le C \frac{\chi_{\{\abs{x} \ge \frac{1}{2}R(t)\}}}{(1 + \abs{x})^2} \abs{\nabla \psi(x, t)}.
\end{align*}
Additionally, by definition of the norm 
$\|\ \|_{\text{out}}$ \eqref{norm_psi}
(where we are using $d_0=\nu_1$),
\begin{align*}
\abs{\nabla \psi(x, t)} \le 
\frac{\log t}{R(t) t^{\nu_1}} \frac{1}{\big(1 + \frac{\abs{x}}{\sqrt{t}}\big)^2}
\|\psi\|_{\mathrm{out},t_1}.
\end{align*}
Then, by the assumption \eqref{eq:claim1b},
\begin{align*}
(1 - \eta) |\LtU{\psi}(x, t)| &\le C     
\frac{\chi_{\{\abs{x} \ge \frac{1}{2}R(t)\}}}{(1 + \abs{x})^2} \frac{\log t}{R(t) t^{\nu_1}} \frac{1}{\big(1 + \frac{\abs{x}}{\sqrt{t}}\big)^2} \|\psi\|_{\mathrm{out},t_1}
\\
&\le  \frac{C}{t_0^{\beta-\theta}}
\frac{\chi_{\{\abs{x} \ge \frac{1}{2}R(t)\}}}{t^{\nu_1+\theta}(1 + \abs{x})^2}  
\frac{1}{\big(1 + \frac{\abs{x}}{\sqrt{t}}\big)^2}\|\psi\|_{\mathrm{out},t_1}
\\
&\le  \frac{C}{t_0^\sigma} \varsigma_3(x, t) \|\psi\|_{\mathrm{out},t_1}
\end{align*}
for some $\sigma>0$.
We deduce that 
\begin{align}
\label{eta_tL_psi}
\big\|(1 - \eta) \LtU{\psi}\big\|_{*,t_1} \le \frac{C}{t_0^\sigma}\|\psi\|_{\mathrm{out},t_1} .
\end{align}

The estimate of the term $(1 - \eta) \LtU{Z}$ follows the same argument as the preceding estimate. Recall that $Z$ solves the homogeneous heat equation with initial data $Z_0$ \eqref{eq:heat-eq-Z}. Since $Z_0$ satisfies \eqref{eq:normZ0}, the corresponding solution $Z$ exhibits the decay described in Lemma~\ref{Lem:est_Z}. In particular, we have
\begin{align*}
    \abs{\nabla Z(x, t)} &\le  C\frac{\|Z_0\|_{0}}{t^{d + \frac{1}{2}}\big(1 + \frac{\abs{x}}{\sqrt{t}}\big)^{2d}} .
\end{align*}
This yields 
\begin{align*}
(1 - \eta) |\LtU{Z}(x, t)| \le \varsigma_3(x, t)\normZ \frac{C}{\sqrt{t_0}}
\end{align*}
It follows that 
\begin{align}
\label{eta_tL_Z}
\big\|(1 - \eta) \LtU{Z}\big\|_{*,t_1} 
\le \frac{C}{t_0^\sigma}\normZ, 
\end{align}
for $0<\sigma<\frac12$.

Let us now estimate $(1 - \eta) \mathcal{K}_0[p, \xi]$, defined in \eqref{eq:K01}.
Using Lemma \ref{lem:estK0} we get 
\begin{align*}
\abs{(1 - \eta) \mathcal{K}_{0}[p, \xi]} 
&\le C\chi_{\{\abs{x} \ge \frac{1}{2}R\}} \frac{1}{(1 + \abs{x})^3} \frac{1}{t\log^2(t)} \frac{1}{\big(1 + \frac{\abs{x}}{\sqrt{t}}\big)^2}
\|\dot p\|_{\dot p,t_1}
\\
&\le C\varsigma_3(x, t) \frac{1}{t_0^\sigma}
\|\dot p\|_{\dot p,t_1}
\end{align*}
for some $\sigma>0$.
This implies that 
\begin{align}
\label{est_K01}
\norm{(1 - \eta)\mathcal{K}_{0}[p, \xi]}_{*,t_1} \le  \frac{C}{t_0^\sigma}
\|\dot p\|_{\dot p,t_1}.
\end{align}

Let us now estimate $(1 - \eta) \mathcal{E}_1[p, \xi]$, defined in \eqref{eq:K1}. We have
\begin{align*}
|(1 - \eta)\mathcal{E}_{1}[p, \xi]| &\le C\frac{\chi_{\{\abs{x} \ge \frac{1}{2}R(t)\}}}{R(t)t^{\nu_1}}\frac{1}{(1 + \abs{x})^2}\|\dot\xi\|_{\dot\xi,t_1}
\\
&\le \frac C{t_0^\sigma}\varsigma_3(x, t) 
\|\dot\xi\|_{\dot\xi,t_1}
\end{align*}
for some $\sigma>0$.
We therefore deduce that
\begin{align}
\label{est_K1_xi}
\norm{(1 - \eta)\mathcal{E}_{1}[p, \xi]}_{*,t_1} \leq \frac C{t_0^\sigma}
\|\dot\xi\|_{\dot\xi,t_1}.
\end{align}

By \eqref{eta_tL_psi}, \eqref{eta_tL_Z}, \eqref{est_K01} and \eqref{est_K1_xi}, we deduce the inequality \eqref{eq:est-g2}.

\medskip\noindent
\textbf{Estimate of $g_3$ \eqref{def_g3}.}
We claim that 
\begin{align}
\label{eq:est-g3}
\|g_3[\psi, p, \xi, Z]\|_{*,t_1} \leq 
\frac{C}{t_0^\sigma}
\left(\|\psi\|_{\mathrm{out},t_1}+\norm{Z_0}_{0}
+\|\dot p\|_{\dot p,t_1}\right).
\end{align}

Let us estimate $\left(\big(\psi + Z +  \Phi^0[p, \xi]\big)  \cdot U\right)U_t$, $\piU{\tilde{\mathcal R}_1}$ is handled similarly. Notice that because of the assumption \eqref{eq:claim1b}, we have
\begin{align}
\label{est_Ut_outer}
\abs{U_t} \le C|\dot \lambda(t)| \abs{\rho w_\rho} + \abs{\dot \omega(t)} \abs{\rho w_\rho} + |\dot \xi(t)| \abs{w_\rho} \le C\frac{1}{t\log^2(t) (1 + \abs{x})}. 
\end{align}
So  
\begin{align*}
\abs{\left(\psi  \cdot U\right)U_t} \le \abs{\psi} \abs{U_t} 
&\le C\frac{1}{t\log^2(t) (1 + \abs{x})} \frac{\log t}{t^{\nu_1}\big(1 + \frac{\abs{x}}{\sqrt{t}}\big)^2}
\|\psi\|_{\mathrm{out},t_1}
\\
&= C\frac{1}{t^{1+\nu_1}\log^2(t) } \frac{\log t}{(1 + \abs{x})\big(1 + \frac{\abs{x}}{\sqrt{t}}\big)^2}
\|\psi\|_{\mathrm{out},t_1}
\\
&\leq \frac{C}{t_0^\sigma} ( \varsigma_1(x,t)+\varsigma_3(x,t))
\|\psi\|_{\mathrm{out},t_1}
\end{align*}
for some $\sigma>0$.
Therefore, we find that 
\begin{align}
\label{est_psiU_Ut}
\norm{\big(\psi\cdot U\big)U_t}_{*,t_1} \le \frac{C}{t_0^\sigma}
\|\psi\|_{\mathrm{out},t_1}.
\end{align}

Let us now estimate the term $(Z \cdot U)U_t$. By Lemma~\ref{Lem:est_Z}, together with the bound for $U_t$ \eqref{est_Ut_outer}, implies that
\begin{align*}
\abs{(Z \cdot U)U_t} 
& =C\frac{1}{t^{1+d}\log^2(t) (1 + \abs{x})} \frac{1}{\big(1 + \frac{\abs{x}}{\sqrt{t}}\big)^{2d}}
\|Z_0\|_{0}.
\end{align*}
Under the assumption
\begin{align}
\label{eq:cond-d-nu1-b}
\nu_1+2\beta<1+d
\end{align}
we get
\begin{align*}
\abs{(Z \cdot U)U_t}  
\leq \frac C{t_0^\sigma}(\varsigma_1(x,t)+\varsigma_3(x,t)) \|Z_0\|_0
\end{align*}
for some $\sigma>0$.
This yields
\begin{align}
\label{est_Z_Ut}
\norm{\big(Z\cdot U\big)U_t}_{*,t_1} \le \frac{C}{t_0^\sigma}\|Z_0\|_{0}.
\end{align}

For $(\Phi^0[p, \xi] \cdot U)U_t$, using Lemma~\ref{lem:phi0} we get
\begin{align*}
|\Phi^0[p, \xi](x, t)| &\le 
C\frac{r}{\log(t_0) t } \frac{1}{\big(1 + \frac{\abs{x}}{\sqrt{t}}\big)^2}
\|\dot p\|_{\dot p,t_1}
\end{align*}
for $t<t_1$.
Additionally, recalling the formulas
\begin{align*}
U(x, t) &= \begin{pmatrix}
e^{i(\theta + \omega)} \sin w(\rho)\\
\cos w(\rho)
\end{pmatrix}, \quad w_\rho(\rho) = - \frac{2}{1 + \rho^2}, 
\\
\sin w(\rho) &= \frac{2\rho}{1 + \rho^2}, \quad \cos w(\rho) = \frac{\rho^2 - 1}{1 + \rho^2},
\end{align*}
and using the estimate \eqref{est_Ut_outer}, there holds
\begin{align*}
\abs{(\Phi^0[p, \xi] \cdot U)U_t} 
&\le  C\frac{1}{t \log^2(t) (1 + \abs{x})} \frac{1}{\log(t_0)t\big(1 + \frac{\abs{x}}{\sqrt{t}}\big)^2}
\|\dot p\|_{\dot p,t_1}
\\
&\leq \frac{C}{t_0^\sigma} (\varsigma_1(x,t)+\varsigma_3(x,t))
\|\dot p\|_{\dot p,t_1},
\end{align*}
for some $\sigma>0$.
It follows that 
\begin{align}
\label{est_Phi0_Ut}
\norm{(\Phi^0[p, \xi] \cdot U)U_t}_{*,t_1} \le 
\frac{C}{t_0^\sigma}\|\dot p\|_{\dot p,t_1}.
\end{align}

For the estimate of $\piU{\tilde{\mathcal R}_1}$, defined in \eqref{def_R1}, it suffices to notice that the main order term is given by 
\begin{align*}
\big|\piU{\tilde{\mathcal R}_1}\big|\le |\dot \xi| \left|\int_{t_0}^t \dot p(s) k(z(r), t - s)ds\right|.
\end{align*}
This integral is precisely the one analyzed in 
Lemma~\ref{lem:phi0}, and we deduce that 
\begin{align*}
\big|\piU{\tilde{\mathcal R}_1}\big| 
&\le \frac{1}{t^{\nu_1+\beta}} \frac{C}{\log(t_0)t\big(1 + \frac{\abs{x}}{\sqrt{t}}\big)^2} 
\|\dot p\|_{\dot p,t_1}
\\
&= C \frac{1}{t^{1+\nu_1+\beta}\log(t_0)} \frac{1}{\big(1 + \frac{\abs{x}}{\sqrt{t}}\big)^2}
\|\dot p\|_{\dot p,t_1}
\\
&\leq \frac{C}{t_0^\sigma}(\varsigma_1(x,t)+\varsigma_3(x,t))
\|\dot p\|_{\dot p,t_1}
\end{align*}
for some $\sigma>0$.
This implies
\begin{align}
\label{est_piR_1new}
\big\|\piU{\tilde{\mathcal R}_1}\big\|_{*,t_1} \le \frac{C}{t_0^\sigma} \|\dot p\|_{\dot p,t_1}.
\end{align}

Using \eqref{est_psiU_Ut}, \eqref{est_Z_Ut}, \eqref{est_Phi0_Ut} and \eqref{est_piR_1new}, we obtain \eqref{eq:est-g3}.

\medskip\noindent
\textbf{Estimate of $g_4$ \eqref{def_g4}.}
We claim that
\begin{align}
\nonumber
\|g_4[\psi, p, \xi, Z]\|_{*,t_1}
&\leq 
\frac{C}{t_0^\sigma} \bigl(
\|\psi\|_{\mathrm{out},t_1}^2
+\|\phi_1\|_{\mathrm{in}1,t_1}^2
+\|\phi_2\|_{\mathrm{in}2,t_1}^2
+\|\phi_3\|_{\mathrm{in}3,t_1}^2
\\
\label{eq:est-g4}
& \qquad 
+\|\dot p\|_{\dot p}^2
+\|Z_0\|_0^2
\bigr).
\end{align}

We recall 
that
 $g_4[\psi, \phi, p, \xi, Z] =N_U\left(\eta  \phi + \piU{\psi +\Phi^0[p, \xi] + Z}\right)$
and 
that the quadratic term is given by
\begin{equation*}
N_U(\zeta) = \left(2\nabla(aU) \cdot \nabla (U + \zeta) + 2 \nabla U \cdot \nabla \zeta + \abs{\nabla \zeta}^2 + \abs{\nabla(aU)}^2 \right) \zeta - aU_t + 2\nabla a \cdot \nabla U,
\end{equation*}
for $\zeta = \piU \varphi = \eta \phi + \piU{\psi + \Phi^0 + Z}$ and $a = a(\zeta) = \sqrt{1 - \abs{\zeta}^2} - 1$.
Notice that $a(\zeta)$ expands, for small $\zeta$, as 
$$a(\zeta) = -\frac{1}{2}\abs{\zeta}^2 + O\big(\abs{\zeta}^4\big) \quad \text{and} \quad \nabla a(\zeta) = -\zeta \cdot \nabla \zeta + O\big(|\nabla \zeta|\abs{\zeta}^3\big).$$
Using this, the last two terms of $N_U$ can be expressed, at main order, as \begin{align*}
-a(\zeta)U_t &= \frac{1}{2}\abs{\zeta}^2U_t + O\big(\abs{\zeta}^4|U_t|\big),\\
2\nabla a(\zeta) \cdot \nabla U &= - \zeta \cdot \nabla \zeta \cdot \nabla U + O\big(|\nabla \zeta|\abs{\zeta}^3|\nabla U|\big).
\end{align*}
It suffices to treat these two terms since all other combinations in $N_U$ can be reduced to these cases. Moreover, for $\zeta = \eta \phi + \piU{\psi + \Phi^0 + Z}$, the contribution of $\piU{\Phi^0 + Z}$ exhibits a faster decay than that of $\eta \phi + \piU{\psi}$. Consequently, it is enough to focus on the latter, as the remaining terms can be handled through analogous computations and yield contributions of lower order. 

We therefore begin with the term involving $\psi$, recalling that $\piU{\psi} = \psi - (\psi \cdot U) U$. We first analyze the quantity $\zeta \cdot \nabla \zeta \cdot \nabla U$, and then the term $\frac{1}{2}\abs{\zeta}^2U_t$.
We recall that 
$$\abs{\nabla U(x, t)} \le C \frac{1}{(1 + \abs{x})^2},$$
and $\nabla \piU{\psi}$ expands as 
$$\nabla \piU{\psi} = \piU{\nabla \psi} -  (\psi \cdot \nabla U) U -  (\psi \cdot U) \nabla U.
$$
Using the definition of $\|\psi\|_{\text{out}}$ \eqref{norm_psi} we see that:
\begin{align}
\nonumber
|\piU{\psi} &~\cdot \nabla \piU{\psi} \cdot \nabla U| 
\\
\nonumber
&\le C \big(|\psi||\nabla \psi| |\nabla U| + |\psi|^2|\nabla U|^2\big)
\|\psi\|_{\mathrm{out},t_1}^2
\\
\nonumber
&\le C\frac{(\log t)^2}{R(t)t^{2\nu_1}\big(1 + \frac{\abs{x}}{\sqrt{t}}\big)^4} \frac{1}{(1 + \abs{x})^2}\|\psi\|_{\mathrm{out},t_1}^2 
\\
&\quad 
\label{eq:nl-psi2}
+  
C\frac{(\log t)^2}{t^{2\nu_1}\big(1 + \frac{\abs{x}}{\sqrt{t}}\big)^4} \frac{1}{(1 + \abs{x})^4}
\|\psi\|_{\mathrm{out},t_1}^2.
\end{align}
Regarding the first term, we have
\begin{align}
\label{eq:nl-psi2a}
\frac{(\log t)^2}{R(t)t^{2\nu_1}\big(1 + \frac{\abs{x}}{\sqrt{t}}\big)^4} \frac{1}{(1 + \abs{x})^2}
&\leq \frac{C}{t_0^\sigma}(\varsigma_1(x,t)+\varsigma_3(x,t))
\end{align}
for some $\sigma>0$.
Regarding the term in \eqref{eq:nl-psi2} we have
\begin{align}
\label{eq:nl-psi2b}
\frac{(\log t)^2}{t^{2\nu_1}\big(1 + \frac{\abs{x}}{\sqrt{t}}\big)^4} \frac{1}{(1 + \abs{x})^4}
\leq \frac{C}{t_0^\sigma}(\varsigma_2(x,t)+\varsigma_3(x,t)).
\end{align}
From \eqref{eq:nl-psi2}, \eqref{eq:nl-psi2a} and \eqref{eq:nl-psi2b} we get
\begin{align}
\label{eq:est-nl1}
\|\piU{\psi} &~\cdot \nabla \piU{\psi} \cdot \nabla U\|_{*,t_1}\leq
\frac{C}{t_0^\sigma} \|\psi\|_{\mathrm{out},t_1}^2.
\end{align}

Next, we discuss the estimate of 
$|\piU{\psi}|^2 U_t$. Using \eqref{est_Ut_outer} it follows that 
\begin{align*}
\abs{\piU{\psi}}^2 \abs{U_t} 
&\le C \frac{\log^2(t)}{t^{2\nu_1}\big(1 + \frac{\abs{x}}{\sqrt{t}}\big)^4}\frac{C}{t\log^2(t) (1 + \abs{x})}
\|\psi\|_{\mathrm{out},t_1}^2
\\
&= \frac C{t_0^\sigma} \frac{1}{t^{\nu_1+1}}\frac{1}{(1 + \abs{x})\big(1 + \frac{\abs{x}}{\sqrt{t}}\big)^4}
\|\psi\|_{\mathrm{out},t_1}^2
\\
&\leq 
\frac C{t_0^\sigma}(\varsigma_1(x,t)+\varsigma_3(x,t))
\|\psi\|_{\mathrm{out},t_1}^2
\end{align*}
for some $\sigma>0$.
This shows that 
\begin{align}
\nonumber
\|\, \abs{\piU{\psi}}^2 \abs{U_t} \,\|_{*,t_1}
\leq \frac C{t_0^\sigma}
\|\psi\|_{\mathrm{out},t_1}^2.
\end{align}

We next estimate the contribution of $\zeta = \eta \phi$ in the terms $\zeta \cdot \nabla \zeta \cdot \nabla U$ and $\frac{1}{2}\abs{\zeta}^2 U_t$. Since the estimate for $\phi_3$ is better than those for $\phi_1$ and $\phi_2$, it is sufficient to restrict our attention to the latter two. The main difficulty associated with these terms arises from their behaviour near $y = 0$. More precisely, one has
$$\phi_1(0, t) \sim \frac{R(t)^{\frac{5 - a_1}{2}}}{t^{\nu_1}} \quad \text{and} \quad \phi_2(0, t) \sim \frac{R(t)^2}{t^{\nu_2}},$$
where 
$$\quad \nu_2 = \nu_1 + \beta.$$

Starting with $\zeta \cdot \nabla \zeta \cdot \nabla U$, first notice that 
$$\nabla (\eta \phi_i) = (\nabla \eta) \phi_i + \eta \nabla \phi_i, \quad i = 1, 2.$$
As $\nabla \eta$ is supported near 
$$y \sim R(t)$$
the term $(\nabla \eta) \phi_i$ is of lower order. We focus our attention on the second term $\eta \nabla \phi_i$ and we get 
\begin{align}
\nonumber
\abs{\eta \phi_1 \cdot \nabla(\eta \phi_1) \cdot \nabla U} &\le \eta \abs{\phi_1} \abs{\nabla  \phi_1} \abs{\nabla U} 
\|\phi_1\|_{\mathrm{in}1,t_1}^2
\\
\nonumber
&\le C\chi_{B_{4R}}\frac{R(t)^{5 - a_1}}{t^{2\nu_1}(1 + |y|)^5}\min\bigg\{1, \frac{R(t)^{\frac{5 - a_1}{2}}}{|y|^2}\bigg\}^2
\|\phi_1\|_{\mathrm{in}1,t_1}^2
\\
\nonumber
&\le C\chi_{B_{4R}} \frac{1}{R(t)^{1 + 2\gamma}t^{\nu_1}(1 + \abs{y})^5} \frac{R(t_0)^{6 - a_1 + 2\gamma}}{t_0^{\nu_1}} 
\|\phi_1\|_{\mathrm{in}1,t_1}^2
\\
\nonumber
&\le \frac{C}{t_0^\sigma} \varsigma_2(x, t)
\|\phi_1\|_{\mathrm{in}1,t_1}^2
\end{align}
for some $\sigma>0$, $\gamma>0$ small under the assumption
\begin{align}
\label{restriction_csts1}
\beta(6 - a_1) < \nu_1.
\end{align}

By a similar calculation, using that $\nu_2=\nu_1+\beta$, 
\begin{align}
\nonumber
\abs{\eta \phi_2 \cdot \nabla(\eta \phi_2) \cdot \nabla U} 
&\le C\chi_{B_{4R}}\frac{R(t)^4}{t^{2\nu_2}(1 + \abs{y})^5} 
\|\phi_2\|_{\mathrm{in}2,t_1}^2
\\
\nonumber
&\le C\frac{\chi_{B_{4R}}}{R(t)^{1 + 2\gamma}t^{\nu_1}(1 + \abs{x})^{5}}\frac{R(t_0)^{3 + 2\gamma }}{t_0^{\nu_1}}
\|\phi_2\|_{\mathrm{in}2,t_1}^2
\\
\nonumber
&\le 
\frac{C}{t_0^\sigma}
\varsigma_2(x, t)
\|\phi_2\|_{\mathrm{in}2,t_1}^2
\end{align}
for some $\gamma>0$, $\sigma>0$ under the assumption
\begin{align}
\label{restriction_csts2}
3\beta <\nu_1,
\end{align}
Similarly, 
\begin{align}
\nonumber
\abs{\eta \phi_3 \cdot \nabla(\eta \phi_3) \cdot \nabla U} 
\leq 
\frac{C}{t_0^\sigma}
\varsigma_2(x, t)
\|\phi_3\|_{\mathrm{in}3,t_1}^2.
\end{align}
Therefore, we deduce that 
\begin{align}
\label{eq:est-nl3}
\norm{\eta \phi \cdot \nabla (\eta \phi) \cdot \nabla U}_{*,t_1} \le \frac{C}{t_0^\sigma}
\left(\|\phi_1\|_{\mathrm{in}1,t_1}^2
+\|\phi_2\|_{\mathrm{in}2,t_1}^2
+\|\phi_3\|_{\mathrm{in}3,t_1}^2\right)
\end{align}

For the term $|\zeta|^2 U_t$, $\zeta = \eta \phi$, we get (using \eqref{est_Ut_outer})
\begin{align*}
|\eta \phi_1|^2 |U_t| &\le C\chi_{B_{4R}}\frac{R(t)^{5 - a_1}}{t^{2\nu_1}(1 + |y|)^2}
\min\bigg\{1, \frac{R(t)^{\frac{5 - a_1}{2}}}{|y|^2}\bigg\}^2
\frac{1}{t \log^2(t)(1 + \abs{y})}\|\phi_1\|_{\mathrm{in}1,t_1}^2
\\
&\le C\frac{\chi_{B_{4R}}}{R(t)^{1+2\gamma} t^{\nu_1}}
\frac{R(t)^{6 - a_1-2\gamma}}{t^{\nu_1}}
\frac{1}{(1 + |y|)^3}\|\phi_1\|_{\mathrm{in}1,t_1}^2
\\
&\le \frac C{t_0^\sigma}\varsigma_2(x, t)\|\phi_1\|_{\mathrm{in}1,t_1}^2
\end{align*}
where again we have used 
\eqref{restriction_csts1}.
We also have
\begin{align*}
|\eta \phi_2|^2 |U_t| &\le C\chi_{B_{4R}}\frac{R(t)^4}{t^{2\nu_2}(1 + |y|)^2}
\frac{1}{t \log^2(t)(1 + \abs{y})}\|\phi_2\|_{\mathrm{in}2,t_1}^2
\\
&\le C\frac{\chi_{B_{4R}}}{t^{\nu_1} R(t)^{1+2\gamma}}
\frac{R(t)^{3-2\gamma}}{t^{1+\nu_1}}
\frac{1}{(1 + \abs{y})^3}\|\phi_2\|_{\mathrm{in}2,t_1}^2
\\
&\le  \frac C{t_0^\sigma}\varsigma_2(x, t)\|\phi_2\|_{\mathrm{in}2,t_1}^2
\end{align*}
where we have used \eqref{restriction_csts2}.
We have a similar estimate for $\phi_3$ and we get
\begin{align}
\label{eq:est-nl4}
\big\||\eta \phi|^2 U_t\big\|_* \le  \frac C{t_0^\sigma}
\left(\|\phi_1\|_{\mathrm{in}1,t_1}^2
+\|\phi_2\|_{\mathrm{in}2,t_1}^2
+\|\phi_3\|_{\mathrm{in}3,t_1}^2\right).
\end{align}

From \eqref{eq:est-nl1}--\eqref{eq:est-nl4} and similar estimate for the remaining
terms we obtain \eqref{eq:est-g4}.

\medskip
Finally, combining \eqref{eq:est-g1}, \eqref{eq:est-g2}, \eqref{eq:est-g3}
and \eqref{eq:est-g4} we deduce \eqref{eq:est-RHS-outer}.
\end{proof}

\begin{lemma}
\label{lem:est-h}
Let $h[p, \xi, \psi,Z]$ be given by \eqref{def_rhs_inner}.
Then
\begin{align}
\label{eq:est-h}
\|h[p, \xi, \psi, Z]\|_{**,a_1, \nu_1 } \le \frac{C}{t_0^\sigma}\big(\|\psi\|_{\mathrm{out}}+\|\dot p\|_{\dot p} + \|\dot \xi\|_{\dot \xi} + \|Z_0\|_{0}\big).    
\end{align}
\end{lemma}
\begin{proof}
We start with $\lambda^2 Q_{-\omega}\tilde{L}_{U}(\psi)\chi_{\D{2}}$. Notice that, for $y \in B_{2R(t)}(0)$ with $t \ge t_0$, we have
$$\xi(t) + \lambda(t) y \in B_{4R(t)}(0) \subset \R^2.$$
By the definition of the norm $\norm{\cdot}_{\text{out}}$ in \eqref{norm_psi}, there holds 
\begin{align*}
|\lambda^2 Q_{-\omega}\tilde{L}_{U}(\psi)(\xi + \lambda y, t)|\chi_{\D{2}}
&\le \frac{C}{(1 + \abs{y})^2} \abs{\nabla \psi(\xi + \lambda y, t))}
\\
&\le  \frac{C}{(1 + \abs{y})^2} \frac{\log t}{R(t)t^{\nu_1}} \|\psi\|_{\text{out}}
\\
&\le \frac{C}{t_0^\sigma}
\frac1{t^{\nu_1}(1 + \abs{y})^{a_1}}
\|\psi\|_{\text{out}}
\end{align*}
where $\sigma>0$ is small and we have used $a_1<3$.
This leads to 
\begin{align}
\label{norm_phi1_tLpsi}
\|\lambda^2 Q_{-\omega}\tilde{L}_{U}(\psi)\chi_{\D{2}}\|_{**,a_1, \nu_1} \le 
\frac{C}{t_0^\sigma}
\|\psi\|_{\text{out}}.
\end{align}

The term $\lambda^2 Q_{-\omega}\tilde{L}_{U}(Z)$ is estimated 
in the same way as the preceding term.
We get that
\begin{align*}
|\lambda^2 Q_{-\omega}\tilde{L}_{U}(Z)(\xi + \lambda y, t)| \chi_{\D{2}}
&\le  \frac{C}{(1 + \abs{y})^2} \frac{\|Z_0\|_{0}}{t^{d + \frac{1}{2}}}
\\
&\le 
\frac{C}{t_0^\sigma}
\frac{C}{t^{\nu_1}(1 + \abs{y})^{a_1}}\|Z_0\|_{0}
\end{align*}
where $\sigma>0$, where we have used \eqref{restriction_nu_a1b}.
This gives us 
\begin{align}
\label{norm_phi1_tLZ}
\|\lambda^2 Q_{-\omega}\tilde{L}_{U}(Z)\chi_{\D{2}}\|_{**,a_1, \nu_1} 
\le \frac{C}{t_0^\sigma}\|Z_0\|_{0} .
\end{align}

To estimate the operator $\lambda^2 Q_{-\omega}\mathcal{K}_0[p, \xi]$ we use Lemma~\ref{lem:estK0} to get
\begin{align*}
|\lambda^2 Q_{-\omega}\mathcal{K}_{0}[p, \xi]| &\le 
C \frac{1}{(1 + |y|)^3} 
\frac{1}{\log(t_0)}\frac1{t(1+\frac{|y|}{\sqrt t})^2}\|\dot p\|_{\dot p}
\\
&\leq 
\frac{C}{t_0^\sigma} \frac{1}{t^{\nu_1}(1 + |y|)^{a_1}} \|\dot p\|_{\dot p}
\end{align*}
where we have used that $\nu_1<1$.
This leads to 
\begin{align}
\label{norm_phi1_K0}
\|\lambda^2 Q_{-\omega}\mathcal{K}_{0}[p, \xi]\|_{**,a_1, \nu_1} \le \frac{C}{t_0^\sigma} \|\dot p\|_{\dot p}.
\end{align}

Finally, the last term $\lambda^2 Q_{-\omega}\mathcal{E}_1[p, \xi]\chi_{\D{2}}$ \eqref{eq:K1} satisfies
\begin{align*}
|\lambda^2 Q_{-\omega}\mathcal{E}_1[p, \xi](y, t)| \chi_{\D{2}}
&\le \frac{C}{(1 + \abs{y})^2} |\dot \xi(t)|
\\
&\le  \frac{C}{(1 + \abs{y})^2} 
\frac{\|\dot\xi\|_{\dot\xi}}{R(t)t^{\nu_1 - \theta}}
\\
&\le \frac{C}{t_0^\sigma}\frac{1}{t^{\nu_1 }(1 + \abs{y})^{a_1}}\|\dot\xi\|_{\dot\xi}
\end{align*}
for some $\sigma>0$ (we use $\theta<\beta$).
This leads to
\begin{align}
\label{norm_phi1_K1}
\|\lambda^2 Q_{-\omega}\mathcal{E}_1[p, \xi]\chi_{\D{2}}\|_{**,a_1, \nu_1} \le \frac{C}{t_0^\sigma}\|\dot\xi\|_{\dot\xi}.
\end{align}

From estimates \eqref{norm_phi1_tLpsi}, \eqref{norm_phi1_tLZ}, \eqref{norm_phi1_K0} and \eqref{norm_phi1_K1}, we deduce the validity of \eqref{eq:est-h}.
\end{proof}

For $t_1>t_0$ let
\[ 
B(t_1) = \|\psi\|_{\mathrm{out},t_1} + \|\phi_1\|_{\mathrm{in}1,t_1} + \|\phi_2\|_{\mathrm{in}2,t_1} 
+ \|\phi_3\|_{\mathrm{in}3,t_1} + \|\dot p\|_{\dot p,t_1} + \|\dot\xi\|_{\dot\xi,t_1}. 
\]

\begin{proposition}
\label{prop:main}
If $t_0$ is sufficiently large, there is $\delta>0$ and $C_{t_0}$ so that if $\|Z_0\|_0<\delta$ 
and $t_1>t_0$ is such that  $B(t_1)\leq 1$, then 
\[ 
B(t_1) \le C_{t_0}\|Z_0\|_0. 
\]
\end{proposition}
\begin{proof}
In the following computations we work with $t\leq t_1$. Since $B(t_1)\le 1$ and $t_0$ is sufficiently large, the bounds in \eqref{assumption_param_outer} hold up to $t_1$.

Let $a_0[p,\xi,\psi,Z]$ denote the right-hand side of \eqref{eq:p-main},
which is defined in \eqref{def_a0}, \eqref{eq:def-a0j}.
By definition of the norm $\|\ \|_{\mathrm{out}}$ \eqref{norm_psi}
and the estimate \eqref{eq:decay-tLUZ-before}
\begin{align*}
|a_0[p,\xi,\psi,Z](t)| \leq C\frac{\log t}{ R(t) t^{\nu_1}}(\|\psi\|_{\mathrm{out},t_1}+\|Z_0\|_{0})
\end{align*}
and similarly
\begin{align*}
|a_0[p,\xi,\psi,Z](t) - a_0[p,\xi,\psi,Z](s)| \leq C \frac{\log t}{R(t) t^{\nu_1}} \frac{|t - s|^\gamma}{R(t)^{2\gamma}}
(\|\psi\|_{\mathrm{out},t_1}+\|Z_0\|_{0}).
\end{align*}
where $0<\gamma<1$.
By Corollary~\ref{cor:localized-p1v2}, the solution $p$ of \eqref{eq:p-main} satisfies
\begin{align*}
|p(t)-1|
& \leq C \frac{1}{t_0^{\nu_1+\beta-1}} \frac{\log^2(t_0)}{\log(t)}
(\|\psi\|_{\mathrm{out},t_1}+\|Z_0\|_{0})
\\
|\dot p(t)|
& \leq C \frac{1}{t_0^{\nu_1+\beta-1}} \frac{\log^2(t_0)}{t\log^2(t)}
(\|\psi\|_{\mathrm{out},t_1}+\|Z_0\|_{0}),
\end{align*}
where we have assumed $\nu_1+\beta>1$ \eqref{eq:d1_plus_beta}.
We deduce that 
\begin{align}
\label{eq:est-dotp-t1}
\|\dot p\|_{\dot p,t_1}
\leq
C \frac{1}{t_0^\sigma}
(\|\psi\|_{\mathrm{out},t_1}+\|Z_0\|_{0}),
\end{align}
for some $\sigma>0$.

Recall $a_1[p,\xi,\psi,Z]$ defined in \eqref{eq:def-a1}. 
By definition of the norm $\|\ \|_{\mathrm{out}}$ \eqref{norm_psi}
and the estimate \eqref{eq:decay-tLUZ-before}
\begin{align*}
|a_1[p,\xi,\psi,Z](t)| \leq C\frac{\log(t)}{ R(t) t^{\nu_1}}(\|\psi\|_{\mathrm{out},t_1}+\|Z_0\|_{0}).
\end{align*}
From the ODE for $\xi$ \eqref{eq:xi-main} we get
\begin{align}
\label{eq:est-dotxi-t1}
\|\dot \xi\|_{\dot \xi,t_1}
\leq (\|\psi\|_{\mathrm{out},t_1}+\|Z_0\|_{0}).
\end{align}

Using Proposition \ref{prop:p2} we also find that
\begin{align}
\nonumber
|\dot p(t)-\dot p(s)|
\leq \left[ C \frac{\log^2(t_0)}{t_0^{\nu-1}}\frac{t-s}{t^2 \log^2(t)}
+ C \frac{(t-s)^\gamma}{t^{\nu_1+\beta+2\gamma}}\right]
(\|\psi\|_{\mathrm{out},t_1}+\|Z_0\|_{0}).
\end{align}
Assuming \eqref{eq:param1}
we get
\begin{align}
\label{eq:est-hatRalphav2}
|\widetilde R_{\alpha j} [p] | \leq \frac C{t_0^\sigma} \frac{1}{R(t)t^{\nu_1}}
(\|\psi\|_{\mathrm{out},t_1}+\|Z_0\|_{0}).
\end{align}
for some $\sigma>0$.
Let $h_2=h_{2,a}+h_{2,b}$ be the right-hand side of equation \eqref{eq_phi2_simp} 
$h_{2,a}$ is defined in \eqref{eq:def-h2a} and $h_{2,b}$ is defined in \eqref{eq:def-h2b}. 
Then \eqref{eq:est-hatRalphav2} implies that 
\begin{align*}
|h_{2,a}|\leq 
\frac C{t_0^\sigma}
\frac{1}{R(t)t^{\nu_1}} \frac{1}{(1 + \rho)^{4}} 
(\|\psi\|_{\mathrm{out},t_1}+\|Z_0\|_{0}).
\end{align*}
Therefore
\[
\|h_{2,a}\|_{**,a_2,\nu_2;t_1}
\leq \frac C{t_0^\sigma}
(\|\psi\|_{\mathrm{out},t_1}+\|Z_0\|_{0}).
\]
Recall that $h_{2,b}$ is given by $h_{2,b}=(1-\eta_{3t_0}(t)) \sum_{j= 1, 2} 
c_{0j}\left[   \lambda^2\big(\tilde L_U(\psi + Z) + \mathcal K_0[p] \big)  \right]
w_\rho^2Z_{0j}$, so it is supported on $[t_0,4t_0]$.
By Proposition~\ref{prop:main1}, we deduce that 
\[
\|h_{2,b}\|_{**,a_2,\nu_2;t_1}
\leq C_{t_0} \|Z_0\|_0.
\]
Therefore
\[
\|h_{2}\|_{**,a_2,\nu_2;t_1}
\leq 
\frac C{t_0^\sigma} \|\psi\|_{\mathrm{out},t_1}
+C_{t_0} \|Z_0\|_0.
\]
By Proposition~\ref{Prop:phi2} we obtain
\begin{align}
\label{eq:est-phi2-t1}
\|\phi_2\|_{\mathrm{in}2,t_1}
\leq \frac C{t_0^\sigma} \|\psi\|_{\mathrm{out},t_1}
+C_{t_0} \|Z_0\|_0.
\end{align}

The proof of Lemma~\ref{lem:est-h}, 
localizing to $[t_0,t_1]$,
gives
\begin{align}
\nonumber
\|h[p, \xi, \psi, Z]\|_{**,a_1, \nu_1;t_1 } \le \frac{C}{t_0^\sigma}\big(\|\psi\|_{\mathrm{out},t_1}+\|\dot p\|_{\dot p,t_1} + \|\dot \xi\|_{\dot \xi,t_1} + \|Z_0\|_{0}\big).    
\end{align}
Recalling $c_{0j}$ \eqref{def:ckj} we find that 
\begin{align*}
|c_{0j} (h[p, \xi, \psi, Z])(t)|\leq \frac{C}{t_0^\sigma} \frac 1{t^{\nu_1}}
\big(\|\psi\|_{\mathrm{out},t_1}+\|\dot p\|_{\dot p,t_1} + \|\dot \xi\|_{\dot \xi,t_1} + \|Z_0\|_{0}\big)
\end{align*}
By Lemma~\ref{lem:est-h} and
Proposition~\ref{Prop:phi1},
\begin{align}
\label{eq:est-phi1-t1}
\|\phi_1\|_{\mathrm{in}1,t_1} \leq \frac{C}{t_0^\sigma}\big(\|\psi\|_{\mathrm{out}}+\|\dot p\|_{\dot p} + \|\dot \xi\|_{\dot \xi} + \|Z_0\|_{0}\big).    
\end{align}

The term $c_{-1 j} (h[p, \xi, \psi, Z])w_\rho^2Z_{-1 j}$ appears in the right-hand side of \eqref{eq:phi3}.
Recalling the definition of  $c_{-1j}$ defined in \eqref{def:ckj}, 
and using $\|\ \|_{\mathrm{out}}$ \eqref{norm_psi}
and the estimate \eqref{eq:decay-tLUZ-before}, we get
\begin{align*}
\left| c_{-1 j} (h[p, \xi, \psi, Z](t)) \right| 
\leq C\frac{\log^2(t)}{ R(t) t^{\nu_1}}(\|\psi\|_{\mathrm{out},t_1}+\|Z_0\|_{0}).
\end{align*}
Since $\nu_1+\beta>1$ we obtain
\[
\left\| c_{-1 j} (h[p, \xi, \psi, Z])w_\rho^2Z_{-1 j} \right\|_{**,a_3,\nu_3;t_1}
\leq \frac C{t_0^\sigma} (\|\psi\|_{\mathrm{out},t_1}+\|Z_0\|_{0}).
\]
Therefore, by Proposition~\ref{Prop:phi3} we get
\begin{align}
\label{eq:est-phi3-t1}
\|\phi_3\|_{\mathrm{in}3,t_1}
\leq \frac{C}{t_0^\sigma}(\|\psi\|_{\mathrm{out},t_1}+\|Z_0\|_{0}).
\end{align}
for some $\sigma>0$.
Collecting \eqref{eq:est-dotp-t1}, \eqref{eq:est-dotxi-t1}, 
\eqref{eq:est-phi2-t1}, \eqref{eq:est-phi1-t1}, \eqref{eq:est-phi3-t1},
we get
\begin{align}
\nonumber
\norm{g[\psi, \phi, p, \xi, Z]}_{*,t_1} 
&\leq
\frac{C}{t_0^\sigma}
\|\psi\|_{\mathrm{out},t_1}
+C_{t_0} \norm{Z_0}_{0}
+\frac{C}{t_0^\sigma} \|\psi\|_{\mathrm{out},t_1}^2
\end{align}
By Proposition~\ref{prop:Estimate_HE_gen}, we find that
\[
\|\psi\|_{\mathrm{out},t_1} \leq 
\frac{C}{t_0^\sigma}
\|\psi\|_{\mathrm{out},t_1}
+\frac{C}{t_0^\sigma} \|\psi\|_{\mathrm{out},t_1}^2
+C_{t_0} \norm{Z_0}_{0}.
\]
A standard argument now shows that if $t_0$ is fixed large, and then $\|Z_0\|_0<\delta$ with $\delta$ sufficiently small,  $\|\psi\|_{\mathrm{out},t_1}\leq 1$, and moreover
\[
\|\psi\|_{\mathrm{out},t_1} \leq 
C_{t_0} \norm{Z_0}_{0}.
\]
Substituting this estimate into the previous bounds for $\phi_1,\phi_2,\phi_3,\dot p,\dot\xi$ yields 
\[ 
\|\phi_1\|_{\mathrm{in}1,t_1} + \|\phi_2\|_{\mathrm{in}2,t_1} 
+ 
\|\phi_3\|_{\mathrm{in}3,t_1} + \|\dot p\|_{\dot p,t_1} + \|\dot\xi\|_{\dot\xi,t_1} \le C_{t_0}\|Z_0\|_0. 
\]
\end{proof}

\begin{proof}[Proof of Theorem~\ref{thm:main}]
Standard arguments
give a solution of the system \eqref{eq:phi1}--\eqref{eq:xi-main} on a maximal interval 
$[t_0,T_{\max})$. 
For $\|Z_0\|_0$ sufficiently small, and $t$ sufficiently close to $t_0$ we have
$B(t)\le \frac12$. 
We claim that for all $t_1\in(t_0,T_{\max})$, $B(t_1)<1$. 
By contradiction, assuming this is false,
let $T\in(t_0,T_{\max})$ be the first time such that  $B(T)=1$. 
Then Proposition~\ref{prop:main} implies 
\[ 
B(T) \le C_{t_0}\|Z_0\|_0. 
\] 
Choosing $\delta$ sufficiently small we obtain 
\[ 
B(T)<\frac12,
\] 
which is not possible.
This proves the claim.

The above argument shows that $T_{\max}=+\infty$. 
Moreover, 
\[ |\dot p(t)|\leq \frac{C\|Z_0\|_0}{t\log^2t}, \qquad |\dot\xi(t)|\leq \frac{C\|Z_0\|_0}{t^{\nu_1+\beta}}. 
\] 
Since $\nu_1+\beta>1$, both right-hand sides are integrable on $(t_0,\infty)$. Hence there exist limits \[ p(t)\longrightarrow p_\infty, \qquad \xi(t)\longrightarrow\xi_\infty. \] Writing $p_\infty=\lambda_\infty e^{i\omega_\infty}$ and using $|p(t)-1|<\frac12$, we obtain \[ \lambda(t)\longrightarrow\lambda_\infty, \qquad \omega(t)\longrightarrow\omega_\infty. \] Finally, the estimates for $\phi_1,\phi_2,\phi_3,\psi,Z$, and $\Phi^0$, together with the decomposition of $u$, imply the decay estimate stated in Theorem~\ref{thm:main}.

\end{proof}

\begin{remark}
\label{rem:1/2}
To adapt the proof of Theorem~\ref{thm:main} to the case $d>\frac12$ the only change needed is to have Proposition~\ref{prop:p1v2} valid for all $\frac12<\alpha<1$. This is possible and explained in Section~\ref{sec:eq-p}, Remark~\ref{rem:adapt-prop:p1v2}. There is also a need to change the parameters, but this is possible, see Lemma~\ref{lem:par1}.
\end{remark}

\subsection{Restrictions on the parameters}
\label{Sect:Choice_of_cst01}

A summary of the parameters is in Section \ref{sect:summary-param}. Here are the collected restrictions on them.
\begin{itemize}
\item
We recall from Section \ref{sect:summary-param} that $d>\frac12$, $\beta\in(\frac14,\frac12)$ and $\alpha\in(\frac12,1)$.
\item 
In order to apply Proposition~\ref{prop:Estimate_HE_gen}
we need 
\begin{align}
\label{hnu1}
&0<\nu_1<1 
\\
\label{hnu1beta}
& \nu_1+\beta>1
\\
\label{hnu1beta2}
& \nu_1+\beta<1+\frac{d}{2}.
\end{align}
We also need $\nu_1+\beta>1$ in order to apply Proposition \ref{prop:p1} (see \eqref{eq:d1_plus_beta}).

\item 
To estimate $\tilde L_U(Z)$ we need
\begin{align}
\label{hnubeta3}
\nu_1< d+\tfrac12-\beta(a_1-2),
\end{align}
see \eqref{restriction_nu_a1b}.

\item 
From \eqref{eq:param1}, which has to do with a gain in the decay of $\mathcal R_\alpha[p]$, we need
\begin{align}
\label{Ralpha1}
&2\beta>\alpha
\\
\label{Ralpha2}
& 2-\alpha \geq \nu_1+\beta
\end{align}

\item 
To estimate nonlinear terms we need
\begin{align}
\label{nl1}
& 3\beta <\nu_1 
\\
\label{nl2}
& \beta (6-a_1)<\nu_1 
\end{align}
The conditions \eqref{nl1} and \eqref{nl2} appear in  \eqref{restriction_csts1} and 
\eqref{restriction_csts2}.

\item Condition 
\begin{align}
\label{eq:cond-d-nu1-b-v2}
\nu_1+2\beta<1+d
\end{align}
appears in the analysis of $(\psi\cdot U)U_t$
\eqref{eq:cond-d-nu1-b}.
\end{itemize}

\begin{lemma}
\label{lem:par1}
For any $d>\frac12$ there are $\beta \in (\frac14,\frac12)$, $\alpha\in(\frac12,1)$,
$\nu_1$ and $a_1\in(2,3)$ such that
\eqref{hnu1}--\eqref{eq:cond-d-nu1-b-v2} hold.
\end{lemma}
\begin{proof}
We first look for conditions so that $\nu_1$ can be found. 
Looking at a lower bound for $\nu_1$ in 
\eqref{hnu1}, \eqref{hnu1beta}, \eqref{nl1}, \eqref{nl2} we see that we need
\[
\nu_1>\max(0,1-\beta,3\beta,\beta(6-a_1)).
\]
Note that $3\beta<\beta(6-a_1)$ if $a_1\in(2,3)$. Also, $1-\beta<\beta(6-a_1)$
for  $\beta \in (\frac14,\frac12)$, $a_1\in(2,3)$.
Thus
\[
\max(0,1-\beta,3\beta,\beta(6-a_1))=
\beta(6-a_1).
\]
An upper bound for $\nu_1$ in \eqref{hnu1}, \eqref{hnu1beta2}, \eqref{hnubeta3}, \eqref{Ralpha2}, and \eqref{eq:cond-d-nu1-b-v2} is
\[
\min(1,1+\frac{d}{2}-\beta,d+\tfrac12-\beta(a_1-2),2-\alpha-\beta,1+d-2\beta).
\]
To find a nonempty interval for $\nu_1$, it suffices to require
\begin{align}
\label{ineq1}
\beta(6-a_1)<\min\left(1,1+\tfrac{d}{2}-\beta,d+\tfrac12-\beta(a_1-2),2-\alpha-\beta,1+d-2\beta\right).
\end{align}
But $\beta(6-a_1)<1$ is equivalent to
\begin{align}
\label{beta1}
\beta<\frac1{6-a_1},
\end{align}
$\beta(6-a_1)<1+\frac{d}{2}-\beta$ is equivalent to 
\begin{align}
\label{beta1b}
\beta<\frac{1+\frac d2}{7-a_1}
,
\end{align}
$\beta(6-a_1)<d+\tfrac12-\beta(a_1-2)$ is equivalent to 
\begin{align}
\label{beta2}
\beta<\frac{d+\tfrac12}{4},
\end{align}
and $\beta(6-a_1)<2-\alpha-\beta$ is equivalent to 
\begin{align}
\label{beta3}
\beta<\frac{2-\alpha}{7-a_1},
\end{align}
and $\beta(6-a_1)<1+d-2\beta$ is equivalent to
\begin{align}
\label{beta4}
\beta<\frac{1+d}{8-a_1}.
\end{align}

Since $d>\frac12$, 
$\frac{d+\tfrac12}{4}>\frac14$.
We fix any $a_1\in(2,3)$.
Then note that 
\[
\frac{1+\frac d2}{7-a_1}>\frac{1}{4},\quad
\frac{1+d}{8-a_1}>\frac14.
\]
Select $\beta>\frac14$ such that \eqref{beta1b}, \eqref{beta2} and \eqref{beta4} hold. Taking $\beta$ close to $\frac14$ and $\alpha\in(\frac12,2\beta)$ we get $\alpha$ close to $\frac12$ and we obtain \eqref{beta3} and \eqref{beta1},
since $\frac{2-\alpha}{7-a_1}$ would be close
to $\frac{\frac32}{7-a_1}>\frac3{10}>\frac14$.

Now choose $\nu_1$ such that
\begin{align}
\label{eq:ineq-nu1b}
\beta(6-a_1)<\nu_1<\min\left(1,1+\tfrac{d}{2}-\beta,d+\tfrac12-\beta(a_1-2),2-\alpha-\beta,1+d-2\beta\right)
\end{align}
which is possible by \eqref{ineq1}. Then \eqref{hnu1}, \eqref{hnu1beta}, \eqref{hnu1beta2}, \eqref{hnubeta3}, \eqref{Ralpha2}, \eqref{nl1}, \eqref{nl2}, and \eqref{eq:cond-d-nu1-b-v2} hold by construction, while \eqref{Ralpha1} holds by the choice of $\alpha$.
\end{proof}

Note that there are more parameters involved (see Section~\ref{sect:summary-param}), but there are no additional restrictions for them. These are $a_2\in (2,3)$, $\gamma\in(0,\frac12)$.

As explained earlier, 
to avoid longer computations in Section~\ref{sec:nonlocal}, it is convenient to assume that $\alpha>\frac58$.
For this we need to assume that $d>\frac34$.

\begin{lemma}
For any $d>\frac34$ there are $\beta \in (\frac5{16},\frac12)$, $\alpha\in(\frac58,1)$,
$\nu_1$ and $a_1\in(2,3)$ such that
\eqref{hnu1}--\eqref{eq:cond-d-nu1-b-v2} hold.
\end{lemma}
\begin{proof}
The proof is similar to that of Lemma~\ref{lem:par1}.

We can assume that $\frac34<d<\frac{15}{16}$.
Choose $a_1<3$ close to $3$ so that
\begin{align}
\label{eq:a1}
\frac{1}{6-a_1} > \frac5{16}
,\quad 
\frac{11}{8}\frac1{7-a_1}> \frac5{16}.
\end{align}
We let $\alpha>\frac58$ be sufficiently close to $\frac58$ such that 
\[
\frac{\alpha}{2}<\min\left(\frac1{6-a_1}, \frac{1+\frac d2}{7-a_1},\frac{d+\tfrac12}{4},
\frac{2-\alpha}{7-a_1},\frac38,\frac{1+d}{8-a_1}\right).
\]
This is possible because of \eqref{eq:a1}.
We also have used that  from \eqref{eq:a1} and $d>\frac34$ we have $\frac{1+\frac d2}{7-a_1}>\frac5{16}$, and 
since $d>\frac34$, we also have $\frac{d+\tfrac12}{4}>\frac5{16}$. 
Now choose $\beta$ such that
\[
\frac{\alpha}{2}<\beta<\min\left(\frac1{6-a_1}, \frac{1+\frac d2}{7-a_1},\frac{d+\tfrac12}{4},
\frac{2-\alpha}{7-a_1},\frac38\right).
\]
These choices guarantee \eqref{beta1}, \eqref{beta1b}, \eqref{beta2}, \eqref{beta3}, \eqref{beta4}.
We then have the inequality \eqref{ineq1} in the proof
of Lemma~\ref{lem:par1} is valid and then it's possible to find $\nu_1$ satisfying \eqref{eq:ineq-nu1b}.
\end{proof}

\subsection{Proof of some of the estimates}
\label{sec:some-estimates}

\begin{proof}[Proof of Proposition~\ref{prop:est-Ralpha}]
Recall that 
$\mathcal R_\alpha=\mathcal R_{\alpha,0}+\mathcal R_{\alpha,1}+\mathcal R_{\alpha,2}$
\eqref{def:Ralpha}
where these terms are defined in \eqref{def:Ralpha0},  \eqref{def:Ralpha1}, and
\eqref{def:Ralpha2}. 
We estimate
\begin{align*}
| \Re{\mathcal R_{\alpha,0}[p]}(t)|
&\leq
\int_{t-t^\alpha}^{t} 
\frac{
|\dot p(t)-\dot p(s)|}{t-s} \Gamma_1\left(\frac{|p(t)|^2}{t - s}\right)
\,d s .
\end{align*}
We use the bounds \eqref{eq:holder-p} to get
\begin{align*}
\int_{t-t^\alpha}^{t} 
\frac{
|\dot p(t)-\dot p(s)|}{t-s} \Gamma_1\left(\frac{|p(t)|^2}{t - s}\right)
\,d s 
&\leq 
I_1+I_2
\end{align*}
where
\begin{align*}
I_1
&=
C \|a_0\|_{\nu,1} \frac{\log^2(t_0)}{t_0^{\nu-1}}
\int_{t-t^\alpha}^t
\frac{t-s}{t^2 \log^2(t)}
\frac{1}{t-s}\,ds
\\
&\leq
C \|a_0\|_{\nu,1} \frac{\log^2(t_0)}{t_0^{\nu-1}}
\frac{1}{t^{2-\alpha} \log^2(t)}
\end{align*}
and
\begin{align*}
I_2
&=
C [a_0]_{\gamma,1,\mu} 
\int_{t-t^\alpha}^t
\frac{(t-s)^\mu}{t^\gamma}
\frac{1}{t-s}\,ds.
\\
&\leq 
C[a]_{\gamma,1,\mu}\frac{t^{\mu \alpha}}{t^\gamma}
=
C[a]_{\gamma,1,\mu}\frac{1}{t^{\gamma-\mu\alpha}}.
\end{align*}
The imaginary part of $\mathcal R_{\alpha,0}[p](t)$ is estimated similarly and we get
\begin{align}
\label{eq:est-Ralpha0}
|\mathcal R_{\alpha,0}[p](t)|
\leq C \|a_0\|_{\nu,1} \frac{\log^2(t_0)}{t_0^{\nu-1}}\frac{1}{t^{2-\alpha} \log^2(t)}
+ C[a]_{\gamma,1,\mu}\frac{1}{t^{\gamma-\mu\alpha}}.
\end{align}

Next we compute, using the bound \eqref{eq:estdotp1} for $\dot p$ and \eqref{eq:bound-gj} for $g_1$ ,
\begin{align*}
| \Re{\mathcal R_{\alpha,1}[p]}(t)|
&\leq |\dot p(t)| g_1\Bigl(\frac{|p(t)|^2}{t^\alpha}\Bigr)
\\
&\leq C \|a_0\|_{\nu,1} \frac{1}{t_0^{\nu-1}} 
\frac{\log^2(t_0)}{t\log^2(t)} \frac{\log(t)}{t^\alpha}
\\
&\leq C \|a_0\|_{\nu,1} \frac{1}{t_0^{\nu-1}} \frac{\log^2(t_0)}{t^{1+\alpha}\log(t)}.
\end{align*}
The imaginary part of $\mathcal R_{\alpha,1}[p](t)$ is estimated similarly
and we get
\begin{align}
\label{eq:est-Ralpha1}
|\mathcal R_{\alpha,1}[p](t)|
\leq C \|a_0\|_{\nu,1} \frac{1}{t_0^{\nu-1}} \frac{\log^2(t_0)}{t^{1+\alpha}\log(t)}.
\end{align}
Finally from the definition \eqref{def:Ralpha2} and  \eqref{eq:estdotp1}
we get
\begin{align*}
| \Re{\mathcal R_{\alpha,2}[p]}(t)|
&\leq 
\int_{t_0}^{t-t^\alpha} \frac{|\dot p(s)|}{t - s}  \left|\Gamma_1\left(\frac{|p(t)|^2}{t - s}\right)-1\right|\, ds.
\end{align*}
Using \eqref{est_Gamma_tau<1} in the form $|\Gamma_1(\tau)-1|\leq C \tau^{1-\sigma}$, $0<\tau<1$,  where $0<\sigma<1$, we obtain
\begin{align*}
| \Re{\mathcal R_{\alpha,2}[p]}(t)|
&\leq 
C \|a_0\|_{\nu,1} \frac{\log^2(t_0)}{t_0^{\nu-1}} 
\int_{t_0}^{t-t^\alpha} \frac{1}{s\log^2(s)}  \frac1{(t-s)^{2-\sigma}}\, ds.
\end{align*}
We split
\begin{align*}
\int_{t_0}^{t-t^\alpha} \frac{1}{s\log^2(s)}  \frac1{(t-s)^{2-\sigma}}\, ds
=\int_{t_0}^{t/2}\ldots + \int_{t/2}^{t-t^\alpha}\ldots
\end{align*}
and estimate
\begin{align*}
\int_{t_0}^{t/2} \frac{1}{s\log^2(s)}  \frac1{(t-s)^{2-\sigma}}\, ds
&\leq \frac{C}{t^{2-\sigma}} 
\\
\int_{t/2}^{t-t^\alpha} \frac{1}{s\log^2(s)}  \frac1{(t-s)^{2-\sigma}}\, ds
&\leq \frac{C}{t^{1+\alpha(1-\sigma)}\log^2(t)}.
\end{align*}
Using that $\alpha>\frac12$ and taking $\sigma>0$ small, we obtain
\begin{align*}
| \Re{\mathcal R_{\alpha,2}[p]}(t)|
&\leq 
C \|a_0\|_{\nu,1} \frac{\log^2(t_0)}{t_0^{\nu-1}}\frac{1}{t^{2-\alpha} \log^2(t)}.
\end{align*}
The imaginary part is estimated similarly, and therefore
\begin{align}
\label{eq:est-Ralpha2}
| \mathcal R_{\alpha,2}[p](t)|
\leq 
C \|a_0\|_{\nu,1} \frac{\log^2(t_0)}{t_0^{\nu-1}}\frac{1}{t^{2-\alpha} \log^2(t)}.
\end{align}

Combining \eqref{eq:est-Ralpha0}, \eqref{eq:est-Ralpha1} and \eqref{eq:est-Ralpha2} 
we get
\begin{align*}
|\mathcal R_\alpha[p](t)|
\leq 
C \|a_0\|_{\nu,1} \frac{\log^2(t_0)}{t_0^{\nu-1}}\frac{1}{t^{2-\alpha} \log^2(t)}
+ C[a]_{\gamma,1,\mu}\frac{1}{t^{\gamma-\mu\alpha}}
\end{align*}
since we assume that $\alpha>\frac12$.
\end{proof}

\begin{proof}[Proof of Lemma~\ref{lem:phi0}]
The expression \eqref{eq:varphi0} can be rewritten as
\[
\varphi^0(z, t) = -\int_{t_0}^t \frac{\dot p(s)}{t-s} K\Bigl(\frac{z^2}{t - s}\Bigr) ds, \quad
K(\zeta)=2\frac{1 - e^{-\frac{\zeta}{4}}}{\zeta}.
\]
Note that $K$ is smooth and bounded on $[0,\infty)$ with $\lim_{\zeta\to0}K(\zeta)=\frac12$ and $\lim_{\zeta\to\infty}\zeta K(\zeta)=2$.
We divide into several cases.

\medskip
\noindent
\textbf{Case $z^2<t-t_0$ and $t<2t_0$.}
\begin{align*}
|\varphi^0(z,t)|
&\leq C\int_{t_0}^{t-z^2} \frac{1}{s\log^2(s)}\frac{1}{t-s}\,ds
+C\int_{t-z^2}^t  \frac{1}{s\log^2(s)}\frac{1}{z^2}\,ds
\\
&\leq 
C \frac{1}{t_0 \log^2(t_0)} \log\Bigl(\frac{t-t_0}{z^2}\Bigr)
+ C  \frac{1}{t_0 \log^2(t_0)} 
\\
&\leq C \frac{1}{t\log(t_0)}
\end{align*}
where we have used that $z\geq\frac12$.

\medskip
\noindent
\textbf{Case $z^2>t-t_0$ and $z^2<\frac t2$.}
\begin{align*}
|\varphi^0(z,t)|
&\leq C\int_{t_0}^{t}  \frac{1}{s\log^2(s)}\frac{1}{z^2}\,ds
\\
&\leq C\int_{t_0}^{t_0+z^2}  \frac{1}{s\log^2(s)}\frac{1}{z^2}\,ds
\\
&\leq C \frac{1}{t_0 \log^2(t_0)}
\\
&\leq C \frac{1}{t\log(t_0)}.
\end{align*}

\medskip
\noindent
\textbf{Case $\frac t2<z^2<t$ and $t<2t_0$.}
Similar to the previous case:
\begin{align*}
|\varphi^0(z,t)|
&\leq C\frac{1}{z^2} \int_{t_0}^{t}  \frac{1}{s\log^2(s)}\,ds
\\
&\leq C\frac{1}{z^2}\int_{t_0}^{2t_0}  \frac{1}{s\log^2(s)}\,ds
\\
&\leq C \frac{1}{z^2 \log^2(t_0)}
\\
&\leq C \frac{1}{t\log(t_0)}
\end{align*}
because $z^2>\frac t2$.

\medskip
\noindent
\textbf{Case $z^2>t$ and $t<2t_0$.}
\begin{align*}
|\varphi^0(z,t)|
&\leq \int_{t_0}^{t}  \frac{1}{s\log^2(s)}\frac{1}{z^2}\,ds
\leq C \frac{1}{z^2 \log(t_0)}.
\end{align*}

\medskip
\noindent
\textbf{Case $z^2<\frac t2$ and $t>2t_0$.}
Note that $z^2<t-t_0$ and so we split
\begin{align}
\label{eq:split-I1-I2}
|\varphi^0(z,t)|
&\leq I_1+I_2
\end{align}
where
\begin{align}
\label{eq:def-I1-I2}
I_1&= \int_{t_0}^{t-z^2} \frac{1}{s\log^2(s)}\frac{1}{t-s}\,ds,
\quad
I_2=\int_{t-z^2}^t  \frac{1}{s\log^2(s)}\frac{1}{z^2}\,ds.
\end{align}
Since $z^2<\frac t2$, 
\begin{align*}
\int_{t_0}^{t-z^2} \frac{1}{s\log^2(s)} \frac{1}{t-s}\,ds
= \int_{t_0}^{t/2} \frac{1}{s\log^2(s)} \frac{1}{t-s}\,ds
+ \int_{t/2}^{t-z^2} \frac{1}{s\log^2(s)} \frac{1}{t-s}\,ds.
\end{align*}
But
\begin{align*}
\int_{t_0}^{t/2} \frac{1}{s\log^2(s)} \frac{1}{t-s}\,ds
&\leq 
\frac Ct\int_{t_0}^{t/2} \frac{1}{s\log^2(s)}\,ds
\leq \frac C{t\log(t_0)}.
\end{align*}
\begin{align*}
\int_{t/2}^{t-z^2} \frac{1}{s\log^2(s)} \frac{1}{t-s}\,ds
&\leq \frac C{t\log^2(t)}\int_{t/2}^{t-z^2}\frac{1}{t-s}\,ds
\leq\frac C{t\log^2(t)} \log\Bigl(\frac{t}{z^2}\Bigr).
\end{align*}
Since $z\geq\frac12$, we get
\[
|I_1|\leq  \frac C{t\log(t_0)}.
\]
Since $z^2<\frac t2$, $t-z^2>\frac t2$ and then
\begin{align*}
I_2&\leq C \frac{1}{t\log^2(t)} \int_{t-z^2}^t \frac{1}{z^2}\,ds
\leq C \frac{1}{t \log(t_0)}.
\end{align*}

\medskip
\noindent
\textbf{Case $\frac t2<z^2<t-t_0$.}
Note that $t>2t_0$. We split as in \eqref{eq:split-I1-I2} with
$I_1$ and $I_2$ given in \eqref{eq:def-I1-I2}.
Since $z^2>\frac t2$, we have
\begin{align*}
I_1
\leq 
\int_{t_0}^{t-z^2} \frac{1}{s\log^2(s)} \frac{1}{t-s}\,ds
\leq
\frac Ct\int_{t_0}^{\infty} \frac{1}{s\log^2(s)}\,ds
\leq \frac C{t\log(t_0)}.
\end{align*}
For $I_2$:
\begin{align*}
I_2&\leq 
C \frac{1}{z^2}\int_{t/2}^\infty \frac{1}{s\log^2(s)}\,ds
\leq C \frac{1}{z^2 \log(t_0)}.
\end{align*}

\medskip
\noindent
\textbf{Case $t-t_0<z^2<t$ and $t>2t_0$.}
\begin{align*}
|\varphi^0(z,t)|
&\leq C \frac1{z^2} \int_{t_0}^\infty \frac{1}{s\log^2(s)}\,ds
\leq C \frac{1}{t \log(t_0)}
\end{align*}
because $t_0<\frac t2$, so $z^2>t-t_0>\frac t2$.

\medskip
\noindent
\textbf{Case $z^2>t$ and $t>2t_0$.}
\begin{align*}
|\varphi^0(z,t)|
&\leq C \frac1{z^2} \int_{t_0}^\infty \frac{1}{s\log^2(s)}\,ds
\leq C \frac{1}{z^2 \log(t_0)}.
\end{align*}

\end{proof}

\begin{proof}[Proof of Lemma~\ref{lem:estK0}]
Let us start with $\mathcal K_{01}$.
From the formula \eqref{eq:K01}
\begin{align*}
|\mathcal{K}_{01}[p, \xi](y,t)| 
&\le \frac{C}{(1 + \abs{y})^3} \int_{t_0}^t |\dot p(s)|  k(z, t - s)\, ds
\\
&\leq \frac{C}{(1 + \abs{y})^3} \int_{t_0}^t \frac{1}{s \log^2(s)}k(z, t - s)\, ds.
\end{align*}
The integral term was estimated in Lemma~\ref{lem:phi0}, which yields
\begin{align*}
|\mathcal{K}_{01}[p, \xi](y,t)| 
&\le C \frac{1}{(1 + \abs{y})^3} 
\frac{1}{\log(t_0)}\frac1{t(1+\frac{|y|}{\sqrt t})^2}.
\end{align*}
The estimate of $\mathcal{K}_{02}[p, \xi]$, defined in \eqref{eq:K02}, follows a very similar computation. It suffices to notice that the main terms appearing in the expression of $\mathcal{K}_{02}[p, \xi]$ are of the form
\[
\int_{t_0}^t \dot p(s)  zk_z(z, t - s) ds \quad \text{and} \quad \int_{t_0}^t \dot p(s)  z^2k_{zz}(z, t - s) ds,
\]
and the kernels $zk_z$ and $z^2k_{zz}$ the same bounds as the kernel $k$:
\begin{align}
\nonumber
\abs{k(z, t)} \le  C \begin{cases}
\frac{1}{t} & \text{if } z^2 < t,\\
\frac{1}{z^2} & \text{if } z^2 > t.
\end{cases} 
\end{align}
\end{proof}

\section{Estimates for the non-local problem}
\label{sec:nonlocal}
In this section we give the proofs of 
Proposition~\ref{prop:p1} and 
\ref{prop:p2}.
We recall the equation given in \eqref{prob12}
\begin{align}
\label{prob12b}
\left\{
\begin{aligned}
&\dot p(t)\log(t^\alpha)
+H(p(t))\dot p(t)
+\eta(t)\int_{t_0}^{t-t^\alpha}\frac{\dot p(s)}{t-s}
\,d s = f(t)
\quad t>t_0
\\
&
p(t_0)=1
\end{aligned}
\right.
\end{align}
where we have changed a bit the notation: $\eta=\eta_{2t_0}$ (defined in \eqref{def:eta2t0})
and $f=a_0$.
Here $H$ is smooth, bounded with bounded gradient in $\R^2$ and with values in $2\times 2$ matrices. It was defined
in Section~\ref{sec:decomp-nonlocal}.
\subsection{Proof of Proposition~\ref{prop:p1}}
\label{sec:eq-p}
We will obtain Proposition~\ref{prop:p1} as a consequence of the
following variant, where the norm is the one defined in \eqref{def:norm-gamma-m}.
\begin{proposition}
\label{prop:p1v2}
Assume $\alpha \in (\frac58,1)$ and $-4<m<-\frac1{2\alpha-1}$.
For $t_0$ sufficiently large, if  $\|f\|_{1,m}\leq \log^{-3-m}(t_0)$,
then the solution $p$ of 
\eqref{prob12b} satisfies
\begin{align}
\label{eq:estp1xx}
|p(t)-1| &\leq C \|f\|_{1,m}   \frac{\log^{m+1}(t_0)}{\log(t)}
\\
\label{eq:estdotp1xx}
|\dot p(t)|&\leq C  \|f\|_{1,m}  \frac{\log^{m+1}(t_0)}{t\log^2(t)}.
\end{align}
\end{proposition}

The restriction $\alpha>\frac58$ can be relaxed to $\alpha>\frac12$, with a
more negative  $m$ and a longer proof.
We choose this regime because it allows the simplest proof (largest $m$) 
such that $\alpha<\frac23$.
The restriction $\|f\|_{1,m}\leq \log^{-3-m}(t_0)$ is so that $t_0$ can be chosen independently of $f$. If  $\|f\|_{1,m}<\infty$, the conclusion is still true for $t_0$ depending on $f$. 
Using Proposition~\ref{prop:p1v2} we can give a proof of Proposition~\ref{prop:p1}.

\begin{proof}[Proof of Proposition~\ref{prop:p1}]
Note that if $\nu>1$ then
\[
\|f\|_{1,m}\leq \frac{\log^{1-m}(t_0)}{t_0^{\nu-1}}\|f\|_{\nu,1}.
\]
So, if $f$ satisfies $\|f\|_{\nu,1}\leq t_0^{\nu-1}\log^{-4}(t_0)$ then $\|f\|_{1,m}\leq \log^{-3-m}(t_0)$ and we are able to apply Proposition~\ref{prop:p1v2} and find $p$ satisfying
\eqref{prob12b} and \eqref{eq:estp1xx} and \eqref{eq:estdotp1xx}.
These estimates imply  \eqref{eq:estp1} and \eqref{eq:estdotp1}.
\end{proof}

Let
\begin{align}
\label{def:M.2}
M[p](t)=
\dot p(t)
\log(t^\alpha) 
+\eta(t)\int_{t_0}^{t-t^\alpha}\frac{\dot p(s)}{t-s}
\,d s 
\end{align}
and
\[
\mathcal S[p](t) = M[p](t) + H(p(t)) \dot p(t) - f(t)
\]
so the equation in \eqref{prob12b} becomes
$\mathcal S[p] =0$.

\begin{lemma}
We have the following formula:
\begin{align*}
M[p](t)
&=
M_0[p](t)
-\frac{\eta(t)}{t}p(t_0)
+M_1[p](t)
\end{align*}
where
\begin{align}
\nonumber
M_0[p](t)
&=
\dot p(t)\log(t)[\alpha (1-\eta(t))+\eta(t)]
+\frac{\eta(t)}{t}p(t)
\\
\nonumber
M_1[p](t)
&=
\eta(t)\Biggl[
\frac{p(t/2)-p(t)}{t}
+\int_{t_0}^{t/2}\dot p(s)\Bigl(\frac{1}{t-s}-
\frac1t\Bigr)\,d s
- \dot p(t)\log(2)
\\
\label{def:M1.2}
&\quad 
-\int_{t/2}^{t-t^\alpha}\frac{\dot p(t)-\dot p(s)}{t-s}\,d s 
\Biggr].
\end{align}
\end{lemma}
\begin{proof}
Let $t_0\ll\tilde t \leq  t-1$ (with the idea to take later $t=t-t^\alpha$ or $\tilde t = \frac t2$).
Let $p$ be differentiable.
Then
\begin{align*}
\int_{t_0}^{t-1}\frac{\dot p(s)}{t-s}\,d s
&=
\int_{t_0}^{\tilde t}\frac{\dot p(s)}{t-s}\,d s
+\int_{\tilde t}^{t-1}\frac{\dot p(s)}{t-s}\,d s
\\
&=
\int_{t_0}^{\tilde t}\frac{\dot p(s)}{t-s}\,d s
+ \dot p(t)\int_{\tilde t}^{t-1}\frac{1}{t-s}\,d s
-\int_{\tilde t}^{t-1}\frac{\dot p(t)-\dot p(s)}{t-s}\,d s
\\
&=\int_{t_0}^{\tilde t}\frac{\dot p(s)}{t-s}\,d s
+ \dot p(t)\log(t-\tilde t)
-\int_{\tilde t}^{t-1}\frac{\dot p(t)-\dot p(s)}{t-s}\,d s .
\end{align*}
Let's rewrite this as
\begin{align}
\label{ident1}
\int_{t_0}^{\tilde t}\frac{\dot p(s)}{t-s}\,d s
+ \dot p(t)\log(t-\tilde t)
=
\int_{t_0}^{t-1}\frac{\dot p(s)}{t-s}\,d s
+\int_{\tilde t}^{t-1}\frac{\dot p(t)-\dot p(s)}{t-s}\,d s .
\end{align}
Using \eqref{ident1} with $\tilde t =t-t^\alpha$ and then with $\tilde t=t/2$, we get
\begin{align*}
\int_{t_0}^{t-t^\alpha}\frac{\dot p(s)}{t-s}\,d s
+ \dot p(t)\log(t^\alpha)
&=
\int_{t_0}^{t-1}\frac{\dot p(s)}{t-s}\,d s
+\int_{t-t^\alpha}^{t-1}\frac{\dot p(t)-\dot p(s)}{t-s}\,d s 
\\
&=
\int_{t_0}^{t/2}\frac{\dot p(s)}{t-s}\,d s
+ \dot p(t)\log(t/2)
-\int_{t/2}^{t-1}\frac{\dot p(t)-\dot p(s)}{t-s}\,d s .
\\
&\quad 
+\int_{t-t^\alpha}^{t-1}\frac{\dot p(t)-\dot p(s)}{t-s}\,d s 
\\
&=
\int_{t_0}^{t/2}\frac{\dot p(s)}{t-s}\,d s
+ \dot p(t)\log(t/2)
-\int_{t/2}^{t-t^\alpha}\frac{\dot p(t)-\dot p(s)}{t-s}\,d s .
\end{align*}
But
\begin{align*}
\int_{t_0}^{t/2}\frac{\dot p(s)}{t-s}\,d s
&=
\frac1t\int_{t_0}^{t/2}\dot p(s)\,d s
+\int_{t_0}^{t/2}\dot p(s)\Bigl(\frac{1}{t-s}-
\frac1t\Bigr)\,d s
\\
&=\frac{p(t/2)-p(t_0)}{t}
+\int_{t_0}^{t/2}\dot p(s)\Bigl(\frac{1}{t-s}-
\frac1t\Bigr)\,d s.
\end{align*}
Thus we arrive at
\begin{align*}
& \int_{t_0}^{t-t^\alpha}\frac{\dot p(s)}{t-s}\,d s
+ \dot p(t)\log(t^\alpha)
\\
&=
\frac{p(t/2)-p(t_0)}{t}
+\int_{t_0}^{t/2}\dot p(s)\Bigl(\frac{1}{t-s}-
\frac1t\Bigr)\,d s
\\
&\quad
+ \dot p(t)\log(t/2)
-\int_{t/2}^{t-t^\alpha}\frac{\dot p(t)-\dot p(s)}{t-s}\,d s 
\\
&=
\dot p(t)\log(t)
+\frac{p(t)}{t}
+\frac{p(t/2)-p(t)}{t}
-\frac{p(t_0)}{t}
+\int_{t_0}^{t/2}\dot p(s)\Bigl(\frac{1}{t-s}-
\frac1t\Bigr)\,d s
\\
&\quad
- \dot p(t)\log(2)
-\int_{t/2}^{t-t^\alpha}\frac{\dot p(t)-\dot p(s)}{t-s}\,d s .
\end{align*}
Therefore
\begin{align*}
&
M[p](t)
\\
&=
(1-\eta(t))\dot p(t)\log(t^\alpha) 
+\eta(t)\Biggl[
\dot p(t)\log(t)
+\frac{p(t)}{t}
+\frac{p(t/2)-p(t)}{t}
-\frac{p(t_0)}{t}
\\
&\quad
+\int_{t_0}^{t/2}\dot p(s)\Bigl(\frac{1}{t-s}-
\frac1t\Bigr)\,d s
- \dot p(t)\log(2)
-\int_{t/2}^{t-t^\alpha}\frac{\dot p(t)-\dot p(s)}{t-s}\,d s 
\Biggr]
\\
&=
\dot p(t)\log(t)[\alpha (1-\eta(t))+\eta(t)]
+\frac{\eta(t)}{t}p(t)
-\frac{\eta(t)}{t}p(t_0)
+M_1[p](t)
\\
&=M_0[p](t)
-\frac{\eta(t)}{t}p(t_0)
+M_1[p](t).
\end{align*}
\end{proof}

We look for $p$ solving \eqref{prob12b} of the form
\[
p=p_0+p_1+p_2+p_3 .
\]
We let $p_0$ solve
\begin{align}
\label{eq:p0}
\left\{
\begin{aligned}
& M_0[p_0] = (a+b) \frac{\eta(t)}{t} , \quad t>t_0
\\
&p_0(t_0)=b
\end{aligned}
\right.
\end{align}
where $a$, $b\in\C$ are parameters.
This is a linear differential equation, so $a$, $b$ can be taken arbitrary, but we expect $a$ and $b$ to be bounded by a fixed constant: $|a|\leq2$, $|b|\leq 2$.

\begin{lemma}
\label{lem:p0}
Let $a$, $b\in\C$ and $p_0$ be the solution to \eqref{eq:p0}.
Then
\begin{align}
\nonumber
|p_0(t)-a-b|&\leq C |a| \frac{\log(t_0)}{\log(t)}
\\
\label{est:dp0.1}
|\dot p_0(t)|&\leq C |a| \frac{\log(t_0)}{t\log^2(t)}
\\
\label{est:dkp0.1}
| p_0^{(k)}(t)|&\leq C_k |a| \frac{\log(t_0)}{t^k\log^2(t)},\quad k=2,3,
\end{align}
Moreover, writing $p_0(t)=p_0(t;a,b)$, we have
\begin{align*}
|p_0(t;a_1,b)-p_0(t;a_2,b)|
&\leq C |a_1-a_2| \frac{\log(t_0)}{\log(t)}
\\
p_0(t;a,b_1)-p_0(t;a,b_2)&=b_1-b_2
\end{align*}
\end{lemma}
\begin{proof}
Note that 
\[
p_0(t;a,b)=a q_0(t) + b
\]
where $q_0$ is the solution to 
\begin{align}
\label{ode-q0}
\left\{
\begin{aligned}
&\dot q_0(t)\log(t)[\alpha(1-\eta(t)) + \eta(t)]+\frac{\eta(t)}{t}q_0(t)=\frac{\eta(t)}t,\quad t>t_0
\\
& q_0(t_0) = 0 .
\end{aligned}
\right.
\end{align}
We claim that
\begin{align}
\label{bq0}
|q_0(t)-1|&\leq C \frac{\log(t_0)}{\log(t)}
\\
\label{bq1}
|\dot q_0(t)|&\leq C \frac{\log(t_0)}{t\log^2(t)}
\\
\label{bq2}
| q_0^{(k)}(t)|&\leq C_k \frac{\log(t_0)}{t^k\log^2(t)},\quad k=2,3,\ldots
\end{align}

We write \eqref{ode-q0} as 
\begin{align*}
&\dot q_0(t)+a(t) q_0(t)=a(t),\quad t>t_0
\end{align*}
where
\begin{align}
\label{def-a}
a(t)&=\frac{\eta(t)}{[\alpha(1-\eta(t)) + \eta(t)]t\log(t)}.
\end{align}
Let 
\begin{align}
\label{def-A}
A(t)=\int_{t_0}^t a(z)\, d z.
\end{align}
Then $q_0$ is given by 
\begin{align*}
q_0(t)&=e^{-A(t)}\int_{t_0}^t e^{A(s)} a(s)\, d s.
\end{align*}
From this we get
\begin{align*}
|q_0(t)|&\leq \int_{t_0}^t e^{-\int_s^ t a} a(s)\, d s
\\
&\leq \frac1\alpha\int_{t_0}^t \frac{1}{s\log(s)}\, d s
\\
&=\frac{1}{\alpha}\log\Bigl(\frac{\log(t)}{\log(t_0)}\Bigr).
\end{align*}
We deduce that 
\begin{align}
\label{q0a}
|q_0(t)|&\leq \frac{C}{\log(t_0)},\quad t_0\leq t\leq 3t_0.
\end{align}
For $t>3t_0$ the ODE is:
\begin{align*}
\dot q_0(t)\log(t)+\frac{q_0(t)}{t}=\frac{1}{t}
\end{align*}
so 
\begin{align}
\label{q0b}
q_0(t)=1-\frac{c_0}{\log(t)}, \quad t>3t_0.
\end{align}
Evaluating at $3t_0$ we find that 
\begin{align}
\label{q0c}
c_0=\log(t_0)+O(1).
\end{align}
From \eqref{q0a}, \eqref{q0b} and \eqref{q0c} we obtain \eqref{bq0}.

Using the equation \eqref{ode-q0}, and \eqref{q0a}, \eqref{q0b}, \eqref{q0c} we obtain \eqref{bq1}.
The estimate \eqref{bq2} follows by differentiating the equation and using the previous estimates.
\end{proof}

We note that $a$ and $p$ in the construction of $p_0$ are such that 
\[
p_0(t_0)=b,\quad\lim_{t\to\infty} p_0(t)=a+b.
\]

In Lemma~\ref{lem:p0}
we considered the equation
\begin{align}
\label{ODE-2.1}
\dot p(t)\log(t)[\alpha (1-\eta(t))+\eta(t)]
+\frac{\eta(t)}{t}p(t)=f(t),
\quad t>t_0
\end{align}
for a particular right-hand side and solved it with initial condition equal to 0.
The next lemma gives a solution and estimate when $p(\infty)=0$.
\begin{lemma}
\label{lem:ODE2.0}
Let $m<-1$. If $\|f\|_{1,m}<\infty$, 
then \eqref{ODE-2.1} has a unique solution $p$ such that $p(\infty)=0$,
and it satisfies
\[
\|p\|_{0,m}+\|\dot p\|_{1,m-1}\leq C \|f\|_{1,m}.
\]
\end{lemma}
\begin{proof}
Write the equation as 
\begin{align*}
\dot p(t)+a(t) p(t)=\tilde f(t),\quad t>t_0
\end{align*}
where $a(t)$ is given in \eqref{def-a} and
\begin{align*}
\tilde f(t)=\frac{f(t)}{\log(t)[\alpha(1-\eta(t)) + \eta(t)]}.
\end{align*}
Let $A(t)$ be given by \eqref{def-A}.
Then $p$ is given by 
\begin{align*}
p(t)=-e^{-A(t)}\int_{t}^\infty e^{A(s)} \tilde f(s)\, d s.
\end{align*}
Note that for $t>3t_0$ we just have
\begin{align*}
p(t)=-\frac{1}{\log(t)}\int_t^\infty f(s)\,ds.
\end{align*}
So for $t>3t_0$
\begin{align*}
|p(t)|&\leq \|f\|_{1,m}
\frac{1}{\log(t)}\int_t^\infty \frac{\log^m(s)}{s}\,ds
\\
&\leq \|f\|_{1,m}\frac{1}{\log(t)} \frac{\log^{m+1}(t)}{|m+1|}
= C \|f\|_{1,m} \log^m(t).
\end{align*}
for which we need $m+1<0$.
In particular at $3t_0$:
\[
|p(3t_0)|\leq C \|f\|_{1,m} \log^m(t).
\]
For $t_0\leq t \leq 3t_0$, we have
\begin{align*}
A(t)&=
\int_{t_0}^t \frac{1}{[\alpha(1-\eta(z))+\eta(z)]}
\frac{1}{z\log(z)}\, d z
\\
&\leq \frac1\alpha\int_{t_0}^{3t_0}
\frac{1}{z\log(z)}\, d z
\\
&\leq C (1 + \log^{-1}(t_0)).
\end{align*}
For $t_0\leq t \leq 3t_0$,
\begin{align*}
p(t)
&= -e^{-A(t)}
\int_{t}^{3t_0} e^{A(s)} \tilde f(s)\, d s
-e^{-A(t)}\int_{3t_0}^\infty e^{A(s)} \tilde f(s)\, d s.
\end{align*}
We estimate for $t_0\leq t \leq 3t_0$,
\begin{align*}
\left|
e^{-A(t)}
\int_{t}^{3t_0} e^{A(s)} \tilde f(s)\, d s
\right|
&\leq C \|f\|_{1,m}
\int_{t_0}^{3t_0}
\frac{\log^{m-1}(s)}{s}\,d s
\\
&\leq C \|f\|_{1,m} \log^m(t_0)
\end{align*}
and
\begin{align*}
\left|
e^{-A(t)}\int_{3t_0}^\infty e^{A(s)} \tilde f(s)\, d s
\right|
&=
\left|
e^{-A(t)+A(t_0)}e^{-A(t_0)}\int_{3t_0}^\infty e^{A(s)} \tilde f(s)\, d s
\right|
\\
&\leq  C \|f\|_{1,m} \log^m(t_0).
\end{align*}
We conclude that for all $t_0$,
\begin{align*}
|p(t)|
&\leq \|f\|_{1,m}\frac{1}{\log(t)} \frac{\log^{m+1}(t)}{|m+1|}
= C \|f\|_{1,m} \log^m(t).
\end{align*}
Using the equation we also get
\begin{align*}
|\dot p(t)|\leq C \|f\|_{1,m} \frac{\log^{m-1}(t)}{t},\quad t>t_0.
\end{align*}
\end{proof}

We give a lemma to estimate $M_1[p]$ \eqref{def:M1.2} that we will apply
for various $p$.
\begin{lemma}
\label{lem:M1-p-gen}
Suppose $\|\dot p\|_{1,-k} +\|\ddot p\|_{2,-k}<\infty$ where $k>0$. 
Then 
\begin{align}
\label{eq:estM1p-gen}
|M_1[p](t)|&\leq C \frac{1}{t\log^k(t)} (\|\dot p\|_{1,-k} +\|\ddot p\|_{2,-k}).
\end{align}
If also $\|p^{(3)}\|_{3,-k}<\infty$, then
\begin{align}
\label{eq:estdotM1p-gen}
|\tfrac d{dt}M_1[p](t)|&\leq C  \frac{1}{t^2\log^k(t)}
(\|\dot p\|_{1,-k} +\|\ddot p\|_{2,-k}+\|p^{(3)}\|_{3,-k}).
\end{align}
\end{lemma}
\begin{proof}
The first term in $M_1[p]$ is estimated as follows
\begin{align*}
|\eta(t)\dot p(t)\log(2) |\leq  C \frac{1}{t\log^k(t)}\|\dot p\|_{1,-k}.
\end{align*}
Next we estimate
\begin{align*}
\eta(t)
\left|\int_{t_0}^{t/2}\dot p(s)\Bigl(\frac{1}{t-s}-\frac 1t\Bigr)\,d s\right|
&\leq C  \eta(t) 
\int_{t_0}^{t/2}\frac1{s\log^k(s)}
\frac{s}{t(t-s)}\,d s \, \|\dot p\|_{1,-k}
\\
&\leq C \frac{1}{t^2} \eta(t) 
\int_{t_0}^{t/2}\frac1{\log^k(s)}\,d s \, \|\dot p\|_{1,-k}
\\
&\leq C \frac{1}{t \log^k(t)} \, \|\dot p\|_{1,-k}.
\end{align*}
Similarly
\begin{align*}
\eta(t)\left|\frac{p(t/2)-p(t)}t \right|
&\leq C \eta(t) \sup_{\frac t2\leq s \leq t}|\dot p(s)|\,\|\dot p\|_{1,-k}
\\
& \leq C \frac{1}{t\log^k(t)}\,\|\dot p\|_{1,-k} 
\end{align*}
and
\begin{align*}
\left| 
\int_{t/2}^{t-t^\alpha}\frac{\dot p(t)-\dot p(s)}{t-s}\,d s 
\right| 
&\leq C  t \sup_{t/2\leq s\leq t-t^\alpha} |\ddot p(s)| 
\\
&\leq 
C \frac{1}{t\log^k(t)} \,\|\ddot p\|_{2,-k} .
\end{align*}
This proves \eqref{eq:estM1p-gen}.

To obtain \eqref{eq:estdotM1p-gen}, let us compute
\begin{align*}
&
\frac{d}{dt}
\int_{t/2}^{t-t^\alpha}\frac{\dot p(t)-\dot p(s)}{t-s}\,d s 
\\
&=
\frac{d}{dt}
\int_{t^\alpha}^{t/2}\frac{\dot p(t)-\dot p(t-u)}{u}\,d u
\\
&=\tfrac12\frac{\dot p(t)-\dot p(\frac t2)}{\frac t2}
-\alpha t^{\alpha-1}\frac{\dot p(t)-\dot p(t-t^\alpha)}{t^\alpha}
+
\int_{t^\alpha}^{t/2}\frac{\ddot p(t)-\ddot p(t-u)}{u}\,d u.
\end{align*}
We estimate each term. In particular for the third one,
using a bound for the third derivative of $p$ we get
\begin{align*}
\left|
\int_{t^\alpha}^{t/2}\frac{\ddot p(t)-\ddot p(t-u)}{u}\,d u
\right|
&\leq 
C\int_{t^\alpha}^{t/2}\frac{1}{\xi_u^3 \log^k(\xi_u) }\,d u \|p^{(3)}\|_{3,-k}
\\
&\leq 
C\int_{t^\alpha}^{t/2}\frac{1}{t^3 \log^k(t) }\,d u\|p^{(3)}\|_{3,-k}
\\
&\leq 
C\frac{1}{t^2 \log^k(t) } \|p^{(3)}\|_{3,-k}.
\end{align*}
The other terms are handled similarly.
\end{proof}

Let $p_0$ be defined by \eqref{eq:p0}.
Using Lemma~\ref{lem:M1-p-gen} and  \eqref{est:dp0.1}, \eqref{est:dkp0.1} 
we get
\begin{align}
\label{eq:estM1p0}
|M_1[p_0](t)|&\leq C  |a| \frac{\log(t_0)}{t\log^2(t)} 
\\
\nonumber
|\tfrac d{dt}M_1[p_0](t)|&\leq C  |a| \frac{\log(t_0)}{t^2\log^2(t)} .
\end{align}

Going back to the equation $\mathcal S[p]=0$, we have
\begin{align*}
\mathcal S[p_0]
&=M_0[p_0]+M_1[p_0] - f(t) + H(p_0) \dot p_0
\\
&=(a+b)\frac{\eta(t)}{t}+M_1[p_0]  + H(p_0) \dot p_0-f(t)
\end{align*}
We introduce a correction $p_1$ to eliminate the terms 
$M_1[p_0]  + H(p_0) \dot p_0$.
So we compute
\begin{align*}
\mathcal S[p_0+p_1]
&=M[p_0]+M[p_1]+ H(p_0+p_1) \tfrac{d}{dt}(p_0+p_1)- f(t)
\\
&=\mathcal S[p_0]+H(p_0+p_1) \tfrac{d}{dt}(p_0+p_1)- H(p_0) \dot p_0 
+ M[p_1]
\\
&=\left(a+b-p_1(t_0)\right)\frac{\eta(t)}{t}
+ M_1[p_0]  + M_0[p_1]+M_1[p_1]
\\
&\quad
+ H(p_0+p_1) \dot p_0 + H(p_0+p_1) \dot p_1 
-f(t).
\end{align*}
Then we solve the following nonlinear equation for $p_1$:
\begin{align}
\label{eq:p1}
\left\{
\begin{aligned}
&M_0[p_1]+M_1[p_0]+ H(p_0+p_1) \dot p_0=0, \quad t>t_0
\\
&p_1(\infty)=0
\end{aligned}
\right.
\end{align}

\begin{lemma}
There is a unique solution $p_1$ of \eqref{eq:p1} satisfying
\begin{align}
\label{est:p1}
|p_1(t)|&\leq  C  |a| \frac{\log(t_0)}{\log^2(t)} 
\\
\label{est:dotp1}
|\dot p_1(t)|&\leq C  |a| \frac{\log(t_0)}{t\log^3(t)}.
\end{align}
Moreover, $p_1$ is Lipschitz with respect to $a$ and $b$ in the sense that
\begin{align}
\label{est:lip-p1}
\left\{
\begin{aligned}
|p_1(t;a_1,b)-p_1(t;a_2,b)|
&\leq  C  |a_1-a_2| \frac{\log(t_0)}{\log^2(t)} 
\\
|\dot p_1(t;a_1,b)-\dot p_1(t;a_2,b)|
&\leq  C  |a_1-a_2| \frac{\log(t_0)}{t\log^3(t)}
\\
|p_1(t;a,b_1)-p_1(t;a,b_2)|
&\leq  C  |a||b_1-b_2| \frac{\log(t_0)}{\log^2(t)}
\\
|\dot p_1(t;a,b_1)-\dot p_1(t;a,b_2)|
&\leq  C  |a||b_1-b_2| \frac{\log(t_0)}{t\log^3(t)} .
\end{aligned}
\right.
\end{align}
\end{lemma}
\begin{proof}
We solve \eqref{eq:p1} by the contraction mapping principle, 
on the space of $C^1$ functions $p_1:[t_0,\infty)\to\C$, with $p_1(t)\to0$ as $t\to\infty$ and norm $\|\dot p_1\|_{1,-3}<\infty$.
Lemma~\ref{lem:ODE2.0} constructs an operator $M_0^{-1}$ and the fixed point problem is written as 
\[
p_1=M_0^{-1}\left[
-M_1[p_0]-H(p_0+p_1)\dot p_0\right].
\]
The fixed point satisfies
\[
\|p_1\|_{0,-2}+\|\dot p_1\|_{1,-3}\leq C |a| \log(t_0),
\]
which are the estimates \eqref{est:p1} and \eqref{est:dotp1}.

To prove the Lipschitz estimate \eqref{est:lip-p1} we take the difference of the fixed point problem associated to two values of $b$.
\end{proof}
Using \eqref{lem:M1-p-gen} we find that 
\begin{align}
\label{est:M1p1}
|M_1[p_1]|\leq C  |a| \frac{\log(t_0)}{t\log^3(t)} .
\end{align}
With this choice of $p_1$ we get
\begin{align*}
\mathcal S[p_0+p_1]
&=\left(a+b-p_1(t_0)\right)\frac{\eta(t)}{t}
+M_1[p_1] + H(p_0+p_1) \dot p_1-f(t).
\end{align*}
We introduce a new correction $p_2$ 
with the idea of eliminating $M_1[p_1] + H(p_0+p_1) \dot p_1 $.
For this we compute
\begin{align*}
&\mathcal S[p_0+p_1+p_2]
\\
&=
M[p_0]+M[p_1]+M[p_2]+H(p_0+p_1+p_2)\tfrac{d}{dt}(p_0+p_1+p_2)
- f(t)
\\
&= \mathcal S[p_0+p_1]+M[p_2]+H(p_0+p_1+p_2) \tfrac{d}{dt}(p_0+p_1+p_2)
\\ 
&\quad 
-H(p_0+p_1) \tfrac{d}{dt}(p_0+p_1)
\\
&=
\left(a+b-p_1(t_0)-p_2(t_0)\right)\frac{\eta(t)}{t}
+ M_0[p_2]+M_1[p_2]
+ M_1[p_1] + H(p_0+p_1) \dot p_1 
\\ 
&\quad 
+ [H(p_0+p_1+p_2)-H(p_0+p_1)]\tfrac{d}{dt}(p_0+p_1)
+H(p_0+p_1+p_2) \dot p_2 
-f(t).
\end{align*}
We solve the following equation for $p_2$:
\begin{align}
\label{eq:p2}
\left\{
\begin{aligned}
&M_0[p_2]+ M_1[p_1] + H(p_0+p_1) \dot p_1=0, \quad t>t_0
\\
&p_2(\infty)=0
\end{aligned}
\right.
\end{align}
\begin{lemma}
There is a unique solution $p_2$ of \eqref{eq:p2}, and it satisfies
\begin{align}
\label{eq:estp2}
|p_2(t)|&\leq  C  |a| \frac{\log(t_0)}{\log^3(t)}
\\
\label{eq:estdotp2}
|\dot p_2(t)|&\leq C  |a| \frac{\log(t_0)}{t\log^4(t)} 
\\
\nonumber
|\ddot p_2(t)|&\leq C  |a| \frac{\log(t_0)}{t^2\log^4(t)} .
\end{align}
We also have Lipschitz estimates
\begin{align}
\label{est:lip-p2}
\left\{
\begin{aligned}
|p_2(t;a_1,b)-p_2(t;a_2,b)|
&\leq  C  |a_1-a_2| \frac{\log(t_0)}{\log^3(t)}
\\
|\dot p_2(t;a_1,b)-\dot p_2(t;a_2,b)|
&\leq  C  |a_1-a_2| \frac{\log(t_0)}{t\log^4(t)}
\\
|p_2(t;a,b_1)-p_2(t;a,b_2)|
&\leq  C  |a||b_1-b_2| \frac{\log(t_0)}{\log^3(t)}
\\
|\dot p_2(t;a,b_1)-\dot p_2(t;a,b_2)|
&\leq  C  |a||b_1-b_2| \frac{\log(t_0)}{t\log^4(t)} .
\end{aligned}
\right.
\end{align}
\end{lemma}
\begin{proof}
Equation \eqref{eq:p2} is linear and we can apply Lemma~\ref{lem:ODE2.0} directly 
using the estimate for $\dot p_1$ \eqref{est:dotp1} and $M_1[p_1]$ \eqref{est:M1p1}.
The proof of the Lipschitz estimates \eqref{est:lip-p2} is standard and we omit it.
\end{proof}
Applying Lemma~\ref{lem:M1-p-gen} we obtain
\begin{align}
\label{eq:M1p2}
|M_1[p_2](t)|\leq C  |a| \frac{\log(t_0)}{t\log^4(t)} 
\end{align}
With this choice of $p_2$ we get
\begin{align*}
\mathcal S[p_0+p_1+p_2]
&=
\left(a+b-p_1(t_0)-p_2(t_0)\right)\frac{\eta(t)}{t}
+M_1[p_2]
\\ 
&\quad 
+[H(p_0+p_1+p_2)-H(p_0+p_1)]\tfrac{d}{dt}(p_0+p_1)
\\ 
&\quad 
+ H(p_0+p_1+p_2)\dot p_2
-f(t).
\end{align*}
Now we add $p_3$:
\begin{align*}
&
\mathcal S[p_0+p_1+p_2+p_3]\\
\\
&=
M[p_0]+M[p_1]+M[p_2]+M[p_3]+H(p_0+p_1+p_2+p_3)\tfrac{d}{dt}(p_0+p_1+p_2+p_3)
\\
&\quad - f(t)
\\
&=\mathcal S[p_0+p_1+p_2]+M[p_3]
\\
&\quad
+H(p_0+p_1+p_2+p_3)\tfrac{d}{dt}(p_0+p_1+p_2+p_3)
\\
&\quad 
-H(p_0+p_1+p_2)\tfrac{d}{dt}(p_0+p_1+p_2)
\\
&=
\left(a+b-p_1(t_0)-p_2(t_0)\right)\frac{\eta(t)}{t}
+M_1[p_2]
\\ 
&\quad 
+ [H(p_0+p_1+p_2)-H(p_0+p_1)]\tfrac{d}{dt}(p_0+p_1)
+ H(p_0+p_1+p_2) \dot p_2
-f(t)
\\
&\quad 
+M[p_3]
\\
&\quad 
+[H(p_0+p_1+p_2+p_3)-H(p_0+p_1+p_2)]\tfrac{d}{dt}(p_0+p_1+p_2)
\\
&\quad 
+H(p_0+p_1+p_2+p_3)\dot p_3.
\end{align*}
We ask that $p_3$ solves 
the nonlinear equation 
\begin{align}
\nonumber
& M[p_3]+M_1[p_2]
+ [H(p_0+p_1+p_2)-H(p_0+p_1)] \tfrac{d}{dt}(p_0+p_1)
\\
\nonumber
&\quad
+ H(p_0+p_1+p_2) \dot p_2
-f(t)
\\
\nonumber
&\quad 
+[H(p_0+p_1+p_2+p_3)-H(p_0+p_1+p_2)] \tfrac{d}{dt}(p_0+p_1+p_2)
\\
\nonumber
&\quad 
+ H(p_0+p_1+p_2+p_3) \dot p_3 
\\
\label{eq:p3}
&=-p_3(t_0)\frac{\eta(t)}{t}
\end{align}
with $p_3(\infty)=0$.
Then
\begin{align*}
\mathcal S[p_0+p_1+p_2+p_3]
=\left(a+b-p_1(t_0)-p_2(t_0)-p_3(t_0)\right)\frac{\eta(t)}{t}.
\end{align*}

For the proof of existence of a solution to \eqref{eq:p3} 
we rewrite the operator $M$ (defined in \eqref{def:M.2})
as follows.
For $t>\frac32t_0$, write
\begin{align*}
\int_{t_0}^{t-t^\alpha}\frac{\dot p(s)}{t-s}ds
&=
\frac{1}{t}\int_{t_0}^{t-t^\alpha}\dot p(s)ds
+\int_{t_0}^{t-t^\alpha}\dot p(s)\Bigl(\frac{1}{t-s}-\frac{1}{t}\Bigr)ds
\\
&=\frac{p(t-t^\alpha)-p(t_0)}{t}
+\int_{t_0}^{t-t^\alpha}\dot p(s)\Bigl(\frac{1}{t-s}-\frac{1}{t}\Bigr)ds
\\
&=
\frac{p(t)}{t}
+\frac{p(t-t^\alpha)-p(t)}{t}
-\frac{p(t_0)}{t}
+\int_{t_0}^{t-t^\alpha}\dot p(s)\Bigl(\frac{1}{t-s}-\frac{1}{t}\Bigr)ds
\end{align*}
and therefore
\begin{align}
\nonumber
M[p](t) = \dot p(t)\log(t^\alpha)+\eta(t)\frac{p(t)}{t}
+\eta(t)\Bigg[&
\frac{p(t-t^\alpha)-p(t)}{t}
-\frac{p(t_0)}{t}
\\
\label{eq:decompM2}
&\quad 
+\int_{t_0}^{t-t^\alpha}\dot p(s)\Bigl(\frac{1}{t-s}-\frac{1}{t}\Bigr)ds
\Bigg].
\end{align}
Thus we are led to the ODE
\begin{align}
\label{eq:ode2}
\dot p(t) \log(t^\alpha)+\frac{\eta(t)}t p(t) = f(t),\quad t>t_0.
\end{align}
\begin{lemma}
\label{lem:ode2}
Assume $m<-\frac1\alpha$. If $\|f\|_{1,m}<\infty$ then 
\eqref{eq:ode2} has a unique solution such that $p(t)\to0$ as $t\to\infty$,
which satisfies
\begin{align}
\nonumber
|\dot p(t)|&\leq \|f\|_{1,m}
\frac{\log^{m-1}(t)}{t}
\left|\frac{m}{1+\alpha m}\right|
\Bigl(1+O\Bigl(\frac{1}{\log(t_0)}\Bigr)\Bigr)
\\
\nonumber
|p(t)|&\leq \frac{\log^m(t)}{|1+\alpha m|}\|f\|_{1,m}(1+O(\log^{-1}(t_0))).
\end{align}
\end{lemma}
\begin{proof}
Let
\[
a(t)=\frac{\eta(t)}{\alpha t \log(t)} ,\quad
u(t)=\exp\Bigl(\int_{t_0}^t a(s)\,ds \Bigr). 
\]
We have
\begin{align*}
u(t)&=u(3t_0)\frac{\log^{\frac{1}{\alpha}}(t)}{\log^{\frac{1}{\alpha}}(3t_0)},\quad t\geq3t_0
\\
u(t)&=1+O\Bigl(\frac{1}{\log(t_0)}\Bigr), \quad t_0\leq t\leq 3t_0.
\end{align*}
In particular
\begin{align*}
u(3t_0)=1+O\Bigl(\frac{1}{\log(t_0)}\Bigr).
\end{align*}
The solution that we want to \eqref{eq:ode2} is given by
\begin{align}
\nonumber
p(t)=-\frac{1}{u(t)}\int_{t}^\infty  \frac{u(s)}{\log(s^\alpha) } f(s)\,ds.
\end{align}
For $t\geq 3t_0$
\begin{align*}
p(t)&=-\frac{1}{\log^{\frac{1}{\alpha}}(t)}
\int_t^\infty  \frac{\log^{\frac1\alpha}(s)}{\log(s^\alpha) } f(s)\,ds .
\end{align*}
But since $m+\frac1\alpha<0$,
\begin{align*}
|p(t)|&\leq
\|f\|_{1,m}
\frac1{\log^{\frac1\alpha}(t)}
\int_t^\infty
\frac{\log^{\frac1\alpha}(s)}{\log(s^\alpha) } \frac{\log^m(s)}{s}\,ds 
\\
&=
\|f\|_{1,m}
\frac1{\alpha\log^{\frac1\alpha}(t)}
\int_t^\infty
\frac{\log^{m+\frac1\alpha-1}(s)}{s}\,ds 
\\
&=
\|f\|_{1,m}
\frac1{|1+\alpha m|}
\log^{m}(t) .
\end{align*}
The calculation for $t_0\leq t\leq 3t_0$ is similar, and we get
\begin{align}
\label{est-p2}
|p(t)|\leq 
\|f\|_{1,m}
\frac{\log^{m}(t)}{|1+\alpha m|}
\Bigl(1+O\Bigl(\frac{1}{\log(t_0)}\Bigr)\Bigr),\quad t\geq t_0.
\end{align}
Using the equation \eqref{eq:ode2} and \eqref{est-p2} we get
\begin{align}
\nonumber
|\dot p(t)| 
&\leq \frac1{\alpha\log(t)}|f(t)|+\frac1{\alpha t\log(t)}|p(t)|
\\
\nonumber
&\leq
\|f\|_{1,m}
\frac{\log^{m-1}(t)}{t}\left[\frac1\alpha+\frac1{\alpha|1+m\alpha|}\right]
\Bigl(1+O\Bigl(\frac{1}{\log(t_0)}\Bigr)\Bigr),\quad t\geq t_0.
\end{align}
\end{proof}

Here is the existence result for \eqref{eq:p3}.
The functions $p_j$,
$j=0,1,2$ that appear in \eqref{eq:p3} depend on the parameters
$a$, $b$. We assume that
\begin{align}
\label{eq:bound-a-b}
|a|\leq2 , \quad |b|\leq 2.
\end{align}

\begin{lemma}
\label{lem:existencep3}
Assume $\alpha \in (\frac58,1)$ and $-4<m<-\frac1{2\alpha-1}$.
Also assume \eqref{eq:bound-a-b} and $\|f\|_{1,m}\leq \log^{-3-m}(t_0)$.
Then for $t_0$ sufficiently large, \eqref{eq:p3} has unique solution $p_3$
of \eqref{eq:p3} with $\|\dot p_3\|_{1,m-1}<\infty$ and it satisfies
\begin{align}
\nonumber
\|\dot p_3\|_{1,m-1} \le C ( \log^{-3-m}(t_0) + \|f\|_{1,m}).
\end{align}
We also have that $p_3$ is Lipschitz with respect to $a$, $b$ and $f$ in the following sense. Let $p=p(t,a,b,f)$. Then
\begin{align}
\nonumber
\| \dot p_3(\cdot,a_1,b_1,f_1)- \dot p_3(\cdot,a_2,b_2,f_2) \|_{1,m-1} 
&\leq 
C \log^{-3-m}(t_0)  ( |a_1-a_2|+|b_1-b_2|) 
\\
\nonumber
&\quad 
+C\|f_1-f_2\|_{1,m}
\end{align}
and
\begin{align}
\nonumber
\| p_3(\cdot,a_1,b_1,f_1)-p_3(\cdot,a_2,b_2,f_2) \|_{0,m} 
&\leq 
C \log^{-3-m}(t_0)  ( |a_1-a_2|+|b_1-b_2|) 
\\
\label{eq:lip-p3}
&\quad 
+C\|f_1-f_2\|_{1,m}.
\end{align}
\end{lemma}
\begin{proof}
Using \eqref{eq:decompM2} we rewrite \eqref{eq:p3} as 
\begin{align}
\label{eq:p3b}
0&=
\dot p(t)\log(t^\alpha)+\eta(t)\frac{p(t)}{t}
+ A(p,\dot p) + E
\end{align}
where
\begin{align*}
A(p,\dot p)
=A_1(p,\dot p)+A_2(p,\dot p)+A_3(p,\dot p)+A_4(p,\dot p)
\end{align*}
\begin{align}
\nonumber
A_1(p,\dot p)(t)
&= \eta(t)\int_{t_0}^{t-t^\alpha}\dot p(s)\Bigl(\frac{1}{t-s}-\frac{1}{t}\Bigr)ds
\\
\nonumber
A_2(p,\dot p)(t)
&=
\eta(t) \frac{p(t-t^\alpha)-p(t)}{t}
\\
\nonumber
A_3(p,\dot p)(t)
&=
[H(p_0+p_1+p_2+p)-H(p_0+p_1+p_2)] \tfrac{d}{dt}(p_0+p_1+p_2)
\\
\nonumber
A_4(p,\dot p)(t)
&= H(p_0+p_1+p_2+p) \dot p 
\end{align}
and
\begin{align}
\nonumber
E&=M_1[p_2]+[H(p_0+p_1+p_2)-H(p_0+p_1)] \tfrac{d}{dt}(p_0+p_1)
\\
\label{def:Ep}
& \quad + H(p_0+p_1+p_2) \dot p_2
-f(t).
\end{align}
Let $X$ be the space of functions $p\in C^1([t_0,\infty))$ such that $p(t)\to0$ as $t\to\infty$ and $\|\dot p\|_{1,m-1}<\infty$. In this space we use the norm $\|p\|_X=\|\dot p\|_{1,m-1}$ defined in \eqref{def:norm-gamma-m}.
Let $T(f)$ denote the solution $p$ of \eqref{eq:ode2} given by Lemma~\ref{lem:ode2}.
We treat \eqref{eq:p3b} as the fixed point problem
\[
p= \mathcal F(p,E), \quad  \mathcal F(p,E)=T(-A(p,\dot p) - E) .
\]
We claim that $\mathcal F$ has a unique fixed point in the closed ball $\overline B_\ell(0) \subset X$, with $\ell>0$ to be determined.
By Lemma~\ref{lem:ode2}, if $p = T(f)$ then
\begin{align*}
\| \dot p\|_{1,m-1} \leq \|f\|_{1,m}
\left|\frac{m}{1+\alpha m}\right|
\Bigl(1+O\Bigl(\frac{1}{\log(t_0)}\Bigr)\Bigr) .
\end{align*}
To prove that $\mathcal F$ is a contraction, we estimate each term of $A$.
Since $A_1$ is linear, we need to compute
$$
\left|
\int_{t_0}^{t-t^\alpha}\dot{p}(s)\Bigl(\frac{1}{t-s}-\frac{1}{t}\Bigr)\,ds\right|
\leq \|\dot p\|_{1,m-1}
\int_{t_0}^{t-t^\alpha}\frac{\log^{m-1}( s)}{s } \frac{s}{(t-s)t}\,ds
$$
and assuming $t\geq 2 t_0$,
\begin{align*}
&\int_{t_0}^{t-t^\alpha} \frac{\log^{m-1} (s)}{(t-s)t}\,ds
\\
&=
\frac1t
\int_{t_0}^{t/2}\frac{\log^{m-1} (s)}{t-s}\,ds
+\frac1t\int_{t/2}^{t-t^\alpha}\frac{\log^{m-1} (s)}{t-s}\,ds
\\
&\leq 
\frac{C}{t^2}\int_{t_0}^{t/2}\log^{m-1} (s)\,ds
+ \log^{m-1}(t)\Bigl(1+O(\frac{1}{\log t})\Bigr) 
\frac{1}{t}\int_{t/2}^{t-t^\alpha}\frac{1}{t-s} \,ds
\\
&\leq C\frac{\log^m(t)}{t \log t}
+ (1-\alpha)\frac{\log^m(t)}{t}\Bigl(1+O(\frac{1}{\log t})\Bigr).
\end{align*}
If $t\leq 2t_0$ the computation is similar and we get
\begin{align*}
\|A_1(\dot p)\|_{1,m} \leq 
\Bigl(1-\alpha+  C \frac{1}{\log(t_0)}\Bigr) \|\dot p\|_{1,m-1}.
\end{align*}
Let's estimate $A_2$:
\begin{align*}
\left|\eta(t)\frac{ p(t-t^\alpha)- p(t)}{t}\right|
&\leq \eta(t) t^{\alpha-1}\sup_{s\in[t-t^\alpha,t]}|\dot p(s)|
\\
&\leq C\eta(t) t^{\alpha-1} t^{-1} \log^{m-1}(t)
\|\dot p\|_{L^\infty}.
\\
&\leq C \frac{\log^m(t)}{t} t_0^{\alpha-1} \log^{-1}(t_0) \|\dot p\|_{L^\infty}.
\end{align*}
so 
\begin{align*}
\|A_2(\dot p)\|_{1,m} \leq 
C t_0^{\alpha-1} \log^{-1}(t_0)  \|\dot p\|_{1,m-1}.
\end{align*}
For $A_3$, the Lipschitz property looks like
\begin{align*}
|A_3(p,\dot p)-A_3(\tilde p,\dot{\tilde p})|
&\leq
C |p(t)-\tilde p(t)| 
|\tfrac{d}{dt}(p_0+p_1+p_2)|
\\
&\leq 
C \log^m(t) \|\dot p-\dot{\tilde p}\|_{1,m-1} 
|a| \frac{\log(t_0)}{t\log^2(t)}
\end{align*}
where we have used the estimates \eqref{est:dp0.1}, \eqref{est:dotp1}, 
\eqref{eq:estdotp2}.
So, using \eqref{eq:bound-a-b},
\begin{align*}
\| A_3(p,\dot p)-A_3(\tilde p,\dot{\tilde p}) \|_{1,m}
\leq C \frac{1}{\log(t_0)}\|\dot p-\dot{\tilde p}\|_{1,m-1} .
\end{align*}
Finally, if $\|\dot p\|_{1,m-1}\leq \ell$, $\|\dot{\tilde p}\|_{1,m-1}\leq \ell$, then
\begin{align*}
|A_4(p,\dot p)-A_4(\tilde p,\dot{\tilde p})|
&\leq |H(p_0+p_1+p_2+p)-H(p_0+p_1+p_2+\tilde p)|\, |\dot p  |
\\
&\quad +  |H(p_0+p_1+p_2+\tilde p)|\, |\dot p(t)-\dot{\tilde p}(t)|
\\
&\leq C |p(t)-\tilde p(t)|\, |\dot p(t)|
+C  |\dot p(t)-\dot{\tilde p}(t)|
\\
&\leq C \log^m(t) \|\dot p-\dot{\tilde p}\|_{1,m-1} 
\frac{\log^{m-1}(t)}{t} \|\dot p\|_{1,m-1} 
\\
&\quad 
+ C \frac{\log^{m-1}(t)}{t} \|\dot p-\dot{\tilde p}\|_{1,m-1} 
\end{align*}
so
\begin{align*}
\| A_4(p,\dot p)-A_4(\tilde p,\dot{\tilde p}) \|_{1,m}
\leq 
\left( 
C \ell \log^{m-1}(t_0)
+ \frac{C}{\log(t_0)} \right) \|\dot p-\dot{\tilde p}\|_{1,m-1} .
\end{align*}
The same computation shows that 
\begin{align*}
|A_4(p,\dot p)|
\leq C\, |\dot p  |
\leq C 
\frac{\log^{m-1}(t)}{t} \|\dot p\|_{1,m-1} 
\end{align*}
so
\begin{align*}
\| A_4(p,\dot p) \|_{1,m}
\leq 
\frac{C}{\log(t_0)}  \|\dot p\|_{1,m-1} .
\end{align*}

It follows that for $p,\tilde p\in \overline{B}_\ell(0)\subset X$,
\begin{align}
\nonumber
& \| \mathcal F(p,E)-\mathcal F(\tilde p,E) \|_X
\\
\label{eq:lipF}
& \quad \leq
\left|\frac{m}{1+\alpha m}\right|
\Bigl(1+O\Bigl(\frac{1}{\log(t_0)}\Bigr)\Bigr) 
\Bigl(1-\alpha+ C \ell \log^{m-1}(t_0) +\frac C{\log(t_0)}\Bigr) \|p-\tilde p\|_X.
\end{align}

Next we compute $\|E\|_{1,m}$ \eqref{def:Ep}.
By \eqref{eq:M1p2}, and the assumption $-4\leq m$, 
\begin{align*}
\|M_1[p_2]\|_{1,m} \leq C \log^{-3-m}(t_0).
\end{align*}
Using \eqref{est:dp0.1}, \eqref{eq:estp2} we get
\begin{align*}
|[H(p_0+p_1+p_2)-H(p_0+p_1)] \tfrac{d}{dt}(p_0+p_1)|
&\leq C |p_2(t)| \,  |a| \frac{\log(t_0)}{t\log^2(t)}
\\
&\leq C  |a| \frac{\log(t_0)}{\log^3(t)} \,  |a| \frac{\log(t_0)}{t\log^2(t)}
\\
&\leq C \frac{\log^2(t_0)}{t \log^5(t)}.
\end{align*}
Therefore
\begin{align*}
\| [H(p_0+p_1+p_2)-H(p_0+p_1)] \tfrac{d}{dt}(p_0+p_1) \|_{1,m}
\leq C \log^{-3-m}(t_0).
\end{align*}
For the next term, using \eqref{eq:estdotp2} we get
\begin{align*}
|H(p_0+p_1+p_2) \dot p_2|\leq C  |a| \frac{\log(t_0)}{t\log^4(t)} 
\end{align*}
so 
\begin{align*}
\| H(p_0+p_1+p_2) \dot p_2 \|_{1,m} \leq C \log^{-3-m}(t_0).
\end{align*}
It follows that 
\begin{align}
\label{eq:estEp}
\|E\|_{1,m} \leq C \log^{-3-m}(t_0)+\|f\|_{1,m}.
\end{align}
We choose $\ell =C_0 \log^{-3-m}(t_0)$ for some big $C_0$.
We are using here the assumption $\|f\|_{1,m}\leq \log^{-3-m}(t_0)$.
Without this assumption we can get existence for $t_0$ large, but
$t_0$ then depends on $f$.

Note that $\ell \ll \log(t_0)$ for $t_0$ large, because $m>-4$.
Also note that $\left|\frac{m}{1+\alpha m}\right|(1-\alpha)<1$ from the assumption
$m<-\frac1{2\alpha-1}$. Therefore for $t_0$ large, from \eqref{eq:lipF}, we see
that $\mathcal F$ is a contraction. Then \eqref{eq:estEp} shows that $\mathcal F$ 
maps $\overline B_\ell(0)$ into itself and therefore it has a unique fixed point in this 
ball. The Lipschitz property follows from a standard argument once we show that
$E$ and $A(p,\dot p)$ are Lipschitz with respect to $a$, $b$ and $f$. 
Let us analyze the Lipschitz dependence of $E$ on $a$.
We have
\begin{align*}
E(t,a_1,b)-E(t,a_2,b)
&= D_1+D_2+D_3
\end{align*}
where
\begin{align*}
D_1&=M_1[p_2(a_1,b)]-M_1[p_2(a_2,b)]
\\
D_2&=[H(p_0+p_1+p_2)-H(p_0+p_1)] \tfrac{d}{dt}(p_0+p_1)\big|_{a=a_1}
\\
& \quad -
[H(p_0+p_1+p_2)-H(p_0+p_1)] \tfrac{d}{dt}(p_0+p_1)\big|_{a=a_2}
\\
D_3 
&= H(p_0+p_1+p_2) \dot p_2\big|_{a=a_1}
-H(p_0+p_1+p_2) \dot p_2\big|_{a=a_2}.
\end{align*}
We have 
\[
|D_1|\leq C  |a_1-a_2| \frac{\log(t_0)}{t\log^4(t)} 
\]
and so
\[
\| D_1 \|_{1,m} \leq C \log^{-3-m}(t_0)\,|a_1-a_2|.
\]
Also
\begin{align*}
|D_2| 
&\leq C |p_2(t,a_1,b)-p_2(t,a_2,b)|\, |\dot p_0(t,a_1,b)|
\\
&\quad + C|p_2(t,a_2,b)|\, |\dot p_0(t,a_1,b)-\dot p_0(t,a_2,b)|
\\
&\leq C |a_1-a_2| \frac{\log(t_0)}{\log^3(t)} \frac{\log(t_0)}{t\log^2(t)}
\end{align*}
and so
\[
\|D_2\|_{1,m}\leq C \log^{-3-m}(t_0)\,|a_1-a_2|.
\]
Next we deal with $D_3$:
\begin{align*}
|D_3|
&\leq
\bigl|H(p_0(t,a_1,b)+p_1(t,a_1,b)+p_2(t,a_1,b)) 
\\
& \qquad
-H(p_0(t,a_2,b)+p_1(t,a_2,b)+p_2(t,a_2,b)) \bigr|
\, |\dot p_2(t,a_1,b)|
\\
&\quad 
+ |H(p_0(t,a_2,b)+p_1(t,a_2,b)+p_2(t,a_2,b))|
\, |\dot p_2(t,a_1,b)-\dot p_2(t,a_2,b)|
\\
&\leq C |a_1-a_2| \frac{\log(t_0)}{\log(t)} \frac{\log(t_0)}{t \log^4(t)}.
\end{align*}
This implies
\[
\|D_3\|_{1,m}\leq C \log^{-3-m}(t_0)\,|a_1-a_2|.
\]
The Lipschitz property of $E$ with respect to $b$ is similar.
The  dependence of $E$ on $f$ is linear.

The analysis of $A(p,\dot p)$ is similar.

\end{proof}

\begin{proof}[Proof of Proposition~\ref{prop:p1v2}]
Given $f$ and parameters $a$, $b$, thanks to Lemma~\ref{lem:existencep3}, we 
have a solution $p_3$ to \eqref{eq:p3}.
Let $p=p_0+p_1+p_2+p_3$ and let us write it as $p(t,a,b)$.
Then
\begin{align*}
\mathcal S[p]
=\left(a+b-p_1(t_0)-p_2(t_0)-p_3(t_0)\right)\frac{\eta(t)}{t}.
\end{align*}
Now we try to adjust $a$, $b$ such that 
$\mathcal S[p(\cdot,a,b)]=0$, which gives the equation
\begin{align}
\label{eq:ab1}
a+b-p_1(t_0,a,b)-p_2(t_0,a,b)-p_3(t_0,a,b)&=0.
\end{align}
We also want that $p(t_0)=1$,
which gives the equation
\begin{align}
\label{eq:ab2}
b+p_1(t_0,a,b)+p_2(t_0,a,b)+p_3(t_0,a,b)&=1.
\end{align}
The system \eqref{eq:ab1}, \eqref{eq:ab2} is equivalent to
\begin{align}
\label{eq:ab1.2}
a+2b-1&=0
\\
\label{eq:ab2.2}
b+p_1(t_0,a,b)+p_2(t_0,a,b)+p_3(t_0,a,b)&=1
\end{align}
so we need to solve \eqref{eq:ab2.2} with $a=1-2b$, that is,
\begin{align}
\label{eq:for-b}
b+p_1(t_0,1-2b,b)+p_2(t_0,1-2b,b)+p_3(t_0,1-2b,b)=1.
\end{align}
But $p_j(t_0,1-2b,b)$ is Lipschitz in $b$ 
with a small Lipschitz constant. Indeed, 
from \eqref{eq:lip-p3} we get
\begin{align}
\nonumber
| p_3(t_0,a_1,b_1,f_1)-p_3(t_0,a_2,b_2,f_2) | 
& \leq 
C \log^{-3}(t_0)  ( |a_1-a_2|+|b_1-b_2|) 
\\
\nonumber
& \quad +C \log^m(t_0)\|f_1-f_2\|_{1,m}.
\end{align}
So the existence of $b$ close to 1 such that \eqref{eq:for-b} holds
follows from the contraction mapping principle. 
We claim that $b$ is Lipschitz with respect to $f$. 
Suppose that $b_1$ is associated to $f_1$ and $b_2$ is associated to $f_2$,
where $f_1$ and $f_2$ satisfy the hypothesis 
$\|f_j\|_{1,m}\leq  \log^{-3-m}(t_0)$.
Then
\begin{align*}
1&=b_1+p_1(t_0,1-2b_1,b_1)+p_2(t_0,1-2b_1,b_1)+p_3(t_0,1-2b_1,b_1,f_1)
\\
1&=b_2+p_1(t_0,1-2b_2,b_2)+p_2(t_0,1-2b_2,b_2)+p_3(t_0,1-2b_2,b_2,f_2).
\end{align*}
Subtracting and using \eqref{est:lip-p1}, \eqref{est:lip-p2}, \eqref{eq:lip-p3} we get
\begin{align}
\label{eq:lip-ab-wrto-f}
|a_1-a_2|+|b_1-b_2|\leq C \log^m(t_0)\|f_1-f_2\|_{1,m}.
\end{align}

We have constructed $p=p_0+p_1+p_2+p_3$ satisfying \eqref{prob12b}, and
it is the unique solution to this problem. It remains to verify that it satisfies the estimates \eqref{eq:estp1xx} and \eqref{eq:estdotp1xx}.
Note that for $f\equiv 0$ there exist $a_0=-1+O(\log^{-1}(t_0))$ and 
$b_0=1+O(\log^{-1}(t_0))$ that solve the system \eqref{eq:ab1.2}, 
\eqref{eq:ab2.2}. But the solution to \eqref{prob12b} when $f\equiv 0$ is
\[
p(t,a_0,b_0,0)=1\quad\forall t\geq t_0.
\]
Then, by \eqref{eq:lip-ab-wrto-f}
\begin{align*}
|p(t,a,b,f)-1| &= 
|p(t,a,b,f)-p(t,a_0,b_0,0)| 
\\
&\leq C (|a-a_0|+|b-b_0|) \frac{\log(t_0)}{\log(t)}
\\
&\leq C  \log^m(t_0)\|f\|_{1,m}\frac{\log(t_0)}{\log(t)}.
\end{align*}
This implies \eqref{eq:estp1xx}. The proof of \eqref{eq:estdotp1xx} is similar.
\end{proof}

\begin{remark}
\label{rem:adapt-prop:p1v2}
Proposition~\ref{prop:p1v2}
can be extended to any $\alpha\in(\frac12,1)$ as follows. In the preceding argument, the solution $p$ was found of the form $p=p_0+p_1+p_2+p_3$. The function $p_0$ is introduced by solving part of the problem and achieves an error of order $O(\frac{\log(t_0)}{t\log^2(t)})$, expressed in \eqref{eq:estM1p0} (and related estimates for the derivatives). The corrections $p_1$ and $p_2$ are introduced to gain a $\log(t)$, $p_3$ is introduced to solve the complete equation. This function needs to be in the weighted space given by $\|\dot p_3\|_{1,m-1}<\infty$ and the error $E$ in the weighted space $\|E\|_{1,m}<\infty$. But there is a compatibility condition  $m<-\frac1{2\alpha-1}$. If $\alpha$ is close to $\frac12$, then $m$ needs to be negative with $|m|$ large.
This can be achieved by improving the approximation a large but finite number of steps,
and estimating a large but finite number of derivatives of the corresponding functions.
\end{remark}

\subsection{Proof of Proposition~\ref{prop:p2}}
\label{sec:holder-p}

Let us define
\begin{align*}
[f]_{\gamma,m,\mu}=\sup \frac{t^{\gamma}}{\log^m(t)}\frac{|f(t)-f(s)|}{(t-s)^\mu}
\end{align*}
where the supremum is taken over $t>t_0$ and
$\max(\frac34t,t_0)<s<t$.
We will deduce Proposition~\ref{prop:p2}
from the following variant.

\begin{proposition}
Assume $\alpha \in (\frac58,1)$ and $-4<m<-\frac1{2\alpha-1}$.
For $t_0$ sufficiently large, if  $\|f\|_{1,m}\leq \log^{-3-m}(t_0)$, and
$[f]_{\gamma,1,\mu}<\infty$,
then the solution $p$ of \eqref{prob12b} satisfies
\begin{align*}
|\dot p(t)-\dot p(s)|
\leq C \|f\|_{1,m} \log^{1+m}(t_0) \frac{t-s}{t^2 \log^2(t)}
+ C [f]_{\gamma,1,\mu} \frac{(t-s)^\mu}{t^\gamma}
\end{align*}
for $t>t_0$ and $\max(t_0,\frac34t)<s<t$.
\end{proposition}
We recall that $\|f\|_{1,m}$ has been defined in \eqref{def:norm-gamma-m}.

\begin{proof}
Recall that from Proposition~\ref{prop:p1v2} we have the estimates
\eqref{eq:estp1xx}, \eqref{eq:estdotp1xx}.
Let 
\[
H=\sup \, 
\frac{|\dot p(t)-\dot p(s)|}{q(t,s)}
\]
where the supremum is taken over $t>t_0$ and
$\max(\frac34t,t_0)<s<t$, and
\[
q(t,s)=\|f\|_{1,m} \log^{1+m}(t_0) \frac{t-s}{t^2 \log^2(t)}
+[f]_{\gamma,1,\mu} \frac{(t-s)^\mu}{t^\gamma}.
\]
We will assume that $f$ is non-trivial so $q(t,s)>0$.

We want to show that for some constant $C$ we have
\begin{align}
\label{eq:claimH}
H\leq C.
\end{align}
For $t\geq t_0$ we rewrite
\[
\eta(t) \int_{t_0}^{t-t^\alpha}\frac{\dot p(s)}{t-s} \,d s 
=\eta(t) \int_{t^\alpha}^{t-t_0}\frac{\dot p(t-u)}{u}\,d u.
\]
Let $t\geq t_0$ and $\max(\frac34t,t_0)<s<t$.
Let's compute \eqref{prob12b} evaluated at $t$ minus for \eqref{prob12b}
evaluated at $s$ and write the result as
\begin{align}
\label{eq:dotpt-dotps}
\alpha \log(t)( \dot p(t) - \dot p(s) ) 
+\sum_{j=1}^6 D_j=0
\end{align}
where
\begin{align*}
D_1&=\dot p(s)(\log(t^\alpha)-\log(s^\alpha))
\\
D_2&=H(p(t)) [\dot p(t)-\dot p(s)]
\\
D_3&=
\left[ H(p(t))- H(p(s))\right] \dot p(s)
\\
D_4&=(\eta(t)-\eta(s)) \int_{t^\alpha}^{t-t_0}\frac{\dot p(t-u)}{u}\,d u
\\
D_5&=
\eta(s)\left[
\int_{t^\alpha}^{t-t_0}\frac{\dot p(t-u)}{u}\,d u
-\int_{s^\alpha}^{s-t_0}\frac{\dot p(s-u)}{u}\,d u\right]
\\
D_6&=f(s)-f(t).
\end{align*}

\medskip\noindent
{\em Estimate of $D_1$.}
\begin{align*}
|D_1|
&\leq C \|f\|_{1,m}\frac{\log^{1+m}(t_0)}{t\log^2(t)} 
(\log(t^\alpha)-\log(s^\alpha))
\\
&\leq C \|f\|_{1,m}\frac{\log^{1+m}(t_0)}{t\log^2(t)} 
\frac{t-s}{t}.
\end{align*}

\medskip\noindent
{\em Estimate of $D_2$.}
\begin{align*}
|D_2|&\leq C |\dot p(t)-\dot p(s)|.
\end{align*}

\medskip\noindent
{\em Estimate of $D_3$.}
\begin{align*}
|D_3|
&\leq C |p(t)- p(s)| \, |\dot p(s)|
\\
&\leq C |\dot p(\tilde t) | \, |\dot p(s)| \, (t-s)
\\
&\leq C \|f\|_{1,m}^2\log^{2(1+m)}(t_0) \frac{t-s}{t^2\log^4(t)} .
\end{align*}
But we are assuming $\|f\|_{1,m}\leq \log^{-3-m}(t_0)$,
so 
\begin{align*}
|D_3|
&\leq C \|f\|_{1,m}\log^{m-1}(t_0) \frac{t-s}{t^2\log^4(t)}.
\end{align*}
\medskip\noindent
{\em Estimate of $D_4$.}
\begin{align*}
|D_4|
&\leq 
C\frac{t-s}t \left|\int_{t^\alpha}^{t-t_0}\frac{\dot p(t-u)}{u}\,d u\right|.
\end{align*}
But, by \eqref{eq:estdotp1xx},
\begin{align*}
\left| \int_{t^\alpha}^{t-t_0}\frac{\dot p(t-u)}{u}\,d u \right|
&\leq 
C \|f\|_{1,m} \log^{1+m}(t_0)
\int_{t^\alpha}^{t-t_0}
\frac1{(t-u)\log^2(t-u)}
\frac{1}{u}\,d u .
\end{align*}
We split
\begin{align*}
\int_{t^\alpha}^{t-t_0}
\frac1{(t-u)\log^2(t-u)}
\frac{1}{u}\,d u 
&= \int_{t^\alpha}^{\frac t2} ... 
+ \int_{\frac t2}^{t-t_0}...
\end{align*}
and estimate
\begin{align*}
\int_{t^\alpha}^{\frac t2} 
\frac1{(t-u)\log^2(t-u)}
\frac{1}{u}\,d u 
&\leq C \frac1{t\log^2(t)}
\int_{t^\alpha}^{\frac t2} 
\frac{1}{u}\,d u 
\\
&\leq C \frac1{t\log(t)}
\end{align*}
\begin{align*}
\int_{\frac t2}^{t-t_0}
\frac1{(t-u)\log^2(t-u)}
\frac{1}{u}\,d u 
& \leq 
\frac Ct
\int_{\frac t2}^{t-t_0}
\frac1{(t-u)\log^2(t-u)} \,d u 
\\
& \leq 
\frac C{t\log(t)}.
\end{align*}
It follows that 
\begin{align*}
|D_4|&\leq 
C \|f\|_{1,m} \log^{1+m}(t_0) \frac{t-s}{t^2\log(t)}.
\end{align*}

\medskip\noindent
{\em Estimate of $D_5$.}
In the case of overlapping intervals, $s-t_0>t^\alpha$,
we have
\begin{align*}
\int_{t^\alpha}^{t-t_0}\frac{\dot p(t-u)}{u}\,d u
-\int_{s^\alpha}^{s-t_0}\frac{\dot p(s-u)}{u}\,d u
&=I_1+I_2+I_3
\end{align*}
where
\begin{align*}
I_1&=
-\int_{s^\alpha}^{t^\alpha}
\frac{\dot p(s-u)}{u}\,d u
\\
I_2&=
\int_{t^\alpha}^{s-t_0}
\left[
\frac{\dot p(t-u)}{u}
-\frac{\dot p(s-u)}{u}
\right] d u
\\
I_3&=
\int_{s-t_0}^{t-t_0} \frac{\dot p(t-u)}{u}\,d u.
\end{align*}
We estimate, using \eqref{eq:estdotp1xx},
\begin{align*}
|I_1| 
&\leq C \|f\|_{1,m}\frac{\log^{1+m}(t_0)}{t\log^2(t)}
\int_{s^\alpha}^{t^\alpha}
\frac{1}{u}\,d u
\\
&\leq C \|f\|_{1,m}\frac{\log^{1+m}(t_0)}{t^2\log^2(t)}
(t-s).
\end{align*}
Similarly
\begin{align*}
|I_3|&\leq
C \|f\|_{1,m}\frac{\log^{1+m}(t_0)}{t\log^2(t)}
\int_{s-t_0}^{t-t_0} \frac{1}{u}\,d u
\\
&\leq C \|f\|_{1,m}\frac{\log^{1+m}(t_0)}{t^2\log^2(t)}
(t-s).
\end{align*}
On the other hand
\begin{align*}
|I_2|\leq 
(1-\alpha)\log(t)
H q(t,s) .
\end{align*}
It follows that
\begin{align}
\label{eq:est-D5}
|D_5|\leq
C \|f\|_{1,m}\frac{\log^{1+m}(t_0)}{t^2\log^2(t)}
(t-s)
+(1-\alpha)\log(t)
H q(t,s) .
\end{align}
In the case that the intervals
$[t^\alpha,t-t_0]$ and $[s^\alpha,s-t_0]$ don't overlap,
with a similar computation
we still
get the estimate \eqref{eq:est-D5}.

So from \eqref{eq:dotpt-dotps} we get
\begin{align}
\nonumber
&
\alpha \log(t)
\Bigl(1-\frac C{\log(t_0)}\Bigr)
|\dot p(t)-\dot p(s)|
\\
\nonumber
&\leq (1-\alpha)\log(t)
H q(t,s)
\\
\nonumber
&\quad 
+ C   \|f\|_{1,m} \log^{1+m}(t_0) \frac{t-s}{t^2 \log(t)}
+C [f]_{\gamma,1,\mu} \log(t) \frac{(t-s)^\mu}{t^\gamma}.
\end{align}
Dividing by $\log(t) q(t,s)$ and taking the supremum, we find that 
\[
H \leq \frac{1-\alpha}{\alpha (1-\frac C{\log(t_0)})} H + C.
\]
Since $\alpha>\frac12$, taking $t_0$ large enough we find that
\[
H\leq \tilde C
\]
for some constant $\tilde C$,
which proves \eqref{eq:claimH}.
\end{proof}

\subsection{Basic existence and uniqueness}
\label{sec:nonlocal3}

\begin{proposition}
\label{prop:nonlocal0}
Let $\alpha\in(\frac12,1)$,
$t_0$ be large,
and  $f\in C([t_0,\infty),\C)$. Then \eqref{prob12b} has a unique solution $p\in C^1([t_0,\infty),\C)$. Moreover, if $\alpha$ and $m$ satisfy
\begin{align}
\label{ineq-m}
m>-\frac1\alpha,\quad m>\frac{2-3\alpha}{\alpha(2\alpha-1)},
\end{align}
and if $\|f\|_{1,m}<\infty$, then
\begin{align}
\nonumber
\|p-1\|_{0,m}+\|\dot p\|_{1,m-1}\leq C  \|f\|_{1,m}.
\end{align}
\end{proposition}

\begin{proof}\

\noindent
\textbf{Uniqueness.}
Let $f:[t_0,\infty)\to\C$ be continuous and suppose that $p_1$, $p_2\in C^1([t_0,\infty))$ satisfy \eqref{prob12b}.
We claim that $p_1(t)=p_2(t)$ for all $t>t_0$. Suppose that this is false and
let 
\[
t_*=\inf\{t\geq t_0\,|\, p_1(t)\not=p_2(t)\}.
\]
Since $\eta(t)=0$ for $t\in[t_0,2t_0]$, equation
\eqref{prob12b} is a Cauchy problem for a nonlinear ODE in this interval, and by the classical Gronwall
inequality $p_1=p_2$ on $[t_0,2t_0]$.
Therefore $t_*\geq 2t_0$.
Select $\delta>0$ so that $(t_*+\delta)-(t_*+\delta)^\alpha<t_*$. 
Then $p_j$, $j=1,2$, satisfy the Cauchy problem for the ODE
\begin{align}
\nonumber
\left\{
\begin{aligned}
&\dot p_j(t)\log(t^\alpha)
+H(p_j(t))\dot p_j(t)
= F_j(t)
\quad t>t_0
\\
&
p_j(t_0)=1
\end{aligned}
\right.
\end{align}
where
\[
F_j(t)=f(t)-\eta(t)\int_{t_0}^{t-t^\alpha}\frac{\dot p_j(s)}{t-s}\,d s.
\]
But note that for $t\in [t_0,t_*+\delta]$, $t-t^\alpha\leq t_*$, 
and hence $F_1(t)=F_2(t)$. Then on $[t_0,t_*+\delta]$, both $p_1$, $p_2$
satisfy the same Cauchy problem, and then by standard uniqueness
they agree on $[t_0,t_*+\delta]$.
This contradicts the definition of $t_*$.

\medskip
\noindent
\textbf{Set up.}
Let $f:[t_0,\infty)\to\C$ be continuous.
We let $T>t_0$ and use the norms \eqref{def:norm-gamma-m} but restricted to $[t_0,T]$:
\begin{align}
\label{eq:def-truncated-norm}
\|f\|_{\gamma,m,T} &= \sup_{t_0<t<T} \frac{t^\gamma}{\log^{m}(t)}|f(t)|.
\end{align}
We will use the Leray-Schauder degree theory, and for this we need
an apriori estimate for fixed points of a family of operators.
First, we rewrite the problem.
Write $q=p-1$, so that the equation becomes
\begin{align}
\label{prob12c2}
\left\{
\begin{aligned}
&\dot q(t)\log(t^\alpha)
+H(1+q(t))\dot q(t)
+\eta(t)\int_{t_0}^{t-t^\alpha}\frac{\dot q(s)}{t-s} \,d s 
= f(t)
\quad t>t_0
\\
&
q(t_0)=0
\end{aligned}
\right.
\end{align}
For $t>2t_0$, write
\begin{align*}
\int_{t_0}^{t-t^\alpha}\frac{\dot q(s)}{t-s}ds
&=
\frac{1}{t}\int_{t_0}^{t-t^\alpha}\dot q(s)ds
+\int_{t_0}^{t-t^\alpha}\dot q(s)\Bigl(\frac{1}{t-s}-\frac{1}{t}\Bigr)ds
\\
&=\frac{q(t-t^\alpha)-q(t_0)}{t}
+\int_{t_0}^{t-t^\alpha}\dot q(s)\Bigl(\frac{1}{t-s}-\frac{1}{t}\Bigr)ds.
\end{align*}
Then \eqref{prob12c2} is equivalent to 
\begin{align}
\nonumber
\left\{
\begin{aligned}
&
\dot q(t) \log(t^\alpha)
+\tilde H(q(t))\dot q(t)
+\eta(t)\frac{q(t)}{t}
+ A[q](t)=f(t)
\\
&q(t_0)=0
\end{aligned}
\right.
\end{align}
where $\tilde H(q)=H(1+q)$, 
\begin{align}
\label{eq:defA}
A[q](t)
=\eta(t)\left[
- \frac{q(t)-q(t-t^\alpha)}{t}
+\int_{t_0}^{t-t^\alpha}\dot q(s)\Bigl(\frac{1}{t-s}-\frac{1}{t}\Bigr)ds
\right]
\end{align}
and we have used that $q(t_0)=0$.
Note that $\tilde H$ is smooth (from $\R^2$ to $2\times 2$ matrices), 
Lipschitz continuous and bounded.

We let $\tau\in[0,1]$ and consider the following family of equations
\begin{align}
\label{prob12e}
\left\{
\begin{aligned}
&
\dot q(t) \log(t^\alpha)
+\tilde H(q(t))\dot q(t)
+\eta(t)\frac{q(t)}{t}
+ \tau A[q](t)=f(t)
\\
&q(t_0)=0
\end{aligned}
\right.
\end{align}

\medskip
\noindent
\textbf{A priori estimate.}
We claim, there is $C$ independent of $f$, $q$ and $T$, and
$\tau\in[0,1]$ such that if $q$ satisfies \eqref{prob12e}
then
\begin{align}
\label{eq:aprioriT}
\|q\|_{0,m,T}+\|\dot q\|_{1,m-1,T}\leq C  \|f\|_{1,m,T}.
\end{align}

For the proof, let $\mathcal T$ be the linear operator that to 
$f$ assigns the solution 
of the ODE
\begin{align}
\label{equ-p2}
\dot q(t) \log(t^\alpha) + \frac{\eta(t)}{t} q(t)=f(t).
\end{align}
Then we write the equation \eqref{prob12e} as the fixed point problem
\begin{align}
\label{eq:fixed3b}
q=\mathcal T[f-\tilde A_\tau[q]]
\end{align}
where
\begin{align}
\label{eq:deftildeA}
\tilde A_\tau&=\tau A+A_3\
\\
\nonumber
A_3[q](t)
&=
\tilde H(q(t))\dot q(t).
\end{align}
We work in the space
$X_T$ of $C^1$ functions $q$ on $[t_0,T)$ such that $q(t_0)=0$ and use the notation
$\|q\|_X=\|\dot q\|_{1,m-1,T}$
We assume that we have a solution of \eqref{eq:fixed3b} in $X_T$.

We claim that 
\begin{align}
\label{est-T2b2}
\|\mathcal T[f] \|_X 
\leq \left( 
\frac1\alpha+\frac{1}{\alpha(1+\alpha m)}+O\Bigl(\frac{1}{\log(t_0)}\Bigr)
\right)
\|f\|_{1,m,T}
\end{align}
To prove this, let $q$ satisfy \eqref{equ-p2} so that $q=\mathcal T[f]$.
Write
$$
\dot q(t) \log(t^\alpha) + \frac{\eta(t)}{t} q(t)
=\log(t^\alpha)\left[ \dot q(t) + \frac{\eta(t)}{\alpha t \log(t)}q(t)\right]
$$
and let
$$
a(t)=\frac{\eta(t)}{\alpha t \log(t)}.
$$
Then
$$
u(t)=\exp\Bigl(\int_{t_0}^t a(s)\,ds \Bigr)
$$
is in the kernel.
We have
\begin{align*}
u(t)&=u(3t_0)\frac{\log^{\frac{1}{\alpha}}(t)}{\log^{\frac{1}{\alpha}}(3t_0)},\quad t\geq3t_0
\\
u(t)&=1+O\Bigl(\frac{1}{\log(t_0)}\Bigr), \quad t_0\leq t\leq 3t_0.
\end{align*}
Then $q$ is given by
\begin{align}
\label{formula-p2b2}
q(t)=\frac{1}{u(t)}\int_{t_0}^t  \frac{u(s)}{\log(s^\alpha) } f(s)\,ds.
\end{align}

Let us estimate $q$ given by \eqref{formula-p2b2}.
For $t\geq 3t_0$
\begin{align*}
q(t)&=\frac1{u(gt_0)}\frac{\log^{\frac{1}{\alpha}}(3t_0)}{\log^{\frac{1}{\alpha}}(t)}
\left[ \int_{t_0}^{3t_0}  \frac{u(s)}{\log(s^\alpha) } f(s)\,ds +\int_{3t_0}^t  \frac{u(s)}{\log(s^\alpha) } f(s)\,ds \right].
\end{align*}
But
\begin{align*}
\left|
\int_{t_0}^{3t_0}  \frac{u(s)}{\log(s^\alpha) } f(s)\,ds
\right|
&\leq
\|f\|_{1,m,T}
\int_{t_0}^{3t_0} 
\Bigl(1+O\Bigl(\frac{1}{\log(t_0)}\Bigr)\Bigr)\frac{\log^{\frac1\alpha}(s)}{\log(s^\alpha) } \frac{\log^{m}(s)}{s}\,ds 
\\
&\leq\frac1\alpha\|f\|_{1,m,T}
\Bigl(1+O\Bigl(\frac{1}{\log(t_0)}\Bigr)\Bigr)
\int_{t_0}^{3t_0} \frac{\log^{\frac1\alpha+m-1}(s)}{s}\,ds
\\
&\leq 
\frac1{1+\alpha m}\|f\|_{1,m,T}
\Bigl(1+O\Bigl(\frac{1}{\log(t_0)}\Bigr)\Bigr)
\alpha\left(\log^{\frac1\alpha+m}(3t_0)-\log^{\frac1\alpha+m}(t_0)\right)
\\
&\leq \left(\frac1{1+\alpha m}+O\Bigl(\frac{1}{\log(t_0)}\Bigr)\right)\|f\|_{1,m,T}.
\end{align*}
under the assumption
$$
m>-\frac1\alpha.
$$
We also have
\begin{align*}
&\left|
\frac1{u(3t_0)}\frac{\log^{\frac{1}{\alpha}}(3t_0)}{\log^{\frac{1}{\alpha}}(t)}\int_{3t_0}^t  \frac{u(s)}{\log(s^\alpha) } f(s)\,ds 
\right|
\\
&\quad \leq \frac1\alpha\|f\|_{1,m,T}
\Bigl(1+O\Bigl(\frac{1}{\log(t_0)}\Bigr)\Bigr)
\frac{\log^{\frac{1}{\alpha}}(t_0)}{\log^{\frac{1}{\alpha}}(t)}
\int_{3t_0}^t  \frac{\log^{\frac1\alpha}(s)}{\log^{\frac1\alpha}(3t_0)\log(s^\alpha) } \frac{\log^m(s)}{s}\,ds 
\\
&\quad \leq \frac1\alpha\|f\|_{1,m,T}
\Bigl(1+O\Bigl(\frac{1}{\log(t_0)}\Bigr)\Bigr)
\frac{1}{\log^{\frac{1}{\alpha}}(t)}
\frac1\alpha
\int_{3t_0}^t  \frac{\log^{\frac1\alpha-1+m}(s)}{s }\,ds 
\\
&\quad \leq \frac1\alpha\|f\|_{1,m,T}
\Bigl(1+O\Bigl(\frac{1}{\log(t_0)}\Bigr)\Bigr)
\frac{1}{\log^{\frac{1}{\alpha}}(t)}
\frac{\log^{\frac1\alpha+m}(t)}{\frac1\alpha+m}
\\
&\quad \leq \frac{1}{1+\alpha m}\log^m(t)\|f\|_{1,m,T}
\Bigl(1+O\Bigl(\frac{1}{\log(t_0)}\Bigr)\Bigr).
\end{align*}
Thus for $t>3t_0$
\begin{align*}
|q(t)|\leq  \frac{\log^m(t)}{1+\alpha m}\|f\|_{1,m,T}
\Bigl(1+O\Bigl(\frac{1}{\log(t_0)}\Bigr)\Bigr).
\end{align*}
The calculation for $t_0\leq t\leq 3t_0$ is similar, and we get
\begin{align}
\label{est-p2b2}
|q(t)|\leq  \frac{\log^m(t)}{1+\alpha m}\|f\|_{1,m,T}
\Bigl(1+O\Bigl(\frac{1}{\log(t_0)}\Bigr)\Bigr), \quad t\in[t_0,T].
\end{align}
Using the equation \eqref{equ-p2} and \eqref{est-p2b2} we get
\begin{align}
\nonumber
|\dot q(t)| 
&\leq \frac1{\alpha\log(t)}|f(t)|+\frac1{\alpha t\log(t)}|q(t)|
\\
\nonumber
&\leq 
\frac1\alpha
\frac{\log^{m-1}(t)}{t}
\|f\|_{1,m,T}
+\frac{1}{\alpha(1+\alpha m)}\frac{\log^{m-1}(t)}t\|f\|_{1,m,T}
\Bigl(1+O\Bigl(\frac{1}{\log(t_0)}\Bigr)\Bigr).
\\
\nonumber
&=
\frac{\log^{m-1}(t)}{t}
\left( 
\frac1\alpha+\frac{1}{\alpha(1+\alpha m)}+O\Bigl(\frac{1}{\log(t_0)}\Bigr)
\right)
\|f\|_{1,m,T}
\end{align}
for all $t>t_0$,
and this implies \eqref{est-T2b2}.

Let $\tilde A_\tau$ be defined by \eqref{eq:deftildeA}. We claim that 
\begin{align}
\label{est-A2b2}
\| \tilde A_\tau[q]\|_{1,m,T} \leq \left(1-\alpha+\frac C{\log(t_0)}\right)
\|q\|_{X}.
\end{align}
Indeed, write $A=A_1+A_2$ where $A$ is defined in \eqref{eq:defA}
and 
\begin{align*}
A_1[q](t)
&=\eta(t)\int_{t_0}^{t-t^\alpha}\dot{q}(s)\Bigl(\frac{1}{t-s}-\frac{1}{t}\Bigr)\,ds
\\
A_2[q](t)
&=\eta(t)\frac{ q(t-t^\alpha)- q(t)}{t}.
\end{align*}
Let's estimate $A_1$:
$$
\left|
\int_{t_0}^{t-t^\alpha}\dot{q}(s)\Bigl(\frac{1}{t-s}-\frac{1}{t}\Bigr)\,ds\right|
\leq \|\dot q\|_{1,m-1}\int_{t_0}^{t-t^\alpha}\frac{\log^{m-1}(s)}{s} \frac{s}{(t-s)t}\,ds.
$$
We have
\begin{align*}
&
\int_{t_0}^{t-t^\alpha}\log^{m-1}(s)\frac{1}{(t-s)t}\,ds
\\
&=
\frac1t
\int_{t_0}^{t/2}\log^{m-1}(s) \frac{1}{t-s}\,ds
+\frac1t\int_{t/2}^{t-t^\alpha}\log^{m-1}(s)\frac{1}{t-s}\,ds
\\
&\leq 
\frac{C}{t^2}\int_{t_0}^{t/2}\log^{m-1}(s)\,ds
+ \log^{m-1}(t)\Bigl(1+O(\frac{1}{\log t})\Bigr) \frac{1}{t}\int_{t/2}^{t-t^\alpha}\frac{1}{t-s} \,ds
\\
&\leq C \frac{\log^{m-1}(t)}{t}
+ (1-\alpha)\frac{\log^{m}(t)}{t}\Bigl(1+O(\frac{1}{\log t})\Bigr).
\end{align*}
Therefore
\begin{align}
\nonumber
\left| A_1[q](t)\right|
\leq 
\frac{\log^{m}(t)}{t}
\left(
1-\alpha + \frac{C}{\log t_0}\right)\|\dot q\|_{1,m-1,T}
\end{align}
and this shows that
\begin{align}
\label{eq:est-A12}
\|A_1[q]\|_{1,m,T}\leq 
\left(
1-\alpha + \frac{C}{\log t_0}\right)\|\dot q\|_{1,m-1,T}
\end{align}
Let's estimate $A_2$:
\begin{align}
\nonumber
\left|\eta(t)\frac{ q(t-t^\alpha)- q(t)}{t}\right|
&\leq \eta(t) t^{\alpha-1}\sup_{s\in[t-t^\alpha,t]}|\dot q(s)|
\\
\nonumber
&\leq C\eta(t) t^{\alpha-1} (t-t^\alpha)^{-1} \log^{m-1}(t-t^\alpha)
\|\dot q\|_{1,m-1,T}
\\
\nonumber
&\leq C \frac{\log^m(t)}{t} t_0^{\alpha-1} \|\dot q\|_{1,m-1,T}
\end{align}
Therefore
\begin{align}
\label{eq:est-A22}
\|A_2[q]\|_{1,m,T}\leq C t_0^{\alpha-1} \|\dot q\|_{1,m-1,T}
\end{align}
Finally,
we have
\begin{align}
\nonumber
|A_3[q](t)|
\leq C |\dot q(t)| \leq C \frac{\log^{m-1}(t)}{t} \|q\|_{1,m-1,T}.
\end{align}
Therefore
\begin{align}
\label{eq:est-A32}
\|A_3[q]\|_{1,m,T} 
\leq  \frac C{\log(t_0)}  \|\dot q\|_{1,m-1,T}.
\end{align}
Combining \eqref{eq:est-A12}, \eqref{eq:est-A22} and \eqref{eq:est-A32},
we get a proof of \eqref{est-A2b2}.

Using \eqref{est-T2b2} and \eqref{est-A2b2} we obtain that 
\begin{align*}
\|q\|_X&=\| \mathcal T[f-\tilde A_\tau[q]]\|_X
\\
&\leq 
C\|f\|_{1,m,T}
+\left( 
(1-\alpha)\Bigl(
\frac1\alpha+\frac{1}{\alpha(1+\alpha m)}\Bigr)+\frac{C}{\log(t_0)}
\right) \|q\|_X.
\end{align*}
Note that the assumptions on $\alpha$ and $m$ imply that 
\[
(1-\alpha)\Bigl(
\frac1\alpha+\frac{1}{\alpha(1+\alpha m)}\Bigr)<1.
\]
Then  choosing $t_0$ large, we get that  
\[
\|q\|_X\leq C \|f\|_{1,m,T}.
\]
This implies the claim \eqref{eq:aprioriT}.

\medskip
\noindent
\textbf{Existence on $[t_0,T]$.}
Consider the operator $\mathcal B$ that to $g\in C([t_0,T])$ associates the solution $q=\mathcal A[g]$, of the Cauchy problem
\begin{align}
\nonumber
\left\{
\begin{aligned}
&\dot q(t)\log(t^\alpha)
+\tilde H(q(t))\dot q(t)
+\eta(t)\frac{q(t)}{t}
= g(t)
\quad t>t_0
\\
&
q(t_0)=0
\end{aligned}
\right.
\end{align}
where we recall, $\tilde H$ is a smooth, bounded and Lipschitz 
continuous function on $\R^2$.
We seek a fixed point of the operator
\[
\mathcal F[q,\tau]=\mathcal B[f-\tau A[q]]
\]
on $X_T$ where $A$ is defined in \eqref{eq:defA}.
We work in the closed ball $\overline B_\ell(0)$ of the space $X_T$ defined in the previous step. 
That same step shows that if we take $\ell=C \|f\|_{1,m,T}$ with $C$ a large fixed constant, then $\mathcal F(\cdot,\tau)$ has no fixed point on $\partial B_\ell(0)$.
We claim that $\mathcal F$ is completely continuous. To show this it is enough to show that $A$ defined in \eqref{eq:defA} (which is linear), is compact from 
$C^1([t_0,T])$ into $C([t_0,T])$, because the nonlinear map $\mathcal B$ is continuous from $C([t_0,T])$ into $C^1([t_0,T])$. This follows from the Arzela-Ascoli theorem.

\medskip
\noindent
\textbf{Existence on $[t_0,\infty)$.}
This follows from the existence on $[t_0,T]$ for all $T>t_0$. 
Uniqueness implies that the solutions on $[t_0,T]$ have a limit
as $T\to\infty$.
\end{proof}

From the proof of Proposition~\ref{prop:nonlocal0} we see that we also have the following corollary, which uses the localized norm \eqref{eq:def-truncated-norm}.

\begin{corollary}
Let $\alpha\in(\frac12,1)$,
and assume $m$ satisfies \eqref{ineq-m}.
Then for $t_0$ large,
and  $f\in C([t_0,\infty),\C)$,
the solution $p$ constructed in Proposition~\ref{prop:nonlocal0} satisfies
\begin{align}
\nonumber
\|p-1\|_{0,m,T}+\|\dot p\|_{1,m-1,T}\leq C  \|f\|_{1,m,T}
\end{align}
for any $T>t_0$.
\end{corollary}

Next we combine the result of Section~\ref{sec:eq-p} with Proposition~\ref{prop:nonlocal0} to obtain a localized version of Proposition~\ref{prop:p1v2}.

\begin{corollary}
\label{cor:localized-p1v2}
Assume $\alpha \in (\frac58,1)$ and $-4<m<-\frac1{2\alpha-1}$.
For $t_0$ sufficiently large, 
and $T>t_0$,
if  $\|f\|_{1,m,T}\leq \log^{-3-m}(t_0)$,
then the solution $p$ of 
\eqref{prob12b} satisfies
\begin{align}
\nonumber
|p(t)-1| &\leq C \|f\|_{1,m,T}   \frac{\log^{m+1}(t_0)}{\log(t)}
\\
\nonumber
|\dot p(t)|&\leq C  \|f\|_{1,m,T}  \frac{\log^{m+1}(t_0)}{t\log^2(t)}
\end{align}
for any $t_0\leq t \leq T$.
\end{corollary}
\begin{proof}
For the proof, let $T>t_0$ and assume that $\|f\|_{1,m,T}\leq \log^{-3-m}(t_0)$. We modify the function $f$ by defining
\[
f_T(t)=
\begin{cases}
f(t),&t_0\leq t \leq T
\\
\frac{f(T)}{T}(2T-t),&T\leq t\leq 2T
\\
0,&t\geq 2T.
\end{cases}
\]
Then $\|f_T\|_{1,m}\leq C \|f\|_{1,m,T} \leq C \log^{-3-m}(t_0)$ and we apply Proposition~\ref{prop:p1v2}.
Although strictly speaking that proposition requires $\|f\|_{1,m,T} \leq \log^{-3-m}(t_0)$, it is valid also under the condition $\|f_T\|_{1,m}\leq C \|f\|_{1,m,T} \leq C \log^{-3-m}(t_0)$. 
Then we find that the solution $p_T$ to \eqref{prob12b} with right hand side $f_T$ satisfies 
\begin{align}
\nonumber
|p_T(t)-1| &\leq C \|f_T\|_{1,m}   \frac{\log^{m+1}(t_0)}{\log(t)}
\\
\nonumber
|\dot p_T(t)|&\leq C  \|f_T\|_{1,m}  \frac{\log^{m+1}(t_0)}{t\log^2(t)}.
\end{align}
Let $p$ be the solution of \eqref{prob12b} with right hand side $f$, whih is unique by Proposition~\ref{prop:nonlocal0}.
That proposition also implies uniqueness up to $T$, and we deduce that $p=p_T$ on $[t_0,T]$.

\end{proof}

\section{Estimates for the outer solution}
\label{Sect:est_HE}
In this section we prove Proposition~\ref{prop:Estimate_HE_gen}
for the solution of the heat equation
\begin{align}
\label{heat_equation_bis}
\left\{
\begin{aligned}
\partial_t \psi &= \Delta \psi + g,
&& \text{in }\R^2 \times (t_0, +\infty),\\
\psi(\cdot, t_0) &= 0, && \text{in } \R^2
\end{aligned}
\right.
\end{align}
by analyzing separately the cases defined by the weights
\eqref{sup_g1}--\eqref{sup_g4}.

We let $R(t)=t^\beta$ with $\beta\in(0,\frac12)$.
The first lemma addresses the weight $\varsigma_1$ defined in \eqref{sup_g1}.
\begin{lemma}
Let $\psi_1$ be the solution to \eqref{heat_equation_bis}
given by Duhamel's formula
with right-hand  side $g_1$ satisfying
\[
\abs{g_1(x, t)} \le \frac{\chi_{B_{4R(t)}(0)}}{R(t)^2t^{d_0}}
\]
where $0<d_0<1$.
Then it satisfies
\begin{align}
\label{est_psi1}
\abs{\psi_1(x, t)}
&\le C \frac{\log t}{t^{d_0}}  \exp\Big(-\sfrac{\abs{x}^2}{16t}\Big),
\\
\label{est_gradpsi1}
\abs{\nabla \psi_1(x, t)} &\le C \frac{\log t}{R(t)t^{d_0}}
\left( 1+ \sfrac{\abs{x}}{\sqrt{t}}\right)
\exp\Big(-\sfrac{\abs{x}^2}{16t}\Big),\quad &&
\end{align}
for all $(x, t) \in \R^2 \times (t_0, +\infty)$
and for any $\gamma \in (0, \frac{1}{2})$, there holds
\begin{align}
\label{est_hold_gradpsi1}
\frac{\abs{\nabla \psi_1(x,t) - \nabla \psi_1 (x', t')}}{\big(\abs{x - x'}^2 + \abs{t-t'}\big)^{\gamma}}
&\le C \frac{ \log t}{R(t)^{1 + 2\gamma}t^{d_0}},
\end{align}
for all $(x, t), (x', t') \in \mathcal{D}_{4R}$ such that
$t' \in \big(\max\{t_0, t - (4R(t))^2\}, t\big)$.
\end{lemma}

We prove \eqref{est_psi1}--\eqref{est_hold_gradpsi1} independently. We let
$$
\upnu_1 = d_0 + 2\beta, 
$$
with $d_0$ and $ \beta$ satisfying \eqref{eq:param-linear-outer} so that $\upnu_1 \in (1, 2)$ and
$\abs{g_1(x, t)} \le \frac{\chi_{B_{4R(t)}(0)}}{t^{\upnu_1}}$.
\vspace{0.3cm}

\noindent
\textit{Proof of }\eqref{est_psi1}. Let $(x, t) \in \R^2 \times (t_0, +\infty)$. Duhamel's formula for the heat equation \eqref{heat_equation_bis} is given by
\begin{align*}
    \psi_1(x, t) = \int_{t_0}^t \int_{\R^2} \frac{1}{4 \pi (t - s)} \exp\left(-\frac{|x - z|^2}{4(t - s)} \right) g_1(z, s) dzds.
\end{align*}
Let us first assume that $t \ge 2t_0$. We get that
\begin{align*}
    \abs{\psi_1(x, t)} &\le C\int_{t_0}^t \frac{1}{s^{\upnu_1}} \int_{\abs{z} \le 4R(s)} \frac{1}{(t - s)} \exp\left(-\frac{|x - z|^2}{4(t - s)} \right)dzds\\
    &\le C\bigg(\int_{t_0}^{t - 1} \frac{1}{s^{\upnu_1}} \frac{1}{(t - s)} \int_{\abs{z} \le 4R(s)} dzds \\
    & \hspace{1cm}+ \int_{t - 1}^t \frac{1}{s^{\upnu_1}} \int_{\abs{\tilde{x} - \xi} \le \frac{4R(s)}{\sqrt{t - s}}}  \exp\left(-\frac{|\xi|^2}{4} \right)d\xi ds\bigg), \quad \tilde{x} = \frac{x}{\sqrt{t - s}}, \\
    &\le C\left(\int_{t_0}^{t - 1} \frac{R(s)^2}{s^{\upnu_1}} \frac{1}{(t - s)} ds + \int_{t - 1}^t \frac{1}{s^{\upnu_1}} ds \right).
\end{align*}
The first integral satisfies
\begin{align*}
    \int_{t_0}^{t - 1} \frac{R(s)^2}{s^{\upnu_1}} \frac{1}{(t - s)} ds &\le \int_{t_0}^{t/2} \frac{R(s)^2}{s^{\upnu_1}} \frac{1}{(t - s)} ds + \int_{t/2}^{t - 1} \frac{R(s)^2}{s^{\upnu_1}} \frac{1}{(t - s)} ds\\
    &\le C\frac{R(t)^2}{t^{\upnu_1}}  \log t,
\end{align*}
provided that $\upnu_1 - 2\beta < 1$. Similarly, the second integral estimates as
\begin{align*}
    \int_{t - 1}^t \frac{1}{s^{\upnu_1}} ds \le \frac{C}{t^{\upnu_1}}.
\end{align*}
The estimates of $\psi_1$ for $t$ satisfying either $t \in (t_0 + 1, 2t_0)$ or $t \in (t_0, t_0 + 1)$ are deduced following the same steps. We infer that
\begin{align*}
    \abs{\psi_1(x, t)} \le C(\upnu_1) \frac{R(t)^2}{t^{\upnu_1}} \log t
\end{align*}
for all $(x, t) \in \R^2 \times (t_0, +\infty)$. \\[0.4cm]
Let us now consider $(x, t)$ such that $|x| > 2\sqrt{t}$ and $t \ge 2t_0$. Notice that the space integral in the expression of $\psi_1$ is such that $\abs{z} \le 4R(s) = 4s^\beta$, $\beta < \frac{1}{2}$. Considering $t_0$ sufficiently large, we deduce
$$4R(t) < \sqrt{t}$$
for all $t \ge t_0$ and, in this region, there holds
\begin{align}
    |x - z| \ge \frac{|x|}{2}, \label{condition_comp_sup}
\end{align}
for all $z \in B_{4R(s)}(0)$. It follows that
\begin{align*}
    \abs{\psi_1(x, t)} &\le C\int_{t_0}^t \frac{1}{s^{\upnu_1}} \int_{\abs{z} \le 4R(s)} \frac{1}{(t - s)} \exp\left(-\frac{|x - z|^2}{4(t - s)} \right)dzds\\
    &\le C\int_{t_0}^t \frac{1}{s^{\upnu_1}} \frac{1}{(t - s)} \exp\left(-\frac{|x|^2}{16(t - s)} \right)\int_{\abs{z} \le 4R(s)} dzds\\
    &\le C\int_{t_0}^t \frac{R(s)^2}{s^{\upnu_1}} \frac{1}{(t - s)} \exp\left(-\frac{|x|^2}{16(t - s)} \right)ds\\
    &\le C\bigg(\int_{t_0}^{3t/4} \frac{R(s)^2}{s^{\upnu_1}} \frac{1}{(t - s)} \exp\left(-\frac{|x|^2}{16(t - s)} \right)ds \\
    &\hspace{3cm} + \int_{3t/4}^t \frac{R(s)^2}{s^{\upnu_1}} \frac{1}{(t - s)} \exp\left(-\frac{|x|^2}{16(t - s)} \right)ds \bigg)\\
    &\le C\bigg(\frac{R(t)^2}{t^{\upnu_1}} \exp\left(-\frac{|x|^2}{16t} \right) +  \frac{R(t)^2}{t^{\upnu_1}} \int_{0}^{t/4} \frac{1}{v}\exp\left(-\frac{|x|^2}{16v} \right)dv\bigg)\\
    &\le C\left(\frac{R(t)^2}{t^{\upnu_1}} \exp\left(-\frac{|x|^2}{16t} \right)+  \frac{R(t)^2}{t^{\upnu_1}} \exp\left(-\frac{|x|^2}{4t} \right)\right)
\end{align*}
where the last inequality comes from the fact that
$$g(v) = \frac{1}{v}\exp\left(-\frac{|x|^2}{16v} \right)$$
is increasing in $v$ provided that $v< \frac{\abs{x}^2}{16}$, which holds true on the interval $(0,\frac{t}{4})$ since $\abs{x} > 2\sqrt{t}$. This estimate in the interval $t \in (t_0, 2t_0)$ follows the same steps. We conclude that
\begin{align*}
    \abs{\psi_1(x, t)} \le C(\upnu_1)\frac{R(t)^2\log t }{t^{\upnu_1}} \exp\left(-\frac{|x|^2}{16t}\right)
\end{align*}
for all $(x, t) \in \R^2 \times (t_0, +\infty)$.
This completes the proof of \eqref{est_psi1}.\hfill $\square$

\medskip
\noindent
\textit{Proof of }\eqref{est_gradpsi1}.
Let $(x, t) \in \R^2 \times (t_0, +\infty)$. The gradient of Duhamel's formula is given by
\begin{align*}
    \nabla\psi_1(x, t) = -\int_{t_0}^t \int_{\R^2} \frac{x - z}{8 \pi (t - s)^2} \exp\left(-\frac{|x - z|^2}{4(t - s)} \right) g_1(z, s) dzds.
\end{align*}
Let us assume that $t \ge 2t_0$. We get that
\begin{align*}
    \abs{\nabla\psi_1(x, t)} &\le C\int_{t_0}^t \frac{1}{s^{\upnu_1}} \int_{\abs{z} \le 4R(s)} \frac{\abs{x - z}}{(t - s)^2} \exp\left(-\frac{|x - z|^2}{4(t - s)} \right)dzds\\
    &\le C\int_{t_0}^t \frac{1}{s^{\upnu_1}} \frac{1}{(t - s)^{3/2}} \int_{\abs{z} \le 4R(s)} \frac{\abs{x - z}}{(t - s)^{1/2}} \exp\left(-\frac{|x - z|^2}{4(t - s)} \right)dzds\\
    &\le C\bigg(\int_{t_0}^{t/2} \frac{R(s)^2}{s^{\upnu_1}} \frac{1}{(t - s)^{3/2}}ds  + \int_{t/2}^{t - 1} \frac{R(s)}{s^{\upnu_1}} \frac{1}{(t - s)}ds + \int_{t - 1}^{t} \frac{1}{s^{\upnu_1}} \frac{1}{(t - s)^{1/2}}ds\bigg)\\
    &\le C(\upnu_1) \frac{R(s)}{t^{\upnu_1}} \log t,
\end{align*}
provided that $\upnu_1 - \beta < \frac{3}{2}.$ The proof for $t \in (t_0, 2t_0)$ follows the same steps.\\[0.3cm]
Let us now assume that $\abs{x} > 2\sqrt{t}$ and $t \ge 2t_0$. Considering $t_0$ sufficiently large and using \eqref{condition_comp_sup}, we deduce
\begin{align*}
    \abs{\nabla\psi_1(x, t)} &\le C\int_{t_0}^t \frac{1}{s^{\upnu_1}}\frac{1}{(t - s)^{3/2}} \int_{\abs{z} \le 4R(s)} \frac{\abs{x - z}}{(t - s)^{1/2}} \exp\left(-\frac{|x - z|^2}{4(t - s)} \right)dzds\\
    &\le  C\int_{t_0}^t \frac{R(s)^2}{s^{\upnu_1}}\frac{1}{(t - s)^{3/2}} \frac{\abs{x}}{(t - s)^{1/2}} \exp\left(-\frac{|x|^2}{16(t - s)} \right)ds\\
    &= C\bigg(\int_{t_0}^{3t/4} \frac{R(s)^2}{s^{\upnu_1}}\frac{1}{(t - s)^{3/2}} \frac{\abs{x}}{(t - s)^{1/2}} \exp\left(-\frac{|x|^2}{16(t - s)} \right)ds \\
    &\hspace{1cm} + \int_{3t/4}^t \frac{R(s)^2}{s^{\upnu_1}}\frac{1}{(t - s)^{3/2}} \frac{\abs{x}}{(t - s)^{1/2}} \exp\left(-\frac{|x|^2}{16(t - s)} \right)ds\bigg) \\
    &\le C\bigg(\frac{1}{t^{3/2}} \frac{\abs{x}}{\sqrt{t}} \exp\left(-\frac{|x|^2}{16t} \right)\int_{t_0}^{3t/4} \frac{R(s)^2}{s^{\upnu_1}}ds + \frac{R(t)^2}{t^{\upnu_1}}\int_0^{t/4} \frac{\abs{x}}{v^2} \exp\left(-\frac{|x|^2}{16v} \right)ds\bigg) \\
    &\le  C\frac{1}{t^{\upnu_1}} \frac{R(t)^2}{t^{1/2}}\frac{\abs{x}}{\sqrt{t}}\exp\left(-\frac{|x|^2}{16t}\right), \\
    &\le  C\frac{R(t)}{t^{\upnu_1}}\frac{\abs{x}}{\sqrt{t}}\exp\left(-\frac{|x|^2}{16t}\right),
\end{align*}
provided that $\upnu_1 - \beta < \frac{3}{2}.$ We deduce that
\begin{align*}
    \abs{\nabla \psi_1(x, t)} \le C\frac{R(t) \log t}{t^{\upnu_1}} \left( 1+ \frac{\abs{x}}{\sqrt{t}}\right)\exp\bigg(-\frac{\abs{x}^2}{16t}\bigg)
\end{align*}
for all $(x, t) \in \R^2 \times (t_0, +\infty).$ The estimate for $t \in (t_0, 2t_0)$ is carried out similarly.
This completes the proof of \eqref{est_gradpsi1}. \hfill $\square$

\medskip
\noindent
\textit{Proof of }\eqref{est_hold_gradpsi1}.
Let $t_1 > t_0$, and $R_1:= 5R(t_1) = 5 t_1^\beta$. Let us consider the rescaled function
\begin{align*}
\widetilde{\psi}_1(z, s) = \psi_1(R_1z, R_1^2s), \quad z \in B_1(0), ~s > \frac{t_0}{R_1^2}.
\end{align*}
We distinguish between two regimes: $t_1 - R_1^2 > t_0$ and $t_1 - R_1^2 < t_0$. In the former case, we apply interior gradient estimates for the parabolic equation, whereas the latter case uses estimates for parabolic equations with initial data. We consider the case $t_1 - R_1^2 > t_0$. We denote, for $a, b \in \R^+$,
$$\widetilde{Q}_a := B_a(0) \times \Big(\sfrac{t_1}{R_1^2} - a^2, \sfrac{t_1}{R_1^2}\Big],  \quad \text{and} \quad Q_{bR_1} := B_{bR_1}(0)\times \big(t_1 - (bR_1)^2, t_1\big].$$
Since $\psi_1$ solves the heat equation with right-hand side $g_1$, the rescaled $\widetilde{\psi}_1$ solves
\begin{equation*}
\partial_s \widetilde{\psi}_1 := \Delta_z  \widetilde{\psi}_1  + \widetilde{g}_1, \quad \text{in }\widetilde{Q}_1,
\end{equation*}
with right-hand side $\widetilde{g}_1$ is
\begin{align*}
\widetilde{g}_1(z, s) = R_1^2g_1(R_1 z, R_1^2 s).
\end{align*}
Since $\widetilde{g}_1 \in L^\infty(\widetilde{Q}_1)$, and using \eqref{sup_g1} and \eqref{est_psi1}, standard interior gradient estimates for the heat equation imply that, for any $\gamma \in (0, \frac{1}{2})$, there holds
\begin{align*}
[\nabla_z \widetilde{\psi}_1]_{\gamma;\widetilde{Q}_{\frac{4}{5}}} &\le C\left(\|\widetilde{\psi}_1\|_{L^{\infty}(\widetilde{Q}_1)} + \norm{\widetilde{g}_1}_{L^{\infty}(\widetilde{Q}_1)} \right)\\
&= C\left(\|\psi_1\|_{L^{\infty}(Q_{R_1})} + R_1^2\norm{g_1}_{L^{\infty}(Q_{R_1})} \right),\\
    &\le C \frac{R_1^2}{t_1^{\upnu_1}} \log t_1.
\end{align*}
Additionally, the H\"older semi-norm satisfies
\begin{align*}
[\nabla_z \widetilde{\psi}_1]_{\gamma;\widetilde{Q}_{\frac{4}{5}}} &= \sup_{\widetilde{Q}_{\frac{4}{5}}} \frac{|\nabla_z\widetilde{\psi}_1(z, s) - \nabla_z\widetilde{\psi}_1(z', s')|}{(\abs{z - z'}^2 + \abs{s-s'})^{\gamma}}
\\
&= R_1^{1+2\gamma } \sup_{Q_{\frac{4}{5}R_1}} \frac{|\nabla_x\psi_1(x, t) - \nabla_x\psi_1(x', t')|}{(\abs{x - x'}^2 + \abs{t-t'})^{\gamma}}.
\end{align*}
In particular, we have that
\begin{align*}
\frac{|\nabla_x\psi_1(x, t_1) - \nabla_x\psi_1(x', t')|}{(\abs{x - x'}^2 + \abs{t_1-t'})^{\gamma}} \le C \frac{1}{t_1^{\upnu_1}R_1^{1+2\gamma}} \log t_1,
\end{align*}
for all
\[
x, x' \in B_{\frac{4}{5}R_1}(0) = B_{4R(t_1)}(0) \quad \text{and} \quad t' \in \Big(t_1 - \Big(\sfrac{4}{5}R_1\Big)^2, t_1\Big) = \big(t_1 - (4R(t_1))^2, t_1\big).
\]
This shows \eqref{est_hold_gradpsi1} when $t_1-R_1^2>t_0$.
The proof in the case $t_1 - R_1^2 < t_0$ follows the same steps but using estimates for parabolic equations with initial data and the fact that $\psi_1(\cdot , t_0) = 0$. This concludes the proof of \eqref{est_hold_gradpsi1}. \hfill $\square$
\vspace{0.3cm}

The next lemma deals with the weight $\varsigma_2$ in \eqref{sup_g2}.
\begin{lemma}
\label{estimate_psi2}
Let $\psi_2$ be the solution to \eqref{heat_equation_bis}
given by Duhamel's formula
with right-hand
side $g_2$ satisfying
\[
|g_2(x,t)|\leq  \frac{\chi_{B_{4R(t)}(0)}}{R(t)^{1 + 2\gamma}t^{d_0}\big(1 + \abs{x}\big)^{a_2}},
\]
where $a_2 > 2$, $0<d_0<1$ and $\gamma\in(0,\frac12)$.
Then it satisfies
\begin{align}
\label{est_psi2}
\abs{\psi_2(x, t)}
&\le C \frac{\log t}{t^{d_0}}  \exp\Big(-\sfrac{\abs{x}^2}{16t}\Big),
\\
\label{est_gradpsi2}
\abs{\nabla \psi_2(x, t)}
&\le C \frac{\log t}{R(t)t^{d_0}} \left( 1+ \sfrac{\abs{x}}{\sqrt{t}}\right)\exp\Big(-\sfrac{\abs{x}^2}{16t}\Big),
\end{align}
for all $(x, t) \in \R^2 \times (t_0, +\infty)$ and
\begin{align}
\label{est_hold_gradpsi2}
\frac{\abs{\nabla \psi_2(x,t) - \nabla \psi_2 (x', t')}}{\big(\abs{x - x'}^2 + \abs{t-t'}\big)^{\gamma}}
&\le C \frac{ \log t}{R(t)^{1 + 2\gamma}t^{d_0}},
\end{align}
for all $(x, t), (x', t') \in \mathcal{D}_{4R}$ such that
$t' \in \big(\max\{t_0, t - (4R(t))^2\}, t\big)$.
\end{lemma}
We prove \eqref{est_psi2}--\eqref{est_hold_gradpsi2} independently. We let
$$
\upnu_2 = d_0 + \beta(1 + 2\gamma), 
$$
with $d_0, \beta$ and $\gamma$ satisfying \eqref{eq:param-linear-outer} so that $\upnu_2 \in (1, 2)$ and
$\abs{g_2(x, t)} \le \frac{\chi_{B_{4R(t)}(0)}}{t^{\upnu_2}\big(1 + \abs{x}\big)^{a_2}}$.

\medskip
\noindent
\textit{Proof of \eqref{est_psi2}}. Let $(x, t) \in \R^2 \times (t_0, +\infty)$.
By Duhamel's formula
\begin{align*}
    \psi_2(x, t) = \int_{t_0}^t \int_{\R^2} \frac{1}{4 \pi (t - s)} \exp\left(-\frac{|x - z|^2}{4(t - s)} \right) g_2(z, s) dzds.
\end{align*}
Let us assume that $t \ge 2t_0$. We get that
\begin{align*}
    \abs{\psi_2(x, t)} &\le C\int_{t_0}^t \frac{1}{s^{\upnu_2}(t - s)}\int_{\abs{z} \le 4R(s)}  \exp\left(-\frac{\abs{x - z}^2}{4(t - s)}\right)\frac{1}{(1 + |z|)^{a_2}}dzds\\
    &\le  C\int_{t_0}^t \frac{1}{s^{\upnu_2}}\int_{\R^2}   \exp(-|\xi|^2/4)\frac{1}{(1 + |x - \sqrt{t - s}\xi|)^{a_2}}d\xi ds.
\end{align*}
Let $\tau = \sqrt{t - s}$ and
\begin{align*}
    A(\tau) := \int_{\R^2}   \exp(-|\xi|^2/4)\frac{1}{(1 + |x - \tau \xi|)^{a_2}}d\xi.
\end{align*}
We claim that
\begin{align*}
    A(\tau)  \le C\frac{1}{1 + \tau^2}.
\end{align*}
For $\tau \le 1$, we get
\begin{align*}
    A(\tau) \le \int_{\R^2}   \exp(-|\xi|^2/4) d\xi \le C
\end{align*}
and for $\tau > 1$,
\begin{align*}
    A(\tau) &\le \int_{\R^2} \frac{1}{(1 + |x - \tau\xi|)^{a_2}}ds \\
    &\le \frac{1}{\tau^2}\int_{\R^2}  \frac{1}{(1 + |x - \zeta|)^{a_2}}d\zeta\\
    &\le \frac{C}{\tau^2}
\end{align*}
Therefore,
\begin{align*}
    \abs{\psi_2(x, t)} &\le C \int_{t_0}^t \frac{1}{s^{\upnu_2}} \frac{1}{1 + t - s} ds \\
    &\le C \frac{\log t}{t^{\upnu_{2, 1}}},
\end{align*}
for all $\upnu_{2, 1}$ such that
$$\upnu_{2, 1} < \min\{1, \upnu_2\}.$$
Let us now consider $(x, t)$ such that $|x| > 2 \sqrt{t}$ and $t \ge 2t_0$. Assuming $t_0$ large enough and using \eqref{condition_comp_sup}, we deduce that
\begin{align*}
    \abs{\psi_2(x, t)} &\le C\int_{t_0}^t \frac{1}{s^{\upnu_2}} \int_{\abs{z} \le 4R(s)} \frac{1}{(t - s)} \exp\left(-\frac{|x - z|^2}{4(t - s)} \right)\frac{1}{(1 + |z|)^{a_2}}dzds\\
    &\le C\int_{t_0}^t \frac{1}{s^{\upnu_2}} \frac{1}{(t - s)} \exp\left(-\frac{|x|^2}{16(t - s)} \right) \int_{\abs{z} \le 4R(s)} \frac{1}{(1 + |z|)^{a_2}} dzds\\
    &\le C \frac{\log t}{t^{\upnu_{2, 1}}} \exp\left(-\frac{|x|^2}{16t} \right),
\end{align*}
for all $\upnu_{2, 1}$ such that $\upnu_{2, 1} < \min\{1, \upnu_2\}$. We deduce that
\begin{align*}
    \abs{\psi_2(x, t)} \le C \frac{\log t}{t^{\upnu_{2, 1}}} \exp\left(-\frac{|x|^2}{16t} \right)
\end{align*}
for all $(x, t) \in \mathbb R^2 \times (t_0, +\infty)$. The lemma follows by taking $\upnu_{2, 1} = d_0 < 1$. \hfill $\square$

\medskip
\noindent
\textit{Proof of \eqref{est_gradpsi2}}. Let $(x, t) \in \R^2 \times (t_0, +\infty)$. The gradient of Duhamel's formula is given by
\begin{align*}
    \nabla\psi_2(x, t) = -\int_{t_0}^t \int_{\R^2} \frac{x - z}{8 \pi (t - s)^2} \exp\left(-\frac{|x - z|^2}{4(t - s)} \right) g_2(z, s) dzds.
\end{align*}
Let us assume that $t \ge 2t_0$. We get that
\begin{align*}
    \abs{\nabla \psi_2(x, t)} &\le C\int_{t_0}^t \frac{1}{s^{\upnu_2}(t - s)^{3/2}}\int_{\abs{z} \le 4R(s)}  \frac{|x - z|}{(t - s)^{1/2}}\exp\left(-\frac{\abs{x - z}^2}{4(t - s)}\right)\frac{1}{(1 + |z|)^{a_2}}dzds\\
    &\le  C\int_{t_0}^t \frac{1}{s^{\upnu_2}(t - s)^{1/2}}\int_{\R^2}   \abs{\xi}\exp(-|\xi|^2/4)\frac{1}{(1 + |x - \sqrt{t - s}\xi|)^{a_2}}d\xi ds.
\end{align*}
Let $\tau = \sqrt{t - s}$ and
\begin{align*}
    B(\tau) := \int_{\R^2}   \abs{\xi}\exp(-|\xi|^2/4)\frac{1}{(1 + |x - \tau \xi|)^{a_2}}d\xi.
\end{align*}
We claim that
\begin{align*}
    B(\tau)  \le C\frac{1}{1 + \tau^2}.
\end{align*}
For $\tau \le 1$, we get
\begin{align*}
    B(\tau) \le \int_{\R^2}   \abs{\xi}\exp(-|\xi|^2/4) d\xi \le C,
\end{align*}
and for $\tau > 1$,
\begin{align*}
    B(\tau) &\le \int_{\R^2} \frac{1}{(1 + |x - \tau\xi|)^{a_2}}ds \\
    &\le \frac{1}{\tau^2}\int_{\R^2}  \frac{1}{(1 + |x - \zeta|)^{a_2}}d\zeta\\
    &\le \frac{C}{\tau^2}.
\end{align*}
Therefore,
\begin{align*}
    \abs{\nabla \psi_2(x, t)} &\le C \int_{t_0}^t \frac{1}{s^{\upnu_2}(t - s)^{1/2}} \frac{1}{1 + t - s} ds \\
    &\le C \frac{\log t}{t^{\upnu_{2, 2}}},
\end{align*}
for all $\upnu_{2, 2}$ such that
$$\upnu_{2, 2} \le \min\Big\{\sfrac{3}{2}, \upnu_2\Big\}.$$
Let us now consider $(x, t)$ such that $|x| > 2 \sqrt{t}$ and $t \ge 2t_0$. We deduce that
\begin{align*}
    \abs{\nabla \psi_2(x, t)} &\le C\int_{t_0}^t \frac{1}{s^{\upnu_2}} \int_{\abs{z} \le 4R(s)} \frac{|x - z|}{(t - s)^2} \exp\left(-\frac{|x - z|^2}{4(t - s)} \right)\frac{1}{(1 + |z|)^{a_2}}dzds\\
    &\le C\int_{t_0}^t \frac{1}{s^{\upnu_2}} \frac{\abs{x}}{(t - s)^2} \exp\left(-\frac{|x|^2}{16(t - s)} \right) \int_{\abs{z} \le 4R(s)} \frac{1}{(1 + |z|)^{a_2}} dzds\\
    &\le C \frac{\log t}{t^{\upnu_{2, 2}}} \exp\left(-\frac{|x|^2}{16t} \right),
\end{align*}
for all $\upnu_{2, 2}$ such that
$$\upnu_{2, 2} \le \min\Big\{\sfrac{3}{2}, \upnu_2\Big\}.$$
We deduce that
\begin{align*}
\abs{\nabla\psi_2(x, t)} \le C \frac{\log t}{t^{\upnu_{2, 2}}}
\exp\left(-\frac{|x|^2}{16t} \right)
\end{align*}
for all $(x, t) \in \mathbb R^2 \times (t_0, +\infty)$. The lemma follows by taking $\upnu_{2, 2} = d_0 + \beta < \frac{3}{2}$. \hfill $\square$

\medskip
\noindent
\textit{Proof of \eqref{est_hold_gradpsi2}.} We establish the decay of the Hölder semi-norm of $\nabla \psi_2$ by using the heat kernel. To this end, we treat the space and time variations separately. More precisely, we show that
\begin{align}
    \frac{\abs{\nabla \psi_2(x_1, t) - \nabla \psi_2(x_2, t)}}{\abs{x_1 - x_2}^{2\gamma}} \le C \frac{\log t}{t^{\upnu_{2, 3}}} \label{est_hold_gradpsi2_space}
\end{align}
for all $x_1, x_2 \in B_{4R(t)}(0)$ and $t> t_0$, and that
\begin{align}
    \frac{\abs{\nabla \psi_2(x, t_1) - \nabla \psi_2(x, t_2)}}{\abs{t_1 - t_2}^{\gamma}} \le C \frac{\log t_2}{t_2^{\upnu_{2, 3}}} \label{est_hold_gradpsi2_time}
\end{align}
for $x \in B_{4R(t_2)}(0)$, $t_2 > t_0$, and $t_1 \in (t_2 - (4R(t_2))^2, t_2)$. The exponent $\upnu_{2, 3}$ is required to satisfy
$$\upnu_{2, 3} < \min \Big\{\sfrac{3}{2}, \upnu_2\Big\}.$$
The proof of the lemma follows by choosing $\upnu_{2, 3} = d_0 + \beta(1 + 2\gamma) < \frac{3}{2}$.\\[0.3cm]
We begin by proving \eqref{est_hold_gradpsi2_space}. Let $t > 2t_0$, $x_1, x_2 \in B_{4R(t)}(0)$ and $L = \abs{x_1 - x_2}$. Let us first assume that $L < \frac{1}{2}$. We have
\begin{align*}
    |\nabla \psi_2(x_1, t) &- \nabla \psi_2(x_2, t)| \\
    &\le \int_{t_0}^t \int_{\abs{z} \le 4R(s)} \abs{\nabla_x G(x_1 - z, t- s) - \nabla_x G(x_2 - z, t- s)} \frac{1}{s^{\upnu_2}(1 + \abs{z})^{a_2}}dzds\\
    &\le \int_{t_0}^{t - L^2} \int_{\abs{z} \le 4R(s)} \abs{\nabla_x G(x_1 - z, t- s) - \nabla_x G(x_2 - z, t- s)} \frac{1}{s^{\upnu_2}(1 + \abs{z})^{a_2}}dzds\\
    &\quad+ \int_{t - L^2}^t \int_{\abs{z} \le 4R(s)} \abs{\nabla_x G(x_1 - z, t- s) - \nabla_x G(x_2 - z, t- s)} \frac{1}{s^{\upnu_2}(1 + \abs{z})^{a_2}}dzds\\
    &=: I_1 + I_2
\end{align*}
where $G(x, t)$ is the heat kernel. We find that
\begin{align*}
    \abs{\nabla_x G(x_1 - z, t- s) - \nabla_x G(x_2 - z, t- s)} \le L \int_0^1 \abs{D^2_x G(x_v -z, t - s)}dv
\end{align*}
for $x_v = vx_1 + (1 - v) x_2$ and
\begin{align*}
    (D^2_x G(x, t))_{ij} = -\frac{\delta_{ij}}{8\pi t^2} \exp\bigg(-\frac{\abs{x}^2}{4t}\bigg) + \frac{x_i x_j}{16\pi t^3}\exp\Big(-\frac{\abs{x}^2}{4t}\Big).
\end{align*}
Therefore, the first integral $I_1$ estimates as
\begin{align*}
    I_1 &\le CL \int_0^1\int_{t_0}^{t - L^2} \frac{1}{s^{\upnu_2}(t - s)^2}\int_{\abs{z} \le 4R(s)} \bigg(1 + \frac{\abs{x_v - z}^2}{(t- s)}\bigg)\exp\bigg(-\frac{\abs{x - z}^2}{4(t -s)}\bigg)\frac{1}{(1 + \abs{z})^{a_2}}dzdsdv\\
    &\le CL \bigg\{\int_0^1\int_{t_0}^{t/2} \frac{1}{s^{\upnu_2}(t - s)^2}\int_{\abs{z} \le 4R(s)} \frac{1}{(1 + \abs{z})^{a_2}}dzdsdv\\
    &\quad + \int_0^1\int_{t/2}^{t - L^2} \frac{1}{s^{\upnu_2}(t - s)}\int_{\abs{z} \le 4R(s)} \big(1 + \abs{\xi}^2\big)\exp\big(-\abs{\xi}^2\big)d\xi dsdv\bigg\}\\
    &\le  CL \bigg\{\int_{t_0}^{t/2} \frac{1}{s^{\upnu_2}(t - s)^2} ds + \int_{t/2}^{t - L^2} \frac{1}{s^{\upnu_2}(t - s)} ds\bigg\}\\
    &\le CL \bigg\{ \frac{1}{t^2} + \frac{1}{t^{\upnu_2}} |\log(L^2)| \bigg\}\\
    &\le C \frac{\log t}{t^{\upnu_{2, 3}}} L |\log L| ,
\end{align*}
for
\begin{align*}
    \upnu_{2, 3} \le \min \{2, \upnu_2\}.
\end{align*}
Furthermore, the second integral estimates as
\begin{align*}
    I_2 &\le C\int_{t - L^2}^t \int_{\abs{z} \le 4R(s)} \big(\abs{\nabla_x G(x_1 - z, t- s)} + \abs{\nabla_x G(x_2 - z, t- s)}\big) \frac{1}{s^{\upnu_2}(1 + \abs{z})^{a_2}}dzds\\
    &\le C\int_{t - L^2}^t \frac{1}{s^{\upnu_2}(t - s)^{1/2}} \int_{\R^2} \abs{\xi} \exp\big(-\abs{\xi}^2\big)d\xi ds\\
    &\le \frac{C}{t^{\upnu_2}} L
\end{align*}
As we assumed $L < \frac{1}{2}$, we deduce
\begin{align*}
    \abs{\nabla \psi_2(x_1, t) - \nabla \psi_2(x_2, t)} \le C\frac{\log t}{t^{\upnu_{2, 3}}} \abs{x_1 - x_2}^{2\gamma}
\end{align*}
for any $\gamma \in (0, \frac{1}{2})$ and $\upnu_{2, 3} < \min \{2, \upnu_2\}$. For $L > \frac{1}{2}$, we find, using \eqref{est_gradpsi2}, that
\begin{align*}
    \abs{\nabla \psi_2(x_1, t) - \nabla \psi_2(x_2, t)} &\le \abs{\nabla \psi_2(x_1, t)} + \abs{\nabla \psi_2(x_2, t)}\\
    &\le C\frac{\log t}{t^{\upnu_{2, 3}}} \\
    &\le C\abs{x_1 - x_2}^{2\gamma}\frac{\log t}{t^{\upnu_{2, 3}}}
\end{align*}
for any $\upnu_{2, 3} \le \min\{\frac{3}{2}, \upnu_2\}.$ This proves \eqref{est_hold_gradpsi2_space}.\\[0.3cm]
Let $t_1, t_2$ be such that $t_0 < t_1 < t_2$ and $t_2 > 2t_0$, and let $x \in B_{4R(t_2)}(0)$. We further assume that $\frac{3t_2}{4} < t_1 < t_2$ and $t_2 - t_1 < 1$. We get
\begin{align*}
    |\nabla \psi_2(x, t_2) &- \nabla \psi_2(x, t_1)| \\
    &\le \int_{t_0}^{t_1 - (t_2 - t_1)} \int_{\abs{z} \le 4R(s)} \abs{\nabla_x G(x - z, t_2- s) - \nabla_x G(x - z, t_1- s)} \frac{1}{s^{\upnu_2}(1 + \abs{z})^{a_2}}dzds\\
    & \quad + \int_{t_1 - (t_2 - t_1)}^{t_1} \int_{\abs{z} \le 4R(s)} \abs{\nabla_x G(x - z, t_1- s)} \frac{1}{s^{\upnu_2}(1 + \abs{z})^{a_2}}dzds\\
    &\quad + \int_{t_1 - (t_2 - t_1)}^{t_2} \int_{\abs{z} \le 4R(s)} \abs{\nabla_x G(x - z, t_2- s)} \frac{1}{s^{\upnu_2}(1 + \abs{z})^{a_2}}dzds\\
    &=: J_1 + J_2 + J_3
\end{align*}
We let $t_v = vt_1 + (1 - v)t_2$, for $v \in (0, 1)$. The integral $J_1$ estimates as
\begin{align*}
    J_1 &\le (t_2 - t_1) \int_0^1 \int_{t_0}^{t_1 - (t_2 - t_1)} \int_{\abs{z} \le 4R(s)} \abs{\partial_t\nabla_x G(x - z, t_v- s)} \frac{1}{s^{\upnu_2}(1 + \abs{z})^{a_2}}dzds\\
    &\le (t_2 - t_1) \int_0^1 \int_{t_0}^{t_1 - (t_2 - t_1)} \frac{1}{s^{\upnu_2}(t_v - s)^{5/2}} \\
    &\hspace{4cm}\int_{\abs{z} \le 4R(s)} \bigg(\frac{\abs{x - z}}{(t_v - s)^{1/2}} +\frac{\abs{x - z}^3}{(t_v - s)^{3/2}}\bigg)e^{-\frac{\abs{x - z}^2}{4(t_v - s)}} \frac{1}{(1 + \abs{z})^{a_2}}dzdsdv\\
    &\le C(t_2 - t_1) \bigg\{\int_{t_0}^{t_2/2} \frac{1}{s^{\upnu_2}(t_2 - s)^{5/2}} \int_{\abs{z} \le 4R(s)} \frac{1}{(1 + \abs{z})^{a_2}}dzds \\
    & \hspace{3cm}  + \int_{t_2/2}^{t_1 - (t_2 - t_1)} \frac{1}{s^{\upnu_2}(t_2 - s)^{3/2}} \int_{\abs{z} \le 4R(s)} \big(\abs{\xi} + \abs{\xi}^3\big) \exp\big(-\abs{\xi}^2\big)d\xi ds\bigg\}\\
    &\le C(t_2 - t_1) \bigg\{\frac{1}{t_2^{5/2}} + \frac{1}{t_2^{\upnu_2}(t_2 - t_1)^{1/2}}\bigg\}\\
    &\le C\frac{(t_2 - t_1)^{1/2}}{t^{\upnu_{2, 3}}}, \quad \upnu_{2, 3} \le \min\{2, \upnu_{2}\}.
\end{align*}
where, in the third inequality, we used that
$$\frac{t_2 - s}{t_v - s} \le C, \quad \text{for all }s \in (t_0, t_1 - (t_2 - t_1))$$
for $C$ independent of $v, t_1, t_2$. For $J_2$, we get
\begin{align*}
    J_2 &\le \int_{t_1 - (t_2 - t_1)}^{t_1} \frac{1}{s^{\upnu_2}(t_1 - s)^{1/2}}\int_{\abs{z} \le 4R(s)} \abs{\xi}\exp\big(-\abs{\xi}^2\big)d\xi ds\\
    &\le C\frac{(t_2 - t_1)^{1/2}}{t_2^{\upnu_2}}\\
    &\le C\frac{(t_2 - t_1)^{\gamma}}{t_2^{\upnu_2}},
\end{align*}
for all $\gamma \in (0, \frac{1}{2})$. The same estimate holds for $J_3$. Finally, for $t_2 - t_1 > 1$ and $\frac{3t_2}{4} < t_1 < t_2$, we get, using \eqref{est_gradpsi2}, that
\begin{align*}
\abs{\nabla \psi_2(x, t_2) - \nabla \psi_2(x, t_1)} &\le \abs{\nabla \psi_2(x, t_2)} + \abs{\nabla \psi_2(x, t_1)}\\
&\le C\frac{(t_2 - t_1)^{\gamma}}{t^{\upnu_{2, 3}}} \log t
\end{align*}
for $\upnu_{2, 3} \le \min\{\frac{3}{2}, \upnu_2\}$. This proves \eqref{est_hold_gradpsi2_time} and concludes the proof of the lemma.\hfill $\square$

The next lemma deals with the weight $\varsigma_3$ in \eqref{sup_g4}.
\begin{lemma}
Let $\psi_3$ be the solution to \eqref{heat_equation_bis}
given by Duhamel's formula
with right-hand  side $g_3$ satisfying
\[
\abs{g_3(x, t)} \leq
\frac{(1 - \chi_{B_{R(t)}(0)})}{t^{d_0 + \theta}(1 + \abs{x})^2},
\]
with $0 < d_0 + \theta < 1$. Then $\psi_3$ satisfies
\begin{align}
\label{est_psi4}
\abs{\psi_3(x, t)} &\le C \frac{(\log t)^2}{t^{d_0 + \theta}\big(1 + \frac{\abs{x}}{\sqrt{t}}\big)^2}  ,
\\
\label{est_gradpsi4}
\abs{\nabla \psi_3(x, t)} &\le C \frac{(\log t)^2}{R(t)t^{d_0 + \theta}\big(1 + \frac{\abs{x}}{\sqrt{t}}\big)^2} ,
\end{align}
for all $(x, t) \in \R^2 \times (t_0, +\infty)$ and for
$\gamma \in (0, \frac{1}{2})$, there holds
\begin{align}
\frac{\abs{\nabla \psi_3(x,t) - \nabla \psi_3 (x', t')}}{\big(\abs{x - x'}^2 + \abs{t-t'}\big)^{\gamma}} &\le C \frac{(\log t)^2}{R(t)^{1 + 2\gamma}t^{d_0 + \theta}}, \label{est_hold_gradpsi4}
\end{align}
for all $(x, t), (x', t') \in \mathcal{D}_{4R}$ such that
$t' \in \big(\max\{t_0, t - (4R(t))^2\}, t\big)$.
\end{lemma}
To prove \eqref{est_psi4}--\eqref{est_gradpsi4}, we divide $g_3$ as
\begin{align*}
    \abs{g_3(x, t)} &\le  \frac{\chi_{\abs{x} \ge R(t)}}{t^{d_0 + \theta}(1 + \abs{x})^2}\\
    &\le \frac{\chi_{\{R(t) \le \abs{x} \le \sqrt{t}\}}}{R(t)^2t^{d_0 + \theta} \big(1 + \frac{\abs{x}}{R(t)}\big)^2}  + \frac{\chi_{\{\abs{x} \ge \sqrt{t}\}}}{t^{d_0 + \theta + 1}\big(1 +  \frac{\abs{x}}{\sqrt{t}}\big)^2} =: \tilde{g}_{31}(x, t) + \tilde{g}_{32}(x, t)
\end{align*}
We denote $\psi_{31}$ and $\psi_{32}$ the solutions to the heat equation \eqref{heat_equation_bis}, given by Duhamel's formula, with forcings $\tilde{g}_{31}$ and $\tilde{g}_{32}$, respectively. We claim that these solutions satisfy
\begin{alignat}{3}
    \abs{\psi_{31}(x, t)} &\le C \frac{(\log t)^2}{t^{d_0 + \theta}}  \exp\Big(-\sfrac{\abs{x}^2}{16t}\Big), \quad &&   \abs{\nabla\psi_{31}(x, t)} &&\le C \frac{(\log t)^2}{R(t)t^{d_0 + \theta}} \Big(1 + \sfrac{\abs{x}}{\sqrt t}\Big)\exp\Big(-\sfrac{\abs{x}^2}{16t}\Big)    , \label{est_psi41}\\
    \abs{\psi_{32}(x, t)} &\le C \frac{(\log t)^2}{t^{d_0 + \theta}}  \frac{1}{\big(1 + \frac{\abs{x}}{\sqrt{t}}\big)^2}, \quad &&   \abs{\nabla\psi_{32}(x, t)} &&\le C \frac{(\log t)^2}{R(t)t^{d_0 + \theta}} \frac{1}{\big(1 + \frac{\abs{x}}{\sqrt{t}}\big)^2}, \label{est_psi42}
\end{alignat}
for all $(x, t) \in \R^2 \times (t_0, +\infty)$. The estimate \eqref{est_hold_gradpsi4} follows from the scaling argument used in the proof of \eqref{est_hold_gradpsi1}.

\medskip
\noindent
\textit{Proof of \eqref{est_psi41}.} We prove the estimate for $\psi_{31}$, the proof for the gradient is similar.
Let $(x, t) \in \R^2 \times (t_0, +\infty)$ and $\psi_{31}$ given by
\begin{align*}
\psi_{31}(x, t) = \int_{t_0}^t \int_{\R^2} \frac{1}{4 \pi (t - s)} \exp\left(-\frac{|x - z|^2}{4(t - s)} \right) \tilde{g}_{31}(z, s) dzds.
\end{align*}
Let us assume that $t \ge 2t_0$. We get that
\begin{align*}
    \abs{\psi_{31}(x, t)} &\le C\int_{t_0}^t \frac{1}{R(s)^2s^{d_0 + \theta}}\int_{\abs{z} \le \sqrt{s}} \frac{1}{t - s} \exp\left(-\frac{|x - z|^2}{4(t - s)} \right)  \frac{1}{\big(1 + \frac{|z|}{R(s)}\big)^2}dzds\\
    &\le C\int_{t_0}^t \frac{1}{R(s)^2s^{d_0 + \theta}}\int_{\abs{x - \zeta} \le \sqrt{s}} \frac{1}{t - s} \exp\left(-\frac{|\zeta|^2}{4(t - s)} \right)  \frac{1}{\big(1 + \frac{\abs{x - \zeta}}{R(s)}\big)^2}dzds
\end{align*}
For $\tau = \sqrt{t - s}$, we claim that
\begin{align*}
    A(\tau, s) := \frac{1}{\tau^2} \int_{\abs{x - \zeta} \le \sqrt{s}} \exp\left(-\frac{|\zeta|^2}{4\tau^2} \right)\frac{1}{\big(1 + \frac{\abs{x - \zeta}}{R(s)}\big)^2}d\zeta \le \frac{R(s)^2}{1 + \tau^2} \log s.
\end{align*}
For $\tau > 1,$ we find
\begin{align*}
    A(\tau, s) &\le \frac{1}{\tau^2} \int_{\abs{x - \zeta} \le \sqrt{s}} \frac{1}{\big(1 + \frac{\abs{x - \zeta}}{R(s)}\big)^2}d\zeta\\
    &= \frac{R(s)^2}{\tau^2} \int_{\abs{\xi} \le \frac{\sqrt{s}}{R(s)}} \frac{1}{\big(1 + \abs{\xi}\big)^2}d\zeta\\
    &\le C\frac{R(s)^2}{\tau^2} \log s
\end{align*}
and for $\tau < 1,$
\begin{align*}
    A(\tau, s) \le C\int_{\R^2} e^{-\abs{\xi}^2/4} d\xi \le C
\end{align*}
Therefore,
\begin{align*}
    \abs{\psi_{31}(x, t)} &\le C\int_{t_0}^t \frac{1}{R(s)^2s^{d_0 + \theta}} \frac{R(s)^2}{1 + t - s} \log(s) ds\\
    &\le C\frac{1}{t^{d_0 + \theta}}(\log t)^2.
\end{align*}
Let us now consider $(x, t)$ such that $\abs{x} \ge 2 \sqrt{t}$. Notice that the space integral is integrated over the domain $\abs{z} \le \sqrt{s}$. In particular, there holds $\abs{z} \le \sqrt{t},$ and
\begin{align*}
    \int_{\abs{z} \le \sqrt{s}} \frac{1}{t - s} \exp\left(-\frac{|x - z|^2}{4(t - s)} \right)  \frac{1}{\big(1 + \frac{|z|}{R(s)}\big)^2}dz &\le C\frac{e^{-\frac{|x|^2}{16(t - s)}}}{t - s} \int_{\abs{z} \le \sqrt{s}} \frac{1}{\big(1 + \frac{|z|}{R(s)}\big)^2}dz\\
    &\le C\frac{R(s)^2}{t - s}e^{-\frac{|x|^2}{16(t - s)}} \log s,
\end{align*}
which implies that
\begin{align*}
    \abs{\psi_{31}(x, t)} &\le C\int_{t_0}^t  \frac{1}{s^{d_0 + \theta}}\frac{1}{t - s} \exp\left(-\frac{|x|^2}{16(t - s)} \right) \log s ds\\
    &\le C\frac{(\log t)^2}{t^{d_0 + \theta}} \exp\left(-\frac{|x|^2}{16t} \right),
\end{align*}
which concludes the proof
of \eqref{est_psi41}.

The proof of the estimate \eqref{est_psi41} for the gradient $\nabla \psi_{31}$ follows along the same lines. It suffices to consider the gradient version of Duhamel’s formula for the heat equation and to repeat the same argument as above. \hfill$\square$

\medskip
\noindent
\textit{Proof of \eqref{est_psi42}.} We prove the estimate for $\psi_{32}$, the proof for the gradient is similar. Let $(x, t) \in \R^2 \times (t_0, +\infty)$ and $\psi_{32}$ given by
\begin{align*}
\psi_{32}(x, t) = \int_{t_0}^t \int_{\R^2} \frac{1}{4 \pi (t - s)} \exp\left(-\frac{|x - z|^2}{4(t - s)} \right) \tilde{g}_{32}(z, s) dzds.
\end{align*}
Let us assume that $t > 2t_0$ and $|x| < 2\sqrt{t}$. We have that
\begin{align*}
    \abs{\psi_{32}(x, t)} &\le C \int_{t_0}^t \frac{1}{s^{{d_0 + \theta}+1}}\int_{\R^2} \frac{1}{(t - s)} \exp\left(-\frac{|x - z|^2}{4(t - s)} \right)  \frac{1}{\big(1 + \frac{\abs{z}}{\sqrt{s}}\big)^2} dzds\\
    &= C \int_{t_0}^t \frac{1}{s^{{d_0 + \theta}+1}}\int_{\R^2} \frac{1}{(t - s)} \exp\left(-\frac{|\zeta|^2}{4(t - s)} \right)  \frac{1}{\big(1 + \frac{\abs{x - \zeta}}{\sqrt{s}}\big)^2} d\zeta ds.
\end{align*}
For $\tau = \sqrt{t - s}$, we consider $A$ as before
\begin{align*}
    A(\tau, s) &:= \frac{1}{\tau^2}\int_{\R^2} \exp\left(-\frac{|\zeta|^2}{4\tau^2} \right)  \frac{1}{\big(1 + \frac{\abs{x - \zeta}}{\sqrt{s}}\big)^2} d\zeta
\end{align*}
We claim that
\begin{align*}
    A(\tau, s) \le  \begin{cases}
         C\displaystyle{\frac{s}{\tau^2}} \log t & \text{for } s \in (t_0, t/2),\\[3pt]
         C & \text{for } s \in (t/2, t).
    \end{cases}
\end{align*}
For $s \in (t_0, t/2)$, let us divide $A$ as
\begin{align*}
    A(\tau, s) &= \frac{1}{\tau^2}\bigg(\int_{\abs{\zeta} < 2 \abs{x}}  \exp\left(-\frac{|\zeta|^2}{4\tau^2} \right)  \frac{1}{\big(1 + \frac{\abs{x - \zeta}}{\sqrt{s}}\big)^2} d\zeta\\
    &\hspace{1cm} + \int_{\abs{\zeta} > 2 \abs{x}}  \exp\left(-\frac{|\zeta|^2}{4\tau^2} \right)  \frac{1}{\big(1 + \frac{\abs{x - \zeta}}{\sqrt{s}}\big)^2} d\zeta\bigg) =: A_1 + A_2
\end{align*}
We compute
\begin{align*}
    A_1(\tau, s) &= \frac{1}{\tau^2}\int_{\abs{\zeta} < 2\abs{x}} \exp\left(-\frac{|\zeta|^2}{4\tau^2} \right)  \frac{1}{\big(1 + \frac{\abs{x - \zeta}}{\sqrt{s}}\big)^2} d\zeta\\
    &\le C\frac{1}{\tau^2}\int_{\abs{\zeta} < 2\abs{x}} \frac{1}{\big(1 + \frac{\abs{x - \zeta}}{\sqrt{s}}\big)^2} d\zeta\\
    &\le C\frac{s}{\tau^2} \log\left(1 + \frac{\abs{x}}{\sqrt{s}}\right)\\
    &\le C\frac{s}{\tau^2} \log t
\end{align*}
For $A_2$, we get
\begin{align*}
    A_2(\tau, s) &= \frac{1}{\tau^2} \bigg( \int_{\{\abs{\zeta} > 2 \abs{x}\} \cap \{\abs{x - \zeta} > \tau\}} \exp\left(-\frac{|\zeta|^2}{4\tau^2} \right)  \frac{1}{\big(1 + \frac{\abs{x - \zeta}}{\sqrt{s}}\big)^2} d\zeta \\
    &\hspace{1cm}+  \int_{\{\abs{\zeta} > 2 \abs{x}\} \cap \{\abs{x - \zeta} < \tau\}} \exp\left(-\frac{|\zeta|^2}{4\tau^2} \right) \frac{1}{\big(1 + \frac{\abs{x - \zeta}}{\sqrt{s}}\big)^2}d\zeta\bigg)=: \tilde{A}_1 + \tilde{A}_2\\
\end{align*}
The integrals $\tilde{A}_1$ and $\tilde{A}_2$ can be estimated as follows
\begin{align*}
    \tilde{A}_1(\tau, s)&\le \frac{s}{\tau^4} \int_{\{\abs{\zeta} > 2 \abs{x}\} \cap \{\abs{x - \zeta} > \tau\}}  \exp\left(-\frac{|\zeta|^2}{4\tau^2} \right) d\zeta\\
    &\le \frac{s}{\tau^2} \int_{\mathbb R^2} e^{-\frac{\abs{\xi}^2}{4}}d\xi\le C\frac{s}{\tau^2}\\
    \tilde{A}_2(\tau, s)&\le \frac{1}{\tau^2}\int_{\{\abs{\zeta} > 2 \abs{x}\} \cap \{\abs{x - \zeta} < \tau\}}  \exp\left(-\frac{|\zeta|^2}{4\tau^2} \right)  \frac{1}{\big(1 + \frac{\abs{x - \zeta}}{\sqrt{s}}\big)^2} d\zeta\\
    &\le C\frac{1}{\tau^2}\int_{\{\abs{x - \zeta} < \tau\}}  \frac{1}{\big(1 + \frac{\abs{x - \zeta}}{\sqrt{s}}\big)^2} d\zeta\\
    &\le C\frac{s}{\tau^2}\log\left(1 + \sfrac{\tau}{\sqrt{s}}\right) \le C \frac{s}{\tau^2}\log t
\end{align*}
For $s\in (t/2, t)$, we find that
\begin{align*}
    A(\tau, s) &\le   \frac{1}{\tau^2}\int_{\R^2}  \exp\left(-\frac{|\zeta|^2}{4\tau^2} \right)  d\zeta\\
    &= \int_{\R^2}  \exp\left(-\frac{|\xi|^2}{4} \right)  d\xi \le C.
\end{align*}

Therefore, we get that
\begin{align*}
    \abs{\psi_{32}(x, t)} &\le C\int_{t_0}^t \frac{1}{s^{{d_0 + \theta}+1}} A(\tau, s)ds\\
    &\le C\left(\int_{t_0}^{t/2} \frac{1}{s^{{d_0 + \theta}+1}} \frac{s}{(t - s)} \log t ~ds + \int_{t/2}^t \frac{1}{s^{{d_0 + \theta}+1}} ds\right)\\
    &\le C\frac{1}{t^{d_0 + \theta}} \log t.
\end{align*}
\vspace{0.2cm}
Let us now assume that $\abs{x} > 2\sqrt{t}$. We divide $A$ in three pieces $A_1, A_2, A_3$ as follows
\begin{align*}
    A(\tau, s, x) &= \frac{1}{\tau^2}\bigg(\int_{\abs{\zeta} < \frac{\abs{x}}{2} } \exp\left(-\frac{|\zeta|^2}{4\tau^2} \right)  \frac{1}{\big(1 + \frac{\abs{x - \zeta}}{\sqrt{s}}\big)^2} d\zeta\\
    &\hspace{1cm} + \int_{\frac{\abs{x}}{2}  < \abs{\zeta} < 2 \abs{x}}  \exp\left(-\frac{|\zeta|^2}{4\tau^2} \right)  \frac{1}{\big(1 + \frac{\abs{x - \zeta}}{\sqrt{s}}\big)^2} d\zeta\\
    &\hspace{1cm} + \int_{\abs{\zeta} > 2 \abs{x}}\exp\left(-\frac{|\zeta|^2}{4\tau^2} \right)  \frac{1}{\big(1 + \frac{\abs{x - \zeta}}{\sqrt{s}}\big)^2} d\zeta\bigg) =: A_1 + A_2 + A_3
\end{align*}
We claim that
\begin{align*}
    A_i(\tau, s, x) &\le  \frac{C}{\big(1 + \frac{\abs{x}}{\sqrt{s}} \big)^2}, \quad i = 1, 3,\\
    A_2(\tau, s, x) &\le C\frac{s}{\tau^2} \exp\left(-\sfrac{|x|^2}{16\tau^2}\right) \log\left(1 + \sfrac{\abs{x}}{\sqrt{s}}\right).
\end{align*}
We estimate $A_i$, for $i = 1, 3$, by
\begin{align*}
    A_1(\tau, s, x) &\le \frac{C}{\tau^2}\frac{1}{\big(1 + \frac{\abs{x}}{\sqrt{s}} \big)^2} \int_{\abs{\zeta} < \frac{\abs{x}}{2}} \exp\left(-\frac{|\zeta|^2}{4\tau^2} \right)d\zeta\\
    &= C\frac{1}{\big(1 + \frac{\abs{x}}{\sqrt{s}} \big)^2} \int_{\abs{\xi} < \frac{\abs{x}}{2\tau}} \exp\left(-\frac{|\xi|^2}{4} \right)d\xi\\
    &\le \frac{C}{\big(1 + \frac{\abs{x}}{\sqrt{s}} \big)^2},\\
    A_3(\tau, s, x) &\le \frac{C}{\tau^2}\frac{1}{\big(1 + \frac{\abs{x}}{\sqrt{s}} \big)^2} \int_{\abs{\zeta} > 2 \abs{x}} \exp\left(-\frac{|\zeta|^2}{4\tau^2} \right) d\zeta\\
    &= \frac{C}{\big(1 + \frac{\abs{x}}{\sqrt{s}} \big)^2} \int_{\abs{\xi} > 2 \abs{x}} \exp\left(-\frac{|\xi|^2}{4} \right)d\xi\\
    &\le \frac{C}{\big(1 + \frac{\abs{x}}{\sqrt{s}} \big)^2},
\end{align*}
and for $A_2$ we get
\begin{align*}
    A_2(\tau, s) &\le \frac{C}{\tau^2} \exp\left(-\frac{|x|^2}{16\tau^2} \right)\int_{\frac{\abs{x}}{2}  < \abs{\zeta} < 2 \abs{x}}   \frac{1}{\big(1 + \frac{\abs{x - \zeta}}{\sqrt{s}}\big)^2} d\zeta\\
    &\le C\frac{s}{\tau^2} \exp\left(-\sfrac{|x|^2}{16\tau^2} \right)\log\left(1 + \sfrac{\abs{x}}{\sqrt{s}}\right).
\end{align*}
Therefore, for $i = 1, 3$, we get
\begin{align*}
    \int_{t_0}^t \frac{1}{s^{{d_0 + \theta}+1}} A_i(\tau, s, x)ds &\le C\int_{t_0}^t \frac{1}{s^{{d_0 + \theta}+1}}  \sfrac{1}{\big(\abs{x}/\sqrt{s} \big)^2}ds\\
    &\le C\frac{1}{\abs{x}^2}\int_{t_0}^t \frac{1}{s^{d_0 + \theta}}ds\\
    &\le C\frac{1}{t^{d_0 + \theta}}\sfrac{1}{\big(\abs{x}/\sqrt{t} \big)^2}
\end{align*}
Finally, for the last term we find
\begin{align*}
    \int_{t_0}^t \frac{1}{s^{{d_0 + \theta}+1}} A_2(\tau, s, x)ds &\le C\int_{t_0}^t \frac{1}{s^{d_0 + \theta}} \frac{1}{(t - s)} \exp\left(-\frac{|x|^2}{16(t - s)} \right)\log\left(1 + \frac{\abs{x}}{\sqrt{s}}\right)ds\\
    &\le C\frac{1}{t^{d_0 + \theta}} \exp\left(-\sfrac{|x|^2}{16t} \right) \log\left(1 + \sfrac{\abs{x}}{\sqrt{t}}\right).
\end{align*}
This shows that
\begin{align*}
     \abs{\psi_{32}(x, t)} \le C \frac{(\log t)^2}{t^{d_0 + \theta}} \sfrac{1}{\big(\abs{x}/\sqrt{t} \big)^2}
\end{align*}
for $\abs{x} > 2\sqrt{t}$. The conclusion follows. The proof of the estimate \eqref{est_psi42} for the gradient $\nabla \psi_{32}$ follows along the same lines. It suffices to consider the gradient version of Duhamel’s formula for the heat equation and to repeat the same argument as above. \hfill $\square$

\section{Estimates for \texorpdfstring{$Z_0$}{Z0}}
\label{App:est_Z0}

\begin{proof}[Proof of  Lemma~\ref{Lem:est_Z}]
The solution of \eqref{eq:heat-eq-Z} given by Duhamel's formula with
initial data $Z_0$ is
\begin{align*}
Z(x, t) = \int_{\R^2} \frac{1}{4\pi(t - t_0)} e^{-\frac{\abs{x - y}^2}{4(t - t_0)}}Z_0(y) dy.
\end{align*}
Let us first assume that $t \ge 2t_0$ and $\abs{x} < \sqrt{t}$. We denote $\tau = \sqrt{t - t_0}$. We get that
\begin{align*}
    \abs{Z(x, t)} &\le C\normZ t_0^{-d} \int_{\R^2} \frac{1}{\tau^2} e^{-\frac{\abs{x - y}^2}{4\tau^2}} \frac{1}{\big(1 + \frac{\abs{y}}{\sqrt{t_0}}\big)^{2d}} dy\\
    &= C\normZ t_0^{-d} \int_{\R^2} \frac{1}{\tau^2} e^{-\frac{\abs{\zeta}^2}{4\tau^2}} \frac{1}{\big(1 + \frac{\abs{x - \zeta}}{\sqrt{t_0}}\big)^{2d}} d\zeta\\
    &= C\normZ t_0^{-d} \bigg(\int_{\abs{\zeta} < 2\abs{x}} \frac{1}{\tau^2} e^{-\frac{\abs{\zeta}^2}{4\tau^2}} \frac{1}{\big(1 + \frac{\abs{x - \zeta}}{\sqrt{t_0}}\big)^{2d}} d\zeta \\
    & \hspace{3cm} + \int_{\abs{\zeta} > 2\abs{x}} \frac{1}{\tau^2} e^{-\frac{\abs{\zeta}^2}{4\tau^2}} \frac{1}{\big(1 + \frac{\abs{x - \zeta}}{\sqrt{t_0}}\big)^{2d}} d\zeta\bigg) := I_1 + I_2.
\end{align*}
For $I_1$, we find
\begin{align*}
    I_1 &\le C\normZ \frac{ t_0^{- d }}{\tau^2} \int_{\abs{\zeta} < 2\abs{x}}\frac{1}{\big(1 + \frac{\abs{x - \zeta}}{\sqrt{t_0}}\big)^{2d}} d\zeta\\
    &\le  C\normZ \frac{t_0^{-d + 1}}{\tau^2} \int_{|\tilde \zeta| < 2|\tilde x|}\frac{1}{(1 + |\tilde x - \tilde \zeta|)^{2d}} d\tilde \zeta, \quad \tilde{x} = \frac{x}{\sqrt{t_0}},\\
    &\le C\normZ \frac{t_0^{-d + 1}}{\tau^2} \bigg(\frac{\abs{x}}{\sqrt{t_0}}\bigg)^{2 - 2d}\\
    &\le C\normZ \frac{1}{t^d},
\end{align*}
where the last inequality follows from the assumptions $t \ge 2t_0$ and $\abs{x} < \sqrt{t}$. We now split the integral $I_2$ into two terms
\begin{align*}
    I_2 &= C\normZ t_0^{-d}\bigg(\int_{\{\abs{\zeta} > 2\abs{x}\} \cap \{|x - \zeta| < \tau\}} \frac{1}{\tau^2} e^{-\frac{\abs{\zeta}^2}{4\tau^2}} \frac{1}{\big(1 + \frac{\abs{x - \zeta}}{\sqrt{t_0}}\big)^{2d}}d\zeta \\
    &\hspace{3cm} + \int_{\{\abs{\zeta} > 2\abs{x}\} \cap \{|x - \zeta| > \tau\}}\frac{1}{\tau^2} e^{-\frac{\abs{\zeta}^2}{4\tau^2}} \frac{1}{\big(1 + \frac{\abs{x - \zeta}}{\sqrt{t_0}}\big)^{2d}} d\zeta \bigg) =: I_{21} + I_{22}.
\end{align*}
These two integrals estimate as follows
\begin{align*}
    I_{21} &\le C\normZ t_0^{-d} \int_{\{|x - \zeta| < \tau\}} \frac{1}{\tau^2} \frac{1}{\big(1 + \frac{\abs{x - \zeta}}{\sqrt{t_0}}\big)^{2d}}d\zeta\\
    &\le C\normZ \frac{1}{t^d},\\
    I_{22} &\le  C\normZ \frac{1}{\tau^{2 + 2d}} \int_{\{\abs{\zeta} > 2\abs{x}\}} e^{-\frac{\abs{\zeta}^2}{4\tau^2}} d\zeta \\
    &\le C\normZ \frac{1}{t^d}.
\end{align*}
This proves the desired estimate inside the self-similar region $|x| < \sqrt{t}$. Let us now assume the $\abs{x} > \sqrt{t}$. We divide the integral into three pieces
\begin{align*}
    \abs{Z(x, t)} &\le C\normZ t_0^{-d} \int_{\R^2} \frac{1}{\tau^2} e^{-\frac{\abs{\zeta}^2}{4\tau^2}} \frac{1}{\big(1 + \frac{\abs{x - \zeta}}{\sqrt{t_0}}\big)^{2d}} d\zeta\\
    &\le  C\normZ t_0^{-d} \bigg( \int_{|\zeta| < \frac{|x|}{2}} \frac{1}{\tau^2} e^{-\frac{\abs{\zeta}^2}{4\tau^2}} \frac{1}{\big(1 + \frac{\abs{x - \zeta}}{\sqrt{t_0}}\big)^{2d}} d\zeta\\
    & ~~~+\int_{|\zeta| > 2|x|} \frac{1}{\tau^2} e^{-\frac{\abs{\zeta}^2}{4\tau^2}} \frac{1}{\big(1 + \frac{\abs{x - \zeta}}{\sqrt{t_0}}\big)^{2d}} d\zeta + \int_{\frac{|x|}{2} <|\zeta| < 2|x|} \frac{1}{\tau^2} e^{-\frac{\abs{\zeta}^2}{4\tau^2}} \frac{1}{\big(1 + \frac{\abs{x - \zeta}}{\sqrt{t_0}}\big)^{2d}} d\zeta\bigg) \\
    &=: J_1 + J_2 + J_3.
\end{align*}
The first two integrals estimate in a similar way. More precisely,
\begin{align*}
    J_1 + J_2 &\le C\normZ  \frac{1}{t_0^{d}\big(\abs{x}/\sqrt{t_0}\big)^{2d}} \\
    &= C\normZ \frac{1}{t^d\big(\abs{x}/\sqrt{t}\big)^{2d}}.
\end{align*}
Finally, for $J_3$ we get that
\begin{align*}
    J_3 &\le C\normZ t_0^{-d} \frac{1}{\tau^2}e^{-\frac{\abs{x}^2}{16\tau^2}} \int_{\frac{|x|}{2} <|\zeta| < 2|x|} \frac{1}{\big(1 + \frac{\abs{x - \zeta}}{\sqrt{t_0}}\big)^{2d}} d\zeta\\
    &\le C\normZ t_0^{-d+ 1} \frac{1}{\tau^2}e^{-\frac{\abs{x}^2}{16\tau^2}}  \bigg(\frac{\abs{x}}{\sqrt{t_0}}\bigg)^{2 - 2d}\\
    &\le C\normZ  \frac{1}{t^{d}} \bigg(\frac{\abs{x}}{\sqrt{t}}\bigg)^{2 - 2d}e^{-\frac{\abs{x}^2}{16t}}
\end{align*}
From the estimate of $J_1, J_2$ and $J_3$, we deduce that
\begin{align*}
    \abs{Z(x, t)} \le C\normZ \frac{1}{t^d\big(\abs{x}/\sqrt{t}\big)^{2d}}
\end{align*}
for $\abs{x} > \sqrt{t}$. \\[0.3cm]
We now consider the case of $t \in (t_0, 2t_0).$ We must prove that
\begin{align}
    \abs{Z(x, t)} &\le C\frac{\normZ}{t_0^{d}\big(1 + \frac{\abs{x}}{\sqrt{t_0}}\big)^{2d}}. \label{est_Z0_t0-2t0}
\end{align}
For $\abs{x} < \sqrt{t}$, we get
\begin{align*}
    \abs{Z(x, t)} \le \normZ t_0^{-d} \int_{\R^2} e^{-\frac{\abs{\xi}^2}{4}}d\xi \le C\normZ t_0^{-d}.
\end{align*}
On the other hand, for $\abs{x} > \sqrt{t}$, we have
\begin{align*}
\abs{Z(x, t)} 
&\le \normZ t_0^{-d} \int_{\R^2} \frac{1}{\tau^2}e^{-\frac{\abs{\zeta}^2}{4\tau^2}} \frac{1}{\big(1 + \frac{\abs{x -  \zeta}}{\sqrt{t_0}}\big)^{2d}}d\zeta
\\
&\le C \normZ t_0^{-d} \bigg( \int_{|\zeta| < \frac{|x|}{2}} \frac{1}{\tau^2} e^{-\frac{\abs{\zeta}^2}{4\tau^2}} \frac{1}{\big(1 + \frac{\abs{x - \zeta}}{\sqrt{t_0}}\big)^{2d}} d\zeta
\\
&\quad 
+\int_{|\zeta| > 2|x|} \frac{1}{\tau^2} e^{-\frac{\abs{\zeta}^2}{4\tau^2}} \frac{1}{\big(1 + \frac{\abs{x - \zeta}}{\sqrt{t_0}}\big)^{2d}} d\zeta 
\\
& \quad
+ \int_{\frac{|x|}{2} <|\zeta| < 2|x|} \frac{1}{\tau^2} e^{-\frac{\abs{\zeta}^2}{4\tau^2}} \frac{1}{\big(1 + \frac{\abs{x - \zeta}}{\sqrt{t_0}}\big)^{2d}} d\zeta\bigg)\\
    &=: L_1 + L_2 + L_3.
\end{align*}
The integrals $L_1$ and $L_2$ estimate as $J_1$ and $J_2$ above, leading to
\begin{align*}
L_1 + L_2 &\le C\normZ t_0^{-d} \frac{1}{\big(\abs{x}/\sqrt{t_0}\big)^{2d}}
\end{align*}
whereas for $L_3$, we obtain
\begin{align*}
L_3 &\le C \normZ t_0^{-d} e^{-\frac{\abs{x}^2}{16\tau^2}} \int_{\frac{|x|}{2\tau} <|\xi| < \frac{2|x|}{\tau}}  d\xi
\\
&\le C \normZ t_0^{-d} \left(\frac{\abs{x}}{\tau}\right)^2e^{-\frac{\abs{x}^2}{16\tau^2}}\\
&\le \frac{C \normZ }{t_0^d\big(\abs{x}/\tau \big)^{2d}} \le \frac{C \normZ }{t_0^d\big(\abs{x}/\sqrt{t_0}\big)^{2d}}.
\end{align*}
This proves that
\begin{align*}
    \abs{Z(x, t)} \le C\frac{\|Z_0\|_{0}}{t^d\big(1 + \frac{\abs{x}}{\sqrt{t}}\big)^{2d}} \quad \text{in } \R^2 \times (t_0, +\infty),
\end{align*}
for some constant $C = C(d)$.\\[0.3cm]
The estimate for $\nabla Z$ follows along the same lines. Recall that
\begin{align*}
\nabla Z(x, t) = \int_{\R^2} \nabla_x G(x - y, t - t_0) Z_0(y) dy, 
\end{align*}
where
\begin{align*}
G(x, t) = \frac{1}{4\pi t} e^{-\frac{\abs{x}^2}{4t}}.
\end{align*}
This yields
\begin{align*}
\abs{\nabla Z(x, t)} &\le C\normZ \int_{\R^2} \frac{\abs{x - y}}{\tau^4}  e^{-\frac{\abs{x - y}^2}{4\tau^2}} \frac{1}{t_0^d\big(1 + \sfrac{\abs{y}}{\sqrt{t_0}}\big)^{2d}}dy,
\end{align*}
for $\tau = \sqrt{t - t_0}$. Considering first $t > 2t_0$ and proceeding as in the estimate for $Z$, we divide the estimate in the two regions $\abs{x} < \sqrt{t}$ and $\abs{x} > \sqrt{t}$, together with the same decomposition of the integrals. This leads to estimate $\nabla Z$ as
\begin{align*}
\abs{\nabla Z(x, t)} &\le  C\frac{\|Z_0\|_{0}}{t^{d + \frac{1}{2}}\big(1 + \frac{\abs{x}}{\sqrt{t}}\big)^{2d}}
\end{align*}
in $\R^2 \times (2t_0, +\infty)$. Additionally, for $t \in (t_0, 2t_0)$, notice that
\begin{align*}
    \nabla Z(x, t) &= \int_{\R^2}  G(x - y, t - t_0) \nabla Z_0(y) dy,
\end{align*}
Proceeding as in \eqref{est_Z0_t0-2t0} for the estimate of $Z(x, t)$, we deduce that
\begin{align*}
\abs{\nabla Z(x, t)} &\le  C\frac{\|Z_0\|_{0}}{t_0^{d + \frac{1}{2}}
\big(1 + \frac{\abs{x}}{\sqrt{t_0}}\big)^{2d}},
\end{align*}
in $\R^2 \times (t_0, 2t_0)$. This yields
\begin{align*}
\sqrt{t}\abs{\nabla Z(x, t)} &\le  C\frac{\|Z_0\|_{0}}{t^{d}\big(1 + \frac{\abs{x}}{\sqrt{t}}\big)^{2d}}
\end{align*}
in $\R^2 \times (t_0, +\infty)$, which completes the proof.
\end{proof}

\bigskip\noindent
{\bf Funding.} J.~D\'avila has been supported  by  a Royal Society  Wolfson Fellowship, UK, Grant RSWF/FT/191007. 
    
\bibliographystyle{siam}
\bibliography{bibliography}

@article {jendrej-lawrie-schlag,
    AUTHOR = {Jendrej, Jacek and Lawrie, Andrew and Schlag, Wilhelm},
     TITLE = {Continuous in time bubble decomposition for the harmonic map
              heat flow},
   JOURNAL = {Forum Math. Pi},
  FJOURNAL = {Forum of Mathematics. Pi},
    VOLUME = {13},
      YEAR = {2025},
     PAGES = {Paper No. e4, 37},
      ISSN = {2050-5086},
   MRCLASS = {58E20 (35K58 58J35)},
  MRNUMBER = {4862091},
MRREVIEWER = {Junichi\ Harada},
       DOI = {10.1017/fmp.2024.15},
       URL = {https://doi.org/10.1017/fmp.2024.15},
}

@article {wang-1996,
    AUTHOR = {Wang, Changyou},
     TITLE = {Bubble phenomena of certain {P}alais-{S}male sequences from
              surfaces to general targets},
   JOURNAL = {Houston J. Math.},
  FJOURNAL = {Houston Journal of Mathematics},
    VOLUME = {22},
      YEAR = {1996},
    NUMBER = {3},
     PAGES = {559--590},
      ISSN = {0362-1588},
   MRCLASS = {58E20 (49Q10)},
  MRNUMBER = {1417632},
MRREVIEWER = {Martin\ Fuchs},
}

@article {lin-wang-1998,
    AUTHOR = {Lin, Fanghua and Wang, Changyou},
     TITLE = {Energy identity of harmonic map flows from surfaces at finite
              singular time},
   JOURNAL = {Calc. Var. Partial Differential Equations},
  FJOURNAL = {Calculus of Variations and Partial Differential Equations},
    VOLUME = {6},
      YEAR = {1998},
    NUMBER = {4},
     PAGES = {369--380},
      ISSN = {0944-2669,1432-0835},
   MRCLASS = {58E20 (35K55 58G11)},
  MRNUMBER = {1624304},
MRREVIEWER = {Joseph\ F.\ Grotowski},
       DOI = {10.1007/s005260050095},
       URL = {https://doi.org/10.1007/s005260050095},
}

@book{Lojasiewicz1965,
  author       = {{\L}ojasiewicz, Stanis{\l}aw},
  title        = {Ensembles semi-analytiques},
  publisher    = {Institut des Hautes Études Scientifiques},
  address      = {Bures-sur-Yvette},
  year         = {1965},
  note         = {LaTeX version prepared by Michel Coste on August 29, 2006,
                  based on mimeographed course notes by S. {\L}ojasiewicz},
  url          = {https://perso.univ-rennes1.fr/michel.coste/Lojasiewicz.pdf}
}

@article {delpino-musso-wei-3d-crit,
    AUTHOR = {del Pino, Manuel and Musso, Monica and Wei, Juncheng},
     TITLE = {Infinite-time blow-up for the 3-dimensional energy-critical
              heat equation},
   JOURNAL = {Anal. PDE},
  FJOURNAL = {Analysis \& PDE},
    VOLUME = {13},
      YEAR = {2020},
    NUMBER = {1},
     PAGES = {215--274},
      ISSN = {2157-5045,1948-206X},
   MRCLASS = {35K91 (35B33 35B40 35B44 35K15 35K58)},
  MRNUMBER = {4047646},
       DOI = {10.2140/apde.2020.13.215},
       URL = {https://doi.org/10.2140/apde.2020.13.215},
}

@article {ddwz-deg2,
    AUTHOR = {D\'avila, Juan and del Pino, Manuel and Wei, Juncheng and
              Zheng, Youquan},
     TITLE = {Singularity formation of two dimensional harmonic map flow:
              the degree 2 case},
   JOURNAL = {Discrete Contin. Dyn. Syst.},
  FJOURNAL = {Discrete and Continuous Dynamical Systems},
    VOLUME = {56},
      YEAR = {2026},
     PAGES = {229--270},
      ISSN = {1078-0947,1553-5231},
   MRCLASS = {35K51 (35B44 35K59)},
  MRNUMBER = {5097206},
       DOI = {10.3934/dcds.2026105},
       URL = {https://doi.org/10.3934/dcds.2026105},
}

@article {Chang,
    AUTHOR = {Chang, Kung-Ching},
     TITLE = {Heat flow and boundary value problem for harmonic maps},
   JOURNAL = {Ann. Inst. H. Poincar\'{e} C Anal. Non Lin\'{e}aire},
  FJOURNAL = {Annales de l'Institut Henri Poincar\'{e} C. Analyse Non
              Lin\'{e}aire},
    VOLUME = {6},
      YEAR = {1989},
    NUMBER = {5},
     PAGES = {363--395},
}

@article {HMF_DdPW,
    AUTHOR = {D\'{a}vila, Juan and del Pino, Manuel and Wei, Juncheng},
    TITLE = {Singularity formation for the two-dimensional harmonic map
              flow into {$S^2$}},
   JOURNAL = {Invent. Math.},
   FJOURNAL = {Inventiones Mathematicae},
   VOLUME = {219},
   YEAR = {2020},
   NUMBER = {2},
}

@article{Simon,
  title={Asymptotics for a class of non-linear evolution equations, with applications to geometric problems},
  author={Leon Simon},
  journal={Annals of Mathematics},
  year={1983},
  volume={118},
  pages={525-571},
  url={https://api.semanticscholar.org/CorpusID:124839705}
}

@incollection {topping_almost_harm,
    AUTHOR = {Topping, Peter},
     TITLE = {Bubbling of almost-harmonic maps between 2-spheres at points
              of zero energy density},
 BOOKTITLE = {Variational problems in {R}iemannian geometry},
    SERIES = {Progr. Nonlinear Differential Equations Appl.},
    VOLUME = {59},
     PAGES = {33--42},
 PUBLISHER = {Birkh\"auser, Basel},
      YEAR = {2004},
      ISBN = {3-7643-2432-5},
   MRCLASS = {58E20},
  MRNUMBER = {2076265},
MRREVIEWER = {Andreas\ Gastel},
}

@article {Eells,
    AUTHOR = {Eells, Jr., James and Sampson, J. H.},
     TITLE = {Harmonic mappings of {R}iemannian manifolds},
   JOURNAL = {Amer. J. Math.},
  FJOURNAL = {American Journal of Mathematics},
    VOLUME = {86},
      YEAR = {1964},
     PAGES = {109--160},
}

@article {CLW,
    AUTHOR = {Chen, Guoyuan and Liu, Yong and Wei, Juncheng},
     TITLE = {Nondegeneracy of harmonic maps from {$\mathbb R^2$} to {$\mathbb
              S^2$}},
   JOURNAL = {Discrete Contin. Dyn. Syst.},
  FJOURNAL = {Discrete and Continuous Dynamical Systems},
    VOLUME = {40},
      YEAR = {2020},
    NUMBER = {6},
     PAGES = {3215--3233},
}

@article {Raphael,
    AUTHOR = {Rapha\"{e}l, Pierre and Schweyer, Remi},
     TITLE = {Stable blowup dynamics for the 1-corotational energy critical
              harmonic heat flow},
   JOURNAL = {Comm. Pure Appl. Math.},
  FJOURNAL = {Communications on Pure and Applied Mathematics},
    VOLUME = {66},
      YEAR = {2013},
    NUMBER = {3},
}

@article {Struwe,
    AUTHOR = {Struwe, Michael},
     TITLE = {On the evolution of harmonic mappings of {R}iemannian
              surfaces},
   JOURNAL = {Comment. Math. Helv.},
  FJOURNAL = {Commentarii Mathematici Helvetici},
    VOLUME = {60},
      YEAR = {1985},
    NUMBER = {4},
     PAGES = {558--581},
}

@article {topping_anals,
    AUTHOR = {Topping, Peter},
     TITLE = {Repulsion and quantization in almost-harmonic maps, and
              asymptotics of the harmonic map flow},
   JOURNAL = {Ann. of Math. (2)},
  FJOURNAL = {Annals of Mathematics. Second Series},
    VOLUME = {159},
      YEAR = {2004},
    NUMBER = {2},
     PAGES = {465--534},
}

@article {Ding-Tian,
    AUTHOR = {Ding, Weiyue and Tian, Gang},
     TITLE = {Energy identity for a class of approximate harmonic maps from
              surfaces},
   JOURNAL = {Comm. Anal. Geom.},
  FJOURNAL = {Communications in Analysis and Geometry},
    VOLUME = {3},
      YEAR = {1995},
    NUMBER = {3-4},
     PAGES = {543--554},
      ISSN = {1019-8385,1944-9992},
   MRCLASS = {58E20 (58G11)},
  MRNUMBER = {1371209},
MRREVIEWER = {Joseph\ F.\ Grotowski},
}

@article {Qing,
    AUTHOR = {Qing, Jie},
     TITLE = {On singularities of the heat flow for harmonic maps from
              surfaces into spheres},
   JOURNAL = {Comm. Anal. Geom.},
  FJOURNAL = {Communications in Analysis and Geometry},
    VOLUME = {3},
      YEAR = {1995},
    NUMBER = {1-2},
     PAGES = {297--315},
      ISSN = {1019-8385,1944-9992},
   MRCLASS = {58G11 (35K55 58E20)},
  MRNUMBER = {1362654},
MRREVIEWER = {Joseph\ F.\ Grotowski},
}

@article {GuanGusTsai,
    AUTHOR = {Guan, Meijiao and Gustafson, Stephen and Tsai, Tai-Peng},
     TITLE = {Global existence and blow-up for harmonic map heat flow},
   JOURNAL = {J. Differential Equations},
  FJOURNAL = {Journal of Differential Equations},
    VOLUME = {246},
      YEAR = {2009},
    NUMBER = {1},
     PAGES = {1--20},
}

@article {GusTsaiNak,
    AUTHOR = {Gustafson, Stephen and Nakanishi, Kenji and Tsai, Tai-Peng},
     TITLE = {Asymptotic stability, concentration, and oscillation in
              harmonic map heat-flow, {L}andau-{L}ifshitz, and
              {S}chr\"odinger maps on {$\mathbb R^2$}},
   JOURNAL = {Comm. Math. Phys.},
  FJOURNAL = {Communications in Mathematical Physics},
    VOLUME = {300},
      YEAR = {2010},
    NUMBER = {1},
     PAGES = {205--242},
}

@article {vdBBHK,
    AUTHOR = {van den Berg, Jan Bouwe and Hulshof, Josephus and King, John
              R.},
     TITLE = {Formal asymptotics of bubbling in the harmonic map heat flow},
   JOURNAL = {SIAM J. Appl. Math.},
  FJOURNAL = {SIAM Journal on Applied Mathematics},
    VOLUME = {63},
      YEAR = {2003},
    NUMBER = {5},
     PAGES = {1682--1717},
      ISSN = {0036-1399,1095-712X},
   MRCLASS = {35K55 (35B35 35B40 53C44 58E20 76A15)},
  MRNUMBER = {2001215},
MRREVIEWER = {Takayoshi\ Ogawa},
}

@article {Qing-Tian,
    AUTHOR = {Qing, Jie and Tian, Gang},
     TITLE = {Bubbling of the heat flows for harmonic maps from surfaces},
   JOURNAL = {Comm. Pure Appl. Math.},
  FJOURNAL = {Communications on Pure and Applied Mathematics},
    VOLUME = {50},
      YEAR = {1997},
    NUMBER = {4},
     PAGES = {295--310},
      ISSN = {0010-3640,1097-0312},
   MRCLASS = {58E20 (58G11)},
  MRNUMBER = {1438148},
MRREVIEWER = {Daniel\ Pollack},
}

@article {Keller_segel_infinite_time,
    AUTHOR = {D\'avila, Juan and del Pino, Manuel and Dolbeault, Jean and
              Musso, Monica and Wei, Juncheng},
     TITLE = {Existence and stability of infinite time blow-up in the
              {K}eller-{S}egel system},
   JOURNAL = {Arch. Ration. Mech. Anal.},
  FJOURNAL = {Archive for Rational Mechanics and Analysis},
    VOLUME = {248},
      YEAR = {2024},
    NUMBER = {4},
     PAGES = {Paper No. 61, 154},
}

@article {skyrmions,
    AUTHOR = {Bernand-Mantel, Anne and Muratov, Cyrill B. and Simon, Theresa
              M.},
     TITLE = {A quantitative description of skyrmions in ultrathin
              ferromagnetic films and rigidity of degree {$\pm 1$} harmonic
              maps from {$\mathbb{R}^2$} to {$\mathbb {S}^2$}},
   JOURNAL = {Arch. Ration. Mech. Anal.},
  FJOURNAL = {Archive for Rational Mechanics and Analysis},
    VOLUME = {239},
      YEAR = {2021},
    NUMBER = {1},
     PAGES = {219--299},
}

@article {HZ,
    AUTHOR = {Hirsch, Jonas and Zemas, Konstantinos},
     TITLE = {A note on a rigidity estimate for degree {$\pm 1$} conformal
              maps on {$\mathbb S^2$}},
   JOURNAL = {Bull. Lond. Math. Soc.},
  FJOURNAL = {Bulletin of the London Mathematical Society},
    VOLUME = {54},
      YEAR = {2022},
    NUMBER = {1},
     PAGES = {256--263},
}

@article {Topping_rigidity,
    AUTHOR = {Topping, Peter M.},
     TITLE = {A rigidity estimate for maps from {$S^2$} to {$S^2$} via the
              harmonic map flow},
   JOURNAL = {Bull. Lond. Math. Soc.},
  FJOURNAL = {Bulletin of the London Mathematical Society},
    VOLUME = {55},
      YEAR = {2023},
    NUMBER = {1},
     PAGES = {338--343},
}

@article {DSW_rigidity,
    AUTHOR = {Deng, Bin and Sun, Liming and Wei, Juncheng},
     TITLE = {Quantitative stability of harmonic maps from {$\mathbb R^2$} to {$\mathbb  S^2$} with a higher degree},
   JOURNAL = {Calc. Var. Partial Differential Equations},
  FJOURNAL = {Calculus of Variations and Partial Differential Equations},
    VOLUME = {63},
      YEAR = {2024},
    NUMBER = {4},
     PAGES = {Paper No. 101, 34},
}

@unpublished{Rupflin,
    AUTHOR = {Rupflin, Melanie},
    TITLE = {Sharp quantitative rigidity results for maps from {$\mathbb S^2$} to {$\mathbb S^2$}},
    YEAR = {2023},
    eprint={2305.17045},
    archivePrefix={arXiv},
    note = {ArXiv 2305.17045},
}

\end{document}